\documentclass[a4paper,10pt]{article}
\usepackage[utf8]{inputenc}
\usepackage{amsthm}
\usepackage{amsfonts}
\usepackage{amssymb}	
\usepackage{amsmath}
\allowdisplaybreaks
\usepackage{mathtools}
\usepackage[english]{babel}
\usepackage{color}
\usepackage{cite}
\usepackage[numbers]{natbib}
\usepackage{slashed}
\usepackage{enumerate}
\usepackage{graphicx}
\usepackage{systeme}
\usepackage{dsfont}
\usepackage{relsize}
\usepackage[margin=0.90in]{geometry}
\usepackage{float}
\usepackage{subcaption}
\usepackage{tikz}
\usepackage[labelformat=simple]{subcaption}

\usepackage{graphicx}
\usepackage{subcaption}
\usepackage{bmpsize}
\usepackage{epstopdf}

\usepackage{bbm}
\usepackage{xcolor,etoolbox}
\usepackage{titling}
\usepackage{authblk}
\newcommand*\samethanks[1][\value{footnote}]{\footnotemark[#1]}

\usepackage{blindtext,graphicx}
\usepackage[absolute]{textpos}
\usepackage{physics}
\usepackage[hang,flushmargin]{footmisc}

\usepackage{xfrac}
\usepackage{nicefrac}

\usepackage{esvect}
\usepackage{cases}
\usepackage{empheq}

\usepackage{stmaryrd}

\theoremstyle{definition}

\makeatletter
\newcommand{\pushright}[1]{\ifmeasuring@#1\else\omit\hfill$\displaystyle#1$\fi\ignorespaces}
\newcommand{\pushleft}[1]{\ifmeasuring@#1\else\omit$\displaystyle#1$\hfill\fi\ignorespaces}
\makeatother

\title{\vspace{-2cm}Analytical solutions of the cylindrical bending problem for the relaxed micromorphic continuum and other generalized continua (including full derivations)}
\author{Gianluca Rizzi\thanks{GEOMAS, INSA-Lyon, Universit\'e de Lyon, 20 avenue Albert Einstein,	69621, Villeurbanne cedex, France} \quad and \quad Geralf Hütter\thanks{TU Bergakademie Freiberg, Institute of Mechanics and Fluid Dynamics, Lampadiusstr. 4, 09596 Freiberg, Germany} \quad and \quad Angela Madeo\samethanks[1] \quad and \quad Patrizio Neff\thanks{Head of Chair for Nonlinear Analysis and Modelling, Fakultät für Mathematik, Universität Duisburg-Essen, \\ \indent Thea-Leymann-Straße 9, 45127 Essen, Germany}}

\thanksmarkseries{arabic}
\date{\today}
\begin{document}
\maketitle
\vspace{-0.75cm}
\begin{center}
	\textit{Dedicated to Holm Altenbach on the occasion of his 65th birthday.}
\end{center}
\begin{abstract}
We consider the cylindrical bending problem for an infinite plate as modelled with a family of generalized continuum models,
including the micromorphic approach.
The models allow to describe length scale effects in the sense that thinner specimens are comparatively stiffer. 
We provide the analytical solution for each case and exhibit the predicted bending stiffness.
The relaxed micromorphic continuum shows bounded bending stiffness for arbitrary thin specimens, while classical micromorphic continuum or gradient elasticity as well as Cosserat models \cite{neff2010stable} exhibit unphysical unbounded bending stiffness for arbitrary thin specimens.
This finding highlights the advantage of using the relaxed micromorphic model, which has a definite limit stiffness for small samples and which aids in identifying the relevant material parameters.
\end{abstract}

\section{Introduction of the cylindrical bending problem}
\label{sec:intro}
The classical Cauchy-Boltzmann theory of continuum mechanics comes to its limits when the wavelengths of deformation fields become comparable to characteristic microstructural length scales of a material, and therefore generalized continuum theories are necessary in this regime. 
Micromorphic continuum theories are established today for this purpose due to their direct connection to concepts of classical continuum mechanics and their relatively simple numerical implementation in existing finite element codes.
The micromorphic theory contains other approaches like the Cosserat continuum \cite{neff2009new,lakes1995experimental,lakes1998elastic}, couple stress continuum \cite{ghiba2017variant,madeo2016new}, microstretch continuum, micro-void continuum, the microstrain continuum, the relaxed micromorphic continuum or the strain-gradient theory as special cases \cite{dell2009generalized}, among others.
An overview can be found in \cite{Forest2018,Forest2019}.
In the very most cases, any of these subclasses is favoured over the general micromorphic theory since the latter requires a large number of constitutive parameters, namely already 18 for isotropic linear-elastic material.
Even more, the effect of most of these parameters on measurable quantities is still rather unclear.

An important milestone was set by Gauthier and Jahsman \cite{gauthier1975quest}, who derived analytical solutions for bending and torsion of a linear isotropic Cosserat continuum. The gradients under bending and torsion problems activate the non-classical terms, which is why the respective non-classical parameters affect the observable stiffness.
Further studies for bending of a linear isotropic Cosserat continuum have been carried out by Altenbach \cite{altenbach2009linear,altenbach2010generalized}.
These solutions have been used in many works to identify the constitutive parameters from real or virtual experiments, e.g. in \cite{gauthier1975quest,Yang1982,Lakes1983,Waseem2013,rueger2019cosserat,tekouglu2008size}
or they could be a useful tool to reduce the computational load in numerical identification or homogenization schemes, e.g. in \cite{arroyo2005continuum,brcic2013estimation,corigliano2005chip,rueger2019cosserat,zhang2017application,renda2020geometric}.

However, certain effects at edges or the dispersion of longitudinal waves could not be explained by the Cosserat theory \cite{lurie2018bending}, but more general micromorphic theories are required.
The authors in \cite{Cowin1983} provided the bending solution for the micro-void continuum.
However, this theory cannot explain experimentally observed size effects in torsion or dispersion of shear waves.
De Cicco and Nappa \cite{DeCicco1997} solved the bending and torsion problem for the microstretch continuum analytically.
Further closed-form solutions of the bending and torsion problems have been found for sub-classes of the micromorphic continuum and certain geometries \cite{Iesan1971,Iesan1994,hutter2016application,lurie2018bending,Lakes2015,Park1987,Taliercio2010}.
However, to the best of the authors' knowledge, no solution is known for the general micromorphic continuum with an unconstrained second-order tensor of microdeformation as kinematic degree of freedom.

The scope of the present paper is to derive and discuss the solution for cylindrical bending for such a general micromorphic continuum and several of its sub-classes, notably the relaxed micromorphic model.
A similar investigation with respect to the simple shear problem has been undertaken in \cite{rizzi2020shear}.

The solid studied here has a finite thickness $h$ along the direction $x_2$ while it is infinite in the directions $x_1$ and $x_3$.
Our aim is to describe a state of uniform cylindrical bending in the infinite plate for various models of generalized elasticity.
The bending of a plate can be thought of as a result of applying couples at the lateral faces of a suitable large portion of it.
\begin{figure}[H]
\begin{subfigure}{0.48\textwidth}
\centering
\includegraphics[width=\textwidth]{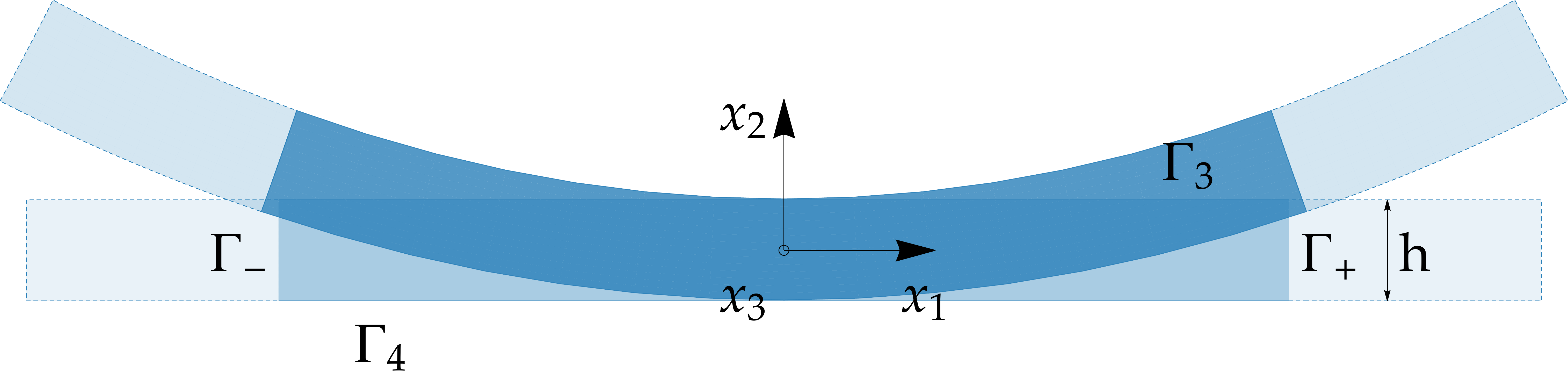}
\caption{}
\end{subfigure}
\hfill
\begin{subfigure}{0.48\textwidth}
\centering
\includegraphics[width=\textwidth]{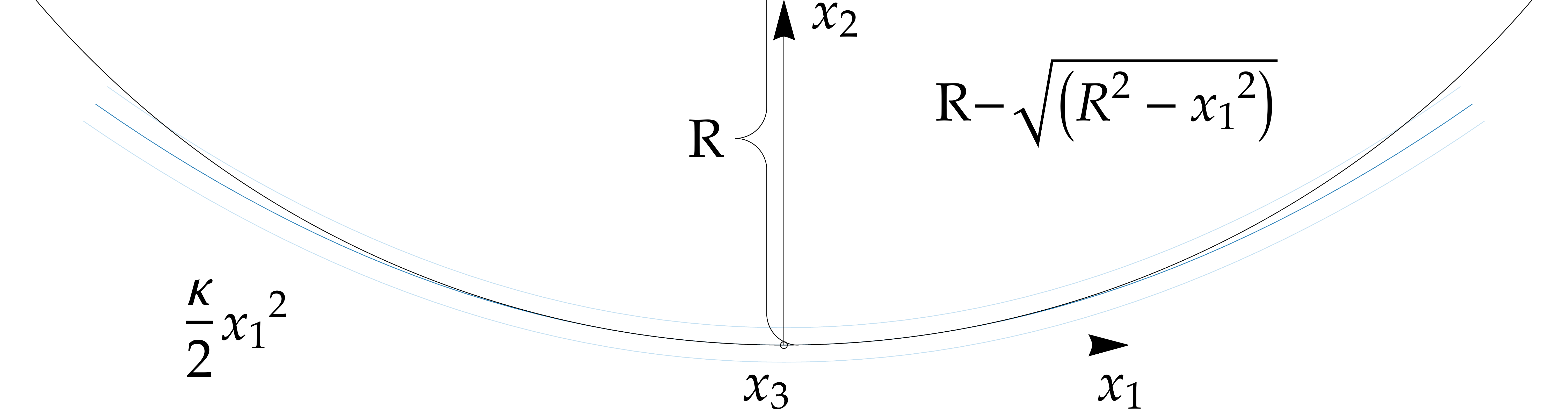}
\caption{}
\end{subfigure}
\centering
\begin{subfigure}{0.48\textwidth}
\centering
\includegraphics[width=\textwidth]{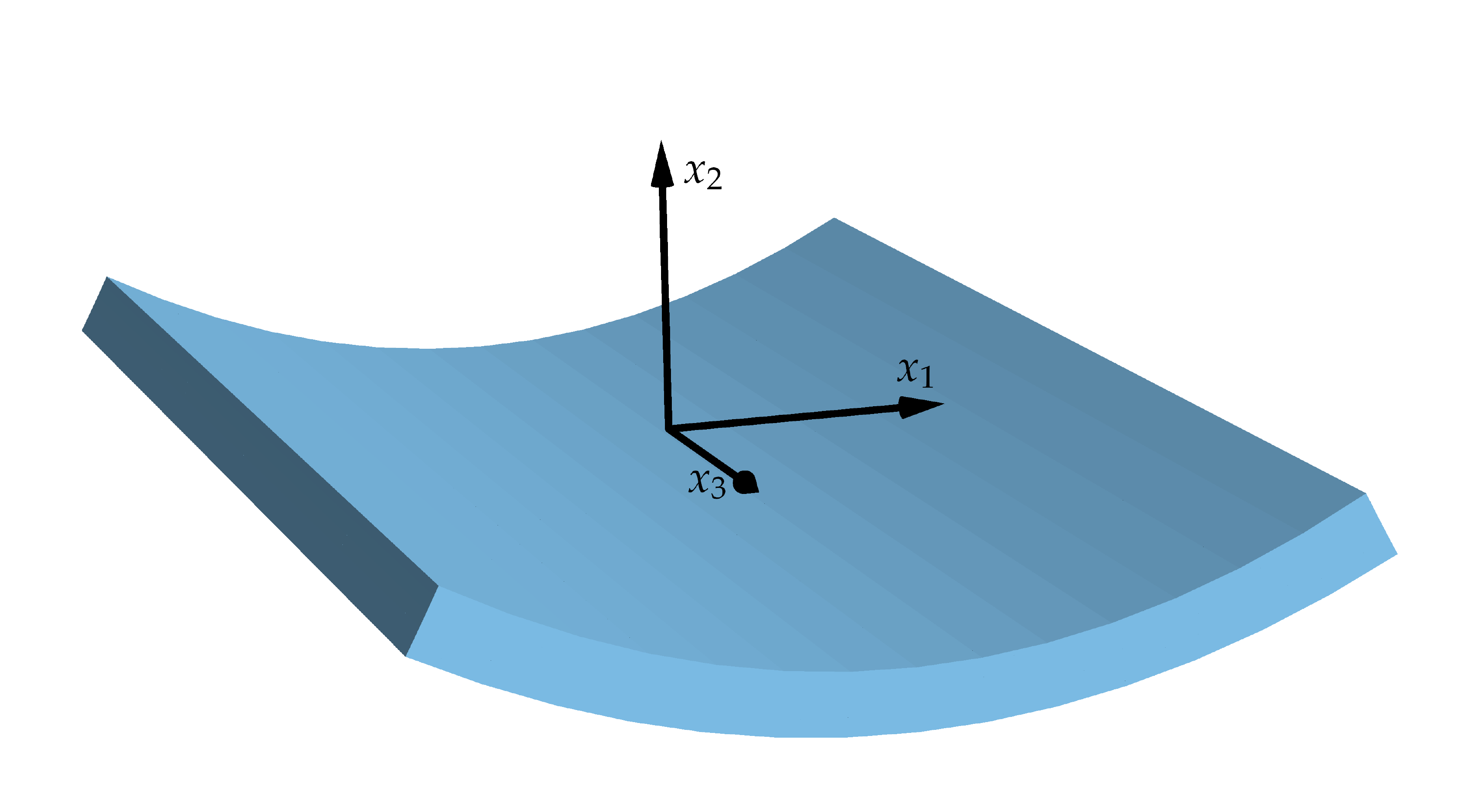}
\caption{}
\end{subfigure}
\caption{
	(a) Transversal section of the plate with the definition of the edges limiting the considered portion ($\Gamma_3$, $\Gamma_4$, $\Gamma_+$, $\Gamma_+$).
	The equation of the bottom half circumference of radius $R$ tangent to the $x_1$ axis in the origin is $x_2 = R - \sqrt{R^2 - x_1^2}$. The second order approximation of the series expansion at $x_1 = 0$ is $x_2\lvert_{x_1 = 0} = R - (R - \frac{x_1^2}{2R} + \dots) \approx \frac{\boldsymbol{\kappa} \, x_1^2}{2}$ where $\boldsymbol{\kappa}=\frac{1}{R}$ is the curvature.
	(b) Mean axis of the transversal section (a).
	(c) Schematic representation of the considered plate's portion.
}
\label{fig:intro}
\end{figure}
\subsection{Notation}
For vectors $a,b\in\mathbb{R}^n$, we consider the scalar product  $\langle \boldsymbol{a},\boldsymbol{b} \rangle \coloneqq \sum_{i=1}^n a_i\,b_i \in \mathbb{R}$, the (squared) norm  $\norm{\boldsymbol{a}}^2\coloneqq\langle \boldsymbol{a},\boldsymbol{a} \rangle$ and  the dyadic product  $\boldsymbol{a}\otimes \boldsymbol{b} \coloneqq \left(a_i\,b_j\right)_{i,j=1,\ldots,n}\in \mathbb{R}^{n\times n}$. Similarly, for tensors  $\boldsymbol{P},\boldsymbol{Q}\in\mathbb{R}^{n\times n}$  with Cartesian coordinates $P_{ij}$ and $Q_{ij}$, we define the scalar product $\langle \boldsymbol{P},\boldsymbol{Q} \rangle \coloneqq\sum_{i,j=1}^n P_{ij}\,Q_{ij} \in \mathbb{R}$ and the (squared) Frobenius-norm $\norm{\boldsymbol{P}}^2\coloneqq\langle \boldsymbol{P},\boldsymbol{P} \rangle$.
Moreover, $\boldsymbol{P}^T\coloneqq (P_{ji})_{i,j=1,\ldots,n}$ denotes the transposition of the matrix $\boldsymbol{P}=(P_{ij})_{i,j=1,\ldots,n}$, which decomposes orthogonally into the symmetric part $\mbox{sym} \boldsymbol{P} \coloneqq \frac{1}{2} (\boldsymbol{P}+\boldsymbol{P}^T)$ and the skew-symmetric part $\mbox{skew} \boldsymbol{P} \coloneqq \frac{1}{2} (\boldsymbol{P}-\boldsymbol{P}^T )$.
The Lie-Algebra of skew-symmetric matrices is denoted by $\mathfrak{so}(3)\coloneqq \{\boldsymbol{A}\in\mathbb{R}^{3\times 3}\mid \boldsymbol{A}^T = -\boldsymbol{A}\}$.
The identity matrix is denoted by $\boldsymbol{\mathbbm{1}}$, so that the trace of a matrix $\boldsymbol{P}$ is given by \ $\tr \boldsymbol{P} \coloneqq \langle \boldsymbol{P},\boldsymbol{\mathbbm{1}} \rangle$.
Using the bijection $\mbox{axl}:\mathfrak{so}(3)\to\mathbb{R}^3$ we have
\begin{align}
\boldsymbol{A} \, \boldsymbol{b} =\mbox{axl}(\boldsymbol{A})\times \boldsymbol{b} \quad \forall\, \boldsymbol{A}\in\mathfrak{so}(3) \, ,
\quad
\boldsymbol{b}\in\mathbb{R}^3.
\label{eq:Aanti_axl}
\end{align}
where $\times$ denotes the standard cross product in $\mathbb{R}^3$.
The inverse is denoted by Anti: $\mathbb{R}^3\to \mathfrak{so}(3)$.
The gradient and the curl for a vector field $\boldsymbol{u}$ are defined as
\begin{equation}
\boldsymbol{\mbox{D}u}
=\!
\left(
\begin{array}{ccc}
u_{1,1} & u_{1,2} & u_{1,3} \\
u_{2,1} & u_{2,2} & u_{2,3} \\
u_{3,1} & u_{3,2} & u_{3,3}
\end{array}
\right)\, ,
\qquad
\mbox{curl} \, \boldsymbol{u} = \boldsymbol{\nabla} \times \boldsymbol{u}
=
\left(
\begin{array}{ccc}
u_{3,2} - u_{2,3}  \\
u_{1,3} - u_{3,1}  \\
u_{2,1} - u_{1,2}
\end{array}
\right)
\, .
\end{equation}
Moreover, we introduce the $\mbox{Curl}$ and the $\mbox{Div}$ operators of the matrix $\boldsymbol{P}$ as
\begin{equation}
\mbox{Curl} \, \boldsymbol{P}
=\!
\left(
\begin{array}{c}
(\mbox{curl}\left( P_{11} , \right. P_{12} , \left. P_{13} \right))^T \\
(\mbox{curl}\left( P_{21} , \right. P_{22} , \left. P_{23} \right))^T \\
(\mbox{curl}\left( P_{31} , \right. P_{32} , \left. P_{33} \right))^T
\end{array}
\right) \!,
\qquad
\mbox{Div}  \, \boldsymbol{P}
=\!
\left(
\begin{array}{c}
\mbox{div}\left( P_{11} , \right. P_{12} , \left. P_{13} \right) \\
\mbox{div}\left( P_{21} , \right. P_{22} , \left. P_{23} \right) \\
\mbox{div}\left( P_{31} , \right. P_{32} , \left. P_{33} \right) 
\end{array}
\right) \, .
\end{equation}
The cross product between a second order tensor and a vector is defined as follow
\begin{equation}
\boldsymbol{m} \times \boldsymbol{b} =
\left(
\begin{array}{ccc}
(b \times (m_{11},m_{12},m_{13}))^T \\
(b \times (m_{21},m_{22},m_{23}))^T \\
(b \times (m_{31},m_{32},m_{33}))^T \\
\end{array}
\right) =
\boldsymbol{m} \cdot \boldsymbol{\epsilon} \cdot \boldsymbol{b} =
m_{ik} \, \epsilon_{kjh} \, b_{h} 
\, ,
\end{equation}
where $\boldsymbol{m} \in \mathbb{R}^{3\times 3}$, $\boldsymbol{b} \in \mathbb{R}^{3}$, and $\boldsymbol{\epsilon}$ is the Levi-Civita tensor.
\section{Cylindrical bending for the isotropic Cauchy continuum}
For comparison we start with the well known classical case.
The expression of the strain energy for an isotropic linear elastic Cauchy continuum is\!\!
\footnote{
Here are reported the macroscopic 3D Poisson's ratio $\nu_{\tiny \mbox{macro}} = \frac{\lambda_{\mbox{\tiny macro}}}{2\left(\lambda_{\mbox{\tiny macro}} + \mu_{\mbox{\tiny macro}}\right)}$, the Young modulus \\$E_{\tiny \mbox{macro}} = \frac{\mu_{\mbox{\tiny macro}} \left(3\lambda_{\mbox{\tiny macro}} + 2 \mu_{\mbox{\tiny macro}}\right)}{\lambda_{\mbox{\tiny macro}} + \mu_{\mbox{\tiny macro}}} = 2\mu_{\mbox{\tiny macro}}\left(1+\nu_{\mbox{\tiny macro}}\right)$, and the bulk modulus $\kappa_{\mbox{\tiny macro}} = \frac{2\mu_{\mbox{\tiny macro}} + 3\lambda_{\mbox{\tiny macro}}}{3}$.
}
\begin{equation}
W \left(\boldsymbol{\mbox{D}u}\right) = 
\mu_{\mbox{\tiny macro}} \left\lVert \mbox{sym} \, \boldsymbol{\mbox{D}u} \right\rVert^{2} + 
\frac{\lambda_{\mbox{\tiny macro}}}{2} \mbox{tr}^2\left(\boldsymbol{\mbox{D}u}\right)
= 
\mu_{\mbox{\tiny macro}} \left\lVert \mbox{dev} \, \mbox{sym} \, \boldsymbol{\mbox{D}u} \right\rVert^{2} + 
\frac{\kappa_{\mbox{\tiny macro}}}{2} \mbox{tr}^2\left(\boldsymbol{\mbox{D}u}\right)
\, ,
\label{eq:energy_Cau}
\end{equation}
where $\mu_{\mbox{\tiny macro}}$ and $\lambda_{\mbox{\tiny macro}}$ are the two classical Lamé constants, and $\kappa_{\mbox{\tiny macro}}$ is the bulk modulus.
Consequently, the equilibrium equations without body forces are
\begin{equation}
\mbox{Div}\, 
\boldsymbol{\sigma} 
= \boldsymbol{0} \, ,
\qquad
\boldsymbol{\sigma} = 
2\,\mu_{\mbox{\tiny macro}}\,\mbox{sym} \, \boldsymbol{\mbox{D}u}
+ \lambda_{\mbox{\tiny macro}}\,\mbox{tr}\left(\boldsymbol{\mbox{D}u}\right)  \boldsymbol{\mathbbm{1}}= 
2\,\mu_{\mbox{\tiny macro}}\,\mbox{dev}\,\mbox{sym} \, \boldsymbol{\mbox{D}u}
+ \kappa_{\mbox{\tiny macro}}\,\mbox{tr}\left(\boldsymbol{\mbox{D}u}\right)  \boldsymbol{\mathbbm{1}} \, ,
\label{eq:equi_Cau}
\end{equation}
where $\boldsymbol{\sigma}$ is the symmetric Cauchy-stress tensor and $\boldsymbol{\epsilon} = \mbox{sym} \, \boldsymbol{\mbox{D}u}$ is the symmetric strain tensor.
The boundary conditions at the upper and lower surface (free surface) are the traction-free forces conditions
\begin{equation}
\boldsymbol{t}(x_2 = \pm \, h/2) = \pm \, \boldsymbol{\sigma}(x_2 = \pm \, h/2) \cdot \boldsymbol{e}_2 = \boldsymbol{0} \, ,
\label{eq:BC_Cau}
\end{equation}
where the expression of $\boldsymbol{\sigma}$ is in eq.(\ref{eq:equi_Cau})$_2$ and $\boldsymbol{e}_2$ is the unit vector aligned to the $x_2$-direction.
We are interested in describing a \textbf{state of uniform curvature} $\boldsymbol{\kappa}$ of the infinite plate.
According to the reference system shown in Fig.~\ref{fig:intro}, the ansatz for the displacement field is
\begin{equation}
\boldsymbol{u}(x_1,x_2)=
\left(
\begin{array}{c}
-x_2 \, w'(x_1) \\
w(x_1) + v(x_2) \\
0 \\
\end{array}
\right).
\label{eq:ansatz_Cau}
\end{equation}
Therein, $w(x_1)$ refers to the deflection of the neutral axis $x_2=0$ and $v(x_2)$ denotes the lateral contraction.
The gradient of the displacement field and its symmetric part (the strain tensor) result to be
\begin{equation}
\boldsymbol{\mbox{D}u} = 
\left(
\begin{array}{ccc}
-x_2 \, w''(x_1) & -w'(x_1) & 0 \\
w'(x_1) & v'(x_2) & 0 \\
0 & 0 & 0 \\
\end{array}
\right),
\qquad
\boldsymbol{\varepsilon} = \mbox{sym} \, \boldsymbol{\mbox{D}u} = 
\left(
\begin{array}{ccc}
-x_2 \, w''(x_1) & 0 & 0 \\
0 & v'(x_2) & 0 \\
0 & 0 & 0 \\
\end{array}
\right) \, .
\label{eq:grad_Cau}
\end{equation}
It is highlighted that the gradient of eq.(\ref{eq:grad_Cau})$_1$ has a symmetric part which has only diagonal components as it can be seen in eq.(\ref{eq:grad_Cau})$_2$.
The equilibrium equations in terms of the ansatz eq.(\ref{eq:ansatz_Cau}) are
\begin{equation}
-x_2 \, w^{(3)}(x_1) \left(\lambda_{\mbox{\tiny macro}} +2 \mu_{\mbox{\tiny macro}} \right) = 0 \, ,
\qquad\qquad
\left(\lambda_{\mbox{\tiny macro}} +2 \mu_{\mbox{\tiny macro}} \right) v''(x_2)-\lambda_{\mbox{\tiny macro}}  w''(x_1) = 0 \, .
\label{eq:equi_equa_Cau}
\end{equation}
Equation.(\ref{eq:equi_equa_Cau})$_2$ requires both $w(x_1)$ and $v(x_2)$ to be a quadratic function in $x_1$ and $x_2$, respectively, and this already satisfies eq.(\ref{eq:equi_equa_Cau})$_1$:
\begin{equation}
w(x_1) = c_3 \, x_1^2 + c_2 \, x_1 + c_1 \, ,
\qquad\qquad
v(x_2) = \frac{c_3 \, \lambda_{\mbox{\tiny macro}} }{\lambda_{\mbox{\tiny macro}}  + 2 \mu_{\mbox{\tiny macro}} } \, x_2^2 + c_5 \, x_2 + c_4 \, .
\label{eq:sol_fun_disp_Cau}
\end{equation}
It is important to highlight that the solution eq.(\ref{eq:sol_fun_disp_Cau}) depends on the elastic parameters of the material.
In addition, we can further simplify the solution eq.(\ref{eq:sol_fun_disp_Cau}) by setting $c_1 = 0$, $c_4 = 0$ and $c_2 = 0$ since they represent rigid body motions (the first two a rigid translation, while the third one a rigid rotation).

The boundary conditions eq.(\ref{eq:BC_Cau}) on the upper and lower surface (stress free surfaces, no tractions), require that $c_5 = 0$. 
Appling these simplifications while substituting the solution eq.(\ref{eq:sol_fun_disp_Cau}) in eq.(\ref{eq:ansatz_Cau}), the  displacement field results to be (with $c_3 = \boldsymbol{\kappa}/2$) \cite{hadjesfandiari2016pure}
\begin{equation}
\boldsymbol{u}(x_1,x_2)=
\frac{\boldsymbol{\kappa}}{2}
\left(
\begin{array}{c}
-2 \, x_1 \, x_2 \\
\frac{\lambda_{\mbox{\tiny macro}}}{\lambda_{\mbox{\tiny macro}}+2 \mu_{\mbox{\tiny macro}}} \, x_2^2 + x_1^2 \\
\end{array}
\right)
=
\frac{1}{2R}
\left(
\begin{array}{c}
- 2\, x_1 \, x_2 \\
\frac{\nu_{\mbox{\tiny macro}}}{1-\nu_{\mbox{\tiny macro}}} \, x_2^2 + x_1^2 \\
\end{array}
\right)
.
\label{eq:sol_disp_Cau}
\end{equation}
Thus, the middle plane of the infinite plate is bent to a cylindrical surface (see Fig.~\ref{fig:intro}) which is approximated by the parabolic cylinder $u_2(x_1,x_2=0) = \frac{1}{2R}x_1^2$.
Since the displacement field solution depends on the elastic coefficients of the material, two extreme cases are shown in the following Fig.~\ref{fig:deformed_shape_Cau} with a plot of $u_{1,1}$ across the thickness $h$.
\begin{figure}[H]
\centering
\begin{subfigure}{.45\textwidth}
\centering
\includegraphics[width=0.75\linewidth]{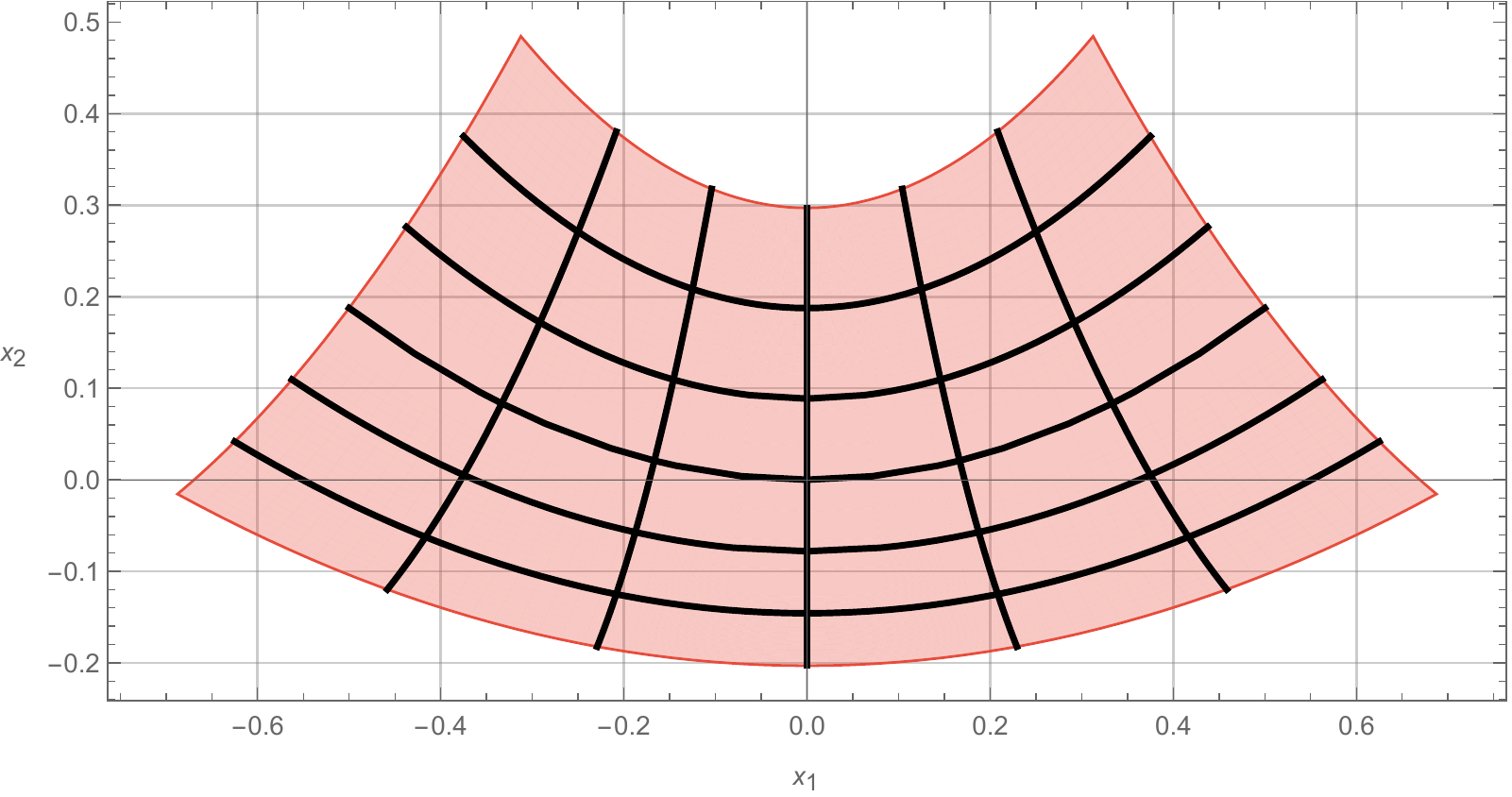}
\caption{}
\label{fig:def_Cau_1}
\end{subfigure}%
\begin{subfigure}{.45\textwidth}
\centering
\includegraphics[width=0.75\linewidth]{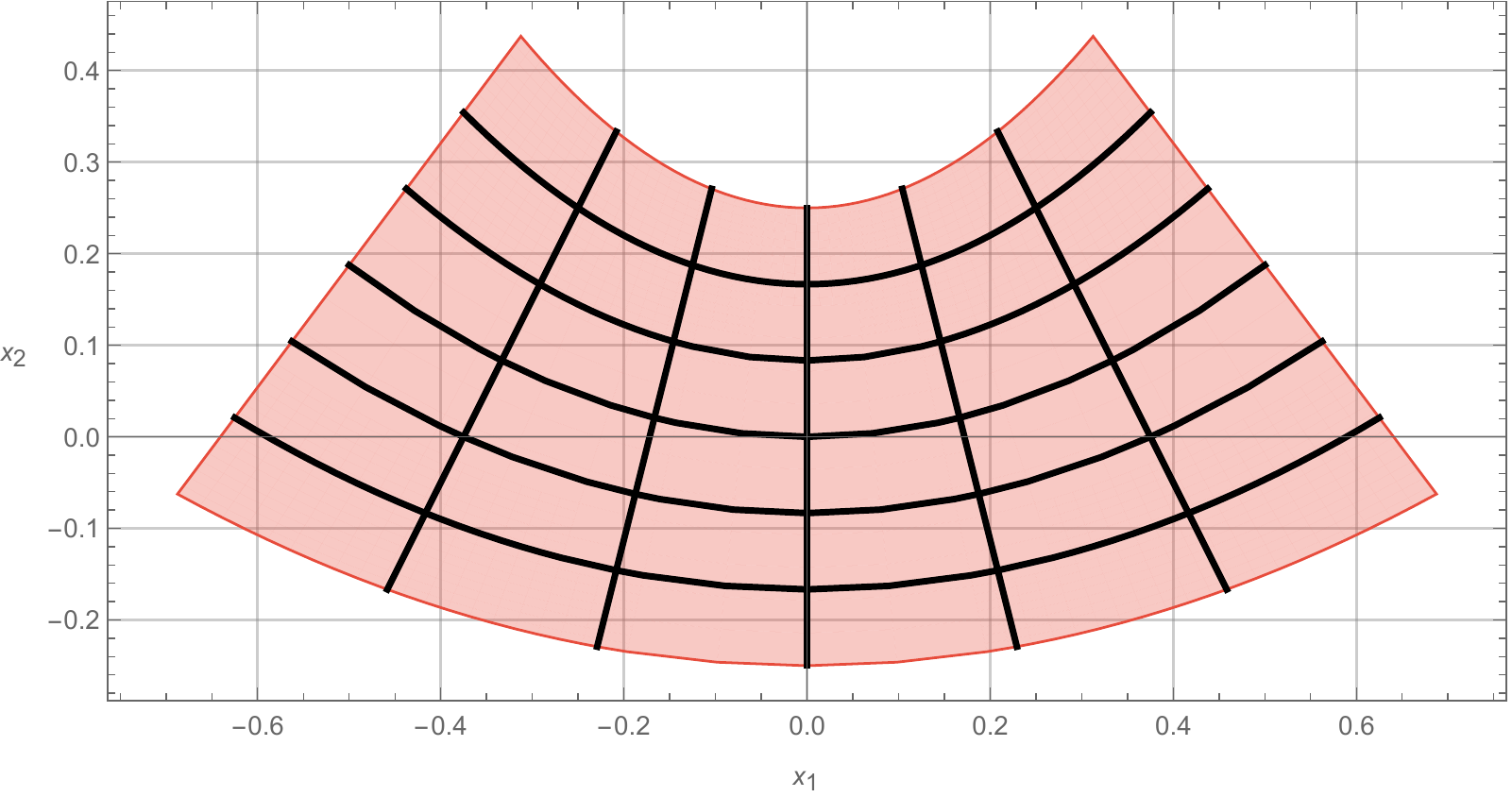}
\caption{}
\label{fig:def_Cau_2}
\end{subfigure}
\begin{subfigure}{.45\textwidth}
\centering
\includegraphics[width=0.75\linewidth]{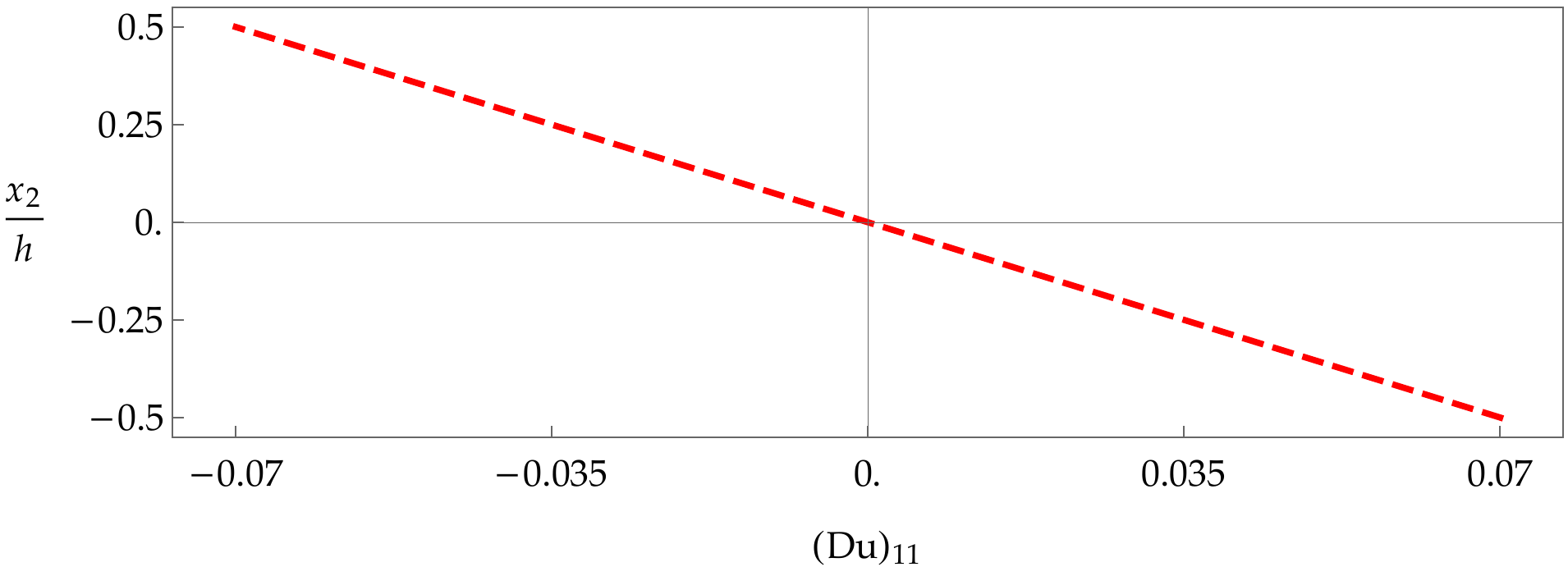}
\caption{}
\label{fig:u11_Cau}
\end{subfigure}
\caption{
	(a) Deformed shape for Poisson's ratio $\nu_{\tiny \mbox{macro}} \to \nicefrac{1}{2}$, and
	(b) deformed shape for Poisson's ratio $\nu_{\tiny \mbox{macro}} \to 0 $;
	(c) plot of the component $\left(\boldsymbol{\mbox{D}u}\right)_{11}$ of the gradient of the displacement across the thickness $h$. The value used for the curvature is $\boldsymbol{\kappa} = 1.5$.}
\label{fig:deformed_shape_Cau}
\end{figure}
The desired bending moment about the $x_3$-axis and energy (per unit area d$x_1$d$x_3$) expressions are
\begin{align}
M_{\mbox{c}} (\boldsymbol{\kappa}) :=&
\displaystyle\int\limits_{-h/2}^{h/2} \langle \boldsymbol{\sigma} \, \boldsymbol{e}_1 , \boldsymbol{e}_1 \rangle \, x_2 \, \mbox{d}x_{2}
=
\frac{2 \, c_3 \, h^3 \mu_{\mbox{\tiny macro}}  (\lambda_{\mbox{\tiny macro}} +\mu_{\mbox{\tiny macro}} )}{3 (\lambda_{\mbox{\tiny macro}} +2 \mu_{\mbox{\tiny macro}} )}
\notag
\\
=
&
\left(\widehat{\lambda}_{\mbox{\tiny macro}} + 2\mu_{\mbox{\tiny macro}}\right) \, J_{x_3} \, \boldsymbol{\kappa}
=
D_{\mbox{\tiny macro}} \, \boldsymbol{\kappa}
\, ,
\label{eq:sigm_ene_dimensionless_Cau}
\\
W_{\mbox{tot}} (\boldsymbol{\kappa}) :=&
\displaystyle\int\limits_{-h/2}^{+h/2} W(\boldsymbol{\mbox{D}u}) \, \mbox{d}x_{2} 
=
\frac{2 \, c_3^2 \, h^3 \mu_{\mbox{\tiny macro}}  (\lambda_{\mbox{\tiny macro}} +\mu_{\mbox{\tiny macro}} )}{3 (\lambda_{\mbox{\tiny macro}} +2 \mu_{\mbox{\tiny macro}} )}
\notag
\\
=
&
\frac{1}{2} \left(\widehat{\lambda}_{\mbox{\tiny macro}} + 2\mu_{\mbox{\tiny macro}}\right) \, J_{x_3} \, \boldsymbol{\kappa}^2
=
\frac{1}{2}D_{\mbox{\tiny macro}} \, \boldsymbol{\kappa}^2
\, ,
\notag
\end{align}
where $\widehat{\lambda}_{\mbox{\tiny macro}}$ is the plane stress first Lamé parameter, $J_{x_3} = \frac{h^3}{12}$ is the moment of inertia for a unit thickness, $\boldsymbol{\kappa} = 2 c_3$ is the curvature, and the quantity 
\begin{equation}
D_{\mbox{\tiny macro}} = (\widehat{\lambda}_{\mbox{\tiny macro}} + 2\mu_{\mbox{\tiny macro}})\, J_{x_3} =
\frac{h^3}{12} \frac{E_{\tiny \mbox{macro}}}{\left(1-\nu_{\tiny \mbox{macro}}^2\right)} =
\frac{h^3}{12} \frac{4 \mu _{\mbox{\tiny macro}} \left(3 \kappa _{\mbox{\tiny macro}}+\mu _{\mbox{\tiny macro}}\right)}{3 \kappa _{\mbox{\tiny macro}}+4 \mu _{\mbox{\tiny macro}}} \, ,
\label{eq:bending_stiffness}
\end{equation}
is the classical cylindrical bending stiffness for a plate (flexural rigidity).
\!\!\!\!\!\!\!
\footnote{
Note  that under the plane stress hypothesis the first Lamé parameter become
$\widehat{\lambda}_{\mbox{\tiny macro}} = \frac{2 \, \lambda_{\mbox{\tiny macro}} \, \mu_{\mbox{\tiny macro}}}{\lambda_{\mbox{\tiny macro}} + 2 \mu_{\mbox{\tiny macro}}}$
, while the shear modulus
$\widehat{\mu}_{\mbox{\tiny macro}} = \mu_{\mbox{\tiny macro}}$
, the Young modulus 
$\widehat{E} = E = \frac{\mu_{\mbox{\tiny macro}}(3\lambda_{\mbox{\tiny macro}} + 2\mu_{\mbox{\tiny macro}})}{\lambda_{\mbox{\tiny macro}} + \mu_{\mbox{\tiny macro}}}$
, and the Poisson's ratio 
$\widehat{\nu} = \mu_{\mbox{\tiny macro}} = \frac{\lambda_{\mbox{\tiny macro}}}{2\lambda_{\tiny \mbox{macro}}+2\mu_{\mbox{\tiny macro}}}$ 
do not change.
It is also reported here the more used bending stiffness expression 
$\widehat{\lambda}_{\mbox{\tiny macro}} + 2\mu_{\mbox{\tiny macro}} = \frac{E_{\mbox{\tiny macro}}}{1- \nu_{\mbox{\tiny macro}}^2}$
.
}
It is also highlighted that 
\begin{equation}
\frac{\mbox{d}}{\mbox{d}\boldsymbol{\kappa}}W_{\mbox{tot}}(\boldsymbol{\kappa}) = M_{\mbox{c}} (\boldsymbol{\kappa}) \, .
\end{equation}

Here and in the remainder of this work, the elastic coefficients $\mu_i,\lambda_i,\kappa_i$ are expressed in [MPa], the coefficients $a_i$ are dimensionless, the lengths $L_c$ and the thickness $h$ in meter [m], the curvature $\boldsymbol{\kappa}$ in [1/m].

\section{Classical Mindlin-Eringen formulation}
The classical micromorphic model couples the displacement $\boldsymbol{u} \in \mathbb{R}^{3}$ with an independent affine field $\boldsymbol{P} \in \mathbb{R}^{3\times3}$, called the microdistortion. Even in the isotropic case the model has 18 parameters (among them 11 for the curvature) and the elastic energy can be represented as 
\begin{align}
W \left(\boldsymbol{\mbox{D}u}, \boldsymbol{P}, \boldsymbol{\mbox{D}P}\right) =
& \,
\widehat{\mu} \, \left \lVert \mbox{sym} \, \boldsymbol{\mbox{D}u} \right \rVert ^2
+ \frac{\widehat{\lambda}}{2} \, \mbox{tr}^2 \left(\boldsymbol{\mbox{D}u}\right)
+ \frac{b_1}{2} \, \mbox{tr}^2 \left(\boldsymbol{\mbox{D}u} - \boldsymbol{P}\right)
\notag
\\*
&
+ \frac{b_2}{2} \, \left \lVert \boldsymbol{\mbox{D}u} - \boldsymbol{P} \right \rVert ^2
+ \frac{b_3}{2} \, \langle \boldsymbol{\mbox{D}u} - \boldsymbol{P},\left(\boldsymbol{\mbox{D}u} - \boldsymbol{P}\right)^{T} \rangle
+ g_1 \, \mbox{tr} \left(\boldsymbol{\mbox{D}u}\right)\mbox{tr} \left(\boldsymbol{\mbox{D}u} - \boldsymbol{P}\right)
\notag
\\*
&
+ g_2 \, \langle \mbox{sym} \, \boldsymbol{\mbox{D}u},\left(\boldsymbol{\mbox{D}u} - \boldsymbol{P}\right)^{T} \rangle
+ \frac{1}{2} \langle \mathbb{A}\boldsymbol{\mbox{D}P},\boldsymbol{\mbox{D}P} \rangle
\notag
\\*
=
&
\, \widehat{\mu} \, \varepsilon_{ij} \, \varepsilon_{ij} 
+ \frac{\widehat{\lambda}}{2} \, \varepsilon_{ii} \, \varepsilon_{jj}
+ \frac{b_1}{2} \, \gamma_{ii} \, \gamma_{jj} 
+ \frac{b_2}{2} \, \gamma_{ij} \, \gamma_{ij} 
+ \frac{b_3}{2} \, \gamma_{ij} \, \gamma_{ji}
\label{eq:energy_MM_Mind}
\\*
&
+ g_1 \, \gamma_{ii} \, \varepsilon_{jj} 
+ g_2 \, \left(\gamma_{ij} + \gamma_{ji} \, \right)\varepsilon_{ij}
+ \widehat{a}_1 \, \chi_{iik} \, \chi_{kjj} 
+ \widehat{a}_2 \, \chi_{iik} \, \chi_{jkj}
\notag
\\*
&
+ \frac{1}{2} \, \widehat{a}_3 \, \chi_{iik} \, \chi_{jjk} 
+ \frac{1}{2} \, \widehat{a}_4 \, \chi_{ijj} \, \chi_{ikk} 
+ \widehat{a}_5 \, \chi_{ijj} \, \chi_{kik}
+ \frac{1}{2} \, \widehat{a}_8 \, \chi_{iji} \, \chi_{kjk}
+ \frac{1}{2} \, \widehat{a}_{10} \, \chi_{ijk} \, \chi_{ijk} 
\notag
\\*
&
+ \widehat{a}_{11} \, \chi_{ijk} \, \chi_{jki}
+ \frac{1}{2} \, \widehat{a}_{13} \, \chi_{ijk} \, \chi_{ikj} 
+ \frac{1}{2} \, \widehat{a}_{14} \, \chi_{ijk} \, \chi_{jik}
+ \frac{1}{2} \, \widehat{a}_{15} \, \chi_{ijk} \, \chi_{kji}
\notag
\end{align}
where $\boldsymbol{\varepsilon} = \mbox{sym} \, \boldsymbol{\mbox{D}u}$ is the symmetric part of the gradient of the displacement field, $\boldsymbol{\gamma} = \boldsymbol{\mbox{D}u} - \boldsymbol{P}$ is the difference between the gradient of the displacement field and the microdistortion tensor, and $\chi_{ijk} = P_{jk,i}$ is the full gradient of the microdistortion and $\mathbb{A}$ is a sixth order tensor.

A general solution of the cylindrical bending problem for the isotropic micromorphic model is not yet known but it could be found by the present ansatz without any fundamental problem, only the formulas get extremely long and complicated.
Therefore, we consider the following simplified energy
\begin{align}
W \left(\boldsymbol{\mbox{D}u}, \boldsymbol{P}, \boldsymbol{\mbox{D}P}\right) = & \,
\frac{\mu _e + \mu_{\mbox{\tiny micro}} + \mu _c}{2} \, \left \lVert \boldsymbol{\mbox{D}u} -\boldsymbol{P} \right \rVert ^2
+ \frac{\mu _e + \mu_{\mbox{\tiny micro}} - \mu _c}{2} \, \langle \boldsymbol{\mbox{D}u} - \boldsymbol{P},\left(\boldsymbol{\mbox{D}u} - \boldsymbol{P}\right)^{T} \rangle 
\notag
\\
&
+ \frac{\lambda _e + \lambda_{\mbox{\tiny micro}}}{2} \, \mbox{tr}^2 \left(\boldsymbol{\mbox{D}u} - \boldsymbol{P}\right) 
+ \mu_{\mbox{\tiny micro}} \, \left \lVert \mbox{sym} \, \boldsymbol{\mbox{D}u} \right \rVert ^2 
+ \frac{\lambda_{\mbox{\tiny micro}}}{2} \, \mbox{tr}^2 \left(\boldsymbol{\mbox{D}u}\right) 
\label{eq:energy_MM_Mind_neff}
\\
&
-2\mu_{\mbox{\tiny micro}}  \, \langle \boldsymbol{\mbox{D}u} - \boldsymbol{P},\mbox{sym} \, \boldsymbol{\mbox{D}u} \rangle
-\lambda_{\mbox{\tiny micro}} \, \mbox{tr} \left(\boldsymbol{\mbox{D}u} - \boldsymbol{P}\right)\mbox{tr} \left(\boldsymbol{\mbox{D}u}\right)
\notag
\\
&
+ \frac{\mu \, L_c^2}{2}
\displaystyle\sum_{i=1}^{3}
\left(
a_1 \, \left\lVert \mbox{dev} \, \mbox{sym} \, \partial_{x_i} \boldsymbol{P} \right\rVert^2
+ a_2 \, \left\lVert \mbox{skew} \, \partial_{x_i} \boldsymbol{P} \right\rVert^2
+ \frac{2}{9} \, a_3 \, \mbox{tr}^2 \Big( \partial_{x_i} \boldsymbol{P} \Big)
\right)
\notag
\end{align}
which is a special case of eq.~(\ref{eq:energy_MM_Mind}), by setting the values of the elastic parameters as follows (see \cite{neff2004material} and Appendix~\ref{app:coeff} for the full derivation for the curvature part)
\begin{align}
\widehat{\mu} &= \mu_{\mbox{\tiny micro}} \, ,
\quad
\widehat{\lambda} = \lambda_{\mbox{\tiny micro}} \, ,
\quad
b_1 = \lambda_e + \lambda_{\mbox{\tiny micro}} \, ,
\quad
b_2 = \mu _e + \mu_{\mbox{\tiny micro}} + \mu _c \, ,
\quad
b_3 = \mu _e + \mu_{\mbox{\tiny micro}} - \mu _c \, ,
\notag
\\
g_1 &= -\lambda_{\tiny \mbox{micro}} \, ,
\quad
g_2 = -2\mu_{\tiny \mbox{micro}} \, ,
\quad
\widehat{a}_{1,2,3,5,8,11,14,15} = 0 \, ,
\quad
\widehat{a}_{4} = \mu \, L_c^2 \frac{2a_3-a_1}{3} \, ,
\quad
\widehat{a}_{10} = \mu \, L_c^2 \frac{a_1 + a_2}{2} \, ,
\label{eq:energy_MM_Mind_neff_coeff}
\\
&
\hspace{6.6cm}
\widehat{a}_{13} = \mu \, L_c^2 \frac{a_1 - a_2}{2}
\, .
\notag
\end{align}
Thus, with this choice of parameters, the former energy can be expressed as
\begin{align}
W \left(\boldsymbol{\mbox{D}u}, \boldsymbol{P}, \boldsymbol{\mbox{D}P}\right)
= &
\, \mu_{e} \left\lVert \mbox{dev} \, \mbox{sym} \left(\boldsymbol{\mbox{D}u} - \boldsymbol{P} \right) \right\rVert^{2}
+ \frac{\kappa_{e}}{2} \mbox{tr}^2 \left(\boldsymbol{\mbox{D}u} - \boldsymbol{P} \right)
+ \mu_{c} \left\lVert \mbox{skew} \left(\boldsymbol{\mbox{D}u} - \boldsymbol{P} \right) \right\rVert^{2}
\notag
\\
&
+ \mu_{\tiny \mbox{micro}} \left\lVert \mbox{dev} \, \mbox{sym}\,\boldsymbol{P} \right\rVert^{2}
+ \frac{\kappa_{\tiny \mbox{micro}}}{2} \mbox{tr}^2 \left(\boldsymbol{P} \right)
\label{eq:energy_MM_Mind_neff_coeff_2}
\\
&
+ \frac{\mu \, L_c^2}{2}
\displaystyle\sum_{i=1}^{3}
\left(
a_1 \, \left\lVert \mbox{dev} \, \mbox{sym} \, \partial_{x_i} \boldsymbol{P} \right\rVert^2
+ a_2 \, \left\lVert \mbox{skew} \, \partial_{x_i} \boldsymbol{P} \right\rVert^2
+ \frac{2}{9} \, a_3 \, \mbox{tr}^2 \Big( \partial_{x_i} \boldsymbol{P} \Big)
\right)
\, ,
\notag
\end{align}
which coincides, apart from the curvature term, with the energy from the relaxed micromorphic model treated in the next paragraph.
Note again carefully that the chosen curvature expression is still isotropic \cite{munch2018rotational} but does not represent the most general isotropic curvature expression in this model.

\section{Cylindrical bending for the isotropic relaxed micromorphic model}
The expression of the strain energy for the isotropic relaxed micromorphic continuum is:\!\!
\footnote{
Are here reported the 3D Poisson's ratio $\nu_{\tiny \mbox{macro}} = \frac{\lambda_{\mbox{\tiny macro}}}{2\left(\lambda_{\mbox{\tiny macro}} + \mu_{\mbox{\tiny macro}}\right)}$, the 3D Young modulus \\$E_{\tiny \mbox{macro}} = \frac{\mu_{\mbox{\tiny macro}} \left(3\lambda_{\mbox{\tiny macro}} + 2 \mu_{\mbox{\tiny macro}}\right)}{\lambda_{\mbox{\tiny macro}} + \mu_{\mbox{\tiny macro}}}$, and the micro and the meso expression of the Poisson's ratio in plane stress $\nu_{\tiny \mbox{micro}} = \frac{\lambda_{\mbox{\tiny micro}}}{2\left(\lambda_{\mbox{\tiny micro}} + \mu_{\mbox{\tiny micro}}\right)}$ and the $\nu_e = \frac{\lambda_e}{2\left(\lambda_e + \mu_e\right)}$, respectively.
}
\begin{align}
W \left(\boldsymbol{\mbox{D}u}, \boldsymbol{P},\mbox{Curl}\,\boldsymbol{P}\right) =
&
\, \mu_{e} \left\lVert \mbox{sym} \left(\boldsymbol{\mbox{D}u} - \boldsymbol{P} \right) \right\rVert^{2}
+ \frac{\lambda_{e}}{2} \mbox{tr}^2 \left(\boldsymbol{\mbox{D}u} - \boldsymbol{P} \right) 
+ \mu_{c} \left\lVert \mbox{skew} \left(\boldsymbol{\mbox{D}u} - \boldsymbol{P} \right) \right\rVert^{2}
\notag
\\*
&
+ \mu_{\tiny \mbox{micro}} \left\lVert \mbox{sym}\,\boldsymbol{P} \right\rVert^{2}
+ \frac{\lambda_{\tiny \mbox{micro}}}{2} \mbox{tr}^2 \left(\boldsymbol{P} \right)
\notag
\\*
&
+ \frac{\mu \,L_c^2 }{2} \,
\left(
a_1 \, \left\lVert \mbox{dev sym} \, \mbox{Curl} \, \boldsymbol{P}\right\rVert^2 +
a_2 \, \left\lVert \mbox{skew} \,  \mbox{Curl} \, \boldsymbol{P}\right\rVert^2 +
\frac{a_3}{3} \, \mbox{tr}^2 \left(\mbox{Curl} \, \boldsymbol{P}\right)
\right)
\notag
\\*
=&
\, \mu_{e} \left\lVert \mbox{dev} \, \mbox{sym} \left(\boldsymbol{\mbox{D}u} - \boldsymbol{P} \right) \right\rVert^{2}
+ \frac{\kappa_{e}}{2} \mbox{tr}^2 \left(\boldsymbol{\mbox{D}u} - \boldsymbol{P} \right) 
+ \mu_{c} \left\lVert \mbox{skew} \left(\boldsymbol{\mbox{D}u} - \boldsymbol{P} \right) \right\rVert^{2}
\label{eq:energy_RM}
\\*
&
+ \mu_{\tiny \mbox{micro}} \left\lVert \mbox{dev} \, \mbox{sym}\,\boldsymbol{P} \right\rVert^{2}
+ \frac{\kappa_{\tiny \mbox{micro}}}{2} \mbox{tr}^2 \left(\boldsymbol{P} \right)
\notag
\\*
&
+ \frac{\mu \,L_c^2 }{2} \,
\left(
a_1 \, \left\lVert \mbox{dev sym} \, \mbox{Curl} \, \boldsymbol{P}\right\rVert^2 +
a_2 \, \left\lVert \mbox{skew} \,  \mbox{Curl} \, \boldsymbol{P}\right\rVert^2 +
\frac{a_3}{3} \, \mbox{tr}^2 \left(\mbox{Curl} \, \boldsymbol{P}\right)
\right)
\, ,
\notag
\end{align}
where $\mu_e$ and $\lambda_e$ are the material parameters related to the meso-scale, $\mu_{\mbox{\tiny micro}}$ and $\lambda_{\mbox{\tiny micro}}$ are the parameters related to the micro-scale, $\mu_c$ is the Cosserat couple modulus, $L_c > 0$ is the characteristic length, and $a_1$, $a_2$, and $a_3$ are the three general isotropic curvature parameters.
The curvature expression is the most general isotropic one in terms of a dependence on the second order tensor $\mbox{Curl} \boldsymbol{P}$ and can be obtained by setting the values of the elastic parameters of the classical micromorphic model eq.(\ref{eq:energy_MM_Mind}) as follows (see Appendix~\ref{app:coeff} for the full derivation)
\begin{align}
	\widehat{\mu} &= \mu_{\mbox{\tiny micro}} \, , \quad \widehat{\lambda} = \lambda_{\mbox{\tiny micro}} \, , \quad b_1 = \lambda_e + \lambda_{\mbox{\tiny micro}} \, ,
	\quad
	b_2 = \mu _e + \mu_{\mbox{\tiny micro}} + \mu _c \, , \quad b_3 = \mu _e + \mu_{\mbox{\tiny micro}} - \mu _c \, , 
	\notag
	\\
	g_1 &= -\lambda_{\tiny \mbox{micro}} \, , \quad g_2 = -2\mu_{\tiny \mbox{micro}} \, ,
	\quad \widehat{a}_{1} = \mu \, L_c^2 \, \frac{2a_1 - a_2}{4} \, ,
	\quad \widehat{a}_{2,5,8} = 0 \, ,
	\quad \widehat{a}_{3} = \widehat{a}_{4} = \mu \, L_c^2 \, \frac{a_2 - a_1}{2} \, ,
	\label{eq:energy_MM_Mind_neff_coeff_b}
	\\
	\widehat{a}_{10} &= - \widehat{a}_{15} = \mu \, L_c^2 \, \frac{2a_1 + a_3}{3} \, ,
	\quad
	\widehat{a}_{11} = -2\widehat{a}_{1} - \widehat{a}_{3} + \frac{\widehat{a}_{10}}{2} \, ,
	\quad
	\widehat{a}_{13} =\widehat{a}_{14} = 4\widehat{a}_{1} + 2 \widehat{a}_{3} - \widehat{a}_{10}
	\, .
	\notag
\end{align}
The most simple isotropic curvature term $\frac{\mu \, L_c^2}{2} \rVert \mbox{Curl} \boldsymbol{P} \lVert^2$ corresponds to $a_1 = a_2 = a_3 = 1$ and would therefore be given by $\widehat{a}_{2,3,4,5,8,11,13,14}=0$, $\widehat{a}_{1}=\frac{\mu \, L_c}{4}$, $\widehat{a}_{10} = -\widehat{a}_{15} = \mu \, L_c$ Due to \cite{lewintan2020korn} the model is well-posed even for $a_1>0$ and $a_2=a_3=0$ if $\mu_{\tiny \mbox{micro}},\kappa_{\tiny \mbox{micro}}>0$.
The equilibrium equations without body forces are
\begin{align}
\mbox{Div}\overbrace{\left[2\mu_{e}\,\mbox{sym} \left(\boldsymbol{\mbox{D}u} - \boldsymbol{P} \right) + \lambda_{e} \mbox{tr} \left(\boldsymbol{\mbox{D}u} - \boldsymbol{P} \right) \boldsymbol{\mathbbm{1}}
+ 2\mu_{c}\,\mbox{skew} \left(\boldsymbol{\mbox{D}u} - \boldsymbol{P} \right)\right]}^{\mathlarger{\widetilde{\sigma}}:=}
&= \boldsymbol{0},
\notag
\\
\widetilde{\sigma}
- 2 \mu_{\mbox{\tiny micro}}\,\mbox{sym}\,\boldsymbol{P} - \lambda_{\tiny \mbox{micro}} \mbox{tr} \left(\boldsymbol{P}\right) \boldsymbol{\mathbbm{1}}
\hspace{8cm}
\label{eq:equi_RM}
\\
- \underbrace{
	\mu \, L_{c}^{2} \, \mbox{Curl}
\left(
a_1 \, \mbox{dev sym} \, \mbox{Curl} \, \boldsymbol{P} +
a_2 \, \mbox{skew} \,  \mbox{Curl} \, \boldsymbol{P} +
\frac{a_3}{3} \, \mbox{tr} \left(\mbox{Curl} \, \boldsymbol{P}\right)\mathbbm{1}
\right)
}_{\boldsymbol{m}}
&= \boldsymbol{0}.
\notag
\end{align}
The boundary conditions at the upper and lower surface (free surface) are 
\begin{align}
\boldsymbol{\widetilde{t}}(x_2 = \pm \, h/2) &= 
\pm \, \boldsymbol{\widetilde{\sigma}}(x_2) \cdot \boldsymbol{e}_2 = 
\boldsymbol{0} \, ,
\label{eq:BC_RM_gen}
\\
\boldsymbol{\eta}(x_2 = \pm \, h/2) &= 
\pm \, \boldsymbol{m} (x_2) \cdot \boldsymbol{\epsilon} \cdot \boldsymbol{e}_2 = 
\pm \, \boldsymbol{m} (x_2) \times \boldsymbol{e}_2 = 
\boldsymbol{0} \, ,
\notag
\end{align}
where the expression of $\boldsymbol{\widetilde{\sigma}}$ and $\boldsymbol{m}$ are in eq.(\ref{eq:equi_RM}),
$\boldsymbol{e}_2$ is the unit vector aligned to the $x_2$-direction,
$\boldsymbol{\epsilon}$ is the Levi-Civita tensor,
and
$\boldsymbol{m} = \mu \, L_c^2 \, \left(
a_1 \, \mbox{dev sym} \, \mbox{Curl} \, \boldsymbol{P} +
a_2 \, \mbox{skew} \,  \mbox{Curl} \, \boldsymbol{P} +
\frac{a_3}{3} \, \mbox{tr} \left(\mbox{Curl} \, \boldsymbol{P}\right)\mathbbm{1}
\right)$
is the generalized second order moment tensor.
The generalised traction vector is $\boldsymbol{\widetilde{t}}(x_2) \in \mathbb{R}^{3}$ and the generalized double traction tensor is $\boldsymbol{\eta}(x_2) \in \mathbb{R}^{3 \times 3}$.

According with the reference system shown in Fig.~\ref{fig:intro}, the ansatz for the displacement field and the microdistortion $\boldsymbol{P}$ is
\begin{equation}
\boldsymbol{u}(x_1,x_2)=
\left(
\begin{array}{c}
- \kappa_1 \, x_1 x_2 \\
v(x_2)+\frac{\kappa_1  x_1^2}{2} \\
0
\end{array}
\right) \, ,
\qquad
\boldsymbol{P}(x_1,x_2) =
\left(
\begin{array}{ccc}
P_{11}(x_2) & -\kappa_2 \, x_1 & 0 \\
\kappa_2 \, x_1 & P_{22}(x_2) & 0 \\
0 & 0 & P_{33}(x_2)
\end{array}
\right)
\, ,
\label{eq:ansatz_RM}
\end{equation}
while the  gradient of the displacement field results to be
\begin{equation}
\boldsymbol{\mbox{D}u} = 
\left(
\begin{array}{ccc}
-\kappa_1 \, x_2 & - \kappa_1 \, x_1 & 0 \\
\kappa_1 \, x_1 & v'(x_2) & 0 \\
0 & 0 & 0 \\
\end{array}
\right) \, .
\label{eq:grad_RM}
\end{equation}
We supply as well the homogenization relations between the macro-parameters and the meso- (with index $(\cdot)_{e}$) and micro-parameters \cite{d2019effective,neff2014unifying,neff2019identification}
\begin{equation}
\begin{array}{c}
\mu _{\tiny \mbox{macro}} = \frac{\mu_e \, \mu_{\tiny \mbox{micro}}}{\mu_e + \mu_{\tiny \mbox{micro}}} \, ,
\qquad\qquad
\kappa _{\tiny \mbox{macro}} = \frac{\kappa_e \, \kappa_{\tiny \mbox{micro}}}{\kappa_e + \kappa_{\tiny \mbox{micro}}} \, ,
\qquad\qquad
\text{with}
\quad
\left\{
\begin{array}{l}
\kappa_{\tiny \mbox{i}} = \frac{2\mu_{i} + 3\lambda_{i}}{3}
\\
i = \left\{e, {\tiny \mbox{micro}}, {\tiny \mbox{macro}}\right\}
\end{array} .
\right.
\end{array}
\label{eq:static_homo_relation}
\end{equation}
\subsection{One curvature parameter and zero Poisson's ratios $\nu_{e}=\nu_{\tiny \mbox{micro}}=0$}
Substituting the ansatz in the following simplified equilibrium eq.(\ref{eq:equi_RM_Poi_zero}) where $\nu_{e}=\nu_{\tiny \mbox{micro}}=0$ and $a_1=a_2=a_3=1$, results in the simplified curvature $\frac{\mu \, L_{c}^{2}}{2} \lVert \mbox{Curl} \, \boldsymbol{P} \rVert^2$ and
\begin{align}
\mbox{Div}\overbrace{\left[
2\mu_{e}\,\mbox{sym} \left(\boldsymbol{\mbox{D}u} - \boldsymbol{P} \right)
+ 2\mu_{c}\,\mbox{skew} \left(\boldsymbol{\mbox{D}u} - \boldsymbol{P} \right)\right]}^{\mathlarger{\widetilde{\sigma}}:=}
&= \boldsymbol{0} \, ,
\label{eq:equi_RM_Poi_zero}
\\
\widetilde{\sigma}
- 2 \mu_{\mbox{\tiny micro}}\,\mbox{sym}\,\boldsymbol{P}
- \mu \, L_{c}^{2} \, \mbox{Curl} \, 
\mbox{Curl} \, \boldsymbol{P}
&= \boldsymbol{0} \, ,
\notag
\end{align}
where the generalized moment tensor is $\boldsymbol{m} = \mu \, L_c^2 \, \mbox{Curl} \, \boldsymbol{P}$.
The equilibrium equations (\ref{eq:equi_RM_Poi_zero}) then are
\begin{align}
2 \, \mu _c (\boldsymbol{\kappa}_1-\boldsymbol{\kappa}_2)+2 \, \mu _e \left(v''\left(x_2\right)-P_{22}'\left(x_2\right)\right) =0 \, ,
\notag
\\
\mu \, L_c^2 \, P_{11}''\left(x_2\right)-2 \, \mu _e \left(P_{11}\left(x_2\right)+\boldsymbol{\kappa}_1 x_2\right)-2 \, \mu _{\tiny \mbox{micro}} \, P_{11}\left(x_2\right) =0 \, ,
\notag
\\
2 \, x_1 \, \mu _c \, (\boldsymbol{\kappa}_2-\boldsymbol{\kappa}_1) =0 \, ,
\label{eq:equi_equa_RM_L0}
\\
2 \, x_1 \, \mu _c \, (\boldsymbol{\kappa}_1-\boldsymbol{\kappa}_2) =0 \, ,
\notag
\\
2 \, \mu _e \, v'\left(x_2\right)-2 \, P_{22}\left(x_2\right) \left(\mu _e+\mu _{\tiny \mbox{micro}}\right) =0 \, ,
\notag
\\
\mu \, L_c^2 \, P_{33}''\left(x_2\right)-2 \, P_{33}\left(x_2\right) \left(\mu _e+\mu _{\tiny \mbox{micro}}\right) =0 \, .
\notag
\end{align}
It is clear that, in order to satisfy eq.(\ref{eq:equi_equa_RM_L0})$_3$ and eq.(\ref{eq:equi_equa_RM_L0})$_4$ either $\mu _c = 0$ or $\kappa_1 = \kappa_2 = \boldsymbol{\kappa}$;
we choose the latter option, which implies that the skew-symmetric part of the gradient of the displacement eq.(\ref{eq:grad_RM}) is the same as the skew-symmetric part of the microdistortion eq.(\ref{eq:ansatz_RM})$_2$.
This also implies that the Cosserat couple modulus $\mu_c$ does not play a role any more.
Consequently, the solution of eq.(\ref{eq:equi_equa_RM_L0}) is
\begin{align}
v(x_2) &= c_{0} \, x_{2} \, \left( 1 + \frac{\mu_{e}}{\mu_{\mbox{\tiny micro}}} \right) \, ,
\qquad
P_{11}(x_2) = 
c_1 \, e^{-\frac{f_{1} x_{2}}{L_c}}
+ c_2 \, e^{\frac{f_{1} x_{2}}{L_c}}
- \frac{\mu _e}{\mu _e+\mu_{\mbox{\tiny micro}}} \, \boldsymbol{\kappa} \, x_{2} \, , 
\label{eq:sol_fun_disp_RM_L0}
\\
P_{22}(x_2) &= c_{0} \, \frac{\mu _e}{\mu_{\mbox{\tiny micro}}} \, ,
\qquad
P_{33}(x_2) =
c_3 \, e^{-\frac{f_{1} x_{2}}{L_c}}
+ c_4 \, e^{\frac{f_{1} x_{2}}{L_c}} \, ,
\qquad
f_1:= \sqrt{\frac{2\left(\mu_{e} + \mu_{\tiny \mbox{micro}}\right)}{\mu}} \, .
\notag
\end{align}

\begin{figure}[H]
\centering
\includegraphics[width=\linewidth]{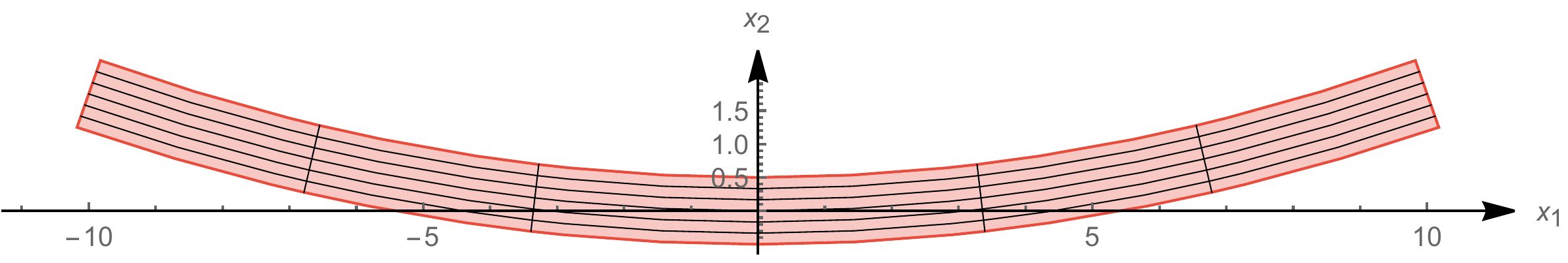}
\caption{Deformed shape of a relaxed micromorphic plate for zero Poisson's ratio $\nu_{e}=\nu_{\tiny \mbox{micro}}=0$ and with one-constant curvature energy $\frac{\mu \, L_{c}^{2}}{2}\lVert\mbox{Curl} \boldsymbol{P}\rVert^2$ for $\boldsymbol{\kappa}=7/200$.}
\label{fig:deformed_shape_RM}
\end{figure}

The boundary conditions eq.(\ref{eq:BC_RM_gen}) on the upper and lower surfaces allow to evaluate the constants $c_i, \, i \in \left\{0,1,2,3,4\right\}$ as shown in the following \!\!\!\!\!
\footnote{
$\text{sech}(x) := \frac{1}{\text{cosh}(x)} = \frac{2}{e^{x}+e^{-x}}$.
}
\begin{equation}
c_0 = c_3 = c_4 = 0 \, ,
\qquad
c_1 = -c_2 = 
\frac{\mu_{\mbox{\tiny micro}} \, \boldsymbol{\kappa}}{2 \left(\mu _e+\mu_{\mbox{\tiny micro}}\right)}
\frac{L_c}{f_{1}} \, 
\text{sech}\left(\frac{f_{1} h}{2 L_c}\right) \, .
\label{eq:BC_RM_L0}
\end{equation}
Finally, the displacement and microdistortion components result in
\begin{align}
u_1(x_1,x_2) &= -\boldsymbol{\kappa} \, x_1 \, x_2 \, ,
\quad
u_2(x_1,x_2) = \frac{\boldsymbol{\kappa} \, x_1^2}{2} \, ,
\quad
P_{12} (x_1) = - P_{21} (x_1) = - \boldsymbol{\kappa} \, x_1 \, ,
\quad
P_{22} (x_2)= 0 \, ,
\label{eq:disp_P_BC_RM}
\\
P_{11} (x_2) &= 
- \frac{\mu_{\mbox{\tiny micro}} \, \boldsymbol{\kappa}}{\mu _e+\mu_{\mbox{\tiny micro}}}
\frac{L_c}{f_{1}}
\text{sech}\left(\frac{f_{1} h}{2 L_c}\right) \sinh \left(\frac{f_{1} x_{2}}{L_c}\right)
- \frac{\mu _e \, \boldsymbol{\kappa}}{\mu _e+\mu_{\mbox{\tiny micro}}} \, x_{2}
\, ,
\qquad
P_{33} (x_2)= 0 \, .
\notag
\end{align}
It is underlined that for this specific case $P_{22}$ and $P_{33}$ turn out to be equal to zero.
In Fig.~\ref{fig:deformed_shape_RM} we show the deformed shape due to the displacement field solution

\begin{figure}[H]
	\centering
	\begin{subfigure}{.45\textwidth}
		\centering
		\includegraphics[width=0.95\linewidth]{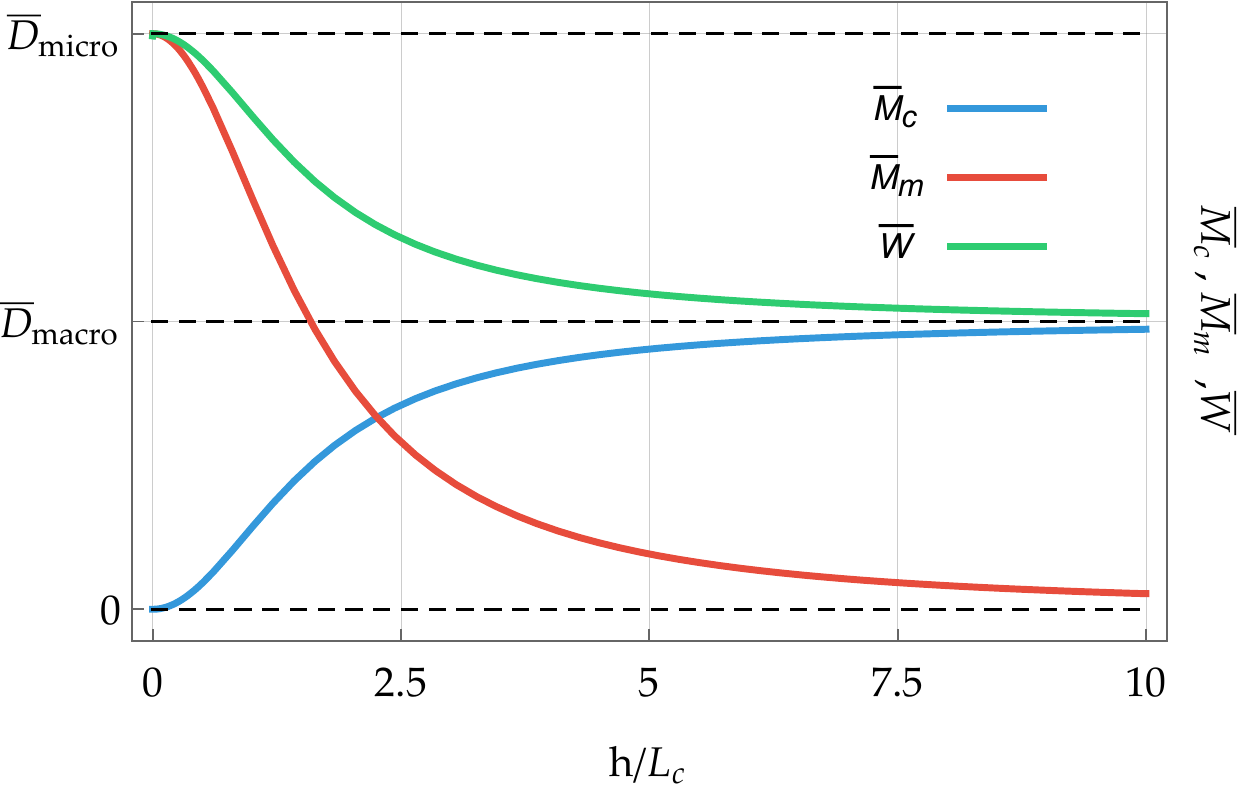}
		\caption{}
		\label{fig:all_plot_RM}
	\end{subfigure}%
	\begin{subfigure}{.45\textwidth}
		\centering
		\includegraphics[width=0.95\linewidth]{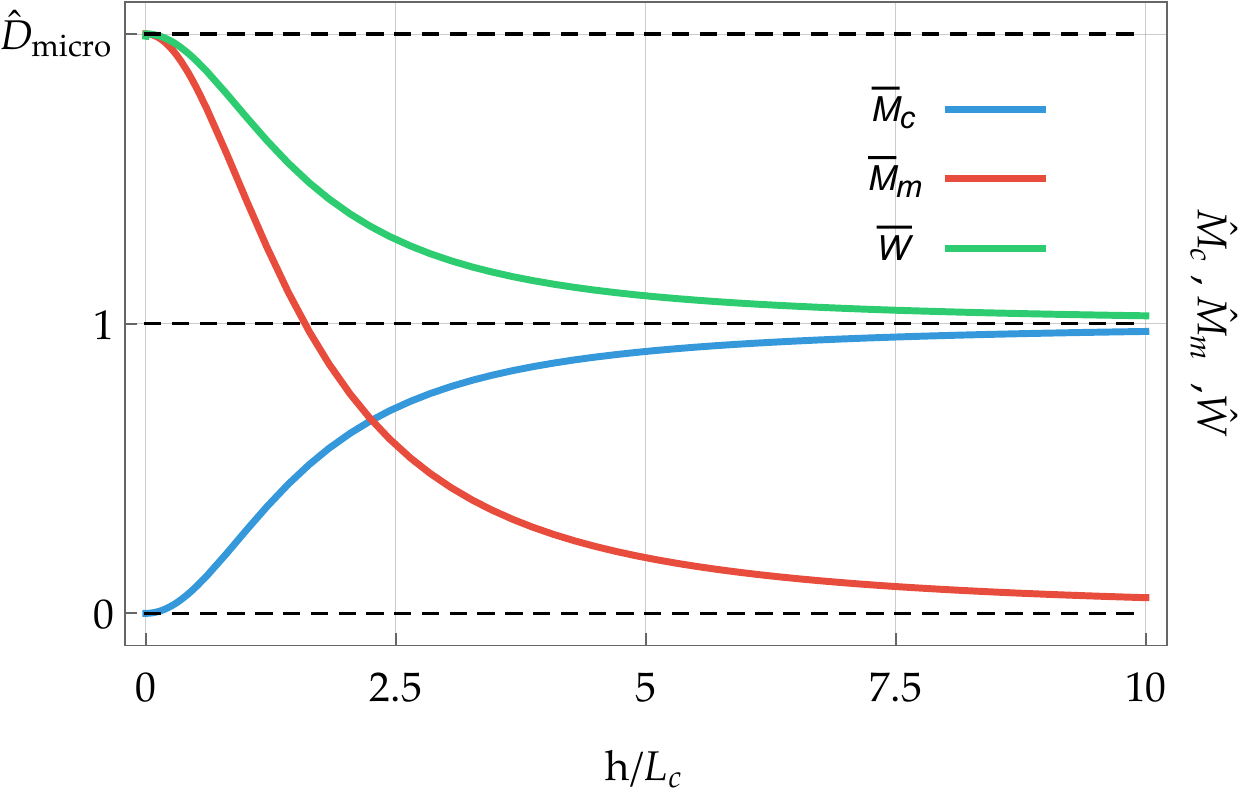}
		\caption{}
		\label{fig:all_plot_RM_2}
	\end{subfigure}
	\caption{(\textbf{Relaxed micromorphic model}, one curvature parameter, zero Poisson) 
		(a) Bending moments and energy while varying $L_c$. Observe that the bending stiffness remains bounded as $L_c \to \infty$ ($h\to 0$). This is one of the distinguishing features of the relaxed micromorphic model;
		(b) dimensionless stiffness where $\widehat{D}_{\tiny \mbox{micro}} = \frac{\overline{D}_{\tiny \mbox{micro}}}{\overline{D}_{\tiny \mbox{macro}}}$ is the stiffness normalized against the classical linear elastic stiffness. Since the two pictures show exactly the same features, in the following we will only show variant (a). The values of the parameters used are: $\mu _e = 1$, $\mu_{\mbox{\tiny micro}} = 1$, $\mu = 1$.
	}
\end{figure}

The classical bending moment, the higher-order bending moment, and energy (per unit area d$x_1$d$x_3$) expressions are reported next in the following eq.(\ref{eq:sigm_ene_dimensionless_RM})
\begin{align}
M_{\mbox{c}} (\boldsymbol{\kappa})&:=
\displaystyle\int\limits_{-h/2}^{h/2}
\langle \boldsymbol{\widetilde{\sigma}} \boldsymbol{e}_1 , \boldsymbol{e}_1 \rangle  \, x_2 \, 
\mbox{d}x_{2}
= 
\frac{h^3}{12}
\frac{2 \, \mu _e \, \mu_{\mbox{\tiny micro}} }{\mu _e+\mu_{\mbox{\tiny micro}}}
\left[
1
- \frac{12}{f_{1}^2}\left(\frac{L_c}{h}\right)^2
+ \frac{24}{f_{1}^3} \left(\frac{L_c}{h}\right)^3 \tanh \left(\frac{f_{1} h}{2 L_c}\right)
\right]
\boldsymbol{\kappa}
\, ,
\notag
\\
M_{\mbox{m}}(\boldsymbol{\kappa}) &:=
\displaystyle\int\limits_{-h/2}^{h/2}
\langle \left(\boldsymbol{m} \times \boldsymbol{e}_1\right) \boldsymbol{e}_2 , \boldsymbol{e}_1 \rangle
\, \mbox{d}x_{2}
= 
\frac{h^3}{12}
\frac{2 \, \mu_{\mbox{\tiny micro}} \, \mu }{\mu _e + \mu_{\mbox{\tiny micro}}}
\left[
6\left(\frac{L_c}{h}\right)^2
-\frac{12}{f_{1}} \left(\frac{L_c}{h}\right)^3
\tanh \left(\frac{f_{1} h}{2 L_c}\right)
\right]
\boldsymbol{\kappa}
\, ,
\label{eq:sigm_ene_dimensionless_RM}
\\
W_{\mbox{tot}} (\boldsymbol{\kappa})&:=
\displaystyle\int\limits_{-h/2}^{+h/2} W \left(\boldsymbol{\mbox{D}u}, \boldsymbol{P}, \mbox{Curl}\boldsymbol{P}\right) \, \mbox{d}x_{2}
= 
\notag
\\*
&
\frac{1}{2}
\frac{h^3}{12}
\frac{2 \, \mu _e \, \mu_{\mbox{\tiny micro}} }{\mu _e+\mu_{\mbox{\tiny micro}}}
\left[
1
- \left(\frac{12}{f_{1}^2} - 6 \frac{\mu}{\mu _e}\right) \left(\frac{L_c}{h}\right)^2
+ \frac{2}{f_{1}}\left(\frac{12}{f_{1}^2} - 6 \frac{\mu}{\mu _e}\right) \left(\frac{L_c}{h}\right)^3 \tanh \left(\frac{f_{1} h}{2 L_c}\right)
\right]
\boldsymbol{\kappa}^2
\, .
\notag
\end{align}
The plot of the bending moments and the strain energy divided by $\frac{h^3}{12}\boldsymbol{\kappa}$ and $\frac{1}{2}\frac{h^3}{12}\boldsymbol{\kappa}^2$, respectively, while changing $L_c$ is shown in Fig.~\ref{fig:all_plot_RM}.
As in the classical Cauchy elastic case, we have
\begin{equation}
\frac{\mbox{d}}{\mbox{d}\boldsymbol{\kappa}}W_{\mbox{tot}}(\boldsymbol{\kappa}) = M_{\mbox{c}} (\boldsymbol{\kappa}) + M_{\mbox{m}} (\boldsymbol{\kappa}) \, .
\end{equation}
\subsubsection{Remarks on the boundary conditions: consistent coupling}
If a finite slice is cut along the $x_1$-direction at distance $\pm \, b/2$ from the origin, a finite solid in the $x_1-x_2$ plane is obtained. The kinematic boundary conditions that arise on the cut surfaces $\Gamma_-$ and $\Gamma_+$ (see Fig.~\ref{fig:intro}) for this solution are
\begin{equation}
\boldsymbol{u} \, \Bigr|_{x_1 = \pm \, b/2} =
\boldsymbol{\kappa}
\left(
\begin{array}{c}
\mp \, b \, x_2/2 \\
b^2/8 \\
0 \\
\end{array}
\right) \, ,
\qquad
\left(\boldsymbol{P}  \cdot \boldsymbol{e}_{2}\right)\Bigr|_{x_1 = \pm \, b/2} =
\left(\boldsymbol{\mbox{D}u} \cdot \boldsymbol{e}_{2}\right)\Bigr|_{x_1 = \pm \, b/2} =
\left(
\begin{array}{c}
- \boldsymbol{\kappa} \, b/2 \\
0 \\
0 \\
\end{array}
\right) \, .
\label{eq:bc_inverse_kinematic_finite}
\end{equation}
Moreover, it holds that $\boldsymbol{P} \, \boldsymbol{e}_3 = \boldsymbol{\mbox{D}u} \, \boldsymbol{e}_3$ at $x_1 = \pm \, b/2$, which implies that for zero Poisson's ratios $\nu_{e}=\nu_{\tiny \mbox{micro}}=0$, the \textbf{consistent coupling condition} 
\begin{equation}
\boldsymbol{P} \times \boldsymbol{e}_1 = \boldsymbol{\mbox{D}u} \times \boldsymbol{e}_1
\end{equation}
at the lateral boundary $\Gamma_-$ and $\Gamma_+$ is exactly verified.
This means that in this exceptional case we could start alternatively with a finite domain boundary value problem and describe the bending condition according to eq.(\ref{eq:bc_inverse_kinematic_finite})$_1$ and in addition require $\boldsymbol{P} \times \boldsymbol{e}_1 = \boldsymbol{\mbox{D}u} \times \boldsymbol{e}_1$ at the lateral boundary $\Gamma_-$ and $\Gamma_+$, giving the ansatz eq.(\ref{eq:ansatz_RM}), since the solution of the problem is unique \cite{neff2014unifying}.

The static boundary conditions that arise on the cut surfaces $\Gamma_-$ and $\Gamma_+$ (see Fig.~\ref{fig:intro}) are
\begin{align}
\widetilde{t}_{1} \, \Bigr|_{x_1 = \pm \, b/2} &= 
\pm
\frac{2 \mu _e \, \mu_{\mbox{\tiny micro}}}{\mu _e+\mu_{\mbox{\tiny micro}}}
\left(
-x_{2} 
+\text{sech}\left(\frac{f_{1} h}{2 L_c}\right) \sinh \left(\frac{f_{1} x_{2}}{L_c}\right)
\frac{L_c}{f_{1}}
\right)
\boldsymbol{\kappa} \, ,
\label{eq:bc_inverse_static_finite}
\\*
\eta_{12}\Bigr|_{x_1 = \pm \, b/2} &=
\pm \, 
\frac{\mu_{\mbox{\tiny micro}}}{\mu _e+\mu_{\mbox{\tiny micro}}}
\mu \,  L_c^2
\left(
1-\text{sech}\left(\frac{f_{1} h}{2 L_c}\right) \cosh \left(\frac{f_{1} x_{2}}{L_c}\right)
\right)
\boldsymbol{\kappa} \, ,
\notag
\end{align}
where $\widetilde{t}_{1}$ and $\eta_{12}$ are the only components different from zero.
It is also highlighted that the compatibility condition expressed by eq.(\ref{eq:bc_inverse_static_finite})$_2$ is satisfied.
The boundary conditions eq.(\ref{eq:BC_RM_gen})  on the upper and lower surface in addition to eq.(\ref{eq:bc_inverse_kinematic_finite}) or eq.(\ref{eq:bc_inverse_static_finite}) (the choice is up to the reader) on the left and right surface are enough to retrieve the full solution eq.(\ref{eq:disp_P_BC_RM}).

\begin{figure}[H]
	\centering
	\begin{subfigure}{.45\textwidth}
		\centering
		\includegraphics[width=0.95\linewidth]{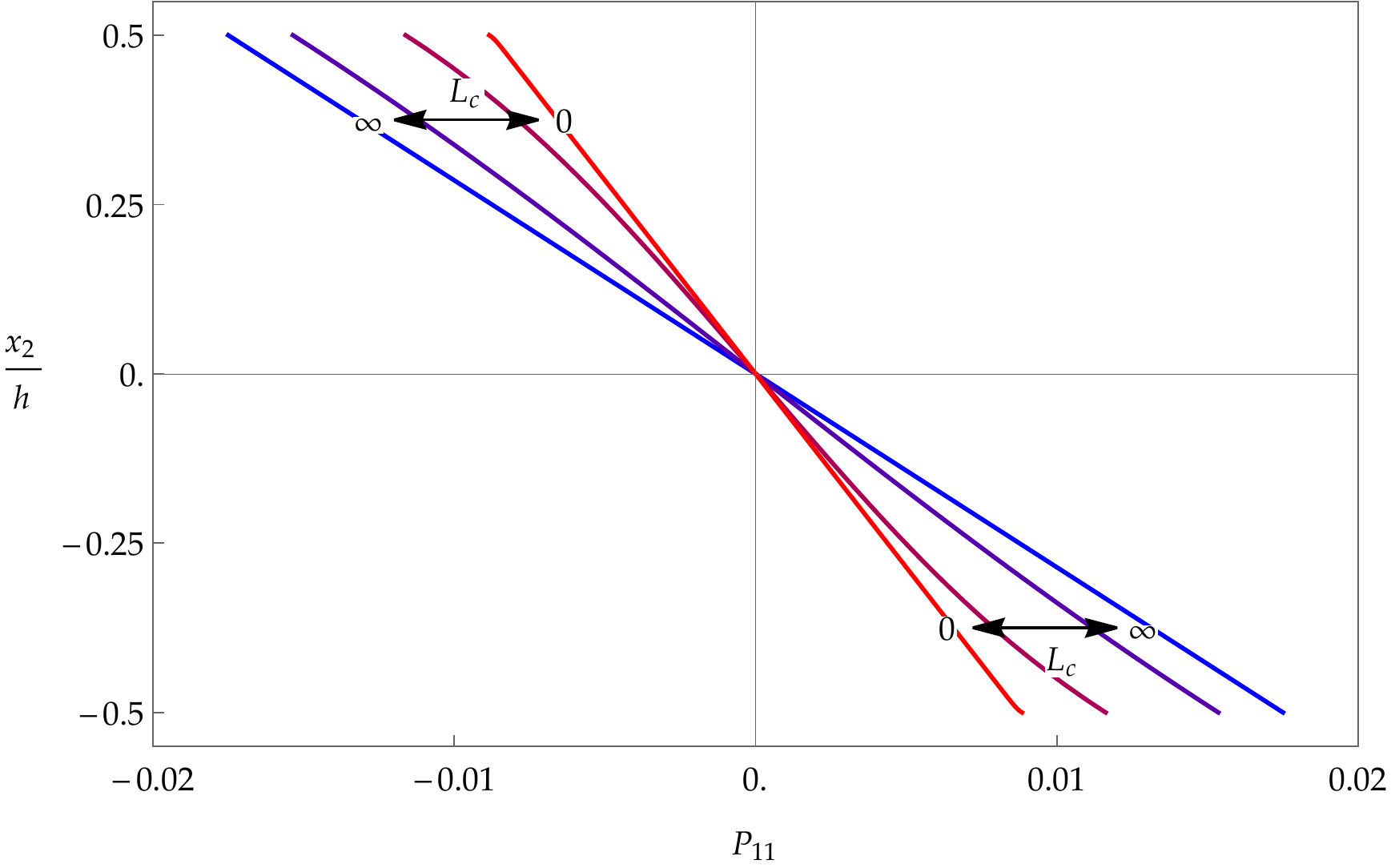}
		\caption{}
		\label{fig:P11_RM_2}
	\end{subfigure}%
	\begin{subfigure}{.45\textwidth}
		\centering
		\includegraphics[width=0.95\linewidth]{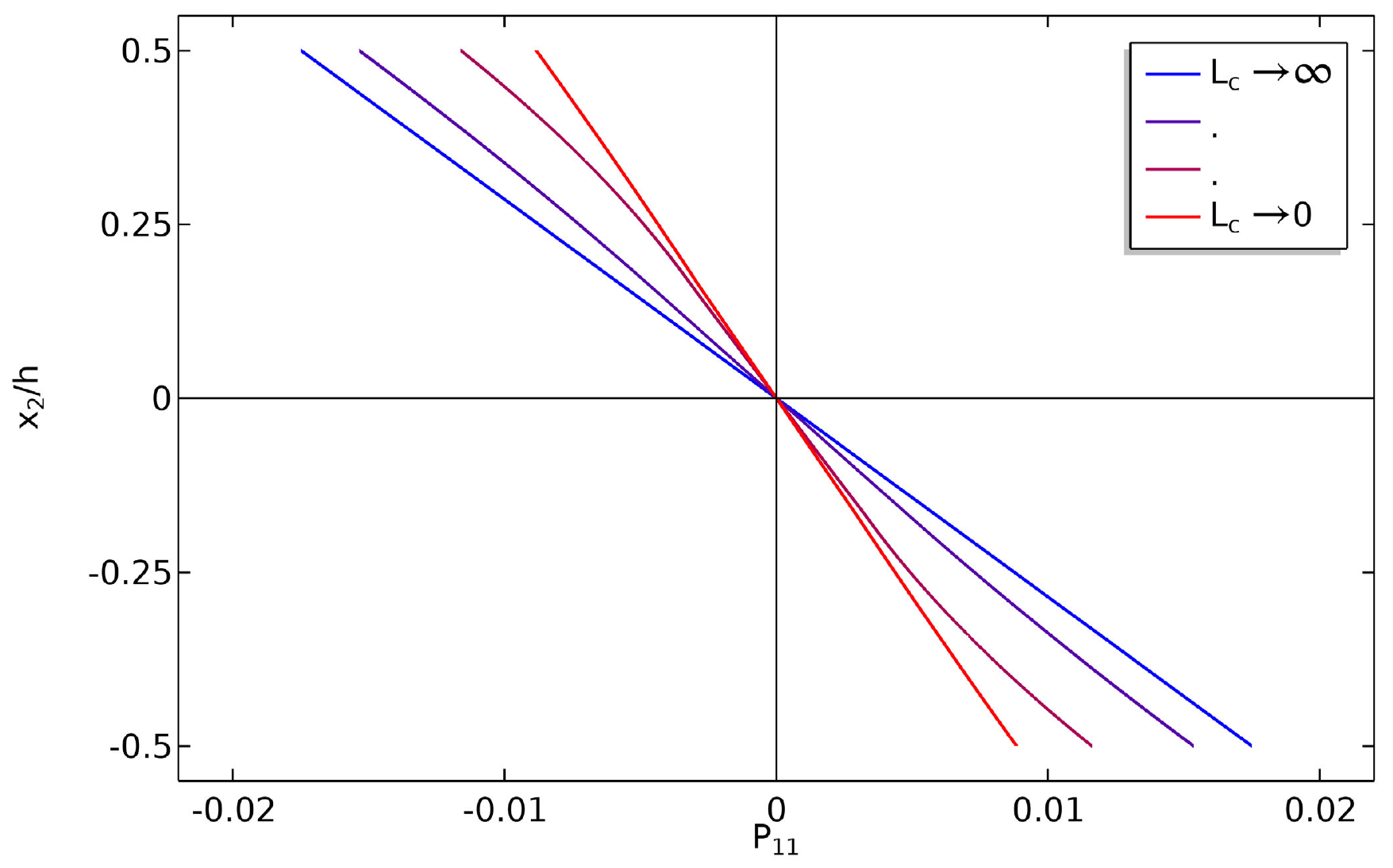}
		\caption{}
		\label{fig:P11_Comsol_RM_2}
	\end{subfigure}
	\caption{(\textbf{Relaxed micromorphic model}, one curvature parameter, zero Poisson)		
		(a) Distribution across the thickness of $P_{11}$ while varying $L_c$.
		On the vertical axis we have the dimensionless thickness while on the horizontal we have the quantity $P_{11}$.
		Notice that the red curve correspond to the limit for $L_c \to 0$ while the blue curve correspond to the limit for $L_c \to \infty$.
		(b) distribution across the thickness of $P_{11}$ while varying $L_c$ obtained via a COMSOL$^{\tiny \mbox{\textregistered}}$ simulation in which a finite cylindrical bending problem is solved by using the consistent coupling kinematic boundary conditions eq.(\ref{eq:bc_inverse_kinematic_finite}). The values of the parameters used are: $\mu _{\tiny \mbox{micro}} = 1$, $\mu _{e} = 1$, $\mu = 1$, $h = 1$, $\boldsymbol{\kappa} = 7/200$.
	}
\end{figure}

\begin{figure}[H]
	\centering
	\begin{subfigure}{.45\textwidth}
		\centering
		\includegraphics[width=0.9\linewidth]{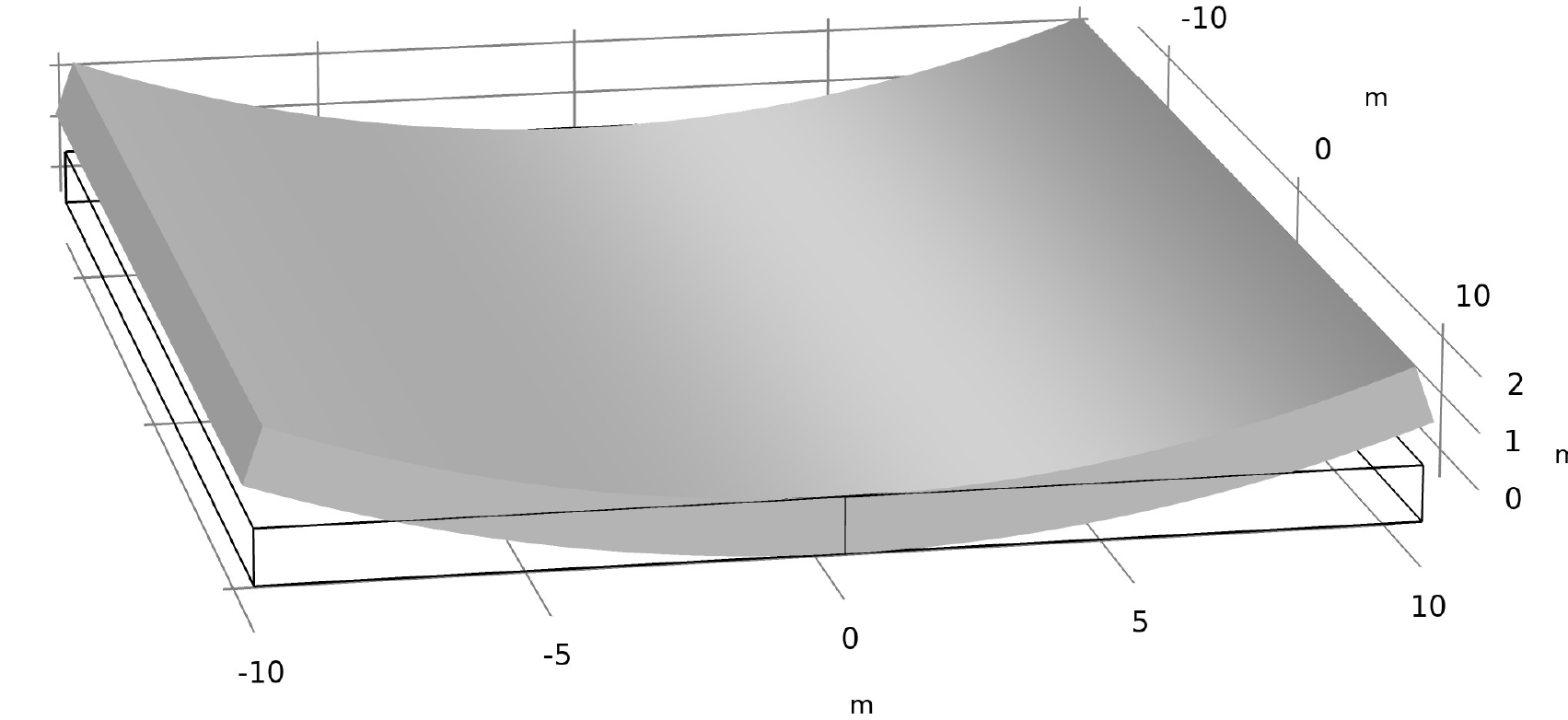}
		\caption{}
		\label{fig:P11_RM_3}
	\end{subfigure}%
	\begin{subfigure}{.45\textwidth}
		\centering
		\includegraphics[width=0.9\linewidth]{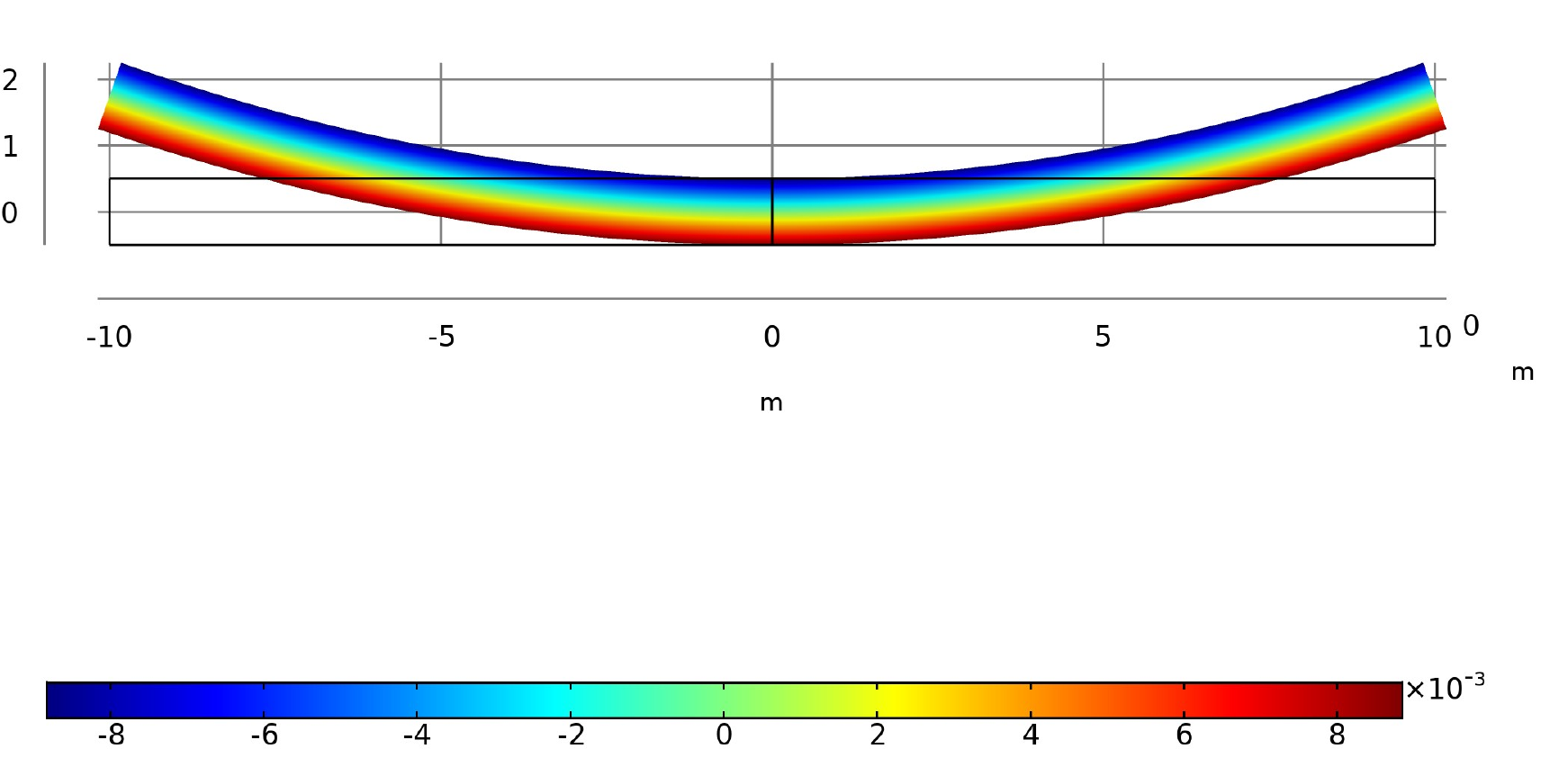}
		\caption{}
		\label{fig:P11_RM_4}
	\end{subfigure}
	\begin{subfigure}{.45\textwidth}
	\centering
	\includegraphics[width=0.9\linewidth]{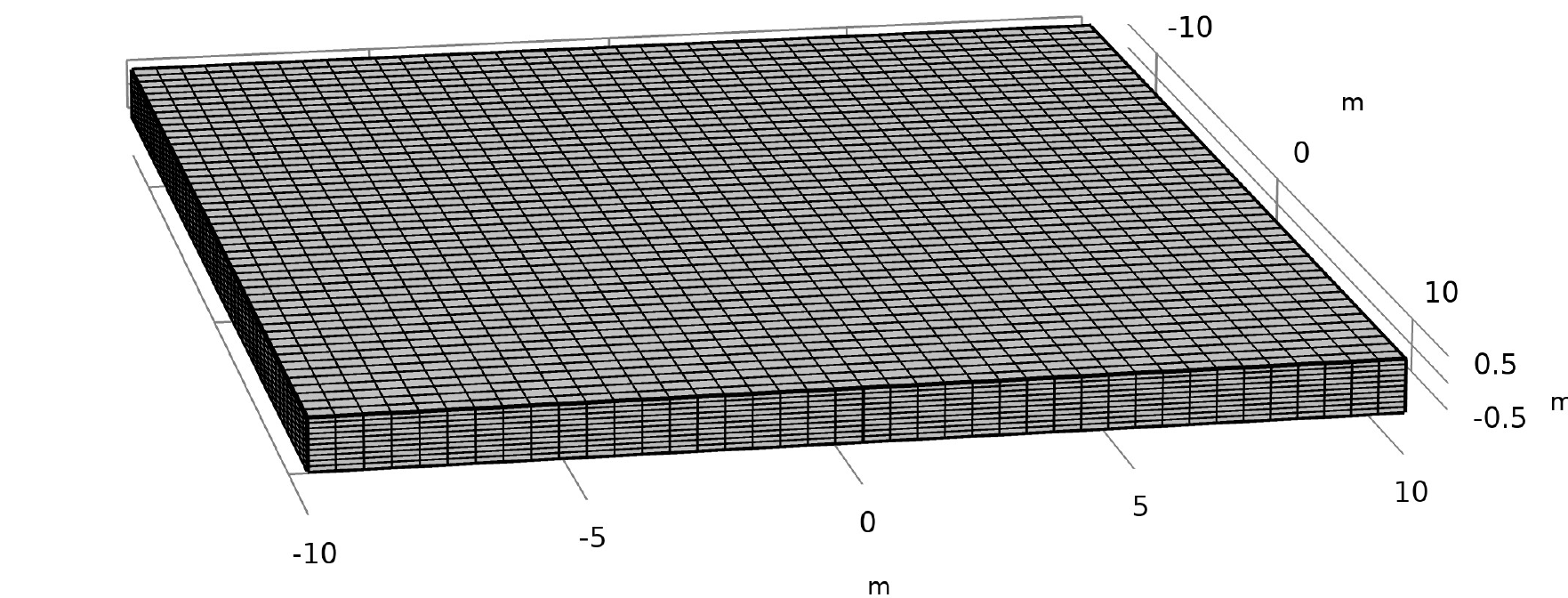}
	\caption{}
	\label{fig:P11_RM_5}
	\end{subfigure}
	\caption{
		(a) 3D-deformed shape of a finite relaxed micromorphic plate for which a finite cylindrical bending problem is solved by using the consistent coupling kinematic boundary conditions;
		(b) plot of the distribution across the thickness of $P_{11}$ on the deformed shape;
		(c) meshed finite plate.
		The values of the parameters used are: $\mu _{\tiny \mbox{micro}} = 1$, $\mu _{e} = 1$, $\mu = 1$, $h = 1$, $\boldsymbol{\kappa} = 7/200$.
	}
\label{fig:P11_RM_tot}
\end{figure}

The plot of $P_{11}$ obtained analytically and numerically (via COMSOL$^{\tiny \mbox{\textregistered}}$) while changing $L_c$ is shown in Fig.~\ref{fig:P11_RM_2} and Fig.~\ref{fig:P11_Comsol_RM_2}, respectively: the two give correctly the exact same results.
In Fig.~\ref{fig:P11_RM_tot} is shown the deformed shape obtained thanks to simulations done via COMSOL$^{\tiny \mbox{\textregistered}}$.
\subsubsection{Limit cases}
We consider in the following the two limit cases: first $L_c \to 0$ corresponds to arbitrary thick specimens, while secondly $L_c \to \infty$ corresponds conceptually to arbitrary thin specimens. We see that in the relaxed micromorphic  model, this corresponds unequivocally to the stiffness $D_{\mbox{\tiny macro}}$ and $D_{\mbox{\tiny micro}}$, respectively, since
\begin{align}
\displaystyle\lim_{L_c\to0} M_{\mbox{c}} (\boldsymbol{\kappa})
=& \frac{h^3}{12} \, \frac{2 \, \mu _e \, \mu _{\tiny \mbox{micro}} }{\mu _e+\mu _{\tiny \mbox{micro}}} \, \boldsymbol{\kappa} 
= D_{\tiny \mbox{macro}} \, \boldsymbol{\kappa} \, , 
\hspace{2.1cm}
\displaystyle\lim_{L_c\to\infty} M_{\mbox{c}} (\boldsymbol{\kappa})
= 0 \, ,
\notag
\\
\displaystyle\lim_{L_c\to0} M_{\mbox{m}} (\boldsymbol{\kappa})
=& 0 \, ,
\hspace{6cm}
\displaystyle\lim_{L_c\to\infty} M_{\mbox{m}} (\boldsymbol{\kappa})
= \frac{h^3}{12} \, 2 \, \mu _{\tiny \mbox{micro}} \, \boldsymbol{\kappa} 
= D_{\tiny \mbox{micro}} \, \boldsymbol{\kappa} \, , 
\\
\displaystyle\lim_{L_c\to0} W_{\mbox{tot}} (\boldsymbol{\kappa})
=& 
\frac{1}{2} \, \frac{h^3}{12} \, \frac{2 \, \mu _e \, \mu_{\tiny \mbox{micro}} }{\mu _e + \mu_{\tiny \mbox{micro}}} \, \boldsymbol{\kappa}^2 
= 
\frac{1}{2} \,D_{\tiny \mbox{macro}} \, \boldsymbol{\kappa}^2 \, , 
\notag
\\
\displaystyle\lim_{L_c\to\infty} W_{\mbox{tot}} (\boldsymbol{\kappa})
=&
\frac{1}{2} \, \frac{h^3}{12} \,  2\mu_{\tiny \mbox{micro}} \, \boldsymbol{\kappa}^2 
= 
\frac{1}{2} \,D_{\tiny \mbox{micro}} \, \boldsymbol{\kappa}^2 \, .
\notag
\end{align}
\subsection{One curvature parameter and arbitrary Poisson's ratios $\nu_{\tiny \mbox{micro}}$ and $\nu_{e}$}
Substituting the ansatz in the following simplified equilibrium eq.(\ref{eq:equi_RM_mono}) where $a_1=a_2=a_3=1$
\begin{align}
\mbox{Div}\overbrace{\left[2\mu_{e}\,\mbox{sym} \left(\boldsymbol{\mbox{D}u} - \boldsymbol{P} \right) + \lambda_{e} \mbox{tr} \left(\boldsymbol{\mbox{D}u} - \boldsymbol{P} \right) \boldsymbol{\mathbbm{1}}
+ 2\mu_{c}\,\mbox{skew} \left(\boldsymbol{\mbox{D}u} - \boldsymbol{P} \right)\right]}^{\mathlarger{\widetilde{\sigma}}:=}
&= \boldsymbol{0},
\label{eq:equi_RM_mono}
\\
\widetilde{\sigma}
- 2 \mu_{\mbox{\tiny micro}}\,\mbox{sym}\,\boldsymbol{P} - \lambda_{\tiny \mbox{micro}} \mbox{tr} \left(\boldsymbol{P}\right) \boldsymbol{\mathbbm{1}}
- \mu \, L_{c}^{2} \, \mbox{Curl} \, 
\mbox{Curl} \, \boldsymbol{P}
&= \boldsymbol{0} \, .
\notag
\end{align}
where the generalized moment tensor is $\boldsymbol{m} = \mu \, L_c^2 \, \mbox{Curl} \, \boldsymbol{P}$.
The equilibrium equation (\ref{eq:equi_RM_mono}) are
\begin{align}
2 \mu _c (\kappa_{1}-\kappa_{2})-\lambda _e \left(\kappa_{1}+P_{11}'(x_{2})+P_{22}'(x_{2})+P_{33}'(x_{2})-v''(x_{2})\right)+2 \mu _e \left(v''(x_{2})-P_{22}'(x_{2})\right) &= 0 \, , 
\notag
\\
\mu \, L_c^2 P_{11}''(x_{2})-P_{11}(x_{2}) \left(2 \left(\mu _e+\mu_{\mbox{\tiny micro}}\right)+\lambda_{\mbox{\tiny micro}}\right)&
\notag \\
-\lambda _e \left(P_{11}(x_{2})+P_{22}(x_{2})+P_{33}(x_{2})-v'(x_{2})+\kappa_{1} x_{2}\right)-2 \kappa_{1} x_{2} \mu _e+\lambda_{\mbox{\tiny micro}} (-P_{22}(x_{2})-P_{33}(x_{2})) &= 0
\notag
\\
2 x_{1} \mu _c (\kappa_{2}-\kappa_{1}) &= 0 \, ,
\notag
\\
2 x_{1} \mu _c (\kappa_{1}-\kappa_{2}) &= 0 \, ,
\label{eq:equi_equa_RM_full}
\\
-\lambda _e \left(P_{11}(x_{2})+P_{22}(x_{2})+P_{33}(x_{2})-v'(x_{2})+\kappa_{1} x_{2}\right)+2 \mu _e \left(v'(x_{2})-P_{22}(x_{2})\right)
\notag
\\
-\lambda_{\mbox{\tiny micro}} P_{11}(x_{2}) -P_{22}(x_{2}) \left(\lambda_{\mbox{\tiny micro}}+2 \mu_{\mbox{\tiny micro}}\right)-\lambda_{\mbox{\tiny micro}} P_{33}(x_{2}) &= 0 \, ,
\notag
\\
\mu \, L_c^2 P_{33}''(x_{2})-\lambda _e \left(P_{11}(x_{2})+P_{22}(x_{2})+P_{33}(x_{2})-v'(x_{2})+\kappa_{1} x_{2}\right)-2 \mu _e P_{33}(x_{2})
\notag
\\
-\lambda_{\mbox{\tiny micro}} P_{11}(x_{2})-\lambda_{\mbox{\tiny micro}} P_{22}(x_{2})-P_{33}(x_{2}) \left(\lambda_{\mbox{\tiny micro}}+2 \mu_{\mbox{\tiny micro}}\right) &= 0 \, .
\notag
\end{align}

In order to satisfy eq.(\ref{eq:equi_equa_RM_full})$_3$ and eq.(\ref{eq:equi_equa_RM_full})$_4$ either $\mu _c = 0$ or $\kappa_1 = \kappa_2 = \boldsymbol{\kappa}$;
we have chosen the latter option, which implies that the skew-symmetric part of the gradient of the displacement eq.(\ref{eq:grad_RM}) is the same as the skew-symmetric part of the microdistortion eq.(\ref{eq:ansatz_RM})$_2$.
This also implies that the Cosserat couple modulus $\mu_c$ does not play a role any more.

From eq.(\ref{eq:equi_equa_RM_full})$_1$ it is possible to evaluate $v''(x_2)$ and consequently $v'(x_2)$
\begin{align}
v''(x_2) &= 
\frac{\lambda _e \left(\boldsymbol{\kappa} +P_{11}'(x_{2})+P_{22}'(x_{2})+P_{33}'(x_{2})\right)+2 \mu _e P_{22}'(x_{2})}{\lambda _e+2 \mu _e} \, ,
\label{eq:vv2_RM_full}
\\
v'(x_2) &= 
\frac{\lambda _e \left(\boldsymbol{\kappa} \, x_2 +P_{11}(x_{2})+P_{22}(x_{2})+P_{33}(x_{2})\right)+2 \mu _e P_{22}(x_{2})}{\lambda _e+2 \mu _e} + c_0 \, .
\notag
\end{align}
By substituting back eq.(\ref{eq:vv2_RM_full}) in eq.(\ref{eq:equi_equa_RM_full}), we can evaluate $P_{22}(x_2)$ and its derivatives from equation eq.(\ref{eq:equi_equa_RM_full})$_5$
\begin{align}
P_{22}(x_2) &= 
- \frac{\lambda_{\mbox{\tiny micro}} (P_{11}(x_{2}) + P_{33}(x_{2}))}{\lambda_{\mbox{\tiny micro}}+2 \mu_{\mbox{\tiny micro}}}
+ \frac{\lambda _e+2 \mu _e}{\lambda_{\mbox{\tiny micro}}+2 \mu_{\mbox{\tiny micro}}}c_0 \, ,
\label{eq:p22_RM_full}
\\
P'_{22}(x_2) &= 
- \frac{\lambda_{\mbox{\tiny micro}} (P'_{11}(x_{2}) + P'_{33}(x_{2}))}{\lambda_{\mbox{\tiny micro}}+2 \mu_{\mbox{\tiny micro}}} \, ,
\qquad\qquad
P''_{22}(x_2) = 
- \frac{\lambda_{\mbox{\tiny micro}} (P''_{11}(x_{2}) + P''_{33}(x_{2}))}{\lambda_{\mbox{\tiny micro}}+2 \mu_{\mbox{\tiny micro}}} \, .
\notag
\end{align}
After substituting eq.(\ref{eq:p22_RM_full}) in eq.(\ref{eq:equi_equa_RM_full}) the following two second order ordinary differential equations in $P_{11}(x_2)$ and $P_{33}(x_2)$ are retrieved
\begin{align}
b_{0} \, c_{0} + \boldsymbol{\kappa} \, b_{3} \, x_{2} - b_{1} \, P_{11}(x_{2}) - b_{2} \, P_{33}(x_{2}) + \mu \, L_c^2 \, P_{11}''(x_{2}) = 0 \, ,
\label{eq:equi_equa_2_RM_full}
\\
b_{0} \, c_{0} + \boldsymbol{\kappa} \, b_{4} \, x_{2} - b_{1} \, P_{33}(x_{2}) - b_{2} \, P_{11}(x_{2}) + \mu \, L_c^2 \, P_{33}''(x_{2}) = 0 \, ,
\notag
\end{align}
where
\begin{align}
b_0 &:= \frac{2 \left(\lambda _e \mu_{\mbox{\tiny micro}}-\mu _e \lambda_{\mbox{\tiny micro}}\right)}{\lambda_{\mbox{\tiny micro}}+2 \mu_{\mbox{\tiny micro}}} \, ,
\qquad\quad
b_3 :=-\frac{4 \mu _e \left(\lambda _e+\mu _e\right)}{\lambda _e+2 \mu _e} \, ,
\qquad\quad
b_4 :=-\frac{2 \mu _e \, \lambda _e}{\lambda _e+2 \mu _e} \, ,
\notag
\\
b_1 &:= \frac{4 \mu _e \left(\lambda _e+\mu _e\right)}{\lambda _e+2 \mu _e}+\frac{4 \mu_{\mbox{\tiny micro}} \left(\lambda_{\mbox{\tiny micro}}+\mu_{\mbox{\tiny micro}}\right)}{\lambda_{\mbox{\tiny micro}}+2 \mu_{\mbox{\tiny micro}}} = \left(\widehat{\lambda} _e + 2 \mu _e\right) + \left(\widehat{\lambda}_{\mbox{\tiny micro}} + 2 \mu_{\mbox{\tiny micro}}\right)\, ,
\label{eq:equi_equa_2_RM_full_const}
\\
b_2 &:= \frac{2 \lambda _e \mu _e}{\lambda _e+2 \mu _e}+\frac{2 \lambda_{\mbox{\tiny micro}} \mu_{\mbox{\tiny micro}}}{\lambda_{\mbox{\tiny micro}}+2 \mu_{\mbox{\tiny micro}}} = \widehat{\lambda}_e + \widehat{\lambda}_{\mbox{\tiny micro}} \, .
\notag
\end{align}
These additional relations between the parameters are satisfied
\begin{align}
b_1 + b_2 &= \frac{2 \mu _e \left(3 \lambda _e+2 \mu _e\right)}{\lambda _e+2 \mu _e} + 
\frac{2 \mu_{\mbox{\tiny micro}} \left(3 \lambda_{\mbox{\tiny micro}}+2 \mu_{\mbox{\tiny micro}}\right)}{\lambda_{\mbox{\tiny micro}}+2 \mu_{\mbox{\tiny micro}}} =
\widehat{K}_{e} + \widehat{K}_{\mbox{\tiny micro}} \, ,
\label{eq:equi_equa_2_RM_full_const_2}
\\
b_1 - b_2 &= 2 \mu_e + 2 \mu_{\mbox{\tiny micro}} \, .
\notag
\end{align}
where $\widehat{K}_{e}$ and $\widehat{K}_{\mbox{\tiny micro}}$ are the plane stress bulk moduli expressions at the micro-and meso-scale, respectively.
Finally, the solution of eq.(\ref{eq:equi_equa_2_RM_full}) is 
\begin{align}
P_{11}(x_2) = \, &
\frac{b_{0}}{b_{1}+b_{2}} c_{0}
- \frac{b_{1} b_{3}-b_{2} b_{4}}{b_{1}^2-b_{2}^2} \, \boldsymbol{\kappa}  \, x_{2}
+ \frac{c_1+c_3}{2} \cosh \left(\frac{f_{1} x_{2}}{L_c}\right)
+ \frac{c_2+c_4}{2}\frac{L_c}{f_{1}}\sinh \left(\frac{f_{1} x_{2}}{L_c}\right)
\notag
\\
&
+ \frac{c_1-c_3}{2} \cosh \left(\frac{f_{2} x_{2}}{L_c}\right)
+ \frac{c_2-c_4}{2}\frac{L_c}{f_{2}}\sinh \left(\frac{f_{2} x_{2}}{L_c}\right)
\, ,
\notag
\\[1.5mm]
P_{33}(x_2) = \, &
\frac{b_{0}}{b_{1}+b_{2}} c_{0}
- \frac{b_{1} b_{4}-b_{2} b_{3}}{b_{1}^2-b_{2}^2} \, \boldsymbol{\kappa}  \, x_{2}
+ \frac{c_1+c_3}{2} \cosh \left(\frac{f_{1} x_{2}}{L_c}\right)
+ \frac{c_2+c_4}{2}\frac{L_c}{f_{1}}\sinh \left(\frac{f_{1} x_{2}}{L_c}\right)
\label{eq:sol_fun_disp_RM_full}
\\
&
- \frac{c_1-c_3}{2} \cosh \left(\frac{f_{2} x_{2}}{L_c}\right)
- \frac{c_2-c_4}{2}\frac{L_c}{f_{2}}\sinh \left(\frac{f_{2} x_{2}}{L_c}\right)
\, ,
\notag
\\[1.5mm]
f_{1} :=\, &  \sqrt{\frac{b_{1} + b_{2}}{\mu}}\, ,
\qquad\qquad
f_{2} := \,   \sqrt{\frac{b_{1} - b_{2}}{\mu}}\,  .
\notag
\end{align}
Given  boundary conditions eq.(\ref{eq:BC_RM_gen}) for this case, the integration constants reduce to
\begin{align}
c_2 &= 
\frac{1}{2} \boldsymbol{\kappa}  \left(\left(\frac{b_{3}+b_{4}}{b_{1}+b_{2}}-1\right) \text{sech}\left(\frac{f_{1} h}{2 L_c}\right)+\left(\frac{b_{3}-b_{4}}{b_{1}-b_{2}}-1\right) \text{sech}\left(\frac{f_{2} h}{2 L_c}\right)\right)
\, ,
\qquad
c_0 = 0 \, , 
\label{eq:BC_RM_full}
\\
c_4 &= 
\frac{1}{2} \boldsymbol{\kappa}  \left(\left(\frac{b_{3}+b_{4}}{b_{1}+b_{2}}-1\right) \text{sech}\left(\frac{f_{1} h}{2 L_c}\right)-\left(\frac{b_{3}-b_{4}}{b_{1}-b_{2}}-1\right) \text{sech}\left(\frac{f_{2} h}{2 L_c}\right)\right)
\, , 
\qquad
c_1 = 0 \, , 
\qquad
c_3 = 0 \, .
\notag
\end{align}
The classical bending moment, the higher-order bending moment, and energy (per unit area d$x_1$d$x_3$) expressions are

\begin{align}
M_{\mbox{c}} (\boldsymbol{\kappa})&:=
\displaystyle\int\limits_{-h/2}^{h/2}
\langle \boldsymbol{\widetilde{\sigma}} \boldsymbol{e}_1 , \boldsymbol{e}_1 \rangle  \, x_2 \, 
\mbox{d}x_{2}
= 
\frac{h^3}{12}
\left( 
p_{1} 
+ p_{2} \, \left(\frac{L_c}{h}\right)^2
+ p_{3} \, \left(\frac{L_c}{h}\right)^3 \tanh \left(\frac{f_{1} h}{2 L_c}\right)
\right.
\notag
\\
&
\pushright{\left.
+ p_{4} \, \left(\frac{L_c}{h}\right)^3 \tanh \left(\frac{f_{2} h}{2 L_c}\right)
\right)
\boldsymbol{\kappa}}
\, ,
\notag
\\
M_{\mbox{m}}(\boldsymbol{\kappa}) &:=
\displaystyle\int\limits_{-h/2}^{h/2}
\langle \left(\boldsymbol{m} \times \boldsymbol{e}_1\right) \boldsymbol{e}_2 , \boldsymbol{e}_1 \rangle
\, \mbox{d}x_{2}
= 
\frac{h^3}{12} 
\left(
q_{1} \, \left(\frac{L_c}{h}\right)^2
-q_{2} \, \left(\frac{L_c}{h}\right)^3 \tanh \left(\frac{f_{1} h}{2 L_c}\right)
\right.
\notag
\\
&
\pushright{\left.
-q_{3} \, \left(\frac{L_c}{h}\right)^3 \tanh \left(\frac{f_{2} h}{2 L_c}\right)
\right)
\boldsymbol{\kappa}}
\, ,
\label{eq:sigm_ene_dimensionless_RM_Full}
\\
W_{\mbox{tot}} (\boldsymbol{\kappa})&:=
\displaystyle\int\limits_{-h/2}^{+h/2} W \left(\boldsymbol{\mbox{D}u}, \boldsymbol{P}, \mbox{Curl}\boldsymbol{P}\right) \, \mbox{d}x_{2}
= 
\frac{1}{2} 
\frac{h^3}{12} 
\left(
p_{1} 
+ \left(p_{2} + q_{1}\right) \, \left(\frac{L_c}{h}\right)^2
\right.
\notag
\\
& \hspace{6.7cm}
+ \left(p_{3} - q_{2}\right) \, \left(\frac{L_c}{h}\right)^3 \tanh \left(\frac{f_{1} h}{2 L_c}\right)
\notag
\\
& \hspace{7.7cm}
\left.
+ \left(p_{4} - q_{3}\right) \, \left(\frac{L_c}{h}\right)^3 \tanh \left(\frac{f_{2} h}{2 L_c}\right)
\right)
\boldsymbol{\kappa}^2
\, ,
\notag
\end{align}
\begin{align}
p_1&:= 
\frac{2 \mu _e}{\lambda _e+2 \mu _e}
\left(
\lambda _e \frac{(b_{1} (2 b_{1}-2 b_{3}-b_{4})-b_{2} (2 b_{2}-b_{3}-2 b_{4}))}{b_{1}^2-b_{2}^2}
+2 \mu _e \frac{(b_{1} (b_{1}-b_{3})-b_{2} (b_{2}-b_{4}))}{b_{1}^2-b_{2}^2}
\right) \, ,
\notag
\\
p_2&:= 
-\frac{12 \mu _e}{\left(\lambda _e+2 \mu _e\right)}
\left(
\left(3 \lambda _e+2 \mu _e\right) \frac{b_{1}+b_{2}-b_{3}-b_{4}}{(b_{1} + b_{2})^2}
+\left(\lambda _e+2 \mu _e\right)   \frac{b_{1}-b_{2}-b_{3}+b_{4}}{(b_{1} - b_{2})^2}
\right) \, \mu \, ,
\notag
\\
p_3&:= 
\frac{12}{f_{1}^3} \left(1-\frac{b_{3}+b_{4}}{b_{1}+b_{2}}\right)
\frac{2\mu _e \left(3 \lambda _e+2 \mu _e\right)}{\left(\lambda _e+2 \mu _e\right)} \, ,
\qquad
p_4:=
\frac{12}{f_{2}^3}
\left(1-\frac{b_{3}-b_{4}}{b_{1}-b_{2}}\right)
2 \mu _e \, ,
\\
q_1&:= 
12 \frac{b_{1} (b_{1}-b_{3})+b_{2} (b_{4}-b_{2})}{b_{1}^2-b_{2}^2} \mu \, ,
\qquad\qquad\quad\;\,
q_2:=
12 \frac{b_{1}+b_{2}-b_{3}-b_{4}}{f_{1}^3} \, ,
\notag
\\
q_3&:=
12 \frac{b_{1}-b_{2}-b_{3}+b_{4}}{f_{2}^3} \, .
\notag
\end{align}

\begin{figure}[H]
\centering
\includegraphics[width=0.5\linewidth]{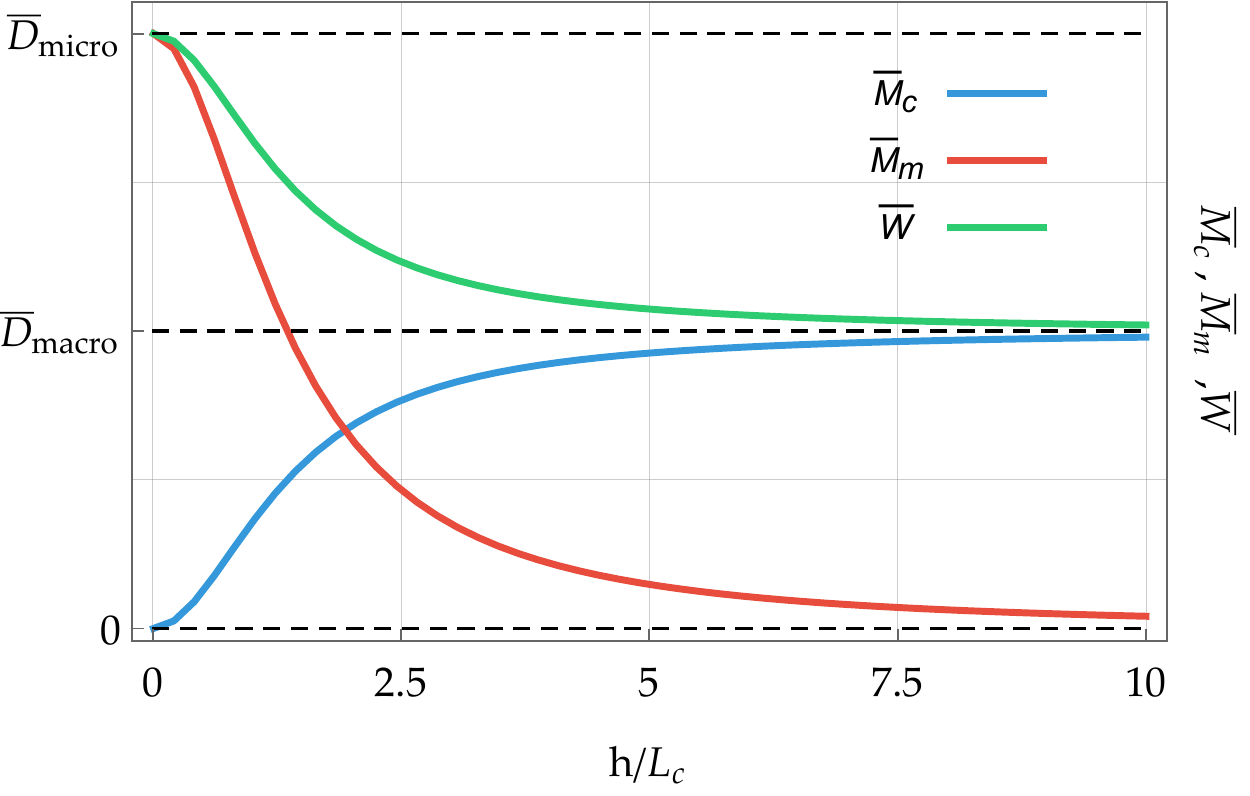}
\caption{(\textbf{Relaxed micromorphic model}, one curvature parameter, arbitrary Poisson) Bending moments and energy while varying $L_c$. Observe that the bending stiffness remains bounded as $L_c \to \infty$ ($h\to 0$). This is a distinguishing feature of the relaxed micromorphic model. The values of the parameters used are: $\mu _e = 1$, $\lambda _e = 1$, $\mu_{\mbox{\tiny micro}} = 1$, $\lambda_{\mbox{\tiny micro}} = 1$, $\mu = 1$.}
\label{fig:all_plot_RM_Full}
\end{figure}

Again, 
$
\frac{\mbox{d}}{\mbox{d}\boldsymbol{\kappa}}W_{\mbox{tot}}(\boldsymbol{\kappa}) = M_{\mbox{c}} (\boldsymbol{\kappa}) + M_{\mbox{m}} (\boldsymbol{\kappa}) \, .
$
The plot of the bending moments and the strain energy divided by $\frac{h^3}{12}\boldsymbol{\kappa}$ and $\frac{1}{2}\frac{h^3}{12}\boldsymbol{\kappa}^2$, respectively, while changing $L_c$ is shown in Fig.~\ref{fig:all_plot_RM_Full}.
\subsubsection{Limit cases}
\begin{align}
\displaystyle\lim_{L_c\to\infty} M_{\mbox{c}} (\boldsymbol{\kappa})
&= 0 \, ,
\qquad
\displaystyle\lim_{L_c\to0} M_{\mbox{m}} (\boldsymbol{\kappa})
= 0 \, ,
\notag
\\*
\displaystyle\lim_{L_c\to0} M_{\mbox{c}} (\boldsymbol{\kappa})
&= \frac{h^3}{12} \, \frac{4 \, \mu_{\mbox{\tiny macro}} \, \left(\lambda_{\mbox{\tiny macro}} + \mu_{\mbox{\tiny macro}}\right) }{\lambda_{\mbox{\tiny macro}} + 2\mu_{\mbox{\tiny macro}}} \, \boldsymbol{\kappa} 
=
\frac{h^3}{12} \frac{4 \mu _{\mbox{\tiny macro}} \left(3 \kappa _{\mbox{\tiny macro}}+\mu _{\mbox{\tiny macro}}\right)}{3 \kappa _{\mbox{\tiny macro}}+4 \mu _{\mbox{\tiny macro}}}
\, \boldsymbol{\kappa}
= D_{\tiny \mbox{macro}} \, \boldsymbol{\kappa} \, , 
\notag
\\*
\displaystyle\lim_{L_c\to\infty} M_{\mbox{m}} (\boldsymbol{\kappa})
&= \frac{h^3}{12} \, \frac{4 \, \mu_{\mbox{\tiny micro}} \, \left(\lambda_{\mbox{\tiny micro}} + \mu_{\mbox{\tiny micro}}\right) }{\lambda_{\mbox{\tiny micro}} + 2\mu_{\mbox{\tiny micro}}} \, \boldsymbol{\kappa} 
=
\frac{h^3}{12} \frac{4 \mu _{\mbox{\tiny micro}} \left(3 \kappa _{\mbox{\tiny micro}}+\mu _{\mbox{\tiny micro}}\right)}{3 \kappa _{\mbox{\tiny micro}}+4 \mu _{\mbox{\tiny micro}}}
\, \boldsymbol{\kappa}
= D_{\tiny \mbox{micro}} \, \boldsymbol{\kappa} \, , 
\\*
\displaystyle\lim_{L_c\to0} W_{\mbox{tot}} (\boldsymbol{\kappa})
&= 
\frac{1}{2} \frac{h^3}{12} \, \frac{4 \, \mu_{\mbox{\tiny macro}} \, \left(\lambda_{\mbox{\tiny macro}} + \mu_{\mbox{\tiny macro}}\right) }{\lambda_{\mbox{\tiny macro}} + 2\mu_{\mbox{\tiny macro}}} \, \boldsymbol{\kappa}^2
=
\frac{1}{2} \frac{h^3}{12} \frac{4 \mu _{\mbox{\tiny macro}} \left(3 \kappa _{\mbox{\tiny macro}}+\mu _{\mbox{\tiny macro}}\right)}{3 \kappa _{\mbox{\tiny macro}}+4 \mu _{\mbox{\tiny macro}}}  \, \boldsymbol{\kappa}^2 
=
\frac{1}{2} \, D_{\tiny \mbox{macro}} \, \boldsymbol{\kappa}^2 \, , 
\notag
\\*
\displaystyle\lim_{L_c\to\infty} W_{\mbox{tot}} (\boldsymbol{\kappa})
&=
\frac{1}{2} \frac{h^3}{12} \, \frac{4 \, \mu_{\mbox{\tiny micro}} \, \left(\lambda_{\mbox{\tiny micro}} + \mu_{\mbox{\tiny micro}}\right) }{\lambda_{\mbox{\tiny micro}} + 2\mu_{\mbox{\tiny micro}}} \, \boldsymbol{\kappa}^2 
=
\frac{1}{2} \frac{h^3}{12} \frac{4 \mu _{\mbox{\tiny micro}} \left(3 \kappa _{\mbox{\tiny micro}}+\mu _{\mbox{\tiny micro}}\right)}{3 \kappa _{\mbox{\tiny micro}}+4 \mu _{\mbox{\tiny micro}}}  \, \boldsymbol{\kappa}^2 
= 
\frac{1}{2} \, D_{\tiny \mbox{micro}} \, \boldsymbol{\kappa}^2 \, .
\notag
\end{align}
\subsection{Full isotropic curvature and zero Poisson's ratios $\nu_{\tiny \mbox{micro}}=\nu_{e}=0$}
The expression of the strain energy for the isotropic relaxed micromorphic continuum is
\begin{align}
W \left(\boldsymbol{\mbox{D}u}, \boldsymbol{P},\mbox{Curl}\,\boldsymbol{P}\right) = &
\, \mu_{e} \left\lVert \mbox{sym} \left(\boldsymbol{\mbox{D}u} - \boldsymbol{P} \right) \right\rVert^{2}
+ \frac{\lambda_{e}}{2} \mbox{tr}^2 \left(\boldsymbol{\mbox{D}u} - \boldsymbol{P} \right) 
+ \mu_{c} \left\lVert \mbox{skew} \left(\boldsymbol{\mbox{D}u} - \boldsymbol{P} \right) \right\rVert^{2}
\notag
\\
&
+ \mu_{\tiny \mbox{micro}} \left\lVert \mbox{sym}\,\boldsymbol{P} \right\rVert^{2}
+ \frac{\lambda_{\tiny \mbox{micro}}}{2} \mbox{tr}^2 \left(\boldsymbol{P} \right)
\label{eq:energy_RM_3_L0}
\\
&
+ \frac{\mu \,L_c^2 }{2} \,
\left(
a_1 \, \left\lVert \mbox{dev sym} \, \mbox{Curl} \, \boldsymbol{P}\right\rVert^2 +
a_2 \, \left\lVert \mbox{skew} \,  \mbox{Curl} \, \boldsymbol{P}\right\rVert^2 +
\frac{a_3}{3} \, \mbox{tr}^2 \left(\mbox{Curl} \, \boldsymbol{P}\right)
\right),
\notag
\end{align}
while the equilibrium equations without body forces are
\begin{align}
\mbox{Div}\overbrace{\left[2\mu_{e}\,\mbox{sym} \left(\boldsymbol{\mbox{D}u} - \boldsymbol{P} \right) + \lambda_{e} \mbox{tr} \left(\boldsymbol{\mbox{D}u} - \boldsymbol{P} \right) \boldsymbol{\mathbbm{1}}
+ 2\mu_{c}\,\mbox{skew} \left(\boldsymbol{\mbox{D}u} - \boldsymbol{P} \right)\right]}^{\mathlarger{\widetilde{\sigma}}:=}
&= \boldsymbol{0},
\notag
\\
\widetilde{\sigma}
- 2 \mu_{\mbox{\tiny micro}}\,\mbox{sym}\,\boldsymbol{P} - \lambda_{\tiny \mbox{micro}} \mbox{tr} \left(\boldsymbol{P}\right) \boldsymbol{\mathbbm{1}}
\hspace{8cm}
\label{eq:equi_RM_3_L0}
\\
- \mu \, L_{c}^{2} \, \mbox{Curl}
\left(
a_1 \, \mbox{dev sym} \, \mbox{Curl} \, \boldsymbol{P} +
a_2 \, \mbox{skew} \,  \mbox{Curl} \, \boldsymbol{P} +
\frac{a_3}{3} \, \mbox{tr} \left(\mbox{Curl} \, \boldsymbol{P}\right)\mathbbm{1}
\right) &= \boldsymbol{0}.
\notag
\end{align}
Substituting the ansatz eq.(\ref{eq:ansatz_RM}) in eq.(\ref{eq:equi_RM_3_L0}) the equilibrium equations result to be
\begin{align}
2 \mu _c (\kappa_{1}-\kappa_{2})+2 \mu _e \left(v''(x_{2})-P_{22}'(x_{2})\right) = 0 \, ,
\notag
\\
\mu \, L_c^2 \left(\left(a_{1}+a_{2}\right) P_{11}''(x_{2})+\left(a_{2}-a_{1}\right) P_{33}''(x_{2})\right)-4 \mu _e (P_{11}(x_{2})+\kappa_{1} x_{2})-4 \mu_{\mbox{\tiny micro}} P_{11}(x_{2})  = 0 \, ,
\notag
\\
2 x_{1} \mu _c (\kappa_{2}-\kappa_{1})  = 0 \, ,
\label{eq:equi_equa_RM_3_L0}
\\
2 x_{1} \mu _c (\kappa_{1}-\kappa_{2})  = 0 \, ,
\notag
\\
2 \mu _e v'(x_{2})-2 P_{22}(x_{2}) \left(\mu _e+\mu_{\mbox{\tiny micro}}\right)  = 0 \, ,
\notag
\\
\mu  L_c^2 \left((a_{2}-a_{1}) P_{11}''(x_{2})+(a_{1}+a_{2}) P_{33}''(x_{2})\right)-4 P_{33}(x_{2}) \left(\mu _e+\mu_{\mbox{\tiny micro}}\right)  = 0 \, .
\notag
\end{align}
It is clear that in order to satisfy eq.(\ref{eq:equi_equa_RM_3_L0})$_3$ and eq.(\ref{eq:equi_equa_RM_3_L0})$_4$ either $\mu _c = 0$ or $\kappa_1 = \kappa_2 = \boldsymbol{\kappa}$;
it has been chosen the latter option, which implies that the skew-symmetric part of the gradient of the displacement eq.(\ref{eq:grad_RM}) is the same as the skew-symmetric part of the microdistortion eq.(\ref{eq:ansatz_RM})$_2$.
This also implies that the Cosserat couple modulus $\mu_c$ does not play a role any more.

\begin{figure}[H]
\centering
\includegraphics[width=0.5\linewidth]{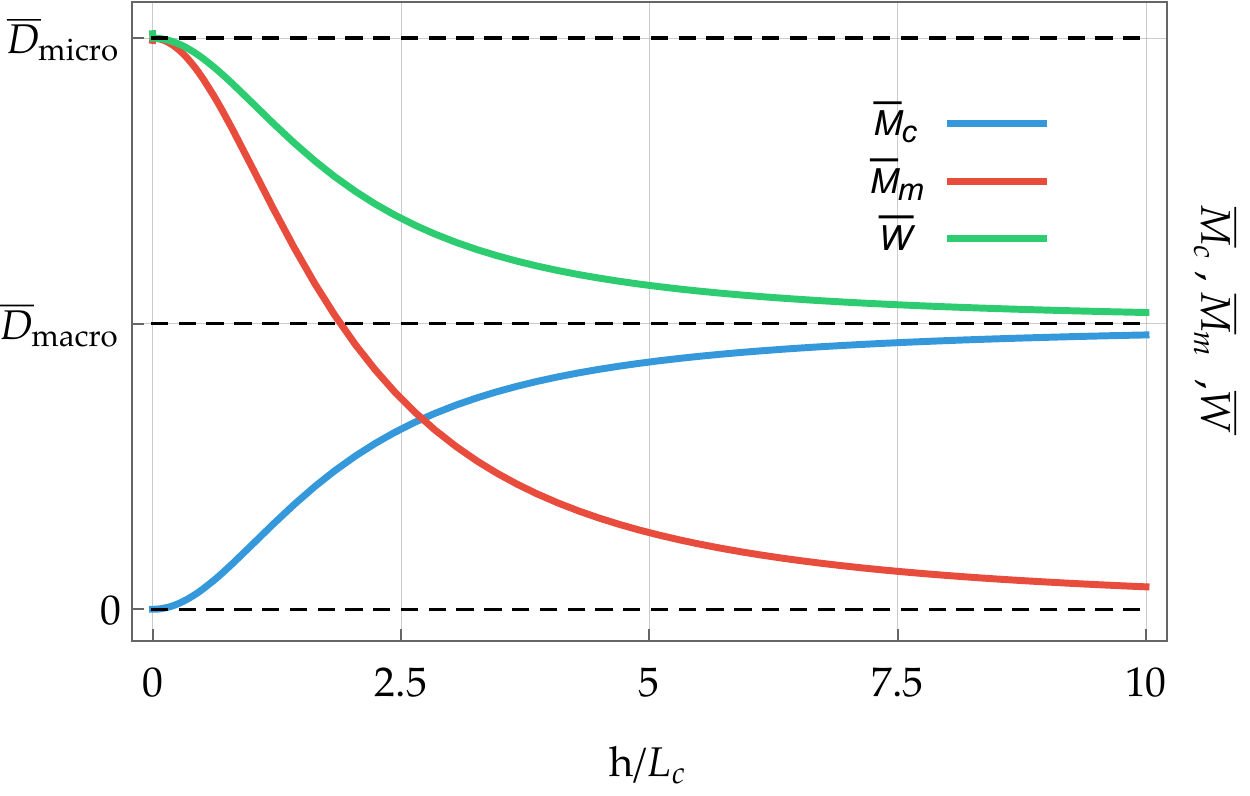}
\caption{(\textbf{Relaxed micromorphic model}, zero Poisson) Bending moments and energy while varying $L_c$. Observe that the bending stiffness remains bounded as $L_c \to \infty$ ($h\to 0$). This is a distinguishing feature of the relaxed micromorphic model. The values of the parameters used are: $\mu _e = 1$, $\mu_{\mbox{\tiny micro}} = 1$, $\mu = 1$, $a _1 = 2$, $a _2 = 1$.}
\label{fig:all_plot_RM_3_L0}
\end{figure}

From eq.(\ref{eq:equi_equa_RM_3_L0})$_1$ it is possible to evaluate $v''(x_2)$ and consequently $v'(x_2)$
\begin{equation}
v''(x_2) = 
P_{22}'(x_{2}) \, ,
\qquad
v'(x_2) = 
P_{22} (x_{2}) + c_0 \, .
\label{eq:vv2_RM_3_L0}
\end{equation}
By substituting back eq.(\ref{eq:vv2_RM_3_L0}) in eq.(\ref{eq:equi_equa_RM_3_L0}), we can evaluate $P_{22}(x_2)$ and its derivatives from equation eq.(\ref{eq:equi_equa_RM_3_L0})$_5$
\begin{equation}
P_{22}(x_2) = 
\frac{\mu_{e}}{\mu_{\mbox{\tiny micro}}}c_{0} \, ,
\qquad
P'_{22}(x_2) = 0 \, ,
\qquad
P''_{22}(x_2) = 0 \, .
\label{eq:p22_RM_3_L0}
\end{equation}
After substituting eq.(\ref{eq:p22_RM_3_L0}) in eq.(\ref{eq:equi_equa_RM_3_L0}) the following two coupled second order ordinary differential equations in $P_{11}(x_2)$ and $P_{33}(x_2)$ are obtained
\begin{align}
\mu \, L_c^2 \left((a_{1}+a_{2}) P_{11}''(x_{2})-(a_{1}-a_{2}) P_{33}''(x_{2})\right)-4 \mu _e (P_{11}(x_{2})+\boldsymbol{\kappa} \, x_{2})-4 \mu_{\mbox{\tiny micro}} P_{11}(x_{2}) = 0 \, ,
\label{eq:equi_equa_2_RM_3_L0}
\\
\mu \, L_c^2 \left((a_{1}+a_{2}) P_{33}''(x_{2})-(a_{1}-a_{2}) P_{11}''(x_{2})\right)-4 P_{33}(x_{2}) \left(\mu _e+\mu_{\mbox{\tiny micro}}\right) = 0 \, .
\notag
\end{align}
Finally, the solution of eq.(\ref{eq:equi_equa_2_RM_3_L0}) is 
\begin{align}
P_{11}(x_2) = \, & 
\frac{c_1-c_3}{2} \cosh \left(\frac{f_{1} x_{2}}{L_c}\right)
+\frac{c_1+c_3}{2} \cosh \left(\frac{f_{2} x_{2}}{L_c}\right)
\notag
\\
&+\frac{c_2-c_4}{2}\frac{L_c}{f_{1}} \sinh \left(\frac{f_{1} x_{2}}{L_c}\right)
+\frac{c_2+c_4}{2}\frac{L_c}{f_{2}} \sinh \left(\frac{f_{2} x_{2}}{L_c}\right)
-\frac{\mu _e}{\mu _e+\mu_{\mbox{\tiny micro}}} \, \boldsymbol{\kappa} \, x_{2}
\, ,
\notag
\\
P_{33}(x_2) = \, & 
\frac{c_3-c_1}{2} \cosh \left(\frac{f_{1} x_{2}}{L_c}\right)
+\frac{c_1+c_3}{2} \cosh \left(\frac{f_{2} x_{2}}{L_c}\right)
\label{eq:sol_fun_disp_RM_3_L0}
\\
&+\frac{c_4-c_2}{2}\frac{L_c}{f_{1}} \sinh \left(\frac{f_{1} x_{2}}{L_c}\right)
+\frac{c_2+c_4}{2}\frac{L_c}{f_{2}} \sinh \left(\frac{f_{2} x_{2}}{L_c}\right)
\, ,
\notag
\\
f_1:= \, &\sqrt{\frac{2\left( \mu _e + \mu_{\mbox{\tiny micro}}\right)}{a_1 \, \mu}}\, ,
\qquad\qquad
f_2:= \, \sqrt{\frac{2\left( \mu _e + \mu_{\mbox{\tiny micro}}\right)}{a_2 \, \mu}}\, .
\notag
\end{align}
Given  boundary conditions eq.(\ref{eq:BC_RM_gen}) for this case, the integration constants reduce to
\begin{align}
c_2 &= 
-\frac{\boldsymbol{\kappa} \, \mu_{\mbox{\tiny micro}}}{2 \left(\mu _e + \mu_{\mbox{\tiny micro}}\right)}
\left(\text{sech}\left(\frac{f_{1} h}{2 L_c}\right) + \text{sech}\left(\frac{f_{2} h}{2 L_c}\right)\right)
\, ,
\qquad\quad
c_0 = 0 \, .
\label{eq:BC_RM_3_L0}
\\
c_4 &= 
\frac{\boldsymbol{\kappa} \, \mu_{\mbox{\tiny micro}}}{2 \left(\mu _e + \mu_{\mbox{\tiny micro}}\right)}
\left(\text{sech}\left(\frac{f_{1} h}{2 L_c}\right) - \text{sech}\left(\frac{f_{2} h}{2 L_c}\right)\right)
\, ,
\qquad\quad\quad
c_1 = 0 \, ,
\qquad\quad
c_3 = 0 \, .
\notag
\end{align}
The classical bending moment, the higher-order bending moment, and energy (per unit area d$x_1$d$x_3$) expressions are

\begin{align}
M_{\mbox{c}} (\boldsymbol{\kappa})&:=
\displaystyle\int\limits_{-h/2}^{h/2}
\langle \boldsymbol{\widetilde{\sigma}} \boldsymbol{e}_1 , \boldsymbol{e}_1 \rangle  \, x_2 \, 
\mbox{d}x_{2}
= 
\frac{h^3}{12}
\frac{2 \mu _e \, \mu_{\mbox{\tiny micro}} }{\mu _e+\mu_{\mbox{\tiny micro}}} 
\left(
1
-6 \, \frac{f_{1}^2+f_{2}^2}{f_{1}^2 f_{2}^2}\left(\frac{L_c}{h}\right)^2
\right.
\notag
\\
& \hspace{2cm}
\left.
+\frac{12}{f_{1}^3 f_{2}^3} \left(\frac{L_c}{h}\right)^3 \left(f_{1}^3 \tanh \left(\frac{f_{2} h}{2 L_c}\right)+f_{2}^3 \tanh \left(\frac{f_{1} h}{2 L_c}\right)\right)
\right)
\boldsymbol{\kappa}
\, ,
\notag
\\
M_{\mbox{m}}(\boldsymbol{\kappa}) &:=
\displaystyle\int\limits_{-h/2}^{h/2}
\langle \left(\boldsymbol{m} \times \boldsymbol{e}_1\right) \boldsymbol{e}_2 , \boldsymbol{e}_1 \rangle
\, \mbox{d}x_{2}
= 
\frac{h^3}{12} \, 
2\mu_{\mbox{\tiny micro}}
\left(
3 \, \frac{a_{1} + a_{2}}{\mu _e + \mu_{\mbox{\tiny micro}}} \mu \left(\frac{L_c}{h}\right)^2
\right.
\notag
\\
& \hspace{2cm}
\left.
-6\frac{a_{1} a_{2}}{\mu _e + \mu_{\mbox{\tiny micro}}} \, \mu \, \left(\frac{L_c}{h}\right)^3
\left(
\frac{1}{a_{1} f_{2}} \tanh \left(\frac{f_{2} h}{2 L_c}\right)
+ \frac{1}{a_{2} f_{1}}\tanh \left(\frac{f_{1} h}{2 L_c}\right)
\right)
\right)
\boldsymbol{\kappa}
\, ,
\label{eq:sigm_ene_dimensionless_RM_3_L0}
\\
W_{\mbox{tot}} (\boldsymbol{\kappa})&:=
\displaystyle\int\limits_{-h/2}^{+h/2} W \left(\boldsymbol{\mbox{D}u}, \boldsymbol{P}, \mbox{Curl}\boldsymbol{P}\right) \, \mbox{d}x_{2}
= 
\frac{1}{2}
\frac{h^3}{12}
\frac{2 \mu _e \, \mu_{\mbox{\tiny micro}} }{\mu _e+\mu_{\mbox{\tiny micro}}} 
\left(
1
+6 \, \frac{\mu_{\mbox{\tiny micro}}}{\mu_e}\frac{f_{1}^2+f_{2}^2}{f_{1}^2 f_{2}^2}\left(\frac{L_c}{h}\right)^2
\right.
\notag
\\
& \hspace{2cm}
\left.
-12 \, \frac{\mu_{\mbox{\tiny micro}}}{\mu_e}\frac{1}{f_{1}^3 f_{2}^3} \left(\frac{L_c}{h}\right)^3 \left(f_{1}^3 \tanh \left(\frac{f_{2} h}{2 L_c}\right)+f_{2}^3 \tanh \left(\frac{f_{1} h}{2 L_c}\right)\right)
\right)
\boldsymbol{\kappa}^2
 \, .
\notag
\end{align}
Again, 
$
\frac{\mbox{d}}{\mbox{d}\boldsymbol{\kappa}}W_{\mbox{tot}}(\boldsymbol{\kappa}) = M_{\mbox{c}} (\boldsymbol{\kappa}) + M_{\mbox{m}} (\boldsymbol{\kappa}) \, .
$
The plot of the bending moments and the strain energy divided by $\frac{h^3}{12}\boldsymbol{\kappa}$ and $\frac{1}{2}\frac{h^3}{12}\boldsymbol{\kappa}^2$, respectively, while changing $L_c$ is shown in Fig.~\ref{fig:all_plot_RM_3_L0}.
\subsubsection{Limit cases}
\begin{align}
\displaystyle\lim_{L_c\to0} M_{\mbox{c}} (\boldsymbol{\kappa})
&= \frac{2 \, \mu _e \, \mu _{\tiny \mbox{micro}} }{\mu _e+\mu _{\tiny \mbox{micro}}} \, \frac{h^3}{12} \, \boldsymbol{\kappa} 
= D_{\tiny \mbox{macro}} \, \boldsymbol{\kappa} \, , 
\qquad\qquad
\displaystyle\lim_{L_c\to\infty} M_{\mbox{c}} (\boldsymbol{\kappa})
= 0 \, ,
\notag
\\
\displaystyle\lim_{L_c\to0} M_{\mbox{m}} (\boldsymbol{\kappa})
&= 0 \, ,
\qquad\qquad\qquad\qquad\qquad\qquad\qquad\quad\,
\displaystyle\lim_{L_c\to\infty} M_{\mbox{m}} (\boldsymbol{\kappa})
= 2 \, \mu _{\tiny \mbox{micro}} \, \frac{h^3}{12} \, \boldsymbol{\kappa} 
= D_{\tiny \mbox{micro}} \, \boldsymbol{\kappa} \, , 
\\
\displaystyle\lim_{L_c\to0} W_{\mbox{tot}} (\boldsymbol{\kappa})
&= 
\frac{1}{2}\frac{2 \, \mu _e \, \mu_{\tiny \mbox{micro}} }{\mu _e + \mu_{\tiny \mbox{micro}}} \, \frac{h^3}{12} \, \boldsymbol{\kappa}^2 
= 
\frac{1}{2} \,D_{\tiny \mbox{macro}} \, \boldsymbol{\kappa}^2 \, , 
\notag
\\
\displaystyle\lim_{L_c\to\infty} W_{\mbox{tot}} (\boldsymbol{\kappa})
&=
\frac{1}{2} 2\mu_{\tiny \mbox{micro}} \, \frac{h^3}{12} \, \boldsymbol{\kappa}^2 
= 
\frac{1}{2} \,D_{\tiny \mbox{micro}} \, \boldsymbol{\kappa}^2 \, .
\notag
\end{align}
\subsection{Full isotropic curvature and arbitrary Poisson's ratios $\nu_{\tiny \mbox{micro}}$ and $\nu_{e}$}
Finally, we are prepared to treat the most general case of an isotropic, linear-elastic relaxed micromorphic continuum..
The expression of the strain energy for the isotropic relaxed micromorphic continuum is:
\begin{align}
W \left(\boldsymbol{\mbox{D}u}, \boldsymbol{P},\mbox{Curl}\,\boldsymbol{P}\right) = &
\, \mu_{e} \left\lVert \mbox{sym} \left(\boldsymbol{\mbox{D}u} - \boldsymbol{P} \right) \right\rVert^{2}
+ \frac{\lambda_{e}}{2} \mbox{tr}^2 \left(\boldsymbol{\mbox{D}u} - \boldsymbol{P} \right) 
+ \mu_{c} \left\lVert \mbox{skew} \left(\boldsymbol{\mbox{D}u} - \boldsymbol{P} \right) \right\rVert^{2}
\notag
\\
&
+ \mu_{\tiny \mbox{micro}} \left\lVert \mbox{sym}\,\boldsymbol{P} \right\rVert^{2}
+ \frac{\lambda_{\tiny \mbox{micro}}}{2} \mbox{tr}^2 \left(\boldsymbol{P} \right)
\label{eq:energy_RM_3}
\\
&
+ \frac{\mu \,L_c^2 }{2} \,
\left(
a_1 \, \left\lVert \mbox{dev sym} \, \mbox{Curl} \, \boldsymbol{P}\right\rVert^2 +
a_2 \, \left\lVert \mbox{skew} \,  \mbox{Curl} \, \boldsymbol{P}\right\rVert^2 +
\frac{a_3}{3} \, \mbox{tr}^2 \left(\mbox{Curl} \, \boldsymbol{P}\right)
\right),
\notag
\end{align}
while the equilibrium equations without body forces are
\begin{align}
\mbox{Div}\overbrace{\left[2\mu_{e}\,\mbox{sym} \left(\boldsymbol{\mbox{D}u} - \boldsymbol{P} \right) + \lambda_{e} \mbox{tr} \left(\boldsymbol{\mbox{D}u} - \boldsymbol{P} \right) \boldsymbol{\mathbbm{1}}
+ 2\mu_{c}\,\mbox{skew} \left(\boldsymbol{\mbox{D}u} - \boldsymbol{P} \right)\right]}^{\mathlarger{\widetilde{\sigma}}:=}
&= \boldsymbol{0},
\notag
\\
\widetilde{\sigma}
- 2 \mu_{\mbox{\tiny micro}}\,\mbox{sym}\,\boldsymbol{P} - \lambda_{\tiny \mbox{micro}} \mbox{tr} \left(\boldsymbol{P}\right) \boldsymbol{\mathbbm{1}}
\hspace{8cm}
\label{eq:equi_RM_3}
\\
- \mu \, L_{c}^{2} \, \mbox{Curl}
\left(
a_1 \, \mbox{dev sym} \, \mbox{Curl} \, \boldsymbol{P} +
a_2 \, \mbox{skew} \,  \mbox{Curl} \, \boldsymbol{P} +
\frac{a_3}{3} \, \mbox{tr} \left(\mbox{Curl} \, \boldsymbol{P}\right)\mathbbm{1}
\right) &= \boldsymbol{0}.
\notag
\end{align}
Substituting the ansatz eq.(\ref{eq:ansatz_RM}) in eq.(\ref{eq:equi_RM_3}) the equilibrium equations are
\begin{align}
2 \mu _c (\kappa_{1} - \kappa_{2}) - \lambda _e \left(\kappa_{1} + P_{11}'(x_{2}) + P_{22}'(x_{2}) + P_{33}'(x_{2}) - v''(x_{2})\right) + 2 \mu _e \left(v''(x_{2}) - P_{22}'(x_{2})\right) &= 0
\, ,
\notag
\\*
\frac{1}{2} \mu  L_c^2 \left((a_{1} + a_{2}) P_{11}''(x_{2}) + (a_{2} - a_{1}) P_{33}''(x_{2})\right) 
\hspace{7cm}
\notag
\\*
- \lambda _e \left(P_{11}(x_{2}) + P_{22}(x_{2}) + P_{33}(x_{2}) - v'(x_{2}) + \kappa_{1} x_{2}\right)
\notag
\\*
- 2 \mu _e (P_{11}(x_{2}) + \kappa_{1} x_{2}) - P_{11}(x_{2}) \left(\lambda_{\mbox{\tiny micro}} + 2 \mu_{\mbox{\tiny micro}}\right) + \lambda_{\mbox{\tiny micro}} ( - P_{22}(x_{2}) - P_{33}(x_{2})) &= 0
\, ,
\label{eq:equi_equa_RM_3}
\\*
2 x_{1} \mu _c (\kappa_{2} - \kappa_{1}) &= 0
\, ,
\notag
\\*
2 x_{1} \mu _c (\kappa_{1} - \kappa_{2}) &= 0
\, ,
\notag
\\*
- \lambda _e \left(P_{11}(x_{2}) + P_{22}(x_{2}) + P_{33}(x_{2}) - v'(x_{2}) + \kappa_{1} x_{2}\right) + 2 \mu _e \left(v'(x_{2}) - P_{22}(x_{2})\right) 
\notag
\\*
+ \lambda_{\mbox{\tiny micro}} ( - P_{11}(x_{2})) - P_{22}(x_{2}) \left(\lambda_{\mbox{\tiny micro}} + 2 \mu_{\mbox{\tiny micro}}\right) - \lambda_{\mbox{\tiny micro}} P_{33}(x_{2}) &= 0
\, ,
\notag
\\*
- \frac{1}{2} \mu  L_c^2 \left((a_{1} - a_{2}) P_{11}''(x_{2}) - (a_{1} + a_{2}) P_{33}''(x_{2})\right) 
\hspace{7cm}
\notag
\\*
- \lambda _e \left(P_{11}(x_{2}) + P_{22}(x_{2}) + P_{33}(x_{2}) - v'(x_{2}) + \kappa_{1} x_{2}\right)
\notag
\\*
- 2 \mu _e P_{33}(x_{2}) - \lambda_{\mbox{\tiny micro}} P_{11}(x_{2}) - \lambda_{\mbox{\tiny micro}} P_{22}(x_{2}) - P_{33}(x_{2}) \left(\lambda_{\mbox{\tiny micro}} + 2 \mu_{\mbox{\tiny micro}}\right) &= 0
\, .
\notag
\end{align}
In order to satisfy eq.(\ref{eq:equi_equa_RM_3})$_3$ and eq.(\ref{eq:equi_equa_RM_3})$_4$ either $\mu _c = 0$ or $\kappa_1 = \kappa_2 = \boldsymbol{\kappa}$;
we have chosen the latter option, which implies that the skew-symmetric part of the gradient of the displacement eq.(\ref{eq:grad_RM}) is the same as the skew-symmetric part of the microdistortion eq.(\ref{eq:ansatz_RM})$_2$.
This also implies that the Cosserat couple modulus $\mu_c$ does not play a role any more.

From eq.(\ref{eq:equi_equa_RM_3})$_1$ it is possible to evaluate $v''(x_2)$ and consequently $v'(x_2)$
\begin{align}
v''(x_2) &= 
\frac{\lambda _e \left(\boldsymbol{\kappa} +P_{11}'(x_{2})+P_{22}'(x_{2})+P_{33}'(x_{2})\right)+2 \mu _e P_{22}'(x_{2})}{\lambda _e+2 \mu _e} \, ,
\label{eq:vv2_RM_3}
\\
v'(x_2) &= 
\frac{\lambda _e \left(\boldsymbol{\kappa} +P_{11}(x_{2})+P_{22}(x_{2})+P_{33}(x_{2})\right)+2 \mu _e P_{22} (x_{2})}{\lambda _e+2 \mu _e} + c_0 \, .
\notag
\end{align}
By substituting back eq.(\ref{eq:vv2_RM_3}) in eq.(\ref{eq:equi_equa_RM_3}), we can evaluate $P_{22}(x_2)$ and its derivatives from equation eq.(\ref{eq:equi_equa_RM_3})$_5$
\begin{align}
P_{22}(x_2) &= 
\frac{\left(\lambda _e+2 \mu _e\right)}{\lambda_{\mbox{\tiny micro}}+2 \mu_{\mbox{\tiny micro}}} c_{0}
- \frac{\lambda_{\mbox{\tiny micro}} (P_{11}(x_{2})+P_{33}(x_{2}))}{\lambda_{\mbox{\tiny micro}}+2 \mu_{\mbox{\tiny micro}}} \, ,
\label{eq:p22_RM_3}
\\
P'_{22}(x_2) &= 
- \frac{\lambda_{\mbox{\tiny micro}} (P_{11}'(x_{2})+P_{33}'(x_{2}))}{\lambda_{\mbox{\tiny micro}}+2 \mu_{\mbox{\tiny micro}}} \, ,
\qquad
P''_{22}(x_2) = 
- \frac{\lambda_{\mbox{\tiny micro}} (P_{11}''(x_{2})+P_{33}''(x_{2}))}{\lambda_{\mbox{\tiny micro}}+2 \mu_{\mbox{\tiny micro}}} \, .
\notag
\end{align}
After substituting eq.(\ref{eq:p22_RM_3}) in eq.(\ref{eq:equi_equa_RM_3}) the following two coupled second order ordinary differential equations in $P_{11}(x_2)$ and $P_{33}(x_2)$ are retrieved
\begin{align}
&\mu  L_c^2 \left((a_{1} + a_{2}) P_{11}''(x_{2}) + (a_{2} - a_{1}) P_{33}''(x_{2})\right) - 2 b_{1} P_{11}(x_{2}) - 2 b_{2} P_{33}(x_{2}) - 2 b_{3} x_{2} + 2 b_{0} c_{0} = 0
\, ,
\notag
\\
&\mu  L_c^2 \left((a_{2} - a_{1}) P_{11}''(x_{2}) + (a_{1} + a_{2}) P_{33}''(x_{2})\right) - 2 b_{1} P_{33}(x_{2}) - 2 b_{2} P_{11}(x_{2}) - 2 b_{4} x_{2} + 2 b_{0} c_{0} = 0
\, ,
\label{eq:equi_equa_2_RM_3}
\\
&b_{0} := \frac{2 \left(\lambda _e \mu_{\mbox{\tiny micro}}-\mu _e \lambda_{\mbox{\tiny micro}}\right)}{\lambda_{\mbox{\tiny micro}}+2 \mu_{\mbox{\tiny micro}}}
\, ,
\hspace{1.5cm}
b_{1} := \frac{4 \mu _e \left(\lambda _e+\mu _e\right)}{\lambda _e+2 \mu _e}+\frac{4 \mu_{\mbox{\tiny micro}} \left(\lambda_{\mbox{\tiny micro}}+\mu_{\mbox{\tiny micro}}\right)}{\lambda_{\mbox{\tiny micro}}+2 \mu_{\mbox{\tiny micro}}}
\, ,
\notag
\\
&b_{2} := \frac{2 \lambda _e \mu _e}{\lambda _e+2 \mu _e}+\frac{2 \lambda_{\mbox{\tiny micro}} \mu_{\mbox{\tiny micro}}}{\lambda_{\mbox{\tiny micro}}+2 \mu_{\mbox{\tiny micro}}}
\, ,
\qquad
b_{3} := \frac{4 \boldsymbol{\kappa}  \mu _e \left(\lambda _e+\mu _e\right)}{\lambda _e+2 \mu _e}
\, ,
\qquad \,\,\,\,\,\,\,\,
b_{4} := \frac{2 \boldsymbol{\kappa}  \lambda _e \mu _e}{\lambda _e+2 \mu _e} \, .
\notag
\end{align}
Finally the solution of eq.(\ref{eq:equi_equa_2_RM_3}) is 
\begin{align}
P_{11}(x_2) = \, & 
\frac{b_{0} c_{0}}{b_{1} + b_{2}}
+ \frac{(b_{2} b_{4} - b_{1} b_{3})x_{2}}{b_{1}^2 - b_{2}^2}
+ \frac{c_1 - c_3}{2} \cosh (\frac{f_{1} x_{2}}{L_c})
+ \frac{c_2 - c_4}{2}\frac{L_c}{f_{1}}\sinh (\frac{f_{1} x_{2}}{L_c})
\notag
\\
&+ \frac{c_1 + c_3}{2} \cosh (\frac{f_{2} x_{2}}{L_c})
+ \frac{c_2 + c_4}{2}\frac{L_c}{f_{2}} \sinh (\frac{f_{2} x_{2}}{L_c})
\, ,
\notag
\\
P_{33}(x_2) = \, & 
\frac{b_{0} c_{0}}{b_{1} + b_{2}}
+ \frac{(b_{2} b_{3} - b_{1} b_{4})x_{2}}{b_{1}^2 - b_{2}^2}
+ \frac{c_3 - c_1}{2} \cosh (\frac{f_{1} x_{2}}{L_c})
+ \frac{c_4 - c_2}{2}\frac{L_c}{f_{1}} \sinh (\frac{f_{1} x_{2}}{L_c})
\label{eq:sol_fun_disp_RM_3}
\\
&+ \frac{c_1 + c_3}{2} \cosh (\frac{f_{2} x_{2}}{L_c})
+ \frac{c_2 + c_4}{2}\frac{L_c}{f_{2}} \sinh (\frac{f_{2}x_{2}}{L_c})
\, ,
\notag
\\
f_1:= \, &\sqrt{\frac{b_{1} - b_{2}}{a_1 \, \mu}}\, ,
\qquad\qquad
f_2:= \, \sqrt{\frac{b_{1} + b_{2}}{a_2 \, \mu}}\, .
\notag
\end{align}
Given boundary conditions eq.(\ref{eq:BC_RM_gen}) for this case, the integration constants reduce to
\begin{align}
c_2 &= 
\frac{1}{2} \left(\left(\frac{b_{3}-b_{4}}{b_{1}-b_{2}}-\boldsymbol{\kappa} \right) \text{sech}\left(\frac{f_{1} h}{2 L_c}\right)+\left(\frac{b_{3}+b_{4}}{b_{1}+b_{2}}-\boldsymbol{\kappa} \right) \text{sech}\left(\frac{f_{2} h}{2 L_c}\right)\right)
\, ,
\quad \quad
c_0 = 0 \, , 
\label{eq:BC_RM_3}
\\
c_4 &= 
\frac{1}{2} \left(\left(\frac{b_{4}-b_{3}}{b_{1}-b_{2}}+\boldsymbol{\kappa} \right) \text{sech}\left(\frac{f_{1} h}{2 L_c}\right)+\left(\frac{b_{3}+b_{4}}{b_{1}+b_{2}}-\boldsymbol{\kappa} \right) \text{sech}\left(\frac{f_{2} h}{2 L_c}\right)\right)
\, , \quad \quad
c_1 = 0 \, , \quad \quad
c_3 = 0 \, .
\notag
\end{align}

\begin{figure}[H]
	\centering
	\begin{subfigure}{.45\textwidth}
		\centering
		\includegraphics[width=0.95\linewidth]{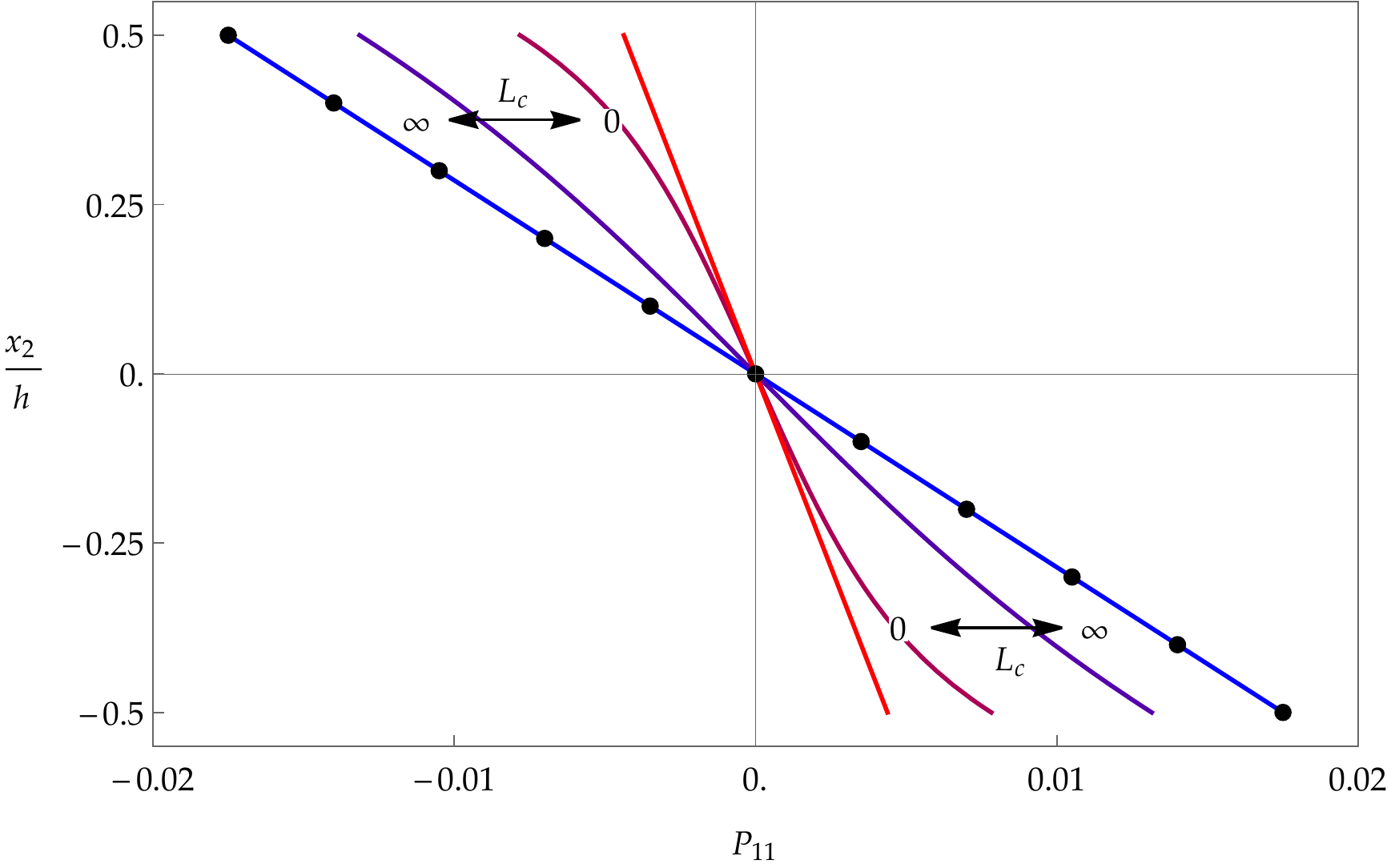}
		\caption{}
		\label{fig:P11_plot_RM_3}
	\end{subfigure}%
	\begin{subfigure}{.45\textwidth}
		\centering
		\includegraphics[width=0.95\linewidth]{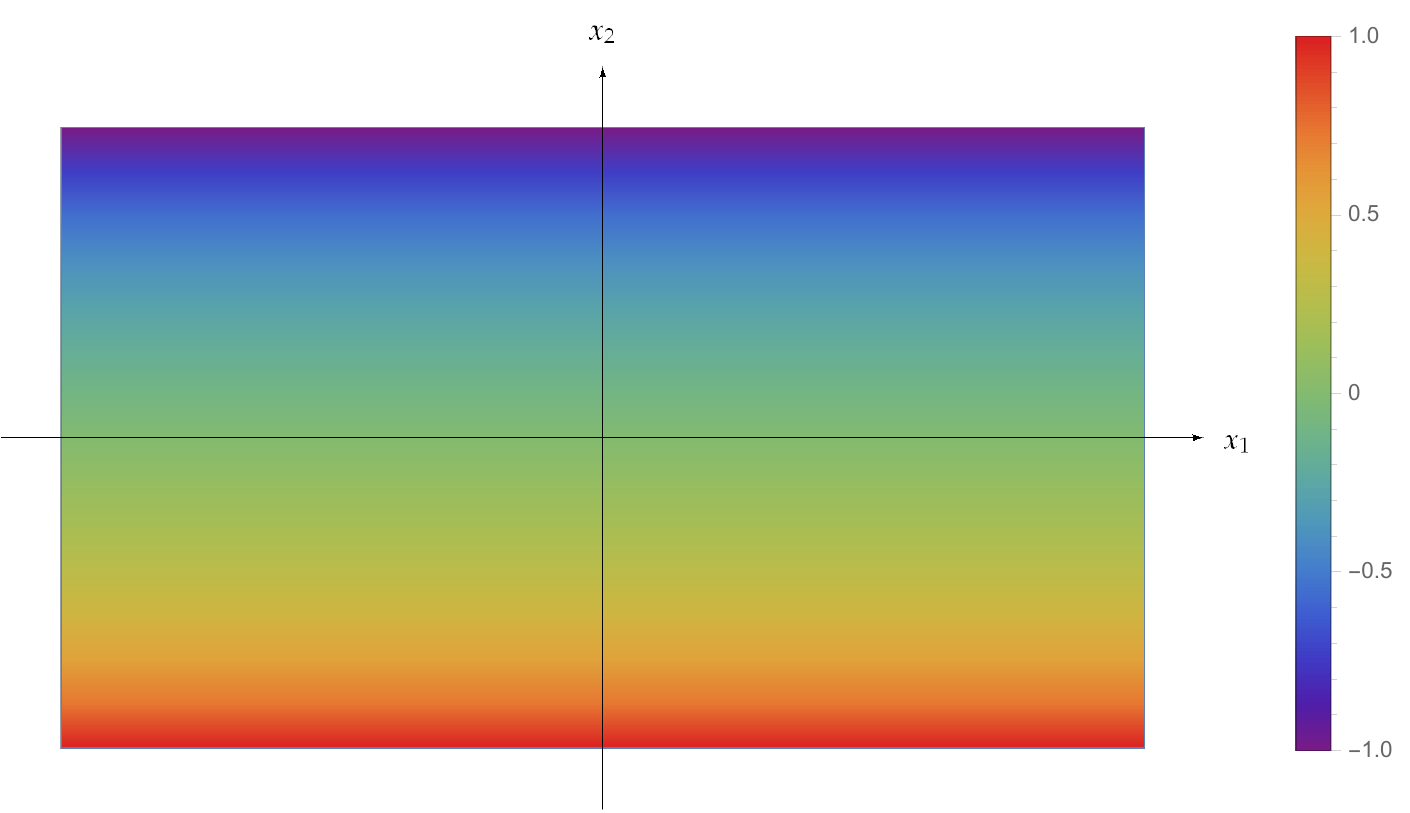}
		\caption{}
		\label{fig:P11_plot_RM_4}
	\end{subfigure}
	\caption{(\textbf{Relaxed micromorphic model}, general case)
		(a) Distribution across the thickness of $P_{11}$ while varying $L_c$ (continuous line) and the distribution across the thickness of $\left(\boldsymbol{\mbox{D}u}\right)_{11}$.
		On the vertical axis we have the dimensionless thickness while on the horizontal axis we have the quantity $P_{11}$.
		Notice that the red curve corresponds to the limit for $L_c \to 0$ while the blue curve correspond to the limit for $L_c \to \infty$.
		(b) dimensionless parametric plot of $P_{11}$ across the thickness for a given $L_c$ for the relaxed micromorphic model. The values of the parameters used are: $\mu _{e} = 2/3$, $\lambda _{e} = 1/3$, $\mu _{\tiny \mbox{micro}} = 2$, $\lambda _{\tiny \mbox{micro}} = 1$, $\mu = 1$, $a _1 = 1$, $a _2 = 1$.
	}
\end{figure}

In Fig.~\ref{fig:P11_plot_RM_3} we present the distribution across the thickness of $P_{11}$ while varying $L_c$.
The increase of $P_{11}$ as $L_c \to \infty$ is counterintuitive in the first instance, since the stiffness increases as $L_c \to \infty$ but it can be explained as follow:
(i) for $L_c \to \infty$ it is possible to demonstrate that $\boldsymbol{P} = \boldsymbol{\mbox{D}u}$ which justifies the fact that $\left(\boldsymbol{\mbox{D}u}\right)_{11}$ is superimposed to $P_{11}$ for $L_c \to \infty$;
(ii) for $L_c \to 0$ it is possible to demonstrate that $\boldsymbol{P}_{11} = \frac{D_{e}}{D_{e} + D_{\tiny \mbox{micro}}} \, \left(\boldsymbol{\mbox{D}u}\right)_{11} = \frac{(\boldsymbol{\mbox{D} u})_{11}}{4}$ and it will be anyway always smaller than $(\boldsymbol{\mbox{D} u})_{11}$.

\begin{figure}[H]
\centering
\includegraphics[width=0.5\linewidth]{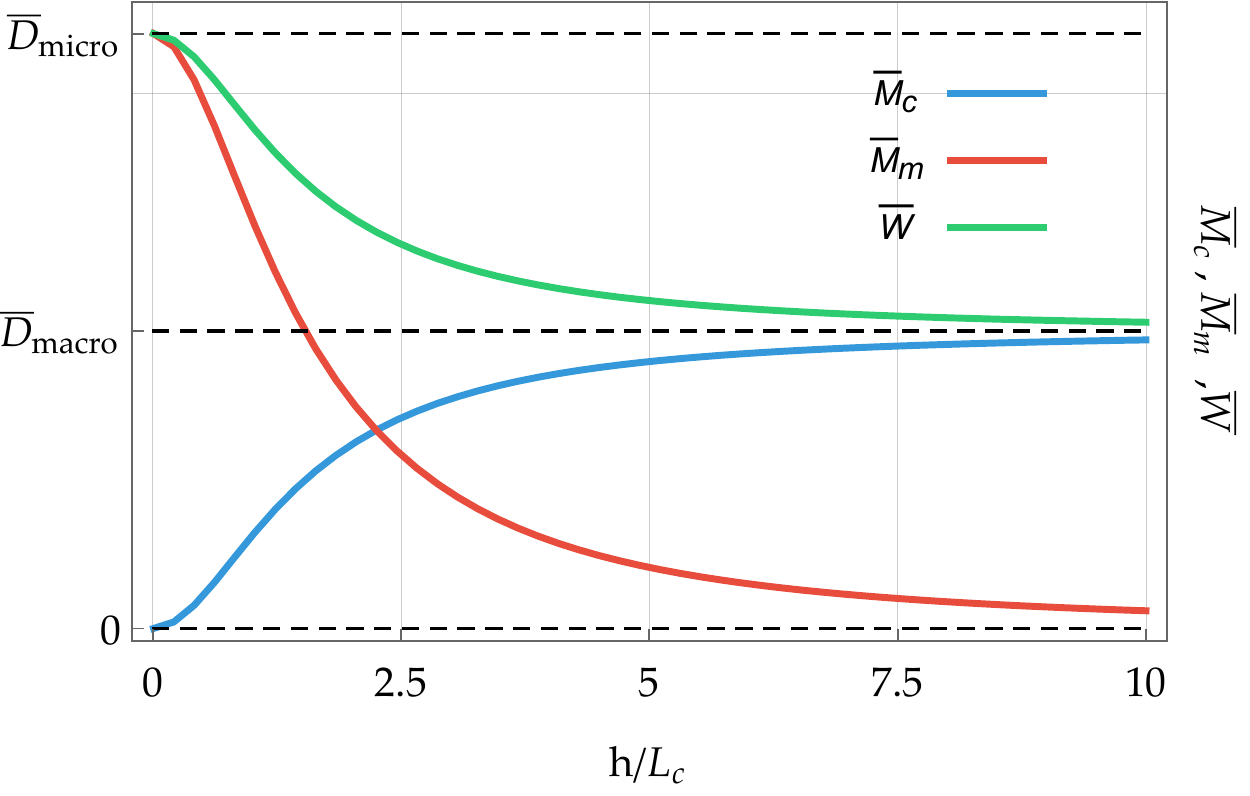}
\caption{(\textbf{Relaxed micromorphic model}, general case) Bending moments and energy while varying $L_c$. Observe that the bending stiffness remains bounded as $L_c \to \infty$ ($h\to 0$). This is a distinguishing feature of the relaxed micromorphic model. The values of the parameters used are: $\mu _e = 1$, $\lambda _e = 1$, $\mu_{\mbox{\tiny micro}} = 1$, $\lambda_{\mbox{\tiny micro}} = 1$, $\mu = 1$, $a _1 = 2$, $a _2 = 1$.}
\label{fig:all_plot_RM_3}
\end{figure}

The classical bending moment, the higher-order bending moment, and energy (per unit area d$x_1$d$x_3$) expressions are
\begin{align}
M_{\mbox{c}} (\boldsymbol{\kappa}) :=&
\displaystyle\int\limits_{-h/2}^{h/2}
\langle \boldsymbol{\widetilde{\sigma}} \boldsymbol{e}_1 , \boldsymbol{e}_1 \rangle  \, x_2 \, 
\mbox{d}x_{2}
=
\notag
\\*
&\hspace{0.5cm}
-z_{0} \left(\lambda _e \left(
p_{1} 
+ p_{2}\left(\frac{L_c}{h}\right)^2 
- p_{3} \left(\frac{L_c}{h}\right)^3 \tanh \left(\frac{f_{1} h}{2 L_c}\right) 
- p_{4}\left(\frac{L_c}{h}\right)^3 \tanh \left(\frac{f_{2} h}{2 L_c}\right) \right) 
\right.
\notag
\\*
&\hspace{0.5cm} \left.
- 2 \mu _e \left( q_{1} 
- q_{2}\left(\frac{L_c}{h}\right)^2 
+ q_{3}\left(\frac{L_c}{h}\right)^3 \tanh \left(\frac{f_{1} h}{2 L_c}\right)
+ q_{4}\left(\frac{L_c}{h}\right)^3 \tanh \left(\frac{f_{2} h}{2 L_c}\right)
\right)\right) \, ,
\notag
\\
M_{\mbox{m}}(\boldsymbol{\kappa}) :=&
\displaystyle\int\limits_{-h/2}^{h/2}
\langle \left(\boldsymbol{m} \times \boldsymbol{e}_1\right) \boldsymbol{e}_2 , \boldsymbol{e}_1 \rangle
\, \mbox{d}x_{2}
= 
\notag
\\*
&\hspace{0.5cm}
-\frac{h^3 \mu}{2 (b_{1}^2-b_{2}^2)}
\left(
r_{1} \left(\frac{L_c}{h}\right)^2
+r_{2} \left(\frac{L_c}{h}\right)^3 \tanh \left(\frac{f_{1} h}{2 L_c}\right) 
-r_{3} \left(\frac{L_c}{h}\right)^3 \tanh \left(\frac{f_{2} h}{2 L_c}\right)
\right) \, ,
\label{eq:sigm_ene_dimensionless_RM_3}
\\
W_{\mbox{tot}} (\boldsymbol{\kappa}) :=&
\displaystyle\int\limits_{-h/2}^{+h/2} W \left(\boldsymbol{\mbox{D}u}, \boldsymbol{P}, \mbox{Curl}\boldsymbol{P}\right) \, \mbox{d}x_{2}
=
\notag
\\*
&\hspace{0.5cm} 
-\frac{z_{0} \, \boldsymbol{\kappa}}{2} \left(\lambda _e \left(
p_{1} 
+ p_{2}\left(\frac{L_c}{h}\right)^2 
- p_{3} \left(\frac{L_c}{h}\right)^3 \tanh \left(\frac{f_{1} h}{2 L_c}\right) 
- p_{4}\left(\frac{L_c}{h}\right)^3 \tanh \left(\frac{f_{2} h}{2 L_c}\right) \right) 
\right.
\notag
\\*
&\hspace{0.5cm} \left.
- 2 \mu _e \left( q_{1} 
- q_{2}\left(\frac{L_c}{h}\right)^2 
+ q_{3}\left(\frac{L_c}{h}\right)^3 \tanh \left(\frac{f_{1} h}{2 L_c}\right)
+ q_{4}\left(\frac{L_c}{h}\right)^3 \tanh \left(\frac{f_{2} h}{2 L_c}\right)
\right)\right)
\notag
\\*
&\hspace{0.5cm}
-\frac{\boldsymbol{\kappa} \, h^3 \mu}{4 (b_{1}^2-b_{2}^2)}
\left(
r_{1} \left(\frac{L_c}{h}\right)^2
+r_{2} \left(\frac{L_c}{h}\right)^3 \tanh \left(\frac{f_{1} h}{2 L_c}\right) 
-r_{3} \left(\frac{L_c}{h}\right)^3 \tanh \left(\frac{f_{2} h}{2 L_c}\right)
\right) \, ,
\notag
\end{align}
where
\begin{align}
z_0 := \, & \frac{h^3 \mu _e}{6 f_{1}^3 f_{2}^3 \left(b_{1}^2-b_{2}^2\right) \left(\lambda _e+2 \mu _e\right)} \, ,
\notag
\\
p_{1} := \, & f_{1}^3 f_{2}^3 \left( - 2 b_{1}^2 \boldsymbol{\kappa}  + b_{1} (2 b_{3} + b_{4}) + b_{2} (2 b_{2} \boldsymbol{\kappa}  - b_{3} - 2 b_{4})\right) \, ,
\notag
\\
p_{2} := \, & 6 f_{1} f_{2} \left(b_{1}^2 \boldsymbol{\kappa}  \left(3 f_{1}^2 + f_{2}^2\right) - 3 b_{1} f_{1}^2 (b_{3} + b_{4}) + b_{1} f_{2}^2 (b_{4} - b_{3}) \right.
\notag
\\
&\left.+ b_{2} b_{3} \left(3 f_{1}^2 - f_{2}^2\right) + b_{2} \left(3 f_{1}^2 + f_{2}^2\right) (b_{4} - b_{2} \boldsymbol{\kappa} )\right) \, ,
\notag
\\
p_{3} := \, & 12 f_{2}^3 (b_{1} + b_{2}) (\boldsymbol{\kappa}  (b_{1} - b_{2}) - b_{3} + b_{4}) \, , \quad
p_{4} := \,  36 f_{1}^3 (b_{1} - b_{2}) (\boldsymbol{\kappa}  (b_{1} + b_{2}) - b_{3} - b_{4}) \, ,
\notag
\\
q_{1} := \, & f_{1}^3 f_{2}^3 \left(b_{1}^2 \boldsymbol{\kappa}  - b_{1} b_{3} + b_{2} (b_{4} - b_{2} \boldsymbol{\kappa} )\right) \, ,
\\
q_{2} := \, & 6 f_{1} f_{2} \left(b_{1}^2 \boldsymbol{\kappa}  \left(f_{1}^2 + f_{2}^2\right) - b_{1} f_{1}^2 (b_{3} + b_{4}) + b_{1} f_{2}^2 (b_{4} - b_{3}) \, , \right.
\notag
\\
&\left. + b_{2} b_{3} (f_{1} - f_{2}) (f_{1} + f_{2}) + b_{2} \left(f_{1}^2 + f_{2}^2\right) (b_{4} - b_{2} \boldsymbol{\kappa} )\right) \, ,
\notag
\\
q_{3} := \, & 12 f_{2}^3 (b_{1} + b_{2}) (\boldsymbol{\kappa}  (b_{1} - b_{2}) - b_{3} + b_{4}) \, , \quad
q_{4} := \,  12 f_{1}^3 (b_{1} - b_{2}) (\boldsymbol{\kappa}  (b_{1} + b_{2}) - b_{3} - b_{4}) \, ,
\notag
\\
r_{1}:= \, & -a_{1} (b_{1}+b_{2}) (b_{1} \boldsymbol{\kappa} -b_{2} \boldsymbol{\kappa} -b_{3}+b_{4})-a_{2} (b_{1}-b_{2}) (\boldsymbol{\kappa}  (b_{1}+b_{2})-b_{3}-b_{4}) \, ,
\notag
\\
r_{2}:= \, & \frac{2 a_{1} (b_{1}+b_{2}) (b_{1} \boldsymbol{\kappa} -b_{2} \boldsymbol{\kappa} -b_{3}+b_{4})}{f_{1}} \, , \quad
r_{3}:= \,  \frac{2 a_{2} (b_{1}-b_{2}) (-\boldsymbol{\kappa}  (b_{1}+b_{2})+b_{3}+b_{4})}{f_{2}} \, .
\notag
\end{align}
Again,
$
\frac{\mbox{d}}{\mbox{d}\boldsymbol{\kappa}}W_{\mbox{tot}}(\boldsymbol{\kappa}) = M_{\mbox{c}} (\boldsymbol{\kappa}) + M_{\mbox{m}} (\boldsymbol{\kappa}) \, .
$
The plot of the bending moments and the strain energy divided by $\frac{h^3}{12}\boldsymbol{\kappa}$ and $\frac{1}{2}\frac{h^3}{12}\boldsymbol{\kappa}^2$, respectively, while changing $L_c$ is shown in Fig.~\ref{fig:all_plot_RM_3}.
\subsubsection{Limit cases}
In addition to the two limits $ 0 \leftarrow L_c \rightarrow \infty$ we consider here as well $\mu_{\mbox{\tiny micro}} \to \infty$ since this makes for the transition to the classical Cosserat model.
\begin{equation}
\begin{array}{rlrl}
\displaystyle\lim_{L_c\to0} M_{\mbox{c}} (\boldsymbol{\kappa})
&= \dfrac{4 \, \mu_{\mbox{\tiny macro}} \, \left(\lambda_{\mbox{\tiny macro}} + \mu_{\mbox{\tiny macro}}\right) }{\lambda_{\mbox{\tiny macro}} + 2\mu_{\mbox{\tiny macro}}} \, \dfrac{h^3}{12} \, \boldsymbol{\kappa} 
= D_{\tiny \mbox{macro}} \, \boldsymbol{\kappa} \, , 
& \hspace{-3.5cm}
\displaystyle\lim_{L_c\to\infty} M_{\mbox{c}} (\boldsymbol{\kappa})
&= 0 \, ,
\\
\displaystyle\lim_{L_c\to\infty} M_{\mbox{m}} (\boldsymbol{\kappa})
&= \dfrac{4 \, \mu_{\mbox{\tiny micro}} \, \left(\lambda_{\mbox{\tiny micro}} + \mu_{\mbox{\tiny micro}}\right) }{\lambda_{\mbox{\tiny micro}} + 2\mu_{\mbox{\tiny micro}}} \, \dfrac{h^3}{12} \, \boldsymbol{\kappa} 
= D_{\tiny \mbox{micro}} \, \boldsymbol{\kappa} \, , 
& \hspace{-3.5cm}
\displaystyle\lim_{L_c\to0} M_{\mbox{m}} (\boldsymbol{\kappa})
&= 0 \, ,
\\
\displaystyle\lim_{L_c\to0} W_{\mbox{tot}} (\boldsymbol{\kappa})
&= 
\dfrac{1}{2}\dfrac{4 \, \mu_{\mbox{\tiny macro}} \, \left(\lambda_{\mbox{\tiny macro}} + \mu_{\mbox{\tiny macro}}\right) }{\lambda_{\mbox{\tiny macro}} + 2\mu_{\mbox{\tiny macro}}} \, \dfrac{h^3}{12} \, \boldsymbol{\kappa}^2 
= 
\dfrac{1}{2} \,D_{\tiny \mbox{macro}} \, \boldsymbol{\kappa}^2 \, ,
\\
\displaystyle\lim_{L_c\to\infty} W_{\mbox{tot}} (\boldsymbol{\kappa})
&=
\dfrac{1}{2} \dfrac{4 \, \mu_{\mbox{\tiny micro}} \, \left(\lambda_{\mbox{\tiny micro}} + \mu_{\mbox{\tiny micro}}\right) }{\lambda_{\mbox{\tiny micro}} + 2\mu_{\mbox{\tiny micro}}} \, \dfrac{h^3}{12} \, \boldsymbol{\kappa}^2 
= 
\dfrac{1}{2} \,D_{\tiny \mbox{micro}} \, \boldsymbol{\kappa}^2 \, ,
\\
\displaystyle\lim_{\mu_{\mbox{\tiny micro}}\to\infty} M_{\small \mbox{c}} (\boldsymbol{\kappa})
&=
\dfrac{4\mu_{\mbox{\tiny macro}} \left(\lambda_{\mbox{\tiny macro}}+\mu_{\mbox{\tiny macro}}\right)}{\lambda_{\mbox{\tiny macro}}+2 \mu_{\mbox{\tiny macro}}} \, 
\dfrac{h^3 }{12} \, 
\boldsymbol{\kappa} \, ,
\hspace{1cm}
\displaystyle\lim_{\mu_{\mbox{\tiny micro}}\to\infty} M_{\small \mbox{m}} (\boldsymbol{\kappa})
=
h \,  \mu \, L_c^2 \, \dfrac{a_1 + a_2}{2} \, \boldsymbol{\kappa} \, .
\\
\displaystyle\lim_{\mu_{\mbox{\tiny micro}}\to\infty} W_{\small \mbox{tot}} (\boldsymbol{\kappa})
&=
\dfrac{1}{2}
\left(
\underbrace{\dfrac{4\mu_{\mbox{\tiny macro}} \left(\lambda_{\mbox{\tiny macro}}+\mu_{\mbox{\tiny macro}}\right)}{\lambda_{\mbox{\tiny macro}}+2 \mu_{\mbox{\tiny macro}}} \, \dfrac{h^3 }{12}}_{D_{\mbox{\tiny macro}}}
+  \mu \, L_c^2 \, \dfrac{a_1 + a_2}{2} \, h
\right)  \, \boldsymbol{\kappa}^2
\, .
\end{array}
\label{eq:limits_RM_3}
\end{equation}
\begin{figure}[H]
\centering
\includegraphics[width=0.5\linewidth]{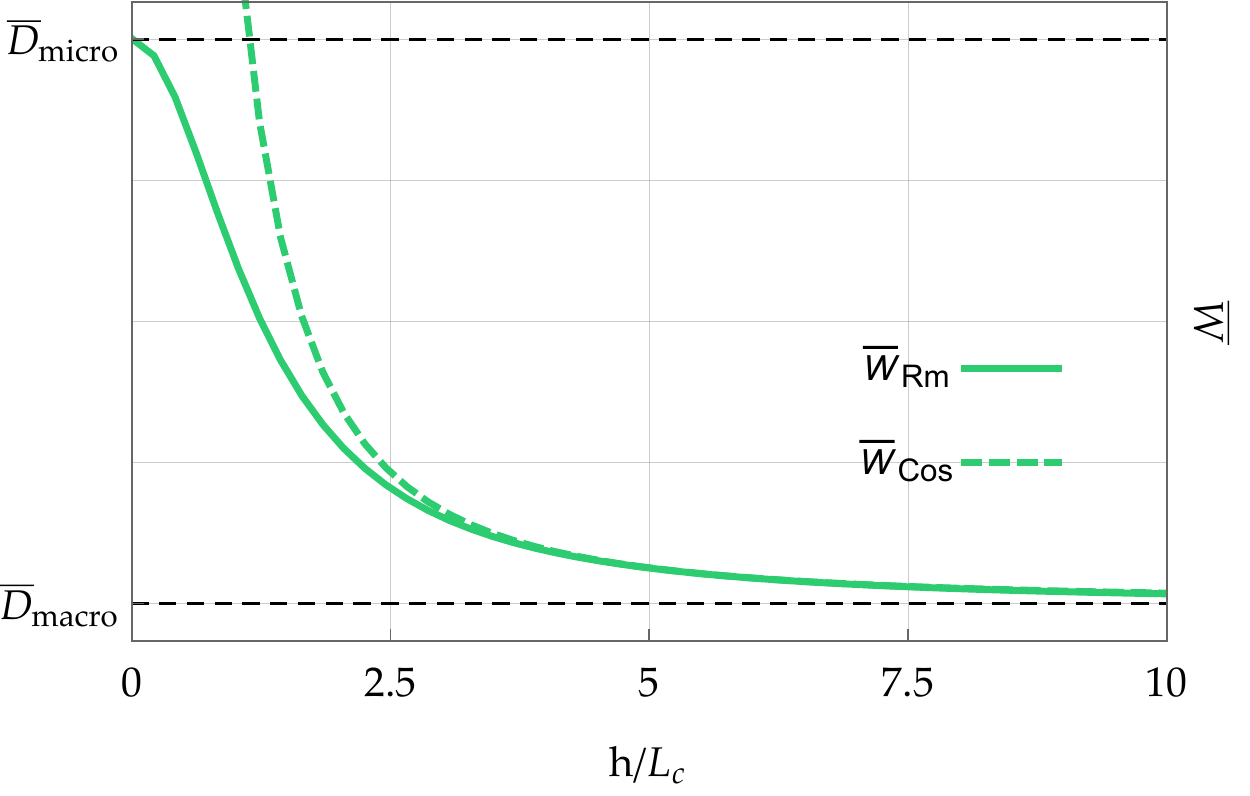}
\caption{Bending stiffness while varying $L_c$ for the relaxed micromorphic model and for the Cosserat model. Observe that the bending stiffness remains bounded for the relaxed micromorphic model while it blows up for the Cosserat model as $L_c \to \infty$ which is connected to the limit $\mu_{\tiny \mbox{micro}} \to \infty$, see Sect.\ref{sec:Cos}. For best comparison, the characteristic length scale of the Cosserat model has been chosen $L_c^{\tiny \mbox{Coss}} := \frac{L_c^{\tiny \mbox{relax}}}{\sqrt{2}}$. The values of the parameters used are: $\mu _{e} = 2$, $\lambda _{e} = 1$, $\mu _{\tiny \mbox{micro}} = 2/3$, $\lambda _{\tiny \mbox{micro}} = 1/3$, $\mu = 1$, $a _1 = 1$, $a _2 = 1$. Of course the micro-parameters are only relevant for the relaxed micromorphic model.}
\label{fig:all_plot_RM_vs_Cos_3}
\end{figure}
\section{The micro-stretch model in dislocation format}
\begin{figure}[H]
	\centering
	\includegraphics[width=0.5\linewidth]{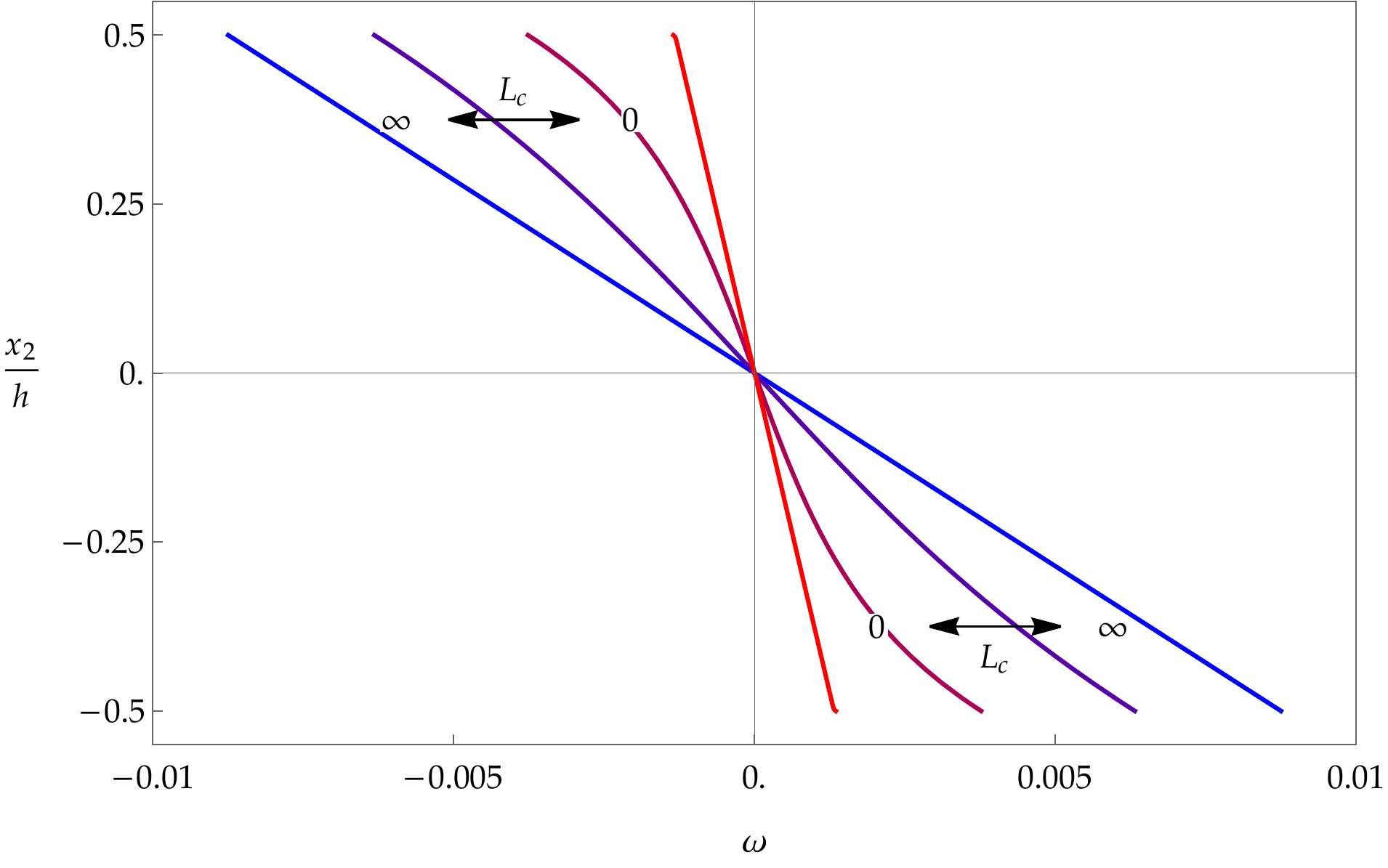}
	\caption{(\textbf{Micro-stretch model}) Distribution across the thickness of $\omega (x_2)$ while varying $L_c$.
		On the vertical axis we have the dimensionless thickness while on the horizontal we have the quantity $\omega(x_2)$.
		Notice that the red curve correspond to the limit for $L_c \to 0$ while the blue curve correspond to the limit for $L_c \to \infty$.
		The same reasoning regarding the increase of $P_{11} \sim \omega(x_2)$ for $L_c \to \infty$ done for the relaxed micromorphic model applies here. The values of the parameters used are: $\mu _e = 2/3$, $\kappa _e = 7/9$, $\mu_{\mbox{\tiny micro}} = 2$, $\kappa_{\mbox{\tiny micro}} = 7/3$, $\mu = 1$, $a _1 = 1/2$, $a _2 = 1/2$, $\boldsymbol{\kappa}=7/200$.}
	\label{fig:P11_plot_stretch}
\end{figure}

The micro-stretch model in dislocation format \cite{neff2014unifying} can be obtained from the relaxed micromorphic model by letting formally $\mu_{\mbox{\tiny micro}}\to\infty$, while $\kappa_{\mbox{\tiny micro}}<\infty$.
For bounded energy, this constrains $\mbox{dev} \, \mbox{sym} \, \boldsymbol{P} = 0 \Leftrightarrow \boldsymbol{P} = \boldsymbol{A} + \omega \boldsymbol{\mathbbm{1}}$, $\boldsymbol{A} \in \mathfrak{so}(3)$, $\omega \in \mathbb{R}$.
Thus the micro-stretch model has 4 additional degrees of freedom.
The expression of the strain energy for the isotropic micro-stretch continuum in dislocation format (i.e. with curvature energy only depending on the dislocation density tensor $\mbox{Curl}\left(\boldsymbol{A}+\omega \boldsymbol{\mathbbm{1}}\right)$) can then be written as \cite{neff2014unifying}:
\begin{align}
W \left(\boldsymbol{\mbox{D}u}, \boldsymbol{A},\omega,\mbox{Curl}\,\left(\boldsymbol{A} - \omega \boldsymbol{\mathbbm{1}}\right)\right) &
\notag
\\*
&
\hspace{-4.5cm}
=
\, \mu_{\tiny \mbox{macro}} \left\lVert \mbox{dev} \, \mbox{sym} \, \boldsymbol{\mbox{D}u} \right\rVert^{2}
+ \frac{\kappa_{e}}{2} \mbox{tr}^2 \left(\boldsymbol{\mbox{D}u} - \omega \boldsymbol{\mathbbm{1}} \right) 
+ \mu_{c} \left\lVert \mbox{skew} \left(\boldsymbol{\mbox{D}u} - \boldsymbol{A} \right) \right\rVert^{2}
+ \frac{\kappa_{\tiny \mbox{micro}}}{2} \mbox{tr}^2 \left(\omega \boldsymbol{\mathbbm{1}} \right)
\notag
\\*
&
\hspace{-4.2cm}
+ \frac{\mu \,L_c^2}{2} \,
\left(
a_1 \, \left\lVert \mbox{dev sym} \, \mbox{Curl} \, \left(\boldsymbol{A} + \omega \boldsymbol{\mathbbm{1}}\right)\right\rVert^2
+ a_2 \, \left\lVert \mbox{skew} \,  \mbox{Curl} \, \left(\boldsymbol{A} + \omega \boldsymbol{\mathbbm{1}}\right)\right\rVert^2 +
\frac{a_3}{3} \, \mbox{tr}^2 \left(\mbox{Curl} \, \left(\boldsymbol{A} + \omega \boldsymbol{\mathbbm{1}}\right)\right)
\right) 
\label{eq:energy_micro_stretch}
\\*
&
\hspace{-4.5cm}
=
\, \mu_{\tiny \mbox{macro}} \left\lVert \mbox{dev} \, \mbox{sym} \, \boldsymbol{\mbox{D}u} \right\rVert^{2}
+ \frac{\kappa_{e}}{2} \mbox{tr}^2 \left(\boldsymbol{\mbox{D}u} - \omega \boldsymbol{\mathbbm{1}} \right) 
+ \mu_{c} \left\lVert \mbox{skew} \left(\boldsymbol{\mbox{D}u} - \boldsymbol{A} \right) \right\rVert^{2}
+ \frac{9}{2} \, \kappa_{\tiny \mbox{micro}} \, \omega^2
\notag
\\*
&
\hspace{-4.2cm}
+ \frac{\mu \,L_c^2}{2} \,
\left(
a_1 \, \left\lVert \mbox{dev sym} \, \mbox{Curl} \, \boldsymbol{A} \right\rVert^2
+ a_2 \, \left\lVert \mbox{skew} \,  \mbox{Curl} \, \left(\boldsymbol{A} + \omega \boldsymbol{\mathbbm{1}}\right) \right\rVert^2
+ \frac{a_3}{3} \, \mbox{tr}^2 \left(\mbox{Curl} \, \boldsymbol{A} \right)
\right) \, ,
\notag
\end{align}
since $\mbox{Curl} \left(\omega \boldsymbol{\mathbbm{1}}\right) \in \mathfrak{so}(3)$.
The equilibrium equations without body forces are then
\begin{align}
\mbox{Div}\overbrace{\left[
2\mu_{\tiny \mbox{macro}}\,\mbox{dev}\,\mbox{sym} \, \boldsymbol{\mbox{D}u}
+ \kappa_{e} \mbox{tr} \left(\boldsymbol{\mbox{D}u} - \omega \boldsymbol{\mathbbm{1}}\right) \boldsymbol{\mathbbm{1}}
+ 2\mu_{c}\,\mbox{skew} \left(\boldsymbol{\mbox{D}u} - \boldsymbol{A}\right) \right]}^{\mathlarger{\widetilde{\sigma}}:=}
&= \boldsymbol{0} \, ,
\notag
\\
2\mu_{c}\,\mbox{skew} \left(\boldsymbol{\mbox{D}u} - \boldsymbol{A}\right)
\hspace{10.1cm}
\notag
\\
-\mu \, L_c^2 \, \mbox{skew} \, \mbox{Curl}\,\left(
a_1 \, \mbox{dev} \, \mbox{sym} \, \mbox{Curl} \, \boldsymbol{A} \, 
+ a_2 \, \mbox{skew} \, \mbox{Curl} \, \left(\boldsymbol{A} + \, \omega \boldsymbol{\mathbbm{1}}\right) \,
+ \frac{a_3}{3} \, \mbox{tr} \left(\mbox{Curl} \, \boldsymbol{A} \right)\boldsymbol{\mathbbm{1}} \, 
\right) &= \boldsymbol{0} \,.
\label{eq:equi_micro_stretch}
\\
\mbox{tr}
\bigg[
2\mu_{\tiny \mbox{macro}}\,\mbox{dev}\,\mbox{sym} \, \boldsymbol{\mbox{D}u}
\hspace{10.25cm}
\notag
\\
+ \kappa_{e} \mbox{tr} \left(\boldsymbol{\mbox{D}u} - \omega \boldsymbol{\mathbbm{1}}\right) \boldsymbol{\mathbbm{1}}
- \kappa_{\mbox{\tiny micro}} \mbox{tr} \left( \omega \boldsymbol{\mathbbm{1}}\right) \boldsymbol{\mathbbm{1}}
-\mu \, L_c^2 \,  a_2 \, \mbox{Curl}\,
\mbox{skew} \, \mbox{Curl} \, \left(\omega \boldsymbol{\mathbbm{1}} + \boldsymbol{A}\right) 
\bigg]
&= \boldsymbol{0} \,.
\notag
\end{align}
The boundary conditions at the upper and lower surface (free surface) are 
\begin{align}
\boldsymbol{\widetilde{t}}(x_2 = \pm \, h/2) &= 
\pm \, \boldsymbol{\widetilde{\sigma}}(x_2) \cdot \boldsymbol{e}_2 = 
\boldsymbol{0} \, ,
\notag
\\
\boldsymbol{\eta}(x_2 = \pm \, h/2) &= 
\pm \, \mbox{skew}\left(\boldsymbol{m} (x_2) \cdot \boldsymbol{\epsilon} \cdot \boldsymbol{e}_2\right) = 
\pm \, \mbox{skew}\left(\boldsymbol{m} (x_2) \times \boldsymbol{e}_2 \right) = 
\boldsymbol{0} \, ,
\label{eq:BC_micro_stretch}
\\
\gamma(x_2 = \pm \, h/2) &= 
\pm \, \frac{1}{3} \mbox{tr}\left(\boldsymbol{m} (x_2) \cdot \boldsymbol{\epsilon} \cdot \boldsymbol{e}_2\right) = 
\pm \, \frac{1}{3} \mbox{tr}\left(\boldsymbol{m} (x_2) \times \boldsymbol{e}_2 \right) = 
\boldsymbol{0} \, .
\notag
\end{align}
According with the reference system shown in Fig.~\ref{fig:intro}, the ansatz for the displacement field and the function $\omega$ is
\begin{align}
\boldsymbol{u}(x_1,x_2) &=
\left(
\begin{array}{c}
-\kappa_1 \, x_1 x_2 \\
v(x_2)+\frac{\kappa_1  x_1^2}{2} \\
0 \\
\end{array}
\right)
\, ,
\qquad
\boldsymbol{A}(x_1) =
\left(
\begin{array}{ccc}
0 & -\kappa_2 \, x_1 & 0 \\
\kappa_2 \, x_1 & 0 & 0 \\
0 & 0 & 0 \\
\end{array}
\right)
\, ,
\qquad
\omega = \omega\left(x_2\right)
\, .
\label{eq:ansatz_micro_stretch}
\end{align}
The equilibrium equations (\ref{eq:equi_micro_stretch}) then result in \!\!
\footnote{Where
$\kappa_e = \frac{2\mu_e + 3\lambda_e}{3}$ and $\kappa_{\mbox{\tiny micro}} = \frac{2\mu_{\mbox{\tiny micro}} + 3\lambda_{\mbox{\tiny micro}}}{3}$ are the meso- and the micro-scale 3D bulk modulus.}
\begin{align}
\kappa _e \left(v''(x_{2}) - \boldsymbol{\kappa} - 3 \omega '(x_{2})\right)
+ \frac{2}{3} \mu _{\tiny \mbox{macro}} \left(\boldsymbol{\kappa} +2 v''(x_{2})\right) = 0
\, ,
\label{eq:equi_equa_micro_stretch}
\\
\frac{2}{3} \text{a2} \, \mu \, L_c^2 \, \omega ''(x_{2})
+ \kappa _e \left(v'(x_{2}) - \boldsymbol{\kappa} \, x_{2}  - 3 \omega (x_{2})\right)
- 3 \kappa _{\mbox{\tiny micro}} \, \omega (x_{2}) = 0
\, ,
\notag
\end{align}
since the second equation eq.(\ref{eq:equi_micro_stretch})$_2$ is already satisfied.
From eq.(\ref{eq:equi_equa_micro_stretch})$_1$ it is possible to evaluate $v''(x_2)$ and consequently $v'(x_2)$ as follows
\begin{align}
v''(x_2) =& \frac{9 \kappa _e}{3 \kappa _e+4 \mu  _{\tiny \mbox{macro}}} \omega '(x_{2})
+ \frac{3 \kappa _e-2 \mu  _{\tiny \mbox{macro}}}{3 \kappa _e+4 \mu  _{\tiny \mbox{macro}}} \boldsymbol{\kappa}
\, ,
\label{eq:d1v_d2v_micro_stretch}
\\*
v'(x_2)  =& \frac{9 \kappa _e}{3 \kappa _e+4 \mu  _{\tiny \mbox{macro}}} \omega (x_{2})
+ \frac{3 \kappa _e-2 \mu  _{\tiny \mbox{macro}}}{3 \kappa _e+4 \mu  _{\tiny \mbox{macro}}} \boldsymbol{\kappa} \, x_{2} 
+ c_0
\, .
\notag
\end{align}
Substituting back the expression of $v''(x_2)$ and $v'(x_2)$ in (\ref{eq:equi_equa_micro_stretch}) it is possible to evaluate $\omega(x_2)$  from eq.(\ref{eq:equi_equa_micro_stretch})$_2$ and consequently $v(x_2)$ which results in
\begin{align}
\omega(x_2) = \,\, & 
c_1 \, e^{-\frac{f_{1} \, x_{2}}{L_c}}
+ c_2 \, e^{\frac{f_{1} \, x_{2}}{L_c}}
- f_{2} \, \boldsymbol{\kappa} \, x_{2}
+ f_{3} \, c_{0}
\, ,
\notag
\\*
v(x_2) = \,\, &
\frac{9 \kappa _e}{3 \kappa _e+4 \mu _{\tiny \mbox{macro}}} \frac{L_c}{f_{1}}
\left( c_2 \, e^{\frac{f_{1} x_{2}}{L_c}} - c_1 \, e^{-\frac{f_{1} x_{2}}{L_c}}\right)
+ \frac{(3-9 f_{2}) \kappa _e-2 \mu _{\tiny \mbox{macro}}}{6 \kappa _e+8 \mu _{\tiny \mbox{macro}}} \, x_{2}^2 \, \boldsymbol{\kappa}
+ \left( 1 + \frac{3 f_{2} \kappa _e}{2 \mu _{\tiny \mbox{macro}}} \right) \, x_{2} \, c_{0}
\, ,
\notag
\\
f_{1} := \,\, & 
\frac{1}{\sqrt{a_{2} \, \mu }}
\sqrt{\frac{9}{2}\left(\kappa_{\tiny \mbox{micro}} + \frac{4\kappa _e \, \mu_e}{3\kappa _e + 4 \mu _{\tiny \mbox{macro}}}\right)}
\, , \qquad\qquad
f_{2} := \,\,
\frac{6 \kappa _e \, \mu _{\tiny \mbox{macro}}}{12 \mu _{\tiny \mbox{macro}} \left(\kappa _e + \kappa _{\mbox{\tiny micro}}\right) + 9 \kappa _e \, \kappa _{\mbox{\tiny micro}}}
\, .
\label{eq:omega_micro_stretch}
\\*
f_{3} := \,\, & \left(\frac{2}{3} + \frac{\kappa _e}{2 \mu _{\tiny \mbox{macro}}} \right) f_{2}
\, .
\notag
\end{align}
In Fig.~\ref{fig:P11_plot_stretch} we present the distribution across the thickness of $P_{11}$ while varying $L_c$.

\begin{figure}[H]
\centering
\includegraphics[width=0.5\linewidth]{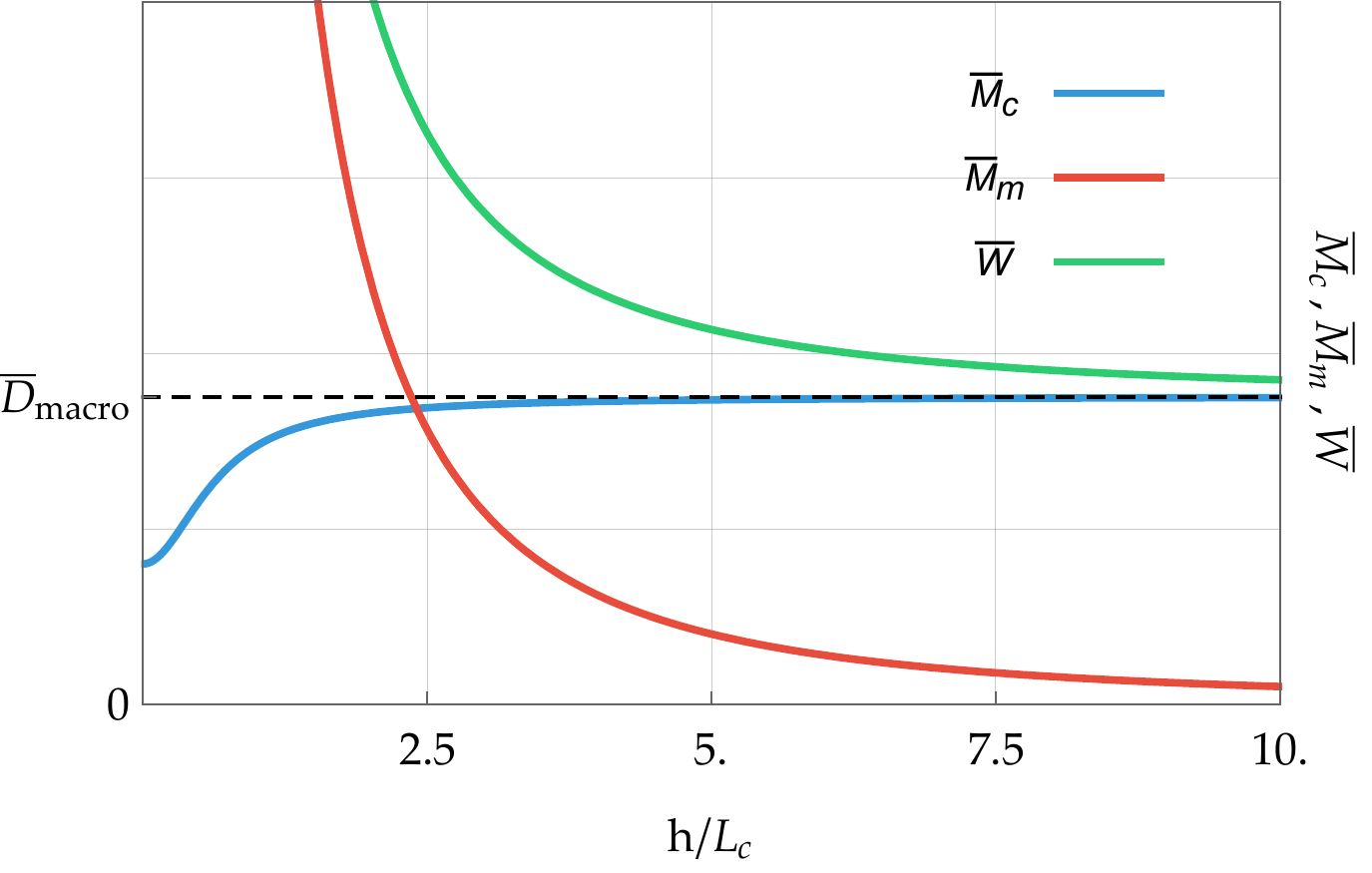}
\caption{(\textbf{Micro-stretch model} in dislocation format) Energy expression while varying $L_c$. Observe that the bending stiffness is unbounded as $L_c \to \infty$ ($h\to 0$), while for $L_c \to 0$ we recover the classical stiffness $D_{\tiny \mbox{macro}}$. The values of the parameters used are: $\mu _e = 2/3$, $\kappa _e = 7/9$, $\mu_{\mbox{\tiny micro}} = 2$, $\kappa_{\mbox{\tiny micro}} = 7/3$, $\mu = 1$, $a _1 = 1$, $a _2 = 1/2$.}
\label{fig:all_plot_micro_stretch}
\end{figure}

The classical bending moment, the higher-order bending moment $ \left(M_{\mbox{m}}(\boldsymbol{\kappa}) = \widehat{M}_{\mbox{m}}(\boldsymbol{\kappa}) + \widetilde{M}_{\mbox{m}}(\boldsymbol{\kappa})\right)$, and energy (per unit area d$x_1$d$x_3$) expressions are
\begin{align}
M_{\mbox{c}} (\boldsymbol{\kappa})&:=
\displaystyle\int\limits_{-h/2}^{h/2}
\langle \boldsymbol{\widetilde{\sigma}} \boldsymbol{e}_1 , \boldsymbol{e}_1 \rangle  \, x_2 \, 
\mbox{d}x_{2}
= 
\frac{h^3}{12}
\mu _{\tiny \mbox{macro}}
\left[
p_{1}
+ p_{2} \, f_{1} \left(\frac{L_c}{h}\right)^2
- 2 p_{2} \left(\frac{L_c}{h}\right)^3  \tanh \left(\frac{f_{1} h}{2 L_c}\right)
\right]
\boldsymbol{\kappa}
\, ,
\notag
\\
\widehat{M}_{\mbox{m}}(\boldsymbol{\kappa}) &:=
\displaystyle\int\limits_{-h/2}^{h/2}
2\langle \mbox{skew} \left(\boldsymbol{m} \times \boldsymbol{e}_1 \right) \boldsymbol{e}_2 , \boldsymbol{e}_1 \rangle
\, \mbox{d}x_{2}
=
\frac{h^3}{12} \,
\mu _{\tiny \mbox{macro}} \,
\left[
q_{1}\left(\frac{L_c}{h}\right)^2 
+ q_{2} \left(\frac{L_c}{h}\right)^3 \tanh \left(\frac{f_{1} h}{2 L_c}\right)
\right]
\boldsymbol{\kappa}
\, ,
\notag
\\
\widetilde{M}_{\mbox{m}}(\boldsymbol{\kappa}) &:=
\displaystyle\int\limits_{-h/2}^{h/2}
\mbox{tr} \left(\boldsymbol{m} \times \boldsymbol{e}_1\right)
\, \mbox{d}x_{2}
=
0
\, ,
\label{eq:sigm_ene_dimensionless_micro_stretch}
\\
W_{\mbox{tot}} (\boldsymbol{\kappa})&:=
\displaystyle\int\limits_{-h/2}^{+h/2} W \left(\boldsymbol{\mbox{D}u}, \boldsymbol{A},\omega,\mbox{Curl}\,\left(\boldsymbol{A} - \omega \boldsymbol{\mathbbm{1}}\right)\right) \, \mbox{d}x_{2}
\notag
\\*
& 
= 
\frac{1}{2}
\frac{h^3}{12} \,
\mu _{\tiny \mbox{macro}}
\Bigg[
p_{1} 
+ ( p_{2} \, f_{1} + \, q_{1}) \left(\frac{L_c}{h}\right)^2
-  (2 p_{2} - q_{2}) \left(\frac{L_c}{h}\right)^3 \tanh \left(\frac{f_{1} h}{2 L_c}\right)
\Bigg]
\boldsymbol{\kappa}^2
\, ,
\notag
\\
p_{1} &:= \frac{2 (6-9 f_{2}) \kappa _e + 4 \mu _{\tiny \mbox{macro}}}{3 \kappa _e+4 \mu _{\tiny \mbox{macro}}} \, ,
\quad
p_{2} := \frac{108 \, (2 f_{2}-1) \kappa _e}{f_{1}^3 \left(3 \kappa _e+4 \mu _{\tiny \mbox{macro}}\right)} \, ,
\notag
\\*
q_{1} &:= \frac{6 \mu \, (a_{1} + a_{2} \left( 1 - 2 f_{2}\right) )}{\mu _{\tiny \mbox{macro}}} \, ,
\quad
q_{2} := \frac{12 a_{2} \, \mu \, (2 f_{2} - 1)}{f_{1} \mu _{\tiny \mbox{macro}}} \, .
\notag
\end{align}
One of the two higher-order bending moment is zero and we have
$
\frac{\mbox{d}}{\mbox{d}\boldsymbol{\kappa}}W_{\mbox{tot}}(\boldsymbol{\kappa}) = M_{\mbox{c}} (\boldsymbol{\kappa}) + \widetilde{M}_{\mbox{m}} (\boldsymbol{\kappa})
$.
The plot of the non zero bending moments and the strain energy divided by $\frac{h^3}{12}\boldsymbol{\kappa}$ and $\frac{1}{2}\frac{h^3}{12}\boldsymbol{\kappa}^2$, respectively, while changing $L_c$ is shown in Fig.~\ref{fig:all_plot_micro_stretch}.
\subsubsection{Limit cases}
\begin{align}
\displaystyle\lim_{L_c\to0} M_{\mbox{c}} (\boldsymbol{\kappa})
=& \, 
\frac{h^3}{12}
\frac{4 \mu _{\mbox{\tiny macro}} \left(3 \kappa _{\mbox{\tiny macro}}+\mu _{\mbox{\tiny macro}}\right)}{3 \kappa _{\mbox{\tiny macro}}+4 \mu _{\mbox{\tiny macro}}}
\boldsymbol{\kappa} 
=D _{\mbox{\tiny macro}} \, \boldsymbol{\kappa}
\, ,
&\displaystyle\lim_{L_c\to0} M_{\mbox{m}} (\boldsymbol{\kappa})
&=  \, 0
\, ,
\notag
\\
\displaystyle\lim_{L_c\to0} W_{\mbox{tot}} (\boldsymbol{\kappa})
= & \, 
\frac{1}{2}
\frac{h^3}{12}
\frac{4 \mu _{\mbox{\tiny macro}} \left(3 \kappa _{\mbox{\tiny macro}}+\mu _{\mbox{\tiny macro}}\right)}{3 \kappa _{\mbox{\tiny macro}}+4 \mu _{\mbox{\tiny macro}}}
\boldsymbol{\kappa}^2
=\frac{1}{2} D _{\mbox{\tiny macro}} \, \boldsymbol{\kappa}^2
,
&\displaystyle\lim_{L_c\to\infty} M_{\mbox{c}} (\boldsymbol{\kappa})
&= \, 
\frac{h^3}{12} \mu _{\tiny \mbox{macro}} \boldsymbol{\kappa}
\, ,
\\
\displaystyle\lim_{L_c\to\infty} M_{\mbox{m}} (\boldsymbol{\kappa})
= & \, \infty
\, ,
&\displaystyle\lim_{L_c\to\infty} W_{\mbox{tot}} (\boldsymbol{\kappa})
&=  \, \infty
\, .
\notag
\end{align}
\section{Cylindrical bending for the isotropic Cosserat continuum}
\label{sec:Cos}
The expression of the strain energy for the isotropic Cosserat continuum in dislocation tensor format can be written \!\!
\footnote{
The equivalent formulation in terms of a rotation vector $\vartheta:=\mbox{axl} (\boldsymbol{A}) \in \mathbb{R}^3$ is given in the appendix.
}
\begin{align}
W \left(\boldsymbol{\mbox{D}u}, \boldsymbol{A},\mbox{Curl}\,\boldsymbol{A}\right) = &
\, \mu_{\mbox{\tiny macro}} \left\lVert \mbox{sym} \, \boldsymbol{\mbox{D}u} \right\rVert^{2}
+ \frac{\lambda_{\mbox{\tiny macro}}}{2} \mbox{tr}^2 \left(\boldsymbol{\mbox{D}u} \right) 
+ \mu_{c} \left\lVert \mbox{skew} \left(\boldsymbol{\mbox{D}u} - \boldsymbol{A} \right) \right\rVert^{2}
\label{eq:energy_Cos}
\\
&+ \frac{\mu \, L_c^2}{2}
\left(
a_1 \, \left \lVert \mbox{dev} \, \mbox{sym} \, \mbox{Curl} \, \boldsymbol{A}\right \rVert^2 \, 
+ a_2 \, \left \lVert \mbox{skew} \, \mbox{Curl} \, \boldsymbol{A}\right \rVert^2 \, 
+ \frac{a_3}{3} \, \mbox{tr}^2 \left(\mbox{Curl} \, \boldsymbol{A} \right)
\right)  \, ,
\notag
\end{align}
where $A \in \mathfrak{so}(3)$.
The Cosserat energy can be obtained from the relaxed micromorphic model by letting formally $\mu_{\mbox{\tiny micro}}, \lambda_{\mbox{\tiny micro}} \to \infty$ or letting $\kappa_{\mbox{\tiny micro}} \to \infty$ in the micro-stretch model.

It is important to underline that, since $\mbox{tr} \left(\mbox{Curl} \, \boldsymbol{A} \right)=0$ and $\left \lVert \mbox{dev} \, \mbox{sym} \, \mbox{Curl} \, \boldsymbol{A}\right \rVert^2 = \left \lVert \mbox{sym} \, \mbox{Curl} \, \boldsymbol{A}\right \rVert^2 = \left \lVert \mbox{skew} \, \mbox{Curl} \, \boldsymbol{A}\right \rVert^2 = \frac{1}{2}\left \lVert \mbox{Curl} \, \boldsymbol{A}\right \rVert^2$ for a plane problem (see Appendix~\ref{Sec:appendix_3D-Curl}), the elastic energy ends up to have only one effective curvature parameter 
\begin{align}
W \left(\boldsymbol{\mbox{D}u}, \boldsymbol{A},\mbox{Curl}\,\boldsymbol{A}\right) = &
\, \mu_{\mbox{\tiny macro}} \left\lVert \mbox{sym} \, \boldsymbol{\mbox{D}u} \right\rVert^{2}
+ \frac{\lambda_{\mbox{\tiny macro}}}{2} \mbox{tr}^2 \left(\boldsymbol{\mbox{D}u} \right) 
+ \mu_{c} \left\lVert \mbox{skew} \left(\boldsymbol{\mbox{D}u} - \boldsymbol{A} \right) \right\rVert^{2}
\label{eq:energy_Cos_2}
\\*
&+ \frac{\mu \, L_c^2}{2}
\,
\frac{a_1 + a_2}{2} \, \left \lVert \mbox{Curl} \, \boldsymbol{A}\right \rVert^2  \, .
\notag
\end{align}
The equilibrium equations without body forces are the following:
\begin{align}
\mbox{Div}\overbrace{\left[2\mu_{\mbox{\tiny macro}}\,\mbox{sym} \, \boldsymbol{\mbox{D}u} + \lambda_{\mbox{\tiny macro}} \mbox{tr} \left(\boldsymbol{\mbox{D}u} \right) \boldsymbol{\mathbbm{1}}
+ 2\mu_{c}\,\mbox{skew} \left(\boldsymbol{\mbox{D}u} - \boldsymbol{A}\right) \right]}^{\mathlarger{\widetilde{\sigma}}:=}
&= \boldsymbol{0} \, ,
\notag
\\
2\mu_{c}\,\mbox{skew} \left(\boldsymbol{\mbox{D}u} - \boldsymbol{A}\right)
\hspace{10cm}
\label{eq:equi_Cos}
\\
-\mu \, L_c^2 \, \mbox{skew} \, \mbox{Curl}\,\left(
a_1 \, \mbox{dev} \, \mbox{sym} \, \mbox{Curl} \, \boldsymbol{A} \, 
+ a_2 \, \mbox{skew} \, \mbox{Curl} \, \boldsymbol{A} \, 
+ \frac{a_3}{3} \, \mbox{tr} \left(\mbox{Curl} \, \boldsymbol{A} \right)\boldsymbol{\mathbbm{1}} \, 
\right) &= \boldsymbol{0} \,.
\notag
\end{align}
The boundary conditions for the upper and lower surface (free surface) are 
\begin{align}
\boldsymbol{\widetilde{t}}(x_2 = \pm \, h/2) &= 
\pm \, \boldsymbol{\widetilde{\sigma}}(x_2) \cdot \boldsymbol{e}_2 = 
\boldsymbol{0} \, ,
\label{eq:BC_Cos_gen}
\\*
\boldsymbol{\eta}(x_2 = \pm \, h/2) &= 
\pm \, \mbox{skew}\left(\boldsymbol{m} (x_2) \times \boldsymbol{e}_2\right) = 
\pm \, \mbox{skew}\left(\boldsymbol{m} (x_2) \cdot \boldsymbol{\epsilon} \cdot \boldsymbol{e}_2\right) = 
\boldsymbol{0} \, ,
\notag
\end{align}
where the expression of $\boldsymbol{\widetilde{\sigma}}$ is in eq.(\ref{eq:equi_Cos}), $\boldsymbol{e}_2$ is defined the unit vector aligned to the $x_2$-direction, $\boldsymbol{\epsilon}$ is the Levi-Civita tensor, and the moment stress tensor $\boldsymbol{m} = \mu \, L_c^2 \, \left(a_1 \, \mbox{dev} \, \mbox{sym} \, \mbox{Curl} \, \boldsymbol{A} + a_2 \, \mbox{skew} \, \mbox{Curl} \, \boldsymbol{A} + a_3/3 \, \mbox{tr} \left(\mbox{Curl} \, \boldsymbol{A} \right)\boldsymbol{\mathbbm{1}} \right)$.
The bending Neumann condition $\mbox{skew}\left(\boldsymbol{m} (x_2) \times \boldsymbol{e}_2\right) = 0$ at the upper and lower surface will be derived in Appendix~\ref{app:neumann_coss}.

According to the reference system shown in Fig.~\ref{fig:intro}, the ansatz for the displacement field and the micro-rotation is
\begin{equation}
\boldsymbol{u}(x_1,x_2)=
\left(
\begin{array}{c}
-\kappa_1 \, x_1 x_2 \\
v(x_2)+\frac{\kappa_1  x_1^2}{2} \\
0 \\
\end{array}
\right) \, ,
\qquad
\boldsymbol{A}(x_1) =
\left(
\begin{array}{ccc}
0 & -\kappa_2 \, x_1 & 0 \\
\kappa_2 \, x_1 & 0 & 0 \\
0 & 0 & 0 \\
\end{array}
\right) \, ,
\label{eq:ansatz_Cos}
\end{equation}
while the  gradient of the displacement field result to be
\begin{equation}
\boldsymbol{\mbox{D}u} = 
\left(
\begin{array}{ccc}
-\kappa_1 \, x_2 & - \kappa_1 \, x_1 & 0 \\
\kappa_1 \, x_1 & v'(x_2) & 0 \\
0 & 0 & 0 \\
\end{array}
\right) \, .
\label{eq:grad_Cos}
\end{equation}
Substituting the ansatz eq.(\ref{eq:ansatz_Cos}) in eq.(\ref{eq:equi_Cos}) the equilibrium equations turn in
\begin{align}
2 \mu _c (\kappa_1-\kappa_2 )-\kappa_1 \lambda_{\mbox{\tiny macro}}+\left(\lambda_{\mbox{\tiny macro}}+2 \mu_{\mbox{\tiny macro}}\right) v''(x_{2}) = 0 \, ,
\label{eq:equi_equa_Cos}
\\
2 x_1 \mu _c (\kappa_2 -\kappa_1) = 0 \, ,
\qquad\qquad\qquad\,\,
2 x_1 \mu _c  (\kappa_1 -\kappa_2) = 0 \, .
\notag
\end{align}
In order to satisfy eq.(\ref{eq:equi_equa_Cos})$_2$ and eq.(\ref{eq:equi_equa_Cos})$_3$ either $\mu _c = 0$ or $\kappa_1 = \kappa_2 = \boldsymbol{\kappa}$;
we have chosen the latter option, which implies that the skew-symmetric part of the gradient of the displacement eq.(\ref{eq:grad_Cos}) is the same as the skew-symmetric part of the micro-rotation eq.(\ref{eq:ansatz_Cos})$_2$.
This also implies that the Cosserat couple modulus $\mu_c$ does not play a role any more.
Consequently, the solution (deprived of the rigid body motion) of eq.(\ref{eq:equi_equa_Cos}) is
\begin{equation}
v(x_2) = 
\frac{\boldsymbol{\kappa} \, \lambda_{\mbox{\tiny macro}}}{2 \lambda_{\mbox{\tiny macro}}+4 \mu_{\mbox{\tiny macro}}} \, x_2^2 + c_2 \, x_2
= 
\frac{\boldsymbol{\kappa} \, \nu_{\mbox{\tiny macro}}}{2 \left(1 - \nu_{\mbox{\tiny macro}}\right)} \, x_2^2 + c_2 \, x_2
\, .
\label{eq:sol_fun_disp_Cos}
\end{equation}

\begin{figure}[H]
\centering
\includegraphics[width=0.5\linewidth]{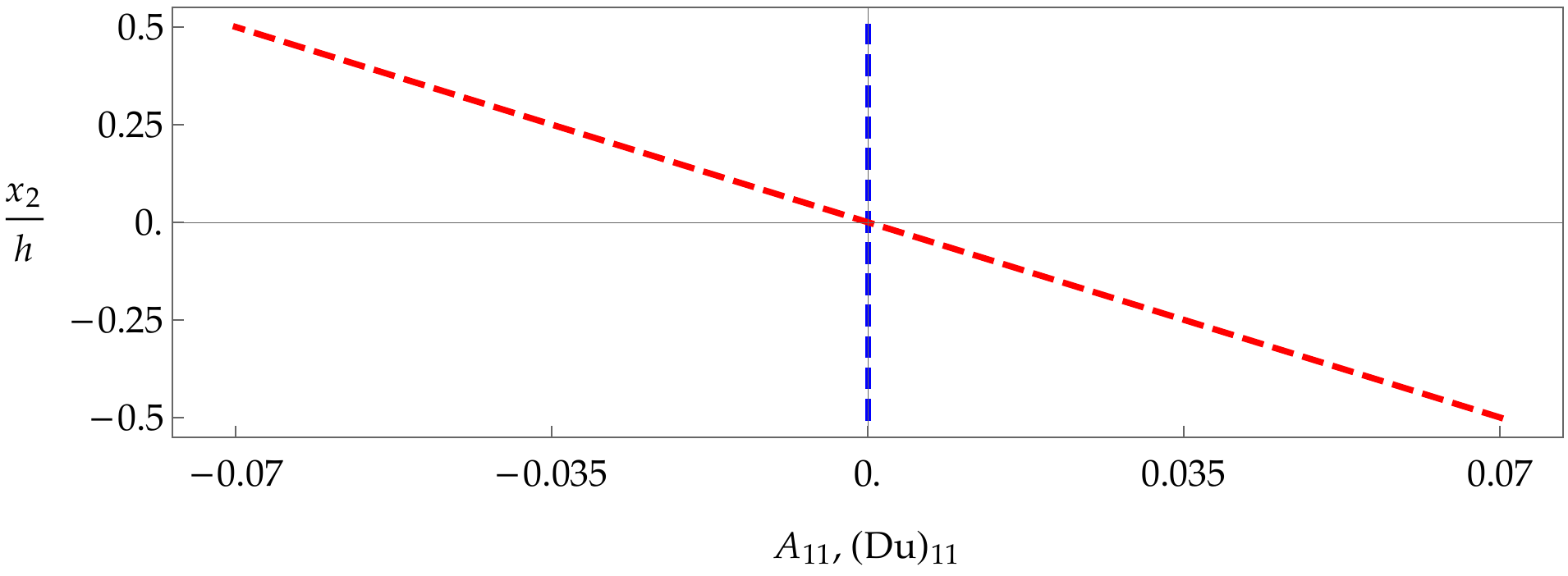}
\caption{(\textbf{Cosserat model}) Plot of $A_{11}$ (blue line) and $\left(\boldsymbol{\mbox{D} u}\right)_{11}$ (red line) across the thickness.
}
\label{fig:P11_plot_Coss}
\end{figure}

\begin{figure}[H]
\centering
\includegraphics[width=0.5\linewidth]{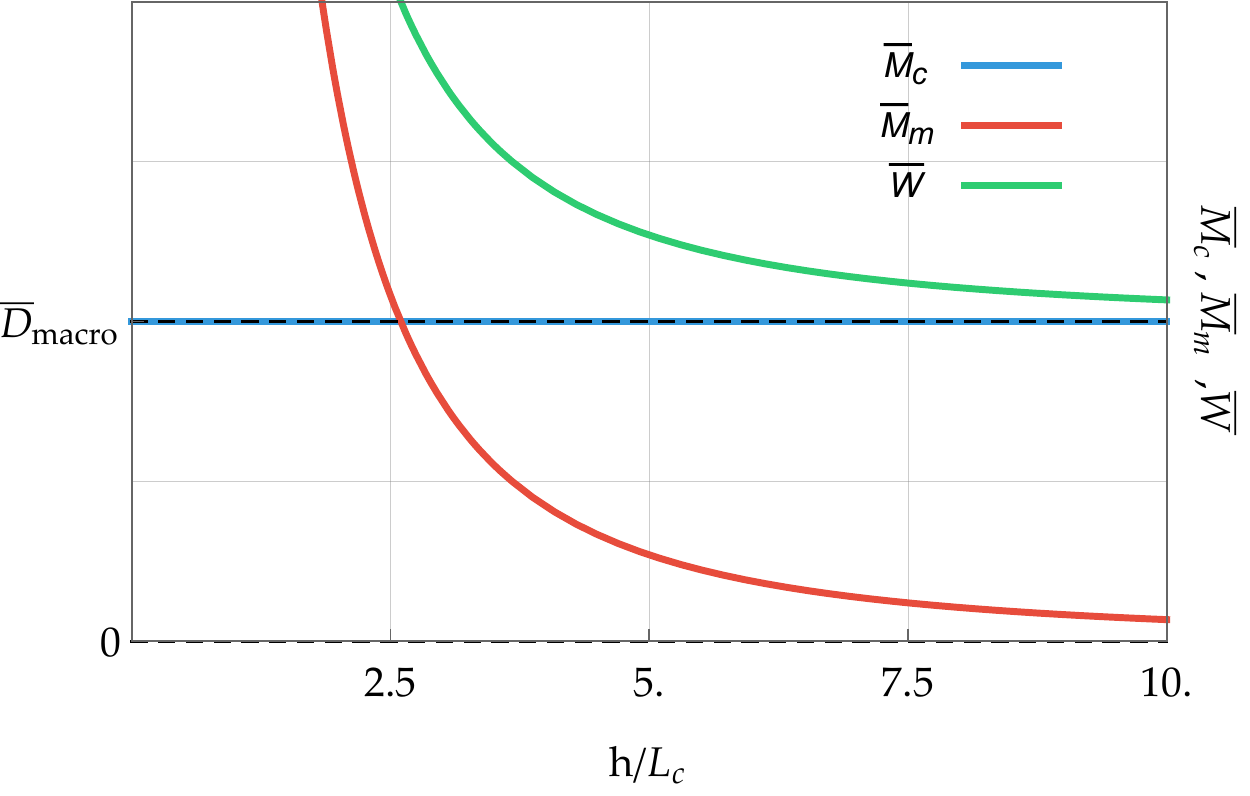}
\caption{(\textbf{Cosserat model}) Bending moments and strain energy while varying $L_c$. Note the singularity of the bending stiffness as $L_c \to \infty$ ($h\to 0$). The values of the parameters used are: $\mu _{\tiny \mbox{macro}} = 1$, $\lambda _{\tiny \mbox{macro}} = 1$, $\mu = 1$, $a _1 = 2$, $a _2 = 1$.}
\label{fig:all_plot_Cos}
\end{figure}

The boundary conditions eq.(\ref{eq:BC_Cos_gen})$_1$ on the upper and lower surfaces constrain $c_2 = 0$, while the second one eq.(\ref{eq:BC_Cos_gen})$_2$ is identically satisfied.
Then, the displacement and micro-rotation fields solution are
\begin{equation}
u_1(x_2) = - \boldsymbol{\kappa} \, x_1 \, x_2 \, ,
\quad
u_2(x_2) = \frac{\boldsymbol{\kappa}}{2} \, \frac{\lambda_{\mbox{\tiny macro}}}{\lambda_{\mbox{\tiny macro}} + 2 \mu_{\mbox{\tiny macro}}} \, x_2^2 
+ \frac{\boldsymbol{\kappa}}{2}x_1^2 \, ,
\quad
A_{12} (x_1) = - A_{21} (x_1) = - \boldsymbol{\kappa} \, x_1 \, .
\label{eq:disp_P_BC_Cos}
\end{equation}
The displacement field eq.(\ref{eq:disp_P_BC_Cos}) is the same as the classical one eq.~(\ref{eq:sol_disp_Cau}).
In Fig.~\ref{fig:P11_plot_Coss} we shown the plot of $A_{11}$ across the thickness.
The classical bending moment, the higher-order bending moment, and energy (per unit area d$x_1$d$x_3$) expressions are reported in the following eq.(\ref{eq:sigm_ene_dimensionless_Cos}).
\begin{align}
M_{\mbox{c}} (\boldsymbol{\kappa})&=
\displaystyle\int\limits_{-h/2}^{h/2}
\langle \boldsymbol{\widetilde{\sigma}} \boldsymbol{e}_1 , \boldsymbol{e}_1 \rangle  \, x_2 \, 
\mbox{d}x_{2}
= 
\frac{h^3 }{12} \, 
\frac{4\mu_{\mbox{\tiny macro}} \left(\lambda_{\mbox{\tiny macro}}+\mu_{\mbox{\tiny macro}}\right)}{\lambda_{\mbox{\tiny macro}}+2 \mu_{\mbox{\tiny macro}}} \, 
\boldsymbol{\kappa}
= D_{\mbox{\tiny macro}} \, \boldsymbol{\kappa} \, ,
\notag
\\*
M_{\mbox{m}}(\boldsymbol{\kappa}) &=
\displaystyle\int\limits_{-h/2}^{h/2}
\Big[\langle \mbox{skew} \left(\boldsymbol{m} \times \boldsymbol{e}_1 \right) \boldsymbol{e}_2 , \boldsymbol{e}_1 \rangle - \langle \mbox{skew} \left(\boldsymbol{m} \times \boldsymbol{e}_1 \right) \boldsymbol{e}_1 , \boldsymbol{e}_2 \rangle \Big] \,
\mbox{d}x_{2}
= 
h \,  \mu \, L_c^2 \, \frac{a_1 + a_2}{2} \, \boldsymbol{\kappa} \, ,
\label{eq:sigm_ene_dimensionless_Cos}
\\*
W_{\mbox{tot}} (\boldsymbol{\kappa})&=
\displaystyle\int\limits_{-h/2}^{+h/2} W \left( \boldsymbol{\mbox{D}u},\boldsymbol{A},\mbox{Curl}\boldsymbol{A} \right) \, \mbox{d}x_{2}
= 
\frac{1}{2}
\left[
\underbrace{
\frac{h^3 }{12}\frac{4\mu_{\mbox{\tiny macro}} \left(\lambda_{\mbox{\tiny macro}}+\mu_{\mbox{\tiny macro}}\right)}{\lambda_{\mbox{\tiny macro}}+2 \mu_{\mbox{\tiny macro}}}}_{D_{\mbox{\tiny macro}}
} +
12 \, \mu \, \left(\frac{L_c}{h}\right)^2 \, \frac{a_1 + a_2}{2} \frac{h^3 }{12}
\right]  \, \boldsymbol{\kappa}^2 
\, .
\notag
\end{align}
Finally, there is only one combination of curvature parameters $\gamma =: \frac{\mu \, L_c^2}{2}\frac{a_1+a_2}{2}$ appearing.
Inverting the formula gives $L_c = \sqrt{\frac{4\gamma}{\mu \left(a_1 + a_2\right)}}$ and this parameter $L_c$ is traditionally called the \textbf{bending length scale} of the Cosserat model (see Appendix \ref{app:lakes}). There is no counterpart to this observation in the relaxed micromorphic model or in the micro-stretch model.
Again, 
$
\frac{\mbox{d}}{\mbox{d}\boldsymbol{\kappa}} W_{\mbox{tot}}(\boldsymbol{\kappa}) = M_{\small \mbox{c}} (\boldsymbol{\kappa}) + M_{\small \mbox{m}} (\boldsymbol{\kappa}) \, .
$
The plot of the bending moments and the strain energy divided by $\frac{h^3}{12}\boldsymbol{\kappa}$ and $\frac{1}{2}\frac{h^3}{12}\boldsymbol{\kappa}^2$, respectively, while changing $L_c$ is shown in Fig.~\ref{fig:all_plot_Cos}.
\section{Micro-void model in dislocation tensor format}
\label{sec:Micro-void}
The expression of the strain energy for the isotropic micro-void continuum with a single curvature parameter in dislocation tensor format (3+1=4 dof's) can be written as:
\begin{align}
W \left(\boldsymbol{\mbox{D}u}, \omega ,\mbox{Curl}\,\left(\omega \boldsymbol{\mathbbm{1}}\right) \right) = &
\, \mu_{\tiny \mbox{macro}} \left\lVert \mbox{dev} \, \mbox{sym} \, \boldsymbol{\mbox{D}u}\right\rVert^{2}
+ \frac{\kappa_{e}}{2} \mbox{tr}^2 \left(\boldsymbol{\mbox{D}u} - \omega \boldsymbol{\mathbbm{1}} \right) 
+ \frac{\kappa_{\tiny \mbox{micro}}}{2} \mbox{tr}^2 \left(\omega \boldsymbol{\mathbbm{1}} \right)
\label{eq:energy_micro_void}
\\
&
+ \frac{\mu \,L_c^2 }{2} \,
a_2 \, \left\lVert \mbox{Curl} \, \left(\omega \boldsymbol{\mathbbm{1}}\right)\right\rVert^2 
\, .
\notag
\end{align}
Here, $\omega : \mathbb{R}^3 \to \mathbb{R}$ describes the additional scalar micro-void degree of freedom.
Since $\mbox{Curl} \left(\omega \boldsymbol{\mathbbm{1}}\right) \in \mathfrak{so}(3)$, the isotropic curvature reduces to 
$\frac{\mu \,L_c^2}{2} \, a_2 \, \left\lVert \mbox{Curl} \, \left(\omega \boldsymbol{\mathbbm{1}}\right)\right\rVert^2$.
The equilibrium equations without body forces are \!\!
\footnote{
Where
$\kappa_e = \frac{2\mu_e + 3\lambda_e}{3}$ and $\kappa_{\mbox{\tiny micro}} = \frac{2\mu_{\mbox{\tiny micro}} + 3\lambda_{\mbox{\tiny micro}}}{3}$ are he meso- and the micro-scale 3D bulk modulus.
}
\begin{align}
\mbox{Div}\overbrace{\left[
2\mu_{\tiny \mbox{macro}} \, \mbox{dev} \, \mbox{sym} \, \boldsymbol{\mbox{D}u}
+ \kappa_{e} \mbox{tr} \left(\boldsymbol{\mbox{D}u} - \omega \boldsymbol{\mathbbm{1}} \right) \boldsymbol{\mathbbm{1}}
\right]}^{\mathlarger{\widetilde{\sigma}}:=}
&= \boldsymbol{0},
\label{eq:equi_micro_void}
\\
\frac{1}{3}\mbox{tr}\left[\widetilde{\sigma}
- \kappa_{\tiny \mbox{micro}} \mbox{tr} \left(\omega \boldsymbol{\mathbbm{1}}\right) \boldsymbol{\mathbbm{1}}
- \mu \, L_{c}^{2} \, a_2 \, \mbox{Curl} \,
\mbox{Curl} \, \left(\omega \boldsymbol{\mathbbm{1}}\right)
\right] &= 0.
\notag
\end{align}
The boundary conditions at the upper and lower surface (free surface) are 
\begin{align}
\boldsymbol{\widetilde{t}}(x_2 = \pm \, h/2) &= 
\pm \, \boldsymbol{\widetilde{\sigma}}(x_2) \cdot \boldsymbol{e}_2 = 
\boldsymbol{0} \, ,
\label{eq:BC_micro_void}
\\
\eta(x_2 = \pm \, h/2) &= 
\pm \, \frac{1}{3}\mbox{tr}\left(\boldsymbol{m} (x_2) \cdot \boldsymbol{\epsilon} \cdot \boldsymbol{e}_2\right) = 
\pm \, \frac{1}{3}\mbox{tr}\left(\boldsymbol{m} (x_2) \times \boldsymbol{e}_2 \right) = 
\boldsymbol{0} \, .
\notag
\end{align}
According with the reference system shown in Fig.~\ref{fig:intro}, the ansatz for the displacement field and the function $\omega$ is
\begin{equation}
\boldsymbol{u}(x_1,x_2)=
\left(
\begin{array}{c}
-\kappa_1 \, x_1 x_2 \\
v(x_2)+\frac{\kappa_1  x_1^2}{2} \\
0 \\
\end{array}
\right) \, ,
\qquad
\omega\left(x_2\right) \boldsymbol{\mathbbm{1}} =
\left(
\begin{array}{ccc}
\omega(x_2) & 0 & 0 \\
0 & \omega(x_2) & 0 \\
0 & 0 & \omega(x_2) \\
\end{array}
\right) \, .
\label{eq:ansatz_micro_void}
\end{equation}
The equilibrium equations (\ref{eq:equi_micro_void}) turn in
\begin{align}
\kappa _e \left(v''(x_{2}) - \boldsymbol{\kappa} - 3 \, \omega '(x_{2})\right)
+ \frac{2}{3} \mu _{\tiny \mbox{macro}} \left(\boldsymbol{\kappa} +2 v''(x_{2})\right) = 0 \, , 
\label{eq:equi_equa_micro_void}
\\*
\frac{2}{3} \, a_{2} \, \mu \, L_c^2 \, \omega ''(x_{2})
+ \kappa _e \left(v'(x_{2}) - \boldsymbol{\kappa} \, x_{2}
- 3 \, \omega (x_{2})\right)
- 3 \, \kappa _{\mbox{\tiny micro}} \, \omega (x_{2}) = 0\, . 
\notag
\end{align}
Form eq.(\ref{eq:equi_equa_micro_void})$_1$ it is possible to evaluate $v''(x_2)$ and consequently $v'(x_2)$ as follows
\begin{align}
v''(x_2) &= \frac{3 \, \kappa _e-2 \mu _{\tiny \mbox{macro}}}{3 \, \kappa _e+4 \mu _{\tiny \mbox{macro}}} \, \boldsymbol{\kappa}
+ \frac{9 \, \kappa _e}{3 \, \kappa _e+4 \mu _{\tiny \mbox{macro}}} \, \omega'(x_{2})
\, ,
\label{eq:d1v_d2v_micro_void}
\\
v'(x_2)  &= \frac{3 \, \kappa _e-2 \mu _{\tiny \mbox{macro}}}{3 \, \kappa _e+4 \mu _{\tiny \mbox{macro}}} \, \boldsymbol{\kappa} \, x_{2}
+ \frac{9 \, \kappa _e}{3 \, \kappa _e+4 \mu _{\tiny \mbox{macro}}} \, \omega(x_{2})
+ c_0
\, .
\notag
\end{align}
Substituting back the expression of $v''(x_2)$ and $v'(x_2)$ in eq.(\ref{eq:equi_equa_micro_void}) it is possible to evaluate $\omega(x_2)$  from eq.(\ref{eq:equi_equa_micro_void})$_2$ and consequently $v(x_2)$ which implies
\begin{align}
\omega(x_2) = \,\, & 
c_1 \, e^{-\frac{f_{1} \, x_{2}}{L_c}} + c_2 \, e^{\frac{f_{1} \, x_{2}}{L_c}} - f_{2} \, \kappa \, x_{2} + f_{3} \, c_{0}
\, ,
\notag
\\
v(x_2) = \,\, &
\left(
c_2 \, e^{\frac{f_{1} x_{2}}{L_c}} - c_1 \, e^{-\frac{f_{1} x_{2}}{L_c}}
\right)
\frac{9 \kappa _e}{3 \kappa _e+4 \mu _{\tiny \mbox{macro}}}
\frac{L_c}{f_{1}}
+ \frac{(3-9 f_{2}) \kappa _e-2 \mu _{\tiny \mbox{macro}}}{6 \kappa _e+8 \mu _{\tiny \mbox{macro}}} \, \boldsymbol{\kappa} \, x_2^2
\label{eq:omega_micro_void}
\\
& + \left( \frac{9 f_{3} \, \kappa _e}{3 \kappa _e+4 \mu _{\tiny \mbox{macro}}}+1 \right) \, c_0 \, x_2
\, ,
\notag
\\
f_{1} := \,\, & \frac{1}{\sqrt{a_{2} \, \mu }}
\sqrt{\frac{9}{2}\left(\frac{4 \kappa _e \, \mu _{\tiny \mbox{macro}}}{3 \kappa _e+4 \mu _{\tiny \mbox{macro}}} + \kappa _{\mbox{\tiny micro}}\right)}
\, , \qquad\qquad
f_{2} := \,\,  \frac{2 \kappa _e \, \mu _{\tiny \mbox{macro}}}{4 \mu _{\tiny \mbox{macro}} \left(\kappa _e + \kappa _{\mbox{\tiny micro}}\right)+3 \kappa _e \,  \kappa _{\mbox{\tiny micro}}}
\, ,
\notag
\\
f_{3} := \,\, & \left(\frac{\kappa _e}{2 \mu _{\tiny \mbox{macro}}}+\frac{2}{3}\right) f_{2}
\, .
\notag
\end{align}

\begin{figure}[H]
	\centering
	\includegraphics[width=0.5\linewidth]{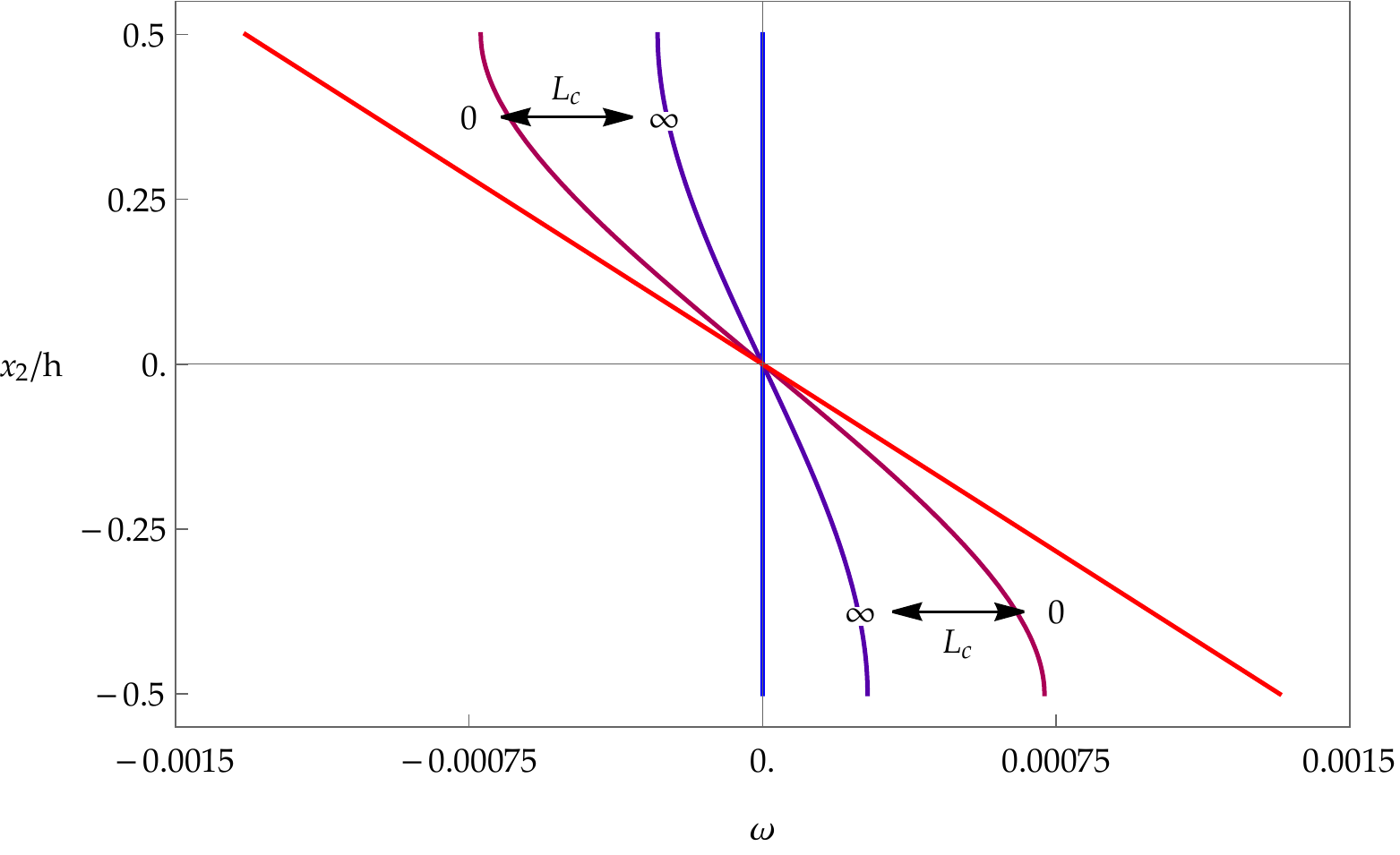}
	\caption{(\textbf{Micro-void model}) Distribution across the thickness of $\omega (x_2)$ while varying $L_c$.
		On the vertical axis we have the dimensionless thickness while on the horizontal we have the quantity $\omega(x_2)$.
		Notice that the red curve correspond to the limit for $L_c \to 0$ while the blue curve correspond to the limit for $L_c \to \infty$. The values of the parameters used are: $\mu _{e} = 2/3$, $\kappa _{e} = 7/9$, $\mu _{\tiny \mbox{micro}} = 2$, $\kappa _{\tiny \mbox{micro}} = 7/3$, $\mu = 1$, $a _1 = 1$, $\kappa=7/200$.}
	\label{fig:P11_plot_micro_void}
\end{figure}

In Fig.~\ref{fig:P11_plot_micro_void} we show the distribution across the thickness of $P_{11}$ while varying $L_c$.

\begin{figure}[H]
\centering
\includegraphics[width=0.5\linewidth]{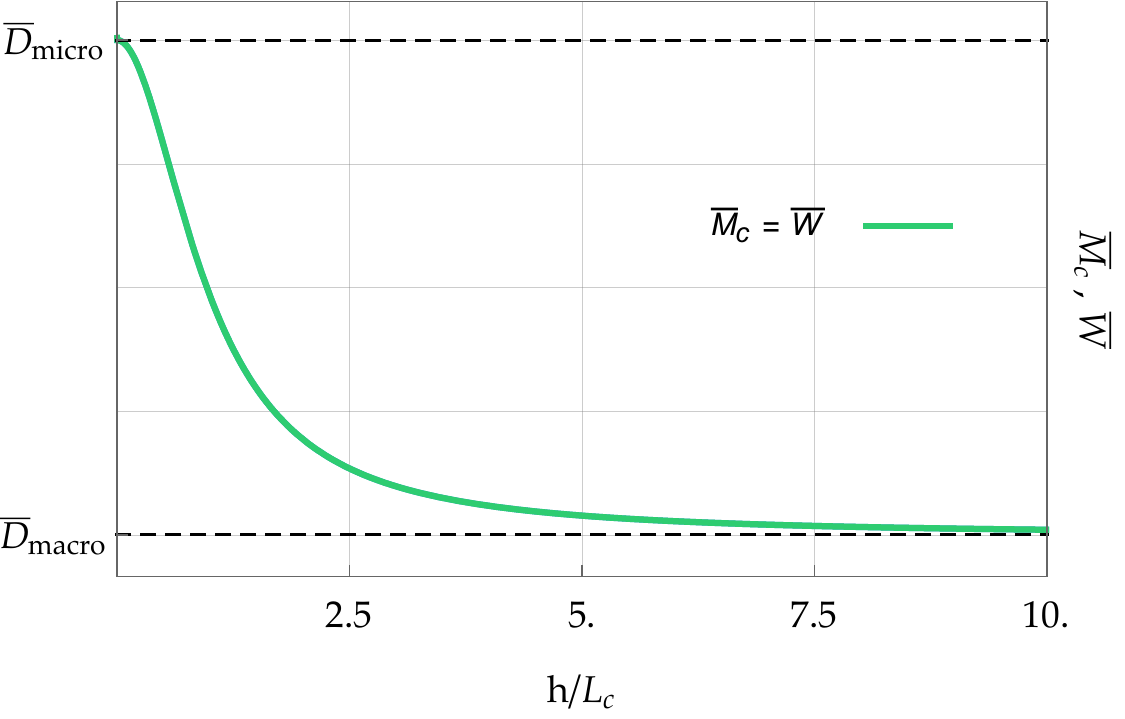}
\caption{(\textbf{Micro-void model}) Bending stiffness while varying $L_c$. Observe that the bending stiffness remains bounded as $L_c \to \infty$ ($h \to 0$). In this sense, the micro-voids model in dislocation format is consistent with the relaxed micromorphic model in cylindrical plate bending. The values of the parameters used are: $\mu _{e} = 1$, $\kappa _{e} = 5/3$, $\mu _{\tiny \mbox{micro}} = 1$, $\kappa _{\tiny \mbox{micro}} = 5/3$, $\mu = 1$, $a _1 = 1$.}
\label{fig:all_plot_micro_void}
\end{figure}

The classical bending moment, the higher-order bending moment, and energy (per unit area d$x_1$d$x_3$) expressions are reported in the following eq.(\ref{eq:sigm_ene_dimensionless_micro_void})
\begin{align}
M_{\mbox{c}} (\boldsymbol{\kappa})&:=
\displaystyle\int\limits_{-h/2}^{h/2}
\langle \boldsymbol{\widetilde{\sigma}} \boldsymbol{e}_1 , \boldsymbol{e}_1 \rangle  \, x_2 \, 
\mbox{d}x_{2}
= 
\frac{h^3}{12}
\frac{4 \mu _{\tiny \mbox{macro}} \left(3 \kappa _e+\mu _{\tiny \mbox{macro}}\right)}{3 \kappa _e+4 \mu _{\tiny \mbox{macro}}}
\left[
1 - \frac{9 f_{2} \, \kappa _e}{2 \left(3 \kappa _e+\mu _{\tiny \mbox{macro}}\right)}
\right.
\notag
\\*
&
\hspace{2cm}
\left.
+ \frac{54 \, \kappa _e}{3 \kappa _e+\mu _{\tiny \mbox{macro}}} \frac{f_{2}}{f_{1}^2} \left(\frac{L_c}{h}\right)^2 \, 
- \frac{108 \, \kappa _e}{3 \kappa _e+\mu _{\tiny \mbox{macro}}} \frac{f_{2}}{f_{1}^3} \left(\frac{L_c}{h}\right)^3 \, 
\tanh \left(\frac{f_{1} h}{2 L_c}\right)
\right]
\boldsymbol{\kappa}
\, ,
\notag
\\
M_{\mbox{m}}(\boldsymbol{\kappa}) &:=
\displaystyle\int\limits_{-h/2}^{h/2}
\mbox{tr} \left(\boldsymbol{m} \times \boldsymbol{e}_1\right)
\, \mbox{d}x_{2}
= 
0
\, ,
\label{eq:sigm_ene_dimensionless_micro_void}
\\
W_{\mbox{tot}} (\boldsymbol{\kappa})&:=
\displaystyle\int\limits_{-h/2}^{+h/2} W \left(\boldsymbol{\mbox{D}u}, \omega ,\mbox{Curl}\,\left(\omega \boldsymbol{\mathbbm{1}}\right) \right) \, \mbox{d}x_{2}
= 
\frac{1}{2}
\frac{h^3}{12}
\frac{4 \mu _{\tiny \mbox{macro}} \left(3 \kappa _e+\mu _{\tiny \mbox{macro}}\right)}{3 \kappa _e+4 \mu _{\tiny \mbox{macro}}}
\left[
1 - \frac{9 f_{2} \, \kappa _e}{2 \left(3 \kappa _e+\mu _{\tiny \mbox{macro}}\right)}
\right.
\notag
\\*
&
\hspace{2cm}
\left.
+ \frac{54 \, \kappa _e}{3 \kappa _e+\mu _{\tiny \mbox{macro}}} \frac{f_{2}}{f_{1}^2} \left(\frac{L_c}{h}\right)^2 \, 
- \frac{108 \, \kappa _e}{3 \kappa _e+\mu _{\tiny \mbox{macro}}} \frac{f_{2}}{f_{1}^3} \left(\frac{L_c}{h}\right)^3 \, 
\tanh \left(\frac{f_{1} h}{2 L_c}\right)
\right]
\boldsymbol{\kappa}^2
\, .
\notag
\end{align}
Since the higher-order bending moment is zero and
$
\frac{\mbox{d}}{\mbox{d}\boldsymbol{\kappa}}W_{\mbox{tot}}(\boldsymbol{\kappa}) = M_{\mbox{c}} (\boldsymbol{\kappa}) + M_{\mbox{m}} (\boldsymbol{\kappa})
= M_{\mbox{c}} (\boldsymbol{\kappa}) \, ,
$
only the plot of energy divided by $\frac{1}{2}\frac{h^3}{12}\boldsymbol{\kappa}^2$ while changing $L_c$ is shown in Fig.~\ref{fig:all_plot_micro_void}.
\subsubsection{Limit cases}
If $\mu_{\tiny \mbox{micro}} \to \infty$ (which is the consistent case for the micro-void model) then $\mu_{\tiny \mbox{macro}} \to \mu_e$ (see the homogenization formulas eq.(\ref{eq:static_homo_relation})) and we obtain
\begin{align}
\displaystyle\lim_{L_c\to0} M_{\mbox{c}} (\boldsymbol{\kappa})
&=
\frac{h^3}{12}
\frac{4 \mu _{\tiny \mbox{macro}} \left(3 \kappa _{\mbox{\tiny macro}}+\mu _{\tiny \mbox{macro}}\right)}{3 \kappa _{\mbox{\tiny macro}}+4 \mu _{\tiny \mbox{macro}}}
\boldsymbol{\kappa}
=
D_{\tiny\mbox{macro}} \, \boldsymbol{\kappa}
\, ,
\notag
\\
\displaystyle\lim_{L_c\to\infty} M_{\mbox{c}} (\boldsymbol{\kappa})
&= 
\frac{h^3}{12}
\frac{4 \mu _{\tiny \mbox{macro}} \left(3 \kappa _e+\mu _{\tiny \mbox{macro}}\right)}{3 \kappa _e+4 \mu _{\tiny \mbox{macro}}}
\boldsymbol{\kappa}
= 
\frac{h^3}{12}
\frac{4 \mu _e \left(3 \kappa _e+\mu _e\right)}{3 \kappa _e+4 \mu _e}
\boldsymbol{\kappa}
= 
D_{e} \, \boldsymbol{\kappa}
\, ,
\\
\displaystyle\lim_{L_c\to0} W_{\mbox{tot}} (\boldsymbol{\kappa})
&= 
\frac{1}{2}
\frac{h^3}{12}
\frac{4 \mu _{\tiny \mbox{macro}} \left(3 \kappa _{\mbox{\tiny macro}}+\mu _{\tiny \mbox{macro}}\right)}{3 \kappa _{\mbox{\tiny macro}}+4 \mu _{\tiny \mbox{macro}}}
\boldsymbol{\kappa}^2
=
\frac{1}{2} D_{\tiny\mbox{macro}} \, \boldsymbol{\kappa}^2 \, ,
\notag
\\
\displaystyle\lim_{L_c\to\infty} W_{\mbox{tot}} (\boldsymbol{\kappa})
&=
\frac{1}{2}
\frac{h^3}{12}
\frac{4 \mu _{\tiny \mbox{macro}} \left(3 \kappa _e+\mu _{\tiny \mbox{macro}}\right)}{3 \kappa _e+4 \mu _{\tiny \mbox{macro}}}
\boldsymbol{\kappa}  ^2
=
\frac{1}{2}
\frac{h^3}{12}
\frac{4 \mu _e \left(3 \kappa _e+\mu _e\right)}{3 \kappa _e+4 \mu _e}
\boldsymbol{\kappa}^2 
= 
\frac{1}{2} \,D_{e} \, \boldsymbol{\kappa}^2 \, .
\notag
\end{align}

\begin{itemize}
\item If $L_c \to 0$, from eq.(\ref{eq:equi_micro_void})$_2$ it is possible to evaluate 
\begin{equation}
\omega = \frac{1}{3}\frac{\kappa_{e}}{\kappa_{e} + \kappa_{\mbox{\tiny micro}}} \mbox{tr}\left( \boldsymbol{\mbox{D}u}\right)\, ,
\end{equation}
and inserting back $\omega \boldsymbol{\mathbbm{1}}$ in eq.(\ref{eq:equi_micro_void})$_1$, the expressions of the macro stiffness coefficients result to be 
\begin{equation}
\mu_{\mbox{\tiny macro}} = \mu_{e} \, , \qquad\qquad \kappa_{\mbox{\tiny macro}} = \frac{\kappa_e \, \kappa_{\mbox{\tiny micro}}}{\kappa_e + \kappa_{\mbox{\tiny micro}}} \, ;
\end{equation}
\item If $L_c \to \infty$, in order for a minimum of the energy to exist, it is required that $\mbox{Curl}\left(\omega \, \boldsymbol{\mathbbm{1}}\right) = 0$ which means that $\omega \, \boldsymbol{\mathbbm{1}}$ has to be the gradient of a scalar function $\zeta$, which implies that $\omega(x_2)$ must be constant.
The boundary conditions eq.(\ref{eq:BC_micro_void})$_2$ requires $\omega'(x_2)$ to be zero on the upper and lower surface and since $\omega \, \boldsymbol{\mathbbm{1}}$ must be constant, this is automatically satisfied.
\end{itemize}

While the micro-void model in dislocation format has bounded stiffness in cylindrical plate bending, the micro-stretch model does not. This may look strange given that the micro-stretch model should be more "flexible" having more independent degrees of freedom (the additional micro-rotation $\boldsymbol{A}\in \mathfrak{so}(3)$).
However, the specific coupling with the curvature terms dictate otherwise. This also shows that the micro-void model is not a simple ``penalization'' of the micro-stretch model.
\section{Cylindrical bending for the isotropic couple stress continuum family}
The indeterminate couple stress models \cite{neff2015correct} appear by letting formally the Cosserat couple modulus $\mu_c \to \infty$. This implies the constraint $\boldsymbol{A}=\mbox{skew} \, \boldsymbol{\mbox{D}u} \in \mathfrak{so}(3)$.
It is highlighted that this constraint is automatically satisfied given the ansatz eq.(\ref{eq:ansatz_Cos}) and the equilibrium equation's requirement that $\kappa_{1} = \kappa_{2}$ (which also implies that the solution is finally independent of the Cosserat couple modulus $\mu_c$), setting the energy for the isotropic Cosserat model equal to the one of the isotropic indeterminate couple stress model.

The indeterminate couple stress elastic energy turns into
\begin{align}
W \left(\boldsymbol{\mbox{D}u},\mbox{Curl}\,\mbox{skew} \, \boldsymbol{\mbox{D}u}\right) = &
\, \mu_{\mbox{\tiny macro}} \left\lVert \mbox{sym} \, \boldsymbol{\mbox{D}u} \right\rVert^{2}
+ \frac{\lambda_{\mbox{\tiny macro}}}{2} \mbox{tr}^2 \left(\boldsymbol{\mbox{D}u} \right)
+ \frac{\mu \, L_c^2}{2}
\left(
a_2 \, \left \lVert \mbox{skew} \, \mbox{Curl} \, \mbox{skew} \, \boldsymbol{\mbox{D}u}\right \rVert^2
\right.
\label{eq:energy_Ind_couple}
\\
&
\left.
+ a_1 \, \left \lVert \mbox{dev} \, \mbox{sym} \, \mbox{Curl} \, \mbox{skew} \, \boldsymbol{\mbox{D}u}\right \rVert^2
+ \frac{a_3}{3} \, \mbox{tr}^2 \left(\mbox{Curl} \, \mbox{skew} \, \boldsymbol{\mbox{D}u} \right)
\right)  \, .
\notag
\end{align}
Since, however, $\mbox{tr}(\mbox{Curl} \, \mbox{skew} \, \boldsymbol{\mbox{D}u}) = 0$ and $\left \lVert \mbox{dev} \, \mbox{sym} \, \mbox{Curl} \, \mbox{skew} \, \boldsymbol{\mbox{D}u}\right \rVert^2 = \left \lVert \mbox{sym} \, \mbox{Curl} \, \mbox{skew} \, \boldsymbol{\mbox{D}u}\right \rVert^2 = $ 
\\$\left \lVert \mbox{skew} \, \mbox{Curl} \, \mbox{skew} \, \boldsymbol{\mbox{D}u}\right \rVert^2 = \frac{1}{2}\left \lVert \mbox{Curl} \, \mbox{skew} \, \boldsymbol{\mbox{D}u}\right \rVert^2$ for a plane problem, we are left with
\begin{align}
W \left(\boldsymbol{\mbox{D}u},\mbox{Curl}\,\mbox{skew} \, \boldsymbol{\mbox{D}u}\right) = &
\, \mu_{\mbox{\tiny macro}} \left\lVert \mbox{sym} \, \boldsymbol{\mbox{D}u} \right\rVert^{2}
+ \frac{\lambda_{\mbox{\tiny macro}}}{2} \mbox{tr}^2 \left(\boldsymbol{\mbox{D}u} \right)
\notag
\\
&
\hspace{1cm}
+ \frac{\mu \, L_c^2}{2}
\left(
a_1 \, \left \lVert \mbox{sym} \, \mbox{Curl} \, \mbox{skew} \, \boldsymbol{\mbox{D}u}\right \rVert^2 \, 
+ a_2 \, \left \lVert \mbox{skew} \, \mbox{Curl} \, \mbox{skew} \, \boldsymbol{\mbox{D}u}\right \rVert^2
\right)  \, ,
\label{eq:energy_Ind_couple_2}
\\
= &
\, \mu_{\mbox{\tiny macro}} \left\lVert \mbox{sym} \, \boldsymbol{\mbox{D}u} \right\rVert^{2}
+ \frac{\lambda_{\mbox{\tiny macro}}}{2} \mbox{tr}^2 \left(\boldsymbol{\mbox{D}u} \right)
+ \frac{\mu \, L_c^2}{2} \, 
\frac{a_1 + a_2}{2} \, \left \lVert \mbox{Curl} \, \mbox{skew} \, \boldsymbol{\mbox{D}u}\right \rVert^2 \, .
\notag
\end{align}
The equilibrium equations without body forces are
\begin{align}
\mbox{Div}\left[2\mu _e \mbox{sym} \, \boldsymbol{\mbox{D}u} + 
\lambda _e \mbox{tr} \left(\boldsymbol{\mbox{D}u}\right) \boldsymbol{\mathbbm{1}}
+ \mu \, L_c^2 \, \mbox{skew} \, \mbox{Curl}\,\Big (
a_1 \, \mbox{dev} \, \mbox{sym} \, \mbox{Curl} \, \mbox{skew} \, \boldsymbol{\mbox{D}u} \, 
\right.
\hspace{3cm}
\label{eq:equi_Ind_couple}
\\
\pushright{\left.\left.
+ a_2 \, \mbox{skew} \, \mbox{Curl} \, \mbox{skew} \, \boldsymbol{\mbox{D}u} \, 
+ \frac{a_3}{3} \, \mbox{tr} \left(\mbox{Curl} \, \mbox{skew} \, \boldsymbol{\mbox{D}u} \right)\boldsymbol{\mathbbm{1}} \, 
\right)
\right] 
= \boldsymbol{0} \, ,}
\notag
\end{align}
while the boundary condition on the upper and lower surface are (for more details see \cite{neff2015correct})
\begin{align}
\pushleft{\boldsymbol{\widetilde{t}}(x_2 = \pm \, h/2) =
\pm \, \bigg\{\left(\boldsymbol{\widetilde{\sigma}} -\frac{1}{2} \mbox{Anti}\left(\mbox{Div} \, \boldsymbol{m} \right)\right)\cdot \boldsymbol{e}_2 
-\frac{1}{2} \boldsymbol{e}_2 \times \boldsymbol{\mbox{D}} \left[\left\langle \boldsymbol{e}_2,\mbox{sym} \, \boldsymbol{m}\cdot \boldsymbol{e}_2 \right\rangle\right]}
\notag
\\
-\frac{1}{2}\boldsymbol{\mbox{D}}\left[\mbox{Anti} \, \left(\left(\boldsymbol{\mathbbm{1}} - \boldsymbol{e}_2 \otimes \boldsymbol{e}_2\right)\cdot \boldsymbol{m}\cdot \boldsymbol{e}_2 \right)\cdot\left(\boldsymbol{\mathbbm{1}} - \boldsymbol{e}_2 \otimes \boldsymbol{e}_2\right)\right]:\left(\boldsymbol{\mathbbm{1}} - \boldsymbol{e}_2 \otimes \boldsymbol{e}_2\right)\bigg\}
=&\,
\boldsymbol{0} \, ,
\label{eq:BC_Ind_couple}
\\
\left(\boldsymbol{\mathbbm{1}} - \boldsymbol{e}_2 \otimes \boldsymbol{e}_2\right) \cdot \boldsymbol{\eta}(x_2 = \pm \, h/2) =
\pm \,  
\left(\boldsymbol{\mathbbm{1}} - \boldsymbol{e}_2 \otimes \boldsymbol{e}_2\right) \cdot 
\mbox{Anti} \left[ \left(\boldsymbol{\mathbbm{1}} - \boldsymbol{e}_2 \otimes \boldsymbol{e}_2\right)\cdot \boldsymbol{m}\cdot \boldsymbol{e}_2  \right]\cdot \boldsymbol{e}_2  =& \,
\boldsymbol{0} \, ,
\notag
\\
\boldsymbol{\pi}(x_2 = \pm \, h/2) = 
\pm \,  
\left(
\mbox{Anti}\left[\left(\boldsymbol{\mathbbm{1}} - \boldsymbol{e}_2 \otimes \boldsymbol{e}_2\right)\cdot \boldsymbol{m}\cdot \boldsymbol{e}_2\right]^+
-
\mbox{Anti}\left[\left(\boldsymbol{\mathbbm{1}} - \boldsymbol{e}_2 \otimes \boldsymbol{e}_2\right)\cdot \boldsymbol{m}\cdot \boldsymbol{e}_2\right]^-
\right)\cdot \boldsymbol{e}_1  =& \,
\boldsymbol{0} \, ,
\notag
\end{align}
where $\boldsymbol{\widetilde{\sigma}} = 2\mu _e \mbox{sym} \, \boldsymbol{\mbox{D}u} + 
\lambda _e \mbox{tr} \left(\boldsymbol{\mbox{D}u}\right) \boldsymbol{\mathbbm{1}}$, $\boldsymbol{e}_2$ is the unit vector aligned to the $x_2$-direction, and the second order moment stress  $\boldsymbol{m} = \mu \, L_c^2 \, \left (a_1 \, \mbox{dev} \, \mbox{sym} \, \mbox{Curl} \, \mbox{skew} \, \boldsymbol{\mbox{D}u} \, \right.$ $+$ $ a_2 \, \mbox{skew} \, \mbox{Curl} \, \mbox{skew} \, \boldsymbol{\mbox{D}u} $ $+$ $\left. \frac{a_3}{3} \, \mbox{tr} \left(\mbox{Curl} \, \mbox{skew} \, \boldsymbol{\mbox{D}u} \right)\boldsymbol{\mathbbm{1}} \right)$.
The operator Anti is defined as the inverse of axl in the context of eq.(\ref{eq:Aanti_axl}).
The term
$(
\mbox{Anti}\left[\left(\boldsymbol{\mathbbm{1}} - \boldsymbol{e}_2 \otimes \boldsymbol{e}_2\right)\cdot \boldsymbol{m}\cdot \boldsymbol{e}_2\right]^+
-
\mbox{Anti}\left[\left(\boldsymbol{\mathbbm{1}} - \boldsymbol{e}_2 \otimes \boldsymbol{e}_2\right)\cdot \boldsymbol{m}\cdot \boldsymbol{e}_2\right]^-
)$
measures the discontinuity of
$\mbox{Anti}\left[\left(\boldsymbol{\mathbbm{1}} - \boldsymbol{e}_2 \otimes \boldsymbol{e}_2\right)\cdot \boldsymbol{m}\cdot \boldsymbol{e}_2\right]$ across the boundary.

According to the reference system shown in Fig.~\ref{fig:intro}, the ansatz for the displacement field and consequently the gradient of the displacement are
\begin{equation}
\boldsymbol{u}(x_1,x_2)=
\left(
\begin{array}{c}
-\boldsymbol{\kappa} \, x_1 x_2 \\
v(x_2)+\frac{\boldsymbol{\kappa}  x_1^2}{2} \\
0 \\
\end{array}
\right) \, ,
\qquad
\boldsymbol{\mbox{D}u} = 
\left(
\begin{array}{ccc}
-\boldsymbol{\kappa} \, x_2 & - \boldsymbol{\kappa} \, x_1 & 0 \\
\boldsymbol{\kappa} \, x_1 & v'(x_2) & 0 \\
0 & 0 & 0 \\
\end{array}
\right) \, .
\label{eq:ansatz_Ind_couple}
\end{equation}

Substituting the ansatz eq.(\ref{eq:ansatz_Ind_couple}) in eq.(\ref{eq:equi_Ind_couple}) the equilibrium equation results to be:
\begin{equation}
-\boldsymbol{\kappa} \lambda_{\mbox{\tiny macro}} + \left(\lambda_{\mbox{\tiny macro}}+2 \mu_{\mbox{\tiny macro}}\right) v''(x_{2}) = 0 \, ,
\label{eq:equi_equa_Ind_couple}
\end{equation}

Consequently, the solution (deprived of the rigid body motion) of eq.(\ref{eq:equi_equa_Ind_couple}) is:
\begin{equation}
v(x_2) = 
\frac{\boldsymbol{\kappa} \, \lambda_{\mbox{\tiny macro}}}{2 \lambda_{\mbox{\tiny macro}}+4 \mu_{\mbox{\tiny macro}}} \, x_2^2 + c_2 \, x_2
= 
\frac{\boldsymbol{\kappa} \, \nu_{\mbox{\tiny macro}}}{2 \left(1 - \nu_{\mbox{\tiny macro}}\right)} \, x_2^2 + c_2 \, x_2
\, .
\label{eq:sol_fun_disp_Ind_couple}
\end{equation}

\begin{figure}[H]
\centering
\includegraphics[width=0.5\linewidth]{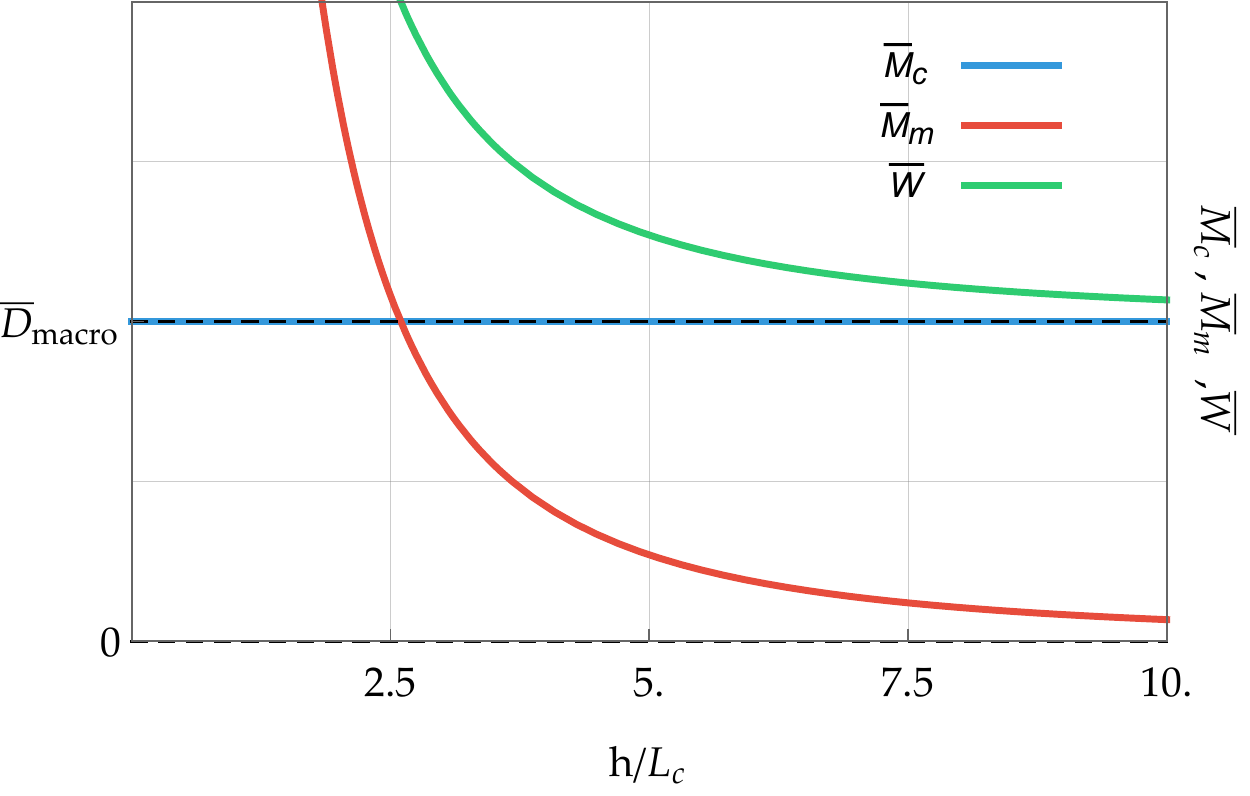}
\caption{(\textbf{Couple stress model}) Bending moments and strain energy while varying $L_c$ with $a_1=a_2=1$. Note the singularity of the bending stiffness as $L_c \to \infty$ ($h\to 0$). The values of the parameters used are: $\mu _{\tiny \mbox{macro}} = 1$, $\lambda _{\tiny \mbox{macro}} = 1$, $\mu = 1$, $a _1 = 2$, $a _2 = 1$.}
\label{fig:all_plot_Ind_couple}
\end{figure}

The boundary conditions eq.(\ref{eq:BC_Ind_couple})$_1$ on the upper and lower surfaces constrain $c_2 = 0$, while eq.(\ref{eq:BC_Ind_couple})$_2$ and eq.(\ref{eq:BC_Ind_couple})$_3$ are identically satisfied.
Then, the displacement field solution results in
\begin{equation}
u_1(x_2) = - \boldsymbol{\kappa} \, x_1 \, x_2 \, ,
\qquad
u_2(x_2) = \frac{\boldsymbol{\kappa}}{2} \, \frac{\lambda_{\mbox{\tiny macro}}}{\lambda_{\mbox{\tiny macro}} + 2 \mu_{\mbox{\tiny macro}}} \, x_2^2 
+ \frac{\boldsymbol{\kappa}}{2}x_1^2  \, .
\label{eq:disp_P_BC_Ind_couple}
\end{equation}
The displacement field eq.(\ref{eq:disp_P_BC_Ind_couple}) is the same as the classical one eq.(\ref{eq:sol_disp_Cau}).
In Fig.~\ref{fig:P11_plot_couple} we show the plot of $\left(\boldsymbol{\mbox{D}u}\right)_{11}$ across the thickness:
\begin{figure}[H]
	\centering
	\includegraphics[width=0.5\linewidth]{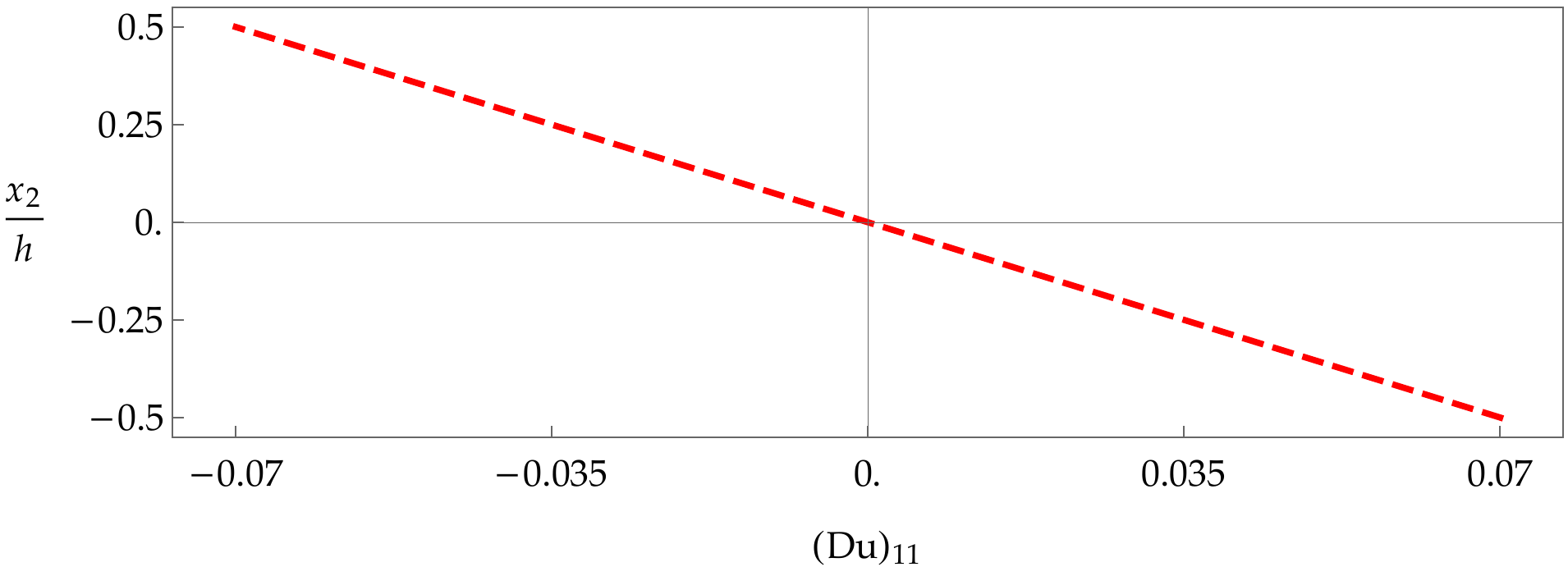}
	\caption{(\textbf{Couple stress model}) Plot of $\left(\boldsymbol{\mbox{D}u}\right)_{11}$ across the thickness.
	}
	\label{fig:P11_plot_couple}
\end{figure}
The classical bending moment, the higher-order bending moment, and energy (per unit area d$x_1$d$x_3$) expressions are reported in the following eq.(\ref{eq:sigm_ene_dimensionless_Ind_couple})
\begin{align}
M_{\mbox{c}} (\boldsymbol{\kappa})
&
=
\displaystyle\int\limits_{-h/2}^{h/2}
\langle \boldsymbol{\widetilde{\sigma}} \boldsymbol{e}_1 , \boldsymbol{e}_1 \rangle  \, x_2 \, 
\mbox{d}x_{2}
= 
\frac{h^3 }{12} \, 
\frac{4\mu_{\mbox{\tiny macro}} \left(\lambda_{\mbox{\tiny macro}}+\mu_{\mbox{\tiny macro}}\right)}{\lambda_{\mbox{\tiny macro}}+2 \mu_{\mbox{\tiny macro}}} \, 
\boldsymbol{\kappa} \, ,
\notag
\\
M_{\mbox{m}}(\boldsymbol{\kappa})
&
=
\displaystyle\int\limits_{-h/2}^{h/2}
\langle \left(\boldsymbol{m} \times \boldsymbol{e}_1 \right) \boldsymbol{e}_2 , \boldsymbol{e}_1 \rangle \,
\mbox{d}x_{2}
= 
\frac{h^3}{12} \,  \mu \, \left(\frac{L_c}{h}\right)^2 \, 12 \, \frac{a_1 + a_2}{2} \, \boldsymbol{\kappa} \, ,
\label{eq:sigm_ene_dimensionless_Ind_couple}
\\
W_{\mbox{tot}} (\boldsymbol{\kappa})
=
&
\displaystyle\int\limits_{-h/2}^{+h/2} W \left( \boldsymbol{\mbox{D}u},\mbox{Curl} \, \mbox{skew} \, \boldsymbol{\mbox{D}u} \right) \, \mbox{d}x_{2}
\notag
\\*
= 
&
\frac{1}{2}
\Bigg(
\underbrace{\frac{h^3 }{12}\frac{4\mu_{\mbox{\tiny macro}} \left(\lambda_{\mbox{\tiny macro}}+\mu_{\mbox{\tiny macro}}\right)}{\lambda_{\mbox{\tiny macro}}+2 \mu_{\mbox{\tiny macro}}}}_{D_{\mbox{\tiny macro}}}
+ 12 \, \mu \, \left(\frac{L_c}{h}\right)^2 \, \frac{a_1 + a_2}{2} \, \frac{h^3 }{12}
\Bigg)
\boldsymbol{\kappa}^2
\, .
\notag
\end{align}
The result coincides with the Cosserat solution eq.(\ref{eq:sigm_ene_dimensionless_Cos}).
Note that
$
\frac{\mbox{d}}{\mbox{d}\boldsymbol{\kappa}} W_{\mbox{tot}}(\boldsymbol{\kappa}) = M_{\small \mbox{c}} (\boldsymbol{\kappa}) + M_{\small \mbox{m}} (\boldsymbol{\kappa}) \, .
$
The plot of the bending moments and the strain energy divided by $\frac{h^3}{12}\boldsymbol{\kappa}$ and $\frac{1}{2}\frac{h^3}{12}\boldsymbol{\kappa}^2$, respectively, while changing $L_c$ is shown in Fig.~\ref{fig:all_plot_Ind_couple}.
\subsection{Cylindrical bending for the isotropic symmetric couple stress continuum}
The \textbf{modified couple stress model} \cite{munch2017modified} consist in choosing $a_1>0$, $a_2=0$ and the \textbf{(``pseudo'')-consistent couple stress model} \cite{hadjesfandiari2011couple} appears for $a_1=0$, $a_2>0$. Since $\mbox{Curl} \, \mbox{sym} \, \boldsymbol{\mbox{D}u} = - \mbox{Curl} \, \mbox{skew} \, \boldsymbol{\mbox{D}u}$ due to $\mbox{Curl} \, \boldsymbol{\mbox{D}u} = 0$ for a plane problem, the form of the energy remains the same.
For the cylindrical bending problem, the higher order Neumann boundary conditions are already identically satisfied for all these models, like for their general case which is the relaxed micromorphic model.
\section{Cylindrical bending for the classical isotropic micromorphic continuum without mixed terms}
The expression of the strain energy for the reduced isotropic micromorphic continuum without mixed terms (like $\langle\mbox{sym}\, \boldsymbol{P}, \mbox{sym} \, \left(\boldsymbol{\mbox{D}u} -\boldsymbol{P}\right)\rangle$, etc.) and simplified isotropic curvature can be written as
\begin{align}
W \left(\boldsymbol{\mbox{D}u}, \boldsymbol{P}, \boldsymbol{\mbox{D}P}\right)
= &
\, \mu_{e} \left\lVert \mbox{dev} \, \mbox{sym} \left(\boldsymbol{\mbox{D}u} - \boldsymbol{P} \right) \right\rVert^{2}
+ \frac{\kappa_{e}}{2} \mbox{tr}^2 \left(\boldsymbol{\mbox{D}u} - \boldsymbol{P} \right)
+ \mu_{c} \left\lVert \mbox{skew} \left(\boldsymbol{\mbox{D}u} - \boldsymbol{P} \right) \right\rVert^{2}
\notag
\\
&
+ \mu_{\tiny \mbox{micro}} \left\lVert \mbox{dev} \, \mbox{sym}\,\boldsymbol{P} \right\rVert^{2}
+ \frac{\kappa_{\tiny \mbox{micro}}}{2} \mbox{tr}^2 \left(\boldsymbol{P} \right)
\notag
\\
&
+ \frac{\mu \, L_c^2}{2}
\displaystyle\sum_{i=1}^{3}
\left[
a_1 \, \left\lVert \mbox{dev} \, \mbox{sym} \Big( \partial_{x_i} \boldsymbol{P} \Big) \right\rVert^2
+ a_2 \, \left\lVert \mbox{skew} \Big( \partial_{x_i} \boldsymbol{P} \Big) \right\rVert^2
+ \frac{2}{9} \, a_3 \, \mbox{tr}^2 \Big( \partial_{x_i} \boldsymbol{P} \Big)
\right]
\notag
\\
= &
\, \mu_{e} \left\lVert \mbox{dev} \, \mbox{sym} \left(\boldsymbol{\mbox{D}u} - \boldsymbol{P} \right) \right\rVert^{2}
+ \frac{\kappa_{e}}{2} \mbox{tr}^2 \left(\boldsymbol{\mbox{D}u} - \boldsymbol{P} \right)
+ \mu_{c} \left\lVert \mbox{skew} \left(\boldsymbol{\mbox{D}u} - \boldsymbol{P} \right) \right\rVert^{2}
\label{eq:energy_MM}
\\
&
+ \mu_{\tiny \mbox{micro}} \left\lVert \mbox{dev} \, \mbox{sym}\,\boldsymbol{P} \right\rVert^{2}
+ \frac{\kappa_{\tiny \mbox{micro}}}{2} \mbox{tr}^2 \left(\boldsymbol{P} \right)
\notag
\\
&
+ \frac{\mu \, L_c^2}{2}
\Bigg(
a_1 \, \left\lVert \mbox{dev} \, \mbox{sym} \, \boldsymbol{\mbox{D}P} \right\rVert^2
+ a_2 \, \left\lVert \mbox{skew} \, \boldsymbol{\mbox{D}P} \right\rVert^2
+ \frac{2}{9} \, a_3 \, \mbox{tr}^2 \Big( \boldsymbol{\mbox{D}P} \Big)
\Bigg)
\, .
\notag
\end{align}
The meaning of $\mbox{dev} \, \mbox{sym} \, \boldsymbol{\mbox{D}P}$, $\mbox{skew} \, \boldsymbol{\mbox{D}P}$, and $\mbox{tr} \left(\boldsymbol{\mbox{D}P}\right)$ for the third order tensor $\boldsymbol{\mbox{D}P}$ can be inferred from eq.(\ref{eq:energy_MM}), i.e. we define
\begin{align}
\left\lVert \mbox{dev} \, \mbox{sym} \, \boldsymbol{\mbox{D}P} \right\rVert^2 :=
\displaystyle\sum_{i=1}^{3}
\left\lVert \mbox{dev} \, \mbox{sym} \Big( \partial_{x_i} \boldsymbol{P} \Big) \right\rVert^2
\, ,
\qquad\qquad
\left\lVert \mbox{skew} \, \boldsymbol{\mbox{D}P} \right\rVert^2 :=
\displaystyle\sum_{i=1}^{3}
\left\lVert \mbox{skew} \Big( \partial_{x_i} \boldsymbol{P} \Big) \right\rVert^2
\, ,
\\
\mbox{tr}^2 \Big( \boldsymbol{\mbox{D}P} \Big) :=
\displaystyle\sum_{i=1}^{3}
\mbox{tr}^2 \Big( \partial_{x_i} \boldsymbol{P} \Big)
\, .
\hspace{4.5cm}
\notag
\end{align}
The equilibrium equations (see the Appendix~\ref{app:full_micro_and_strain_eq_eqa}) without body forces are the following 
\begin{align}
\mbox{Div}\overbrace{\left[2\mu_{e} \, \mbox{dev} \, \mbox{sym} \left(\boldsymbol{\mbox{D}u} - \boldsymbol{P} \right) + \kappa_{e} \mbox{tr} \left(\boldsymbol{\mbox{D}u} - \boldsymbol{P} \right) \boldsymbol{\mathbbm{1}}
+ 2\mu_{c}\,\mbox{skew} \left(\boldsymbol{\mbox{D}u} - \boldsymbol{P} \right)\right]}^{\mathlarger{\widetilde{\sigma}}:=}
= \boldsymbol{0} \, ,
\notag
\\
\widetilde{\sigma}
- 2 \mu_{\mbox{\tiny micro}} \, \mbox{dev} \,\mbox{sym}\,\boldsymbol{P}
- \kappa_{\tiny \mbox{micro}} \mbox{tr} \left(\boldsymbol{P}\right) \boldsymbol{\mathbbm{1}}
\hspace{9cm}
\label{eq:equiMic_MM_three}
\\
+\mu L_{c}^{2} \,
\left[
a_1 \, \mbox{dev} \, \mbox{sym} \, \boldsymbol{\Delta P}
+ a_2 \, \mbox{skew} \, \boldsymbol{\Delta P}
+ \frac{2}{9} \, a_3 \, \mbox{tr} \left(\boldsymbol{\Delta P}\right)\boldsymbol{\mathbbm{1}}
\right]
= \boldsymbol{0} \, ,
\notag
\end{align}
where $\boldsymbol{\Delta P} \in \mathbb{R}^{3\times3}$ is taken component-wise.
The boundary conditions (see the Appendix~\ref{app:full_micro_and_strain_eq_eqa}) at the upper and lower surface (free surface) are 
\begin{align}
\boldsymbol{\widetilde{t}}(x_2 = \pm \, h/2) &= 
\pm \, \boldsymbol{\widetilde{\sigma}}(x_2) \cdot \boldsymbol{e}_2 = 
\boldsymbol{0}_{\mathbb{R}^{3}}
\, ,
\label{eq:BC_MM_mono}
\\
\boldsymbol{\eta}(x_2 = \pm \, h/2) &= 
\pm \, \sum_{i=1}^{3}\boldsymbol{m}_{i} (x_2) \, (\boldsymbol{e}_{2})_{i} =
\boldsymbol{0}_{\mathbb{R}^{3\times3}}
\, ,
\notag
\end{align}
where $\boldsymbol{m}_{i} = a_1 \, \mbox{dev} \, \mbox{sym} \left( \partial_{x_i} \boldsymbol{P} \right)$ $+$ $a_2 \, \mbox{skew} \left( \partial_{x_i} \boldsymbol{P} \right)$ $+$ $\frac{2}{9} \, a_3 \, \mbox{tr} \left( \partial_{x_i} \boldsymbol{P} \right) \boldsymbol{\mathbbm{1}}$, $i=1,2,3$, is a second order tensor, and $(\boldsymbol{e}_{2})_{i}$ is the \textit{i}th component of the unit vector $\boldsymbol{e}_{2}$ (see the Appendix~\ref{app:full_micro_and_strain_eq_eqa}).
According with the reference system shown in Fig.~\ref{fig:intro}, the ansatz for the displacement field and the microdistortion is
\begin{equation}
\boldsymbol{u}(x_1,x_2)=
\left(
\begin{array}{c}
-\kappa_1 \, x_1 x_2 \\
v(x_2)+\frac{\kappa_1  x_1^2}{2} \\
0 \\
\end{array}
\right) \, ,
\qquad
\boldsymbol{P}\left(x_1,x_2\right) =
\left(
\begin{array}{ccc}
P_{11}(x_2) & -\boldsymbol{\kappa} \, x_1 & 0 \\
\boldsymbol{\kappa} \, x_1 & P_{22}(x_2) & 0 \\
0 & 0 & P_{33}(x_2) \\
\end{array}
\right) \, .
\label{eq:ansatz_MM_three}
\end{equation}
Substituting the ansatz eq.(\ref{eq:ansatz_MM_three}) in eq.(\ref{eq:equiMic_MM_three}) the equilibrium equation are
\begin{align}
2 \mu _e \left(v''(x_{2})-P_{22}'(x_{2})\right)-\lambda _e \left(\boldsymbol{\kappa} +P_{11}'(x_{2})+P_{22}'(x_{2})+P_{33}'(x_{2})-v''(x_{2})\right) &= 0 \, ,
\notag
\\*
\frac{1}{9} \mu \, L_c^2 \left(2 (3 a_{1}+a_{3}) P_{11}''(x_{2})-(3 a_{1}-2 a_{3}) \left(P_{22}''(x_{2})+P_{33}''(x_{2})\right)\right)&
\notag
\\*
-\lambda _e \left(P_{11}(x_{2})+P_{22}(x_{2})+P_{33}(x_{2})-v'(x_{2})+\boldsymbol{\kappa}  x_{2}\right)-2 \mu _e (P_{11}(x_{2})+\boldsymbol{\kappa}  x_{2})&
\notag
\\*
+\lambda _{\mbox{\tiny micro}} (-P_{11}(x_{2})-P_{22}(x_{2})-P_{33}(x_{2}))-2 \mu _{\mbox{\tiny micro}} P_{11}(x_{2}) &= 0
\, ,
\notag
\\*
-\frac{1}{9} \mu \, L_c^2 \left((3 a_{1}-2 a_{3}) P_{11}''(x_{2})-2 (3 a_{1}+a_{3}) P_{22}''(x_{2})+(3 a_{1}-2 a_{3}) P_{33}''(x_{2})\right)&
\label{eq:equi_equa_MM_three}
\\*
-\lambda _e \left(P_{11}(x_{2})+P_{22}(x_{2})+P_{33}(x_{2})-v'(x_{2})+\boldsymbol{\kappa}  x_{2}\right)+2 \mu _e \left(v'(x_{2})-P_{22}(x_{2})\right)&
\notag
\\*
+\lambda _{\mbox{\tiny micro}} (-P_{11}(x_{2})-P_{22}(x_{2})-P_{33}(x_{2}))-2 \mu _{\mbox{\tiny micro}} P_{22}(x_{2}) &= 0
\, ,
\notag
\\*
-\frac{1}{9} \mu \, L_c^2 \left((3 a_{1}-2 a_{3}) P_{11}''(x_{2})+(3 a_{1}-2 a_{3}) P_{22}''(x_{2})-2 (3 a_{1}+a_{3}) P_{33}''(x_{2})\right)&
\notag
\\*
-\lambda _e \left(P_{11}(x_{2})+P_{22}(x_{2})+P_{33}(x_{2})-v'(x_{2})+\boldsymbol{\kappa}  x_{2}\right)-2 \mu _e P_{33}(x_{2})&
\notag
\\*
+\lambda _{\mbox{\tiny micro}} (-P_{11}(x_{2})-P_{22}(x_{2})-P_{33}(x_{2}))-2 \mu _{\mbox{\tiny micro}} P_{33}(x_{2}) &= 0 \, .
\notag
\end{align}
It is possible to evaluate $v''(x_2)$ from eq.(\ref{eq:equi_equa_MM_three})$_1$ and then $P_{22}(x_2)$ form a linear combination of the remaining equation.
After substituting the expressions of $v''(x_2)$ and $P_{22}(x_2)$ in eq.(\ref{eq:equi_equa_MM_three}), we are left with two coupled ordinary differential equations of fourth order in $P_{11}(x_2)$ and $P_{33}(x_2)$.
The solution is again of hyperbolic type and the unknown coefficients are determined by the boundary conditions eq.(\ref{eq:BC_MM_mono}).

\begin{figure}[H]
	\centering
	\includegraphics[width=0.5\linewidth]{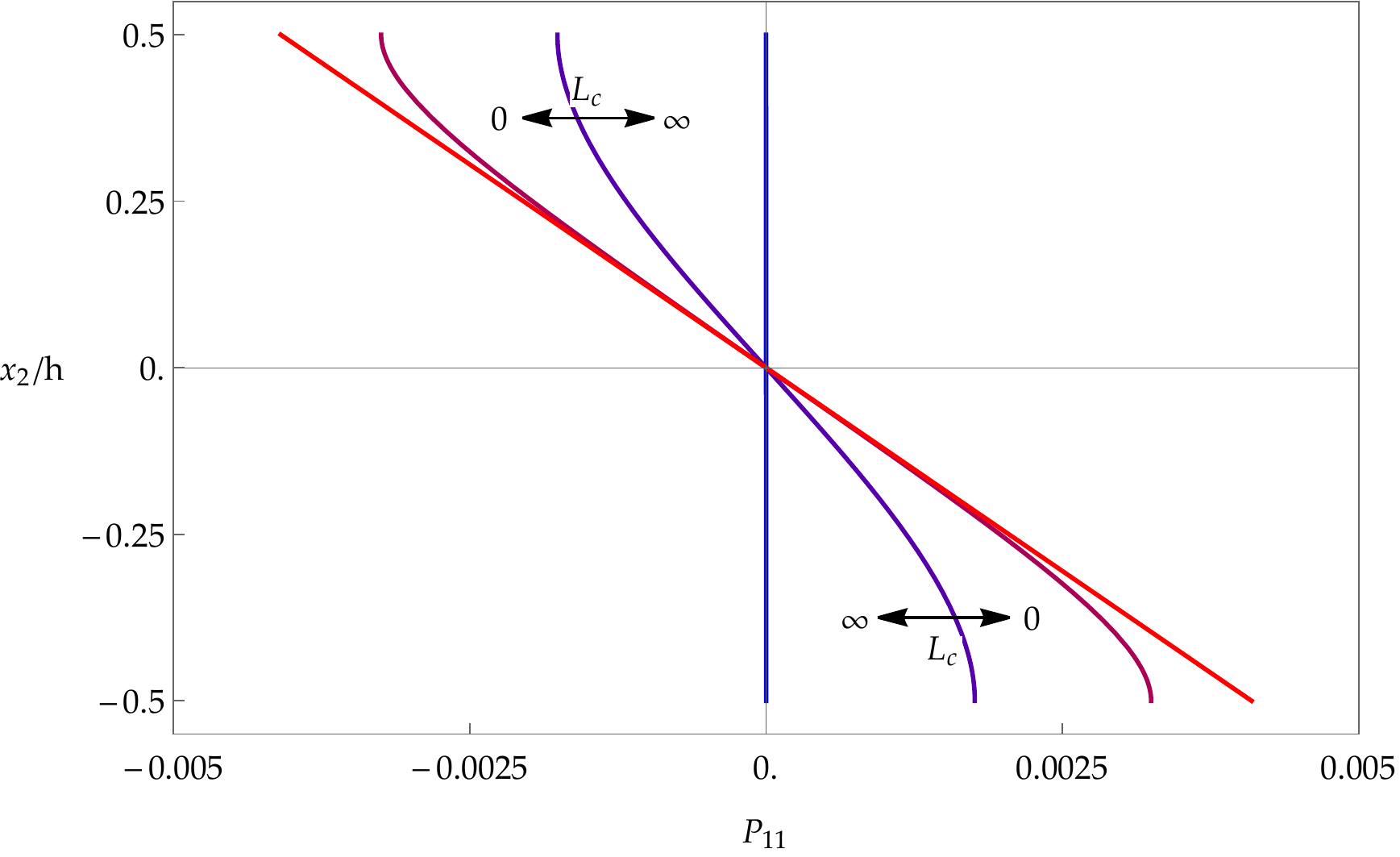}
	\caption{(\textbf{Micromorphic model}, classical case) Distribution across the thickness of $P_{11}$ while varying $L_c$.
		On the vertical axis we have the dimensionless thickness while on the horizontal we have the quantity $P_{11}$.
		Notice that the red curve correspond to the limit for $L_c \to 0$ while the blue curve correspond to the limit for $L_c \to \infty$.
		Again the response is opposite to the relaxed micromorphic model. $L_c \to \infty$ implies $P_{11} \to 0$ (blue line). The values of the parameters used are: $\mu _{e} = 1/3$, $\lambda _{e} = 1/8$, $\mu _{\tiny \mbox{micro}} = 2$, $\lambda _{\tiny \mbox{micro}} = 1$, $\mu = 1$, $a _1 = 1/2$, $a _2 = 1/2$, $a _3 = 3/2$, $\boldsymbol{\kappa}=7/200$..}
	\label{fig:P11_plot_MM_three}
\end{figure}

In Fig.~\ref{fig:P11_plot_MM_three} we show the distribution across the thickness of $P_{11}$ while varying $L_c$.

Subsequently, the bending moments and the strain-energy per unit area are computed by
\begin{align}
M_{\mbox{c}} (\boldsymbol{\kappa})
=
\displaystyle\int\limits_{-h/2}^{h/2}
\langle \boldsymbol{\widetilde{\sigma}} \boldsymbol{e}_1 , \boldsymbol{e}_1 \rangle  \, x_2 \, 
\mbox{d}x_{2}
\, ,
\qquad\qquad
M_{\mbox{m}}(\boldsymbol{\kappa})
=
\displaystyle\int\limits_{-h/2}^{h/2}
\langle \left(\sum_{i=1}^{3}\boldsymbol{m}_{i} (x_2) \cdot (\boldsymbol{e}_{1})_{i} \right) \boldsymbol{e}_2 , \boldsymbol{e}_1 \rangle \,
\mbox{d}x_{2}
\, ,
\label{eq:sigm_ene_dimensionless_MM_three}
\\
W_{\mbox{tot}} (\boldsymbol{\kappa})
=
\displaystyle\int\limits_{-h/2}^{+h/2} W \left( \boldsymbol{\mbox{D}u},\boldsymbol{P},\boldsymbol{\mbox{D}P} \right) \, \mbox{d}x_{2}
\, ,
\hspace{4cm}
\notag
\end{align}
where $(\boldsymbol{e}_{1})_{i}$ is the scalar \textit{i}th component of the unit vector $\boldsymbol{e}_{1}$, similarly to eq.(\ref{eq:BC_MM_mono}) .

\begin{figure}[H]
\centering
\includegraphics[width=0.5\linewidth]{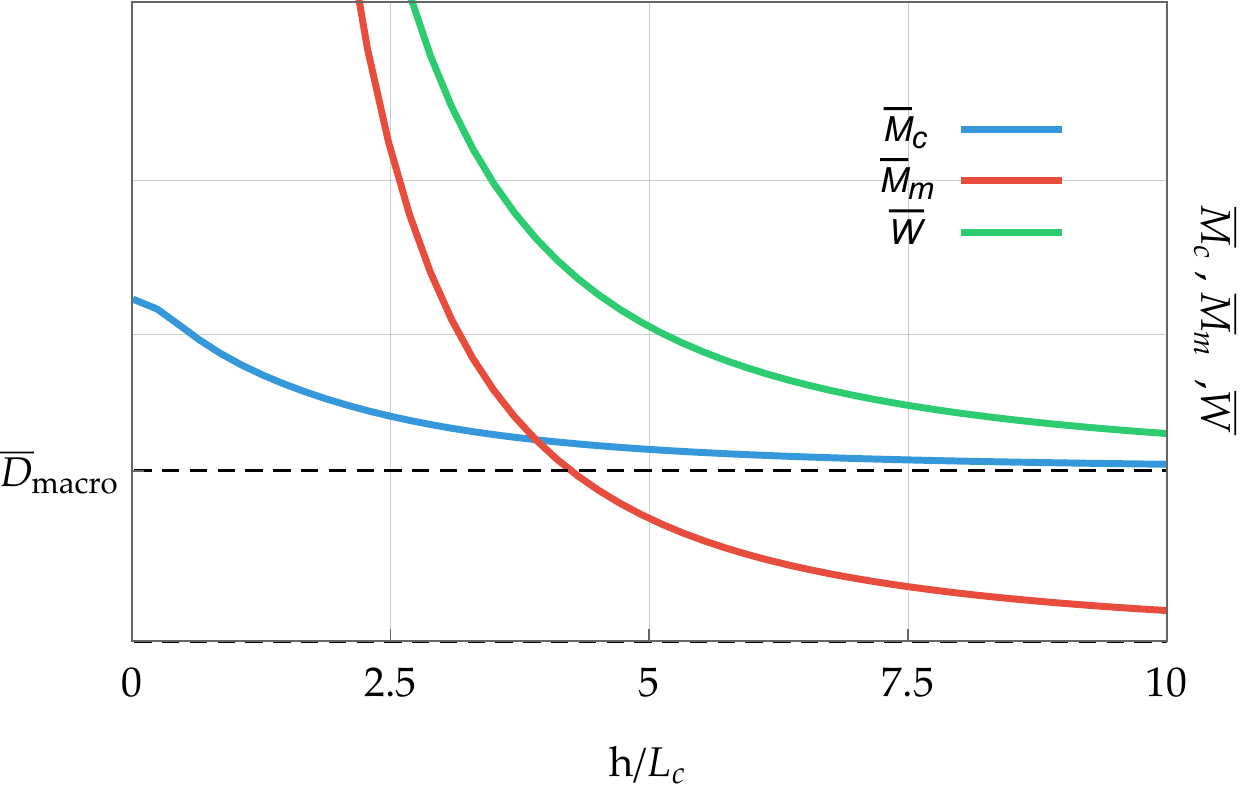}
\caption{(\textbf{Micromorphic model}, classical case) Bending moments and energy while varying $L_c$. Observe that the bending stiffness is unbounded as $L_c \to \infty$ ($h\to 0$). This is a major difference with respect the relaxed micromorphic model. The values of the parameters used are: $\mu _{e} = 1$, $\lambda _{e} = 1$, $\mu _{\tiny \mbox{micro}} = 1$, $\lambda _{\tiny \mbox{micro}} = 1$, $\mu = 1$, $a _1 = 2$, $a _2 = 1$, $a _3 = 1/2$.}
\label{fig:all_plot_MM_three}
\end{figure}

The symbolic expressions are too long to be reported here, but we provide a plot of the bending moments and the strain energy divided by $\frac{h^3}{12}\boldsymbol{\kappa}$ and $\frac{1}{2}\frac{h^3}{12}\boldsymbol{\kappa}^2$, respectively, in Fig.~\ref{fig:all_plot_MM_three} for selected parameter sets while changing $L_c$.
Still we have
$
\frac{\mbox{d}}{\mbox{d}\boldsymbol{\kappa}}W_{\mbox{tot}}(\boldsymbol{\kappa}) = M_{\mbox{c}} (\boldsymbol{\kappa}) + M_{\mbox{m}} (\boldsymbol{\kappa})
= M_{\mbox{c}} (\boldsymbol{\kappa}) \, .
$
\subsection{Penalized second gradient elasticity}
The classical Mindlin-Eringen micromorphic model can also be interpreted as a penalty formulation of second gradient elasticity. Indeed, letting $\mu_e \to \infty,\kappa_e \to \infty,\mu_c \to \infty$ imposes the constraint $\boldsymbol{P}= \boldsymbol{\mbox{D}u}$ and the remaining minimization problem is of the type
\begin{align}
	\min\limits_{\boldsymbol{u}}
	\left(
	\int_{\Omega}
	\mu_{\tiny \mbox{micro}} \left\lVert \mbox{sym}\,\boldsymbol{\mbox{D}u} \right\rVert^{2}
	+ \frac{\lambda_{\tiny \mbox{micro}}}{2} \mbox{tr}^2 \left(\boldsymbol{\mbox{D}u} \right)
	+ \frac{\mu \, L_c^2}{2}
	\lVert
	\boldsymbol{\mbox{D}^2 u}
	\rVert^2
	\right) \, .
\end{align}
We consider the stiffness generated for this limit in the simple case $\nu_{\tiny \mbox{micro}}=\nu_{e}=0$ and we investigate this specific limit case since it remains analytically treatable.
\begin{figure}[H]
\centering
\begin{subfigure}{.45\textwidth}
\centering
\includegraphics[width=0.95\linewidth]{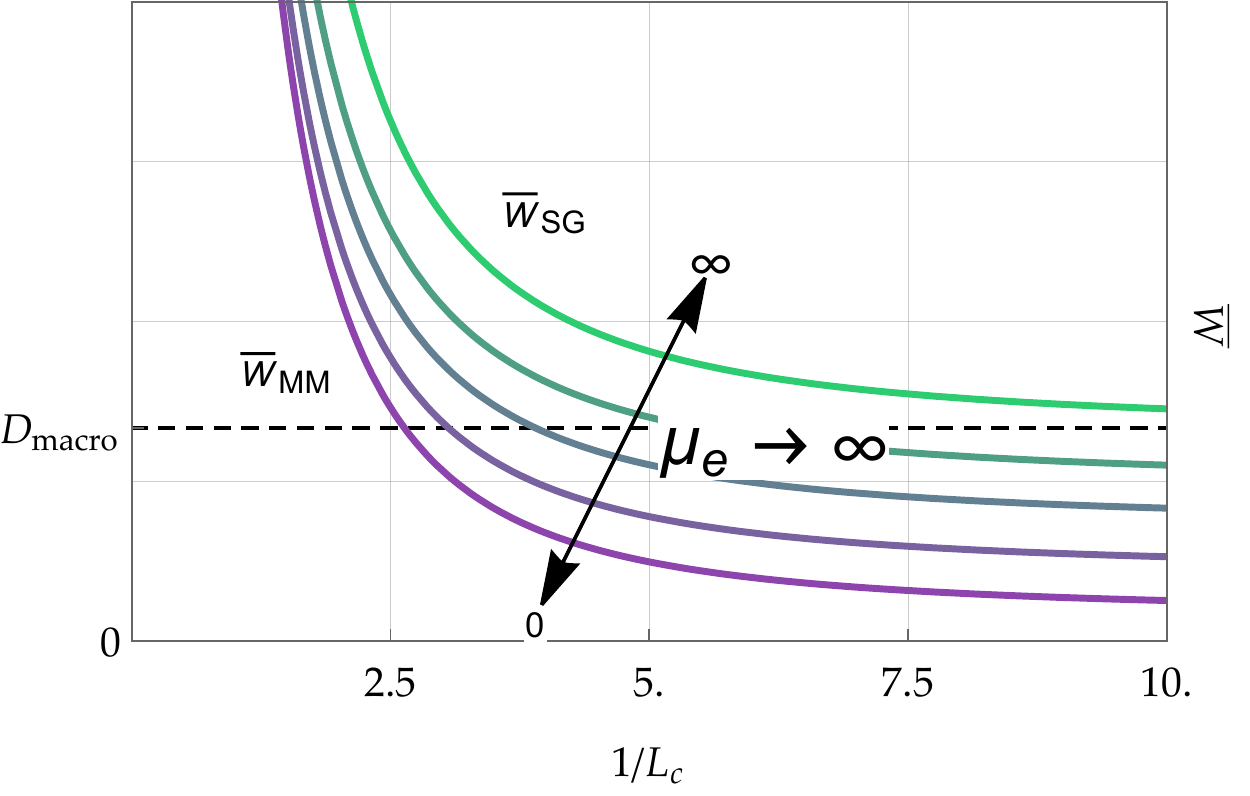}
\caption{}
\label{fig:limit_MM_SG_and_RM_1}
\end{subfigure}%
\begin{subfigure}{.45\textwidth}
\centering
\includegraphics[width=0.95\linewidth]{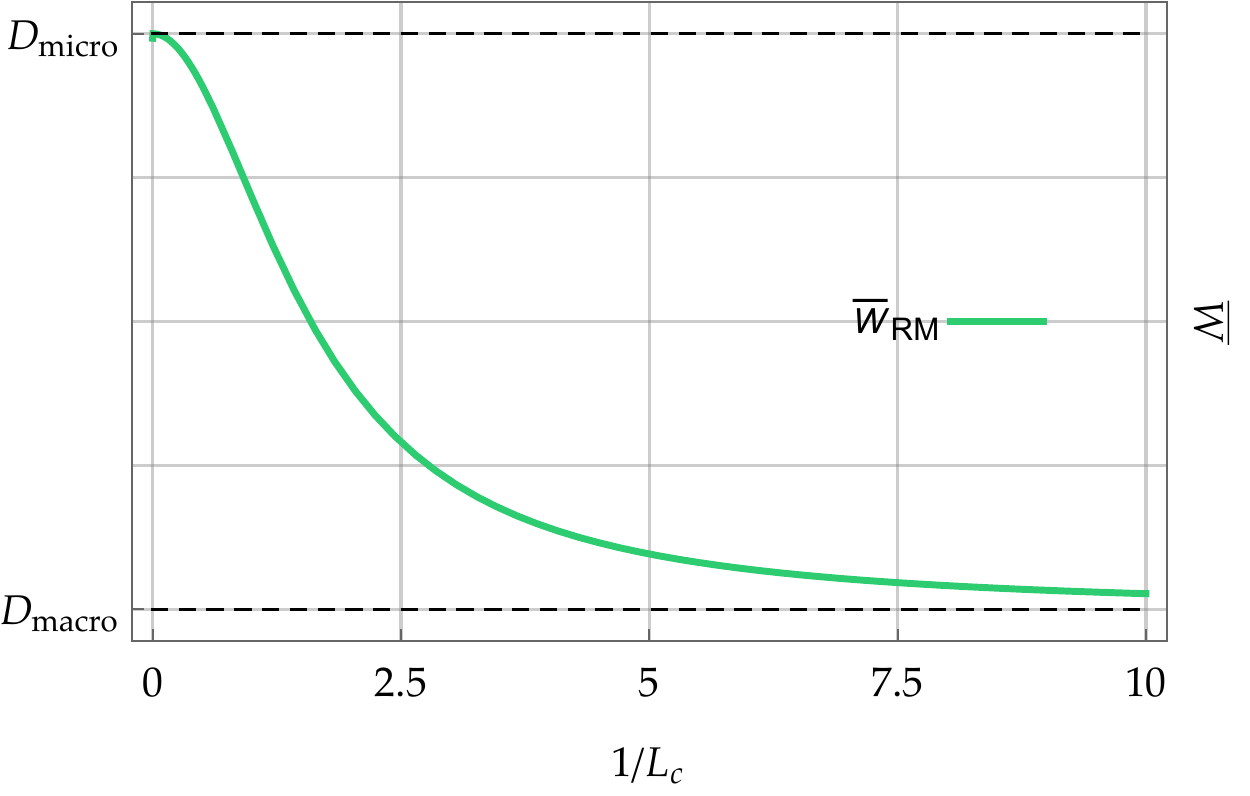}
\caption{}
\label{fig:limit_MM_SG_and_RM_2}
\end{subfigure}
\caption{
(a) (\textbf{Classical micromorphic model}) Limit for $\mu_{e} \to \infty$ of the unbounded stiffness of the Classical micromorphic model. The limit is given by the corresponding unbounded stiffness of the second gradient model $\overline{W}_{\mbox{SG}}$.
The purple line represent the unbounded stiffness for $\mu_{e} = 0$.
(b) (\textbf{Relaxed micromorphic model}) Limit for $\mu_{e} \to \infty$ of the bounded stiffness of the relaxed micromorphic model. The limit $\mu_{e} \to \infty$ (since $\nu_{e}=\lambda_{e}=0$) gives a linear elastic solution with stiffness $D_{\tiny \mbox{micro}}$ which remains bounded. Letting $\mu_{e} \to 0$ simply decouples the problem.
}
\end{figure}
\section{Cylindrical bending for the micro-strain model without mixed terms}

The micro-strain model \cite{forest2006nonlinear,hutter2015micromorphic} can be obtained from the classical Mindlin-Eringen model by assuming a priori that the microdistortion remains symmetric, $\boldsymbol{P}=\boldsymbol{S}\in \mbox{Sym}(3)$.

A bending solution for a particular model of this type has been derived in \cite{hutter2016application} disregarding of the lateral contraction.
This simplification shall be overcome here, whereby we employ a reduced isotropic curvature expression to make the calculations manageable.

Note that the micro-strain model cannot be obtained from the relaxed micromorphic model, although there are certain similarities. The free energy which we consider is given by
\begin{align}
W \left(\boldsymbol{\mbox{D}u}, \boldsymbol{S}, \boldsymbol{\mbox{D}S}\right) 
=& \, \mu_{e} \left\lVert \mbox{dev} \left(\mbox{sym} \, \boldsymbol{\mbox{D}u} - \boldsymbol{S}\right) \right\rVert^{2}
+ \frac{\kappa_{e}}{2} \mbox{tr}^2 \left(\boldsymbol{\mbox{D}u} - \boldsymbol{S} \right)
+ \mu_{\tiny \mbox{micro}} \left\lVert \mbox{dev} \, \boldsymbol{S} \right\rVert^{2} 
+ \frac{\kappa_{\tiny \mbox{micro}}}{2} \mbox{tr}^2 \left(\boldsymbol{S} \right)
\notag
\\*
&
+ \frac{\mu \, L_c^2}{2} \, 
\displaystyle\sum_{i=1}^{3}
\left(
a_1 \, \left\lVert \mbox{dev} \Big( \partial_{x_i} \boldsymbol{S} \Big) \right\rVert^2_{\mathbb{R}^{3\times 3}}
+ \frac{2}{9} \, a_3 \, \mbox{tr}^2 \Big( \partial_{x_i} \boldsymbol{S} \Big)
\right)
\notag
\\*
= &
\, \mu_{e} \left\lVert \mbox{dev} \left(\mbox{sym} \, \boldsymbol{\mbox{D}u} - \boldsymbol{S}\right) \right\rVert^{2}
+ \frac{\kappa_{e}}{2} \mbox{tr}^2 \left(\boldsymbol{\mbox{D}u} - \boldsymbol{S} \right)
+ \mu_{\tiny \mbox{micro}} \left\lVert \mbox{dev} \, \boldsymbol{S} \right\rVert^{2} 
+ \frac{\kappa_{\tiny \mbox{micro}}}{2} \mbox{tr}^2 \left(\boldsymbol{S} \right)
\label{eq:energy_FR}
\\*
&
+ \frac{\mu \, L_c^2}{2} \,
\left(
a_1 \, \left\lVert \mbox{dev} \, \boldsymbol{\mbox{D}S} \right\rVert^2
+ \frac{2}{9} \, a_3 \, \mbox{tr}^2 \Big( \boldsymbol{\mbox{D}S} \Big)
\right)
\, .
\notag
\end{align}
The meaning of $\mbox{dev} \, \mbox{sym} \, \boldsymbol{\mbox{D}S}$ and $\mbox{tr} \left(\boldsymbol{\mbox{D}S}\right)$ for the third order tensor $\boldsymbol{\mbox{D}S}$ can be inferred from eq.(\ref{eq:energy_FR}), i.e. we define
\begin{equation}
\left\lVert \mbox{dev} \, \mbox{sym} \, \boldsymbol{\mbox{D}S} \right\rVert^2 :=
\displaystyle\sum_{i=1}^{3}
\left\lVert \mbox{dev} \, \mbox{sym} \Big( \partial_{x_i} \boldsymbol{S} \Big) \right\rVert^2
\, , \qquad
\mbox{tr}^2 \Big( \boldsymbol{\mbox{D}S} \Big) :=
\displaystyle\sum_{i=1}^{3}
\mbox{tr}^2 \Big( \partial_{x_i} \boldsymbol{S} \Big)
\, .
\end{equation}
The chosen 2-parameter curvature expression represents a simplified isotropic curvature (the full isotropic curvature for the micro-strain model still counts 8 parameters \cite{barbagallo2016transparent}).
If $a_{1},a_{2}>0$ the chosen curvature energy provides a complete control of $\left\lVert \boldsymbol{\mbox{D}S} \right\rVert^2$.

If we assume that $\mu_{\mbox{\tiny micro}}\to\infty$, while $\kappa_{\mbox{\tiny micro}}<\infty$, then the model turns formally also into the micro-void model (see Section \ref{sec:Micro-void}) i.e. $\boldsymbol{S} = \omega \boldsymbol{\mathbbm{1}}$.
Then the curvature turns into
\footnote{
Note that $\left\lVert \mbox{Curl} \left(\omega \boldsymbol{\mathbbm{1}}\right) \right\rVert^2_{\mathbb{R}^{3\times 3}}
=
\left\lVert \mbox{Anti} \left( \boldsymbol{\mbox{D}} \omega \right) \right\rVert^2_{\mathbb{R}^{3\times 3}}
=
2 \left\lVert \mbox{axl} \, \left(\mbox{Anti} \left( \boldsymbol{\mbox{D}} \omega \right)\right) \right\rVert^2_{\mathbb{R}^3}
=
2 \left\lVert \boldsymbol{\mbox{D}} \omega \right\rVert^2_{\mathbb{R}^3}$
}
\begin{align}
\frac{\mu \,L_c^2}{2} \,
\displaystyle\sum_{i=1}^{3} \frac{2}{9} \, a_3 \, \mbox{tr}^2 \left( \partial_{x_i} \left( \omega \boldsymbol{\mathbbm{1}} \right) \right)
&=
\frac{\mu \,L_c^2}{2} \,
\frac{2}{9} \, a_3 \, \displaystyle\sum_{i=1}^{3} \mbox{tr}^2 \left( \partial_{x_i} \omega \cdot \boldsymbol{\mathbbm{1}}  \right)
=
\frac{\mu \,L_c^2}{2} \,
\frac{2}{9} \, a_3 \, \displaystyle\sum_{i=1}^{3} \left\lvert \partial_{x_i} \omega \right\rvert^2_{\mathbb{R}} \cdot 9
\label{eq:energy_FR_micro-void}
\\
&= \,
\frac{\mu \,L_c^2}{2} \,
2 \, a_3 \, \left\lVert \boldsymbol{\mbox{D}} \omega \right\rVert^2_{\mathbb{R}^3}
=
\frac{\mu \,L_c^2}{2} \,
a_3 \, \left\lVert \mbox{Curl} \left(\omega \boldsymbol{\mathbbm{1}} \right) \right\rVert^2_{\mathbb{R}^{3\times 3}}
\, .
\notag
\end{align}
The equilibrium equations without body forces are the following (see Appendix~\ref{app:full_micro_and_strain_eq_eqa})
\begin{align}
\mbox{Div}\overbrace{\left[
2\mu_{e} \, \mbox{dev} \left(\mbox{sym} \, \boldsymbol{\mbox{D}u} - \boldsymbol{S}\right)
+ \kappa_{e} \, \mbox{tr} \left(\boldsymbol{\mbox{D}u} - \boldsymbol{S} \right) \boldsymbol{\mathbbm{1}}
\right]}^{\mathlarger{\widetilde{\sigma}}:=}
= \boldsymbol{0},
\notag
\\*
2\mu_{e} \, \mbox{dev} \left(\mbox{sym} \, \boldsymbol{\mbox{D}u} - \boldsymbol{S}\right)
+ \kappa_{e} \, \mbox{tr} \left(\boldsymbol{\mbox{D}u} - \boldsymbol{S} \right) \boldsymbol{\mathbbm{1}}
- 2 \mu_{\tiny \mbox{micro}} \, \mbox{dev} \,\boldsymbol{S}
- \kappa_{\tiny \mbox{micro}} \, \mbox{tr} \left(\boldsymbol{S}\right) \boldsymbol{\mathbbm{1}} \, \, 
\hspace{2.5cm}
\label{eq:equiMic_FR_three}
\\*
+ \, \mu \, L_{c}^{2}\,
\mbox{sym} \, 
\left[
a_1 \, \mbox{dev} \, \boldsymbol{\Delta S}
+ \frac{2}{9} \, a_3 \, \mbox{tr} \left(\boldsymbol{\Delta S}\right) \boldsymbol{\mathbbm{1}}
\right]
= \boldsymbol{0} \, ,
\notag
\end{align}
where $\boldsymbol{\Delta S} \in \mathbb{R}^{3\times3}$ is taken component-wise.
The boundary conditions (see the Appendix~\ref{app:full_micro_and_strain_eq_eqa}) at the upper and lower surface (free surface) are 
\begin{align}
\boldsymbol{\widetilde{t}}(x_2 = \pm \, h/2) &= 
\pm \, \boldsymbol{\widetilde{\sigma}}(x_2) \cdot \boldsymbol{e}_2 = 
\boldsymbol{0} \, ,
\label{eq:BC_FR_three}
\\
\boldsymbol{\eta}(x_2 = \pm \, h/2) &= 
\pm \, \sum_{i=1}^{3} \mbox{sym} \left(\boldsymbol{m}_{i} (x_2) \, (\boldsymbol{e}_{2})_{i}\right) =
\boldsymbol{0} \, ,
\notag
\end{align}
where $\boldsymbol{m}_{i} = a_1 \, \mbox{dev} \left( \partial_{x_i} \boldsymbol{S} \right)$ $+$ $\frac{2}{9} \, a_3 \, \mbox{tr} \left( \partial_{x_i} \boldsymbol{S} \right) \boldsymbol{\mathbbm{1}}$, $i=1,2,3$ is a second order tensor and $(\boldsymbol{e}_{2})_{i}$ is the scalar \textit{i}th component of the unit vector $\boldsymbol{e}_{2}$ (see the Appendix~\ref{app:full_micro_and_strain_eq_eqa}).
According with the reference system shown in Fig.~\ref{fig:intro}, the ansatz for the displacement field and the microdistortion is
\begin{equation}
\boldsymbol{u}(x_1,x_2)=
\left(
\begin{array}{c}
-\kappa_1 \, x_1 x_2 \\
v(x_2)+\frac{\kappa_1  x_1^2}{2} \\
0 \\
\end{array}
\right) \, ,
\qquad
\boldsymbol{S}\left(x_1,x_2\right) =
\left(
\begin{array}{ccc}
S_{11}(x_2) & 0 & 0 \\
0 & S_{22}(x_2) & 0 \\
0 & 0 & S_{33}(x_2) \\
\end{array}
\right) \, .
\label{eq:ansatz_FR_three}
\end{equation}
Substituting the ansatz eq.(\ref{eq:ansatz_FR_three}) in eq.(\ref{eq:equiMic_FR_three}) the equilibrium equation results in
\begin{align}
2 \mu _e \left(v''(x_{2})-S_{22}'(x_{2})\right)-\lambda _e \left(\boldsymbol{\kappa} +S_{11}'(x_{2})+S_{22}'(x_{2})+S_{33}'(x_{2})-v''(x_{2})\right) &= 0 \, ,
\notag
\\*
\frac{1}{9} \mu \, L_c^2 \left(2 (3 a_{1}+a_{3}) S_{11}''(x_{2})-(3 a_{1}-2 a_{3}) \left(S_{22}''(x_{2})+S_{33}''(x_{2})\right)\right)&
\notag
\\*
-\lambda _e \left(S_{11}(x_{2})+S_{22}(x_{2})+S_{33}(x_{2})-v'(x_{2})+\boldsymbol{\kappa}  x_{2}\right)-2 \mu _e (S_{11}(x_{2})+\boldsymbol{\kappa}  x_{2})&
\notag
\\*
+\lambda _{\mbox{\tiny micro}} (-S_{11}(x_{2})-S_{22}(x_{2})-S_{33}(x_{2}))-2 \mu _{\mbox{\tiny micro}} S_{11}(x_{2}) &= 0
\, ,
\notag
\\*
-\frac{1}{9} \mu \, L_c^2 \left((3 a_{1}-2 a_{3}) S_{11}''(x_{2})-2 (3 a_{1}+a_{3}) S_{22}''(x_{2})+(3 a_{1}-2 a_{3}) S_{33}''(x_{2})\right)&
\label{eq:equi_equa_FR_three}
\\*
-\lambda _e \left(S_{11}(x_{2})+S_{22}(x_{2})+S_{33}(x_{2})-v'(x_{2})+\boldsymbol{\kappa}  x_{2}\right)+2 \mu _e \left(v'(x_{2})-S_{22}(x_{2})\right)&
\notag
\\*
+\lambda _{\mbox{\tiny micro}} (-S_{11}(x_{2})-S_{22}(x_{2})-S_{33}(x_{2}))-2 \mu _{\mbox{\tiny micro}} S_{22}(x_{2}) &= 0
\, ,
\notag
\\*
-\frac{1}{9} \mu \, L_c^2 \left((3 a_{1}-2 a_{3}) S_{11}''(x_{2})+(3 a_{1}-2 a_{3}) S_{22}''(x_{2})-2 (3 a_{1}+a_{3}) S_{33}''(x_{2})\right)&
\notag
\\*
-\lambda _e \left(S_{11}(x_{2})+S_{22}(x_{2})+S_{33}(x_{2})-v'(x_{2})+\boldsymbol{\kappa}  x_{2}\right)-2 \mu _e S_{33}(x_{2})&
\notag
\\*
+\lambda _{\mbox{\tiny micro}} (-S_{11}(x_{2})-S_{22}(x_{2})-S_{33}(x_{2}))-2 \mu _{\mbox{\tiny micro}} S_{33}(x_{2}) &= 0 \, .
\notag
\end{align}
It is possible to evaluate $v''(x_2)$ from eq.(\ref{eq:equi_equa_FR_three})$_1$ and then $S_{22}(x_2)$ form a linear combination of the remaining equation.
After substituting the expressions of $v''(x_2)$ and $S_{22}(x_2)$ in eq.(\ref{eq:equi_equa_FR_three}), we are left with two coupled ordinary differential equations of fourth order in $S_{11}(x_2)$ and $S_{33}(x_2)$.

\begin{figure}[H]
	\centering
	\includegraphics[width=0.5\linewidth]{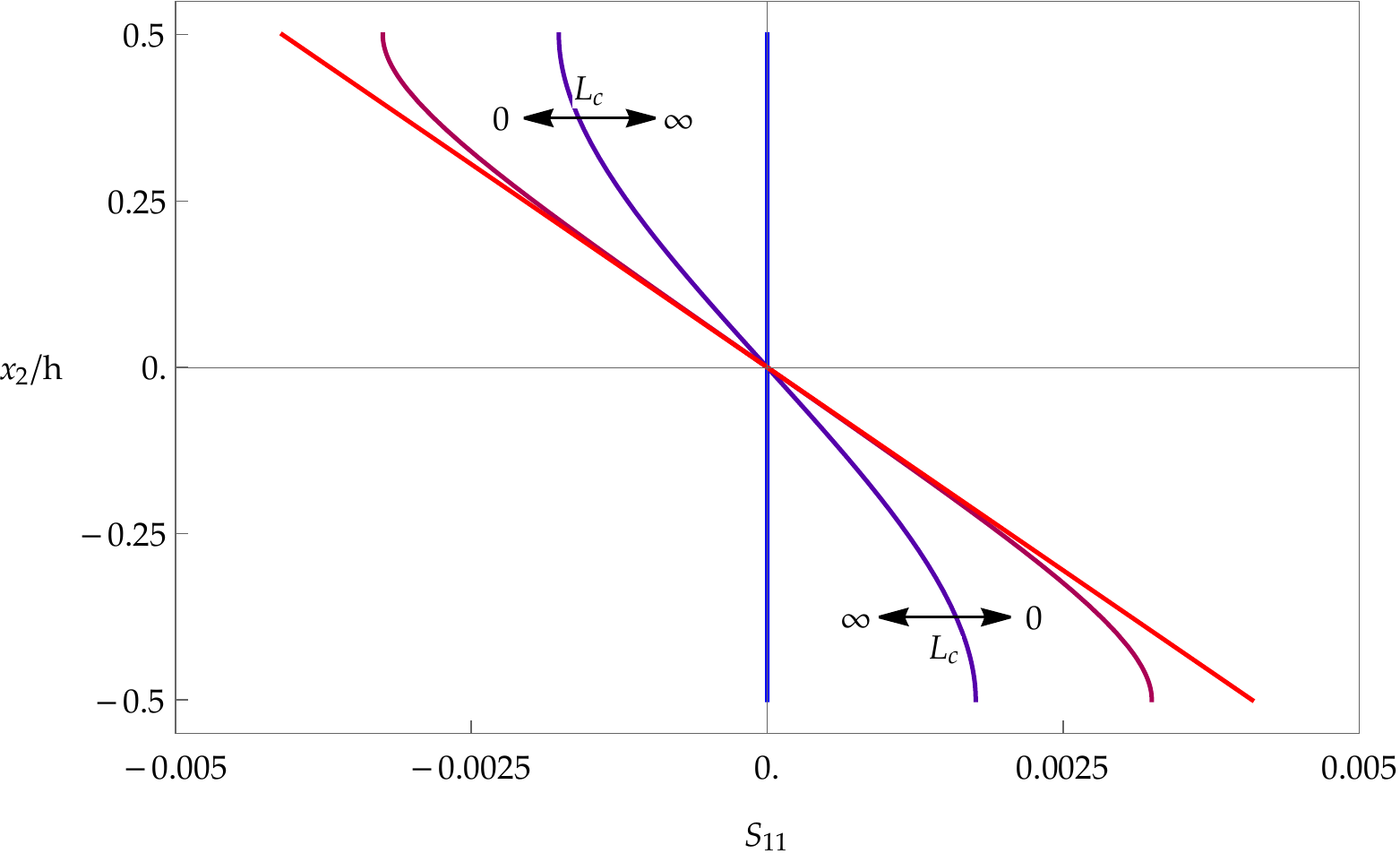}
	\caption{(\textbf{Micro-strain model}) Distribution across the thickness of $S_{11}$ while varying $L_c$.
		On the vertical axis we have the dimensionless thickness while on the horizontal we have the quantity $S_{11}$.
		Notice that the red curve correspond to the limit for $L_c \to 0$ while the blue curve correspond to the limit for $L_c \to \infty$.
		Observe the similar behaviour on $L_c \to \infty$ with the classical micromorphic model, as opposed to the relaxed micromorphic model. The values of the parameters used are: $\mu _{e} = 1/3$, $\lambda _{e} = 1/8$, $\mu _{\tiny \mbox{micro}} = 2$, $\lambda _{\tiny \mbox{micro}} = 1$, $\mu = 1$, $a _1 = 1/2$, $a _2 = 1/2$, $a _3 = 3/2$, $\boldsymbol{\kappa}=7/200$.
		}
	\label{fig:P11_plot_FR_three}
\end{figure}

In Fig.~\ref{fig:P11_plot_FR_three} we show the distribution across the thickness of $P_{11}$ while varying $L_c$.

The classical bending moment, the higher-order bending moment, and energy (per unit area d$x_1$d$x_3$) definitions are reported in the following eq.(\ref{eq:sigm_ene_dimensionless_FR_three})
\begin{align}
M_{\mbox{c}} (\boldsymbol{\kappa})
=
\displaystyle\int\limits_{-h/2}^{h/2}
\langle \boldsymbol{\widetilde{\sigma}} \boldsymbol{e}_1 , \boldsymbol{e}_1 \rangle  \, x_2 \, 
\mbox{d}x_{2}
\, ,
\qquad\qquad
M_{\mbox{m}}(\boldsymbol{\kappa})
=
\displaystyle\int\limits_{-h/2}^{h/2}
\langle \left(\sum_{i=1}^{3}\boldsymbol{m}_{i} (x_2) \cdot (\boldsymbol{e}_{1})_{i} \right) \boldsymbol{e}_2 , \boldsymbol{e}_1 \rangle \,
\mbox{d}x_{2}
\, ,
\label{eq:sigm_ene_dimensionless_FR_three}
\\
W_{\mbox{tot}} (\boldsymbol{\kappa})
=
\displaystyle\int\limits_{-h/2}^{+h/2} W \left( \boldsymbol{\mbox{D}u},\boldsymbol{S},\boldsymbol{\mbox{D}S} \right) \, \mbox{d}x_{2}
\, ,
\hspace{4cm}
\notag
\end{align}
where $(\boldsymbol{e}_{1})_{i}$ is the scalar \textit{i}th component of the unit vector $\boldsymbol{e}_{1}$.

\begin{figure}[H]
\centering
\includegraphics[width=0.5\linewidth]{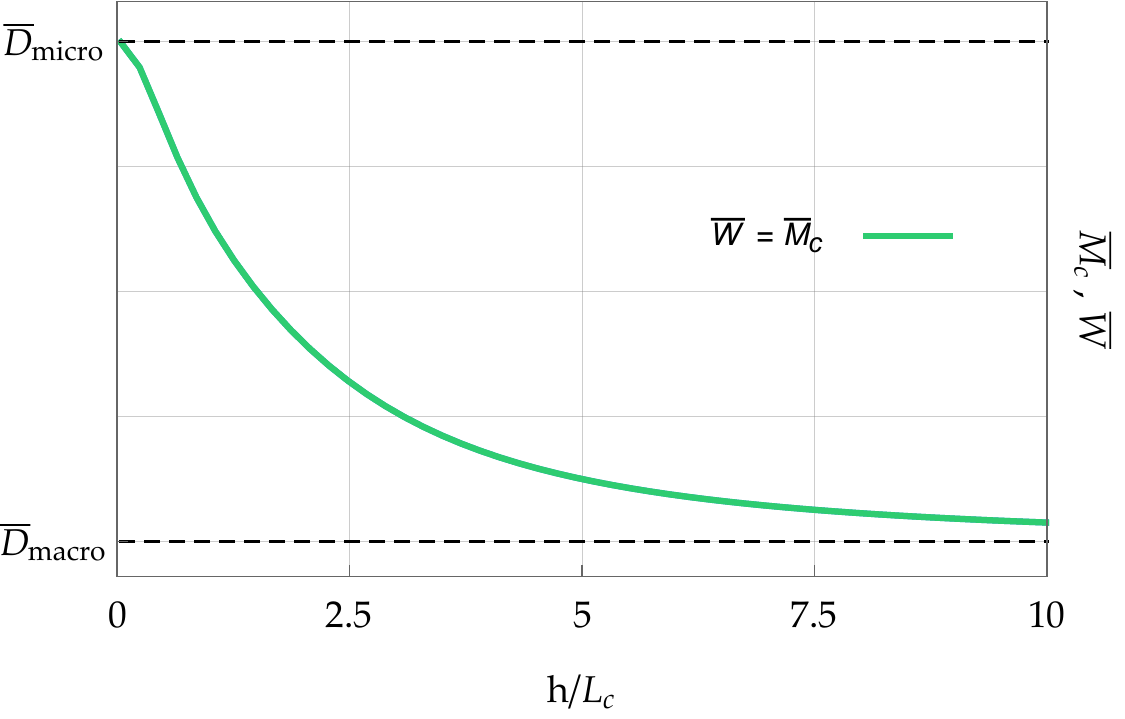}
\caption{(\textbf{Micro-strain model}) Energy expression while varying $L_c$. Observe that the bending stiffness remains bounded as $L_c \to \infty$ ($h \to 0$). Here this is due to the fact that the skew-symmetric part of the micro-distorsion is absent in the curvature energy, similarly to the micro-void curvature energy. The values of the parameters used are: $\mu _{e} = 1$, $\lambda _{e} = 1$, $\mu _{\tiny \mbox{micro}} = 1$, $\lambda _{\tiny \mbox{micro}} = 1$, $\mu = 1$, $a _1 = 2$, $a _2 = 1$, $a _3 = 1/2$.}
\label{fig:all_plot_FR_three}
\end{figure}

Again, the symbolic expressions are lengthy and are thus not reported here in detail. Fig.\ref{fig:all_plot_FR_three} provides a graphical representation of the final result.
Since the higher-order bending moment is zero and the following relation holds
\begin{equation}
\frac{\mbox{d}}{\mbox{d}\boldsymbol{\kappa}}W_{\mbox{tot}}(\boldsymbol{\kappa}) = M_{\mbox{c}} (\boldsymbol{\kappa}) + M_{\mbox{m}} (\boldsymbol{\kappa})
= M_{\mbox{c}} (\boldsymbol{\kappa}) \, ,
\end{equation}
only the plot of energy \cite{shaat2018reduced} (per unit area d$x_1$d$x_3$) while changing $L_c$ is shown in Fig.~\ref{fig:all_plot_FR_three}
The energy of the model remains bounded, as for the micro-void model, since for both models the higher-order bending moment are zero, and this does not create a conflict with the boundary condition as $L_c \to \infty$.
\section{Cylindrical bending for the second gradient continuum}
The expression of the most general isotropic strain energy for the second gradient continuum is \cite{mindlin1964micro}
\begin{align}
W \left(\boldsymbol{\mbox{D}u}, \boldsymbol{\mbox{D}^2 u}\right)
= &
\, \mu_{\tiny \mbox{macro}} \left\lVert \mbox{sym} \, \boldsymbol{\mbox{D}u} \right\lVert^2
+ \frac{\lambda_{\tiny \mbox{macro}}}{2} \, \mbox{tr}^2 \left(\boldsymbol{\mbox{D}u}\right)
\label{eq:energy_SG_gen}
\\
&
+ \widehat{a}_1 \, \chi_{iik} \, \chi_{kjj}
+ \widehat{a}_2 \, \chi_{ijj} \, \chi_{ikk}
+ \widehat{a}_3 \, \chi_{iik} \, \chi_{jjk}
+ \widehat{a}_4 \, \chi_{ijk} \, \chi_{ijk}
+ \widehat{a}_5 \, \chi_{ijk} \, \chi_{kji}
\, ,
\notag
\end{align}
where $\boldsymbol{\chi} = \boldsymbol{\mbox{D}}^2 \boldsymbol{u}$ ($\chi_{ijk} = \frac{\partial^2 u_k}{\partial x_i \, \partial x_j}$).
The expression we are going to use in the following is a simplified isotropic strain energy with three curvature parameters
\begin{align}
W \left(\boldsymbol{\mbox{D}u}, \boldsymbol{\mbox{D}^2 u}\right)
= &
\, \mu_{\tiny \mbox{macro}} \left\lVert \mbox{sym}\,\boldsymbol{\mbox{D}u} \right\rVert^{2}
+ \frac{\lambda_{\tiny \mbox{macro}}}{2} \mbox{tr}^2 \left(\boldsymbol{\mbox{D}u} \right)
\notag
\\
&
+ \frac{\mu \, L_c^2}{2}
\displaystyle\sum_{i=1}^{3}
\left(
a_1 \, \left\lVert \mbox{dev} \, \mbox{sym} \Big( \partial_{x_i} \boldsymbol{\mbox{D}u} \Big) \right\rVert^2_{\mathbb{R}^{3\times 3}}
+ a_2 \, \left\lVert \mbox{skew} \Big( \partial_{x_i} \boldsymbol{\mbox{D}u} \Big) \right\rVert^2_{\mathbb{R}^{3\times 3}}
+ \frac{2}{9} \, a_3 \, \mbox{tr}^2 \Big( \partial_{x_i} \boldsymbol{\mbox{D}u} \Big)
\right)
\notag
\\
= &
\, \mu_{\tiny \mbox{macro}} \left\lVert \mbox{sym}\,\boldsymbol{\mbox{D}u} \right\rVert^{2}
+ \frac{\lambda_{\tiny \mbox{macro}}}{2} \mbox{tr}^2 \left(\boldsymbol{\mbox{D}u} \right)
\label{eq:energy_SG}
\\
&
+ \frac{\mu \, L_c^2}{2}
\left(
a_1 \, \left\lVert \mbox{dev} \, \mbox{sym} \, \boldsymbol{\mbox{D}^2 u} \right\rVert^2
+ a_2 \, \left\lVert \mbox{skew} \, \boldsymbol{\mbox{D}^2 u} \right\rVert^2
+  \frac{2}{9} \, a_3 \, \mbox{tr}^2 \Big( \boldsymbol{\mbox{D}^2 u} \Big)
\right)
\, .
\notag
\end{align}
The meaning of $\mbox{dev} \, \mbox{sym} \, \boldsymbol{\mbox{D}^2 u}$, $\mbox{skew} \, \boldsymbol{\mbox{D}^2 u}$, and $\mbox{tr} \left(\boldsymbol{\mbox{D}^2 u}\right)$ for the third order tensor $\boldsymbol{\mbox{D}^2 u}$ can be inferred from eq.(\ref{eq:energy_SG}), i.e. we define
\begin{align}
\left\lVert \mbox{dev} \, \mbox{sym} \, \boldsymbol{\mbox{D}^2 u} \right\rVert^2 :=
\displaystyle\sum_{i=1}^{3}
\left\lVert \mbox{dev} \, \mbox{sym} \Big( \partial_{x_i} \boldsymbol{\mbox{D}u} \Big) \right\rVert^2
\, , \qquad
\left\lVert \mbox{skew} \, \boldsymbol{\mbox{D}^2 u} \right\rVert^2 :=
\displaystyle\sum_{i=1}^{3}
\left\lVert \mbox{skew} \Big( \partial_{x_i} \boldsymbol{\mbox{D}u} \Big) \right\rVert^2
\, ,
\\
\mbox{tr}^2 \Big( \boldsymbol{\mbox{D}^2 u} \Big) :=
\displaystyle\sum_{i=1}^{3}
\mbox{tr}^2 \Big( \partial_{x_i} \boldsymbol{\mbox{D}u} \Big)
\, .
\notag
\hspace{4cm}
\end{align}
The equilibrium equations (see the Appendix~\ref{app:second_gradient_eq_eqa}) without body forces are
\begin{align}
\mbox{Div}\bigg[
2 \mu_{\mbox{\tiny macro}} \,\mbox{sym}\,\boldsymbol{\mbox{D}u}
+ \lambda_{\tiny \mbox{macro}} \mbox{tr} \left(\boldsymbol{\mbox{D}u}\right) \boldsymbol{\mathbbm{1}}
\hspace{8cm}
\label{eq:equiMic_SG_three}
\\*
- \mu L_{c}^{2} \,
\left(
a_1 \, \mbox{dev} \, \mbox{sym} \, \boldsymbol{\Delta} \left(\boldsymbol{\mbox{D}u}\right)
+ a_2 \, \mbox{skew} \, \boldsymbol{\Delta} \left(\boldsymbol{\mbox{D}u}\right)
+ \frac{2}{9} \, a_3 \, \mbox{tr} \left(\boldsymbol{\Delta} \left(\boldsymbol{\mbox{D}u}\right)\right)\boldsymbol{\mathbbm{1}}
\right)
\bigg]
= \boldsymbol{0} \, ,
\notag
\end{align}
where $\boldsymbol{\Delta} \left(\boldsymbol{\mbox{D}u}\right) \in \mathbb{R}^{3\times3}$, is taken component-wise.
According with the reference system shown in Fig.~\ref{fig:intro}, the ansatz for the displacement field and the microdistortion is
\begin{equation}
\boldsymbol{u}(x_1,x_2)=
\left(
\begin{array}{c}
-\boldsymbol{\kappa} \, x_1 x_2 \\
v(x_2)+\frac{\boldsymbol{\kappa} \, x_1^2}{2} \\
0 \\
\end{array}
\right) \, .
\label{eq:ansatz_SG_three}
\end{equation}
\subsection{One curvature parameter and zero Poisson's ratio $\nu_{\tiny \mbox{macro}}=0$}
\begin{figure}[H]
	\centering
	\includegraphics[width=0.5\linewidth]{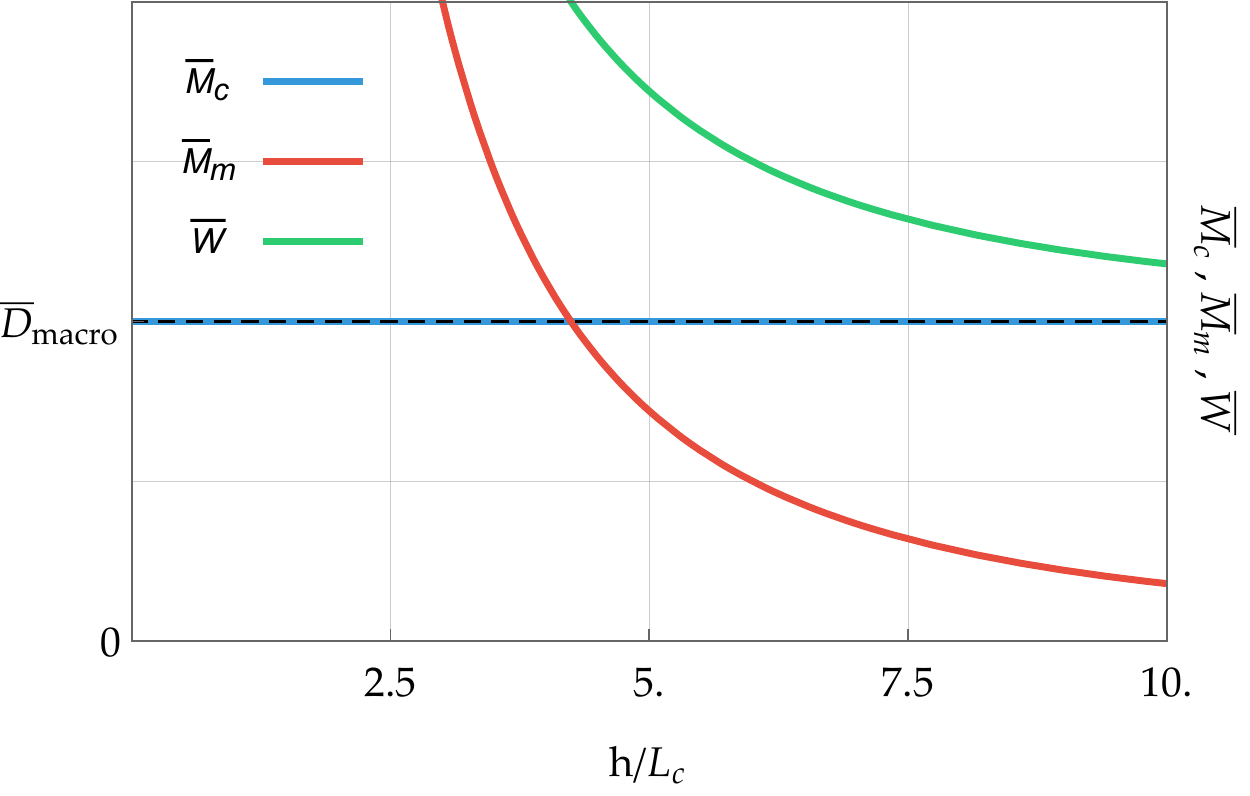}
	\caption{(\textbf{Second gradient model}, one curvature parameter, zero Poisson's ratio $\nu_{\mbox{\tiny macro}}$=0) Bending moments and energy while varying $L_c$. Observe that the bending stiffness is unbounded as $L_c \to \infty$ ($h\to 0$). This is a major difference with respect the relaxed micromorphic model. The values of the parameters used are: $\mu _{\tiny \mbox{macro}} = 1$, $\mu = 1$.}
	\label{fig:all_plot_SG_three_L0}
\end{figure}

Substituting the ansatz eq.(\ref{eq:ansatz_SG_three}) in eq.(\ref{eq:equiMic_SG_three}) while choosing $a_1=a_2=1, a_3=\frac{3}{2}$ and the Poisson's ratio $\nu_{\tiny \mbox{macro}}=0$, the equilibrium equation are
\begin{equation}
	- \mu \, L_c^2 \, v^{(4)}(x_{2}) + 2 \mu _{\tiny \mbox{macro}} \, v''(x_{2}) = 0 \, .
	\label{eq:equi_equa_SG_three_L0}
\end{equation}
It is possible to see that, for the second gradient model there is just one cumulative higher order parameter.
The boundary conditions (are reported here just the non zero terms), in the classical Mindlin formulation \cite{mindlin1964micro}, at the upper and lower surface (free surface) are
\begin{align}
	\widetilde{t}_{k} (x_2 = \pm \, h/2) &= 
	\pm \, \left(\widetilde{\sigma}_{jk} \, n_{j} - n_{i} \, n_{j} \, D\left(\mathfrak{m}_{ijk}\right)-2n_{j} \, D_{i}\left(\mathfrak{m}_{ijk}\right)\right) = 
	0 \, ,
	\label{eq:BC_SG_mono_L0}
	\\
	\eta_{k}(x_2 = \pm \, h/2) &= 
	\pm \, \left(n_{i} \, n_{j} \, \mathfrak{m}_{ijk}\right) =
	0 \, ,
	\notag
\end{align}
where $\widetilde{\boldsymbol{\sigma}} = 2 \mu_{\mbox{\tiny macro}} \, \mbox{sym}\,\boldsymbol{\mbox{D}u}$, the third-order moment stress tensor $\boldsymbol{\mathfrak{m}}=\mu \, L_c \, \boldsymbol{\mbox{D}^2 u}$, $\boldsymbol{n} = \boldsymbol{e}_2$ is the normal to the upper or lower surface and
\begin{equation}
	D_{j}\left(\cdot \right) =\left(\delta_{jl}-n_{j}n_{\ell}\right)\left(\cdot \right)_{,\ell}
	\, ,
	\qquad
	D \left(\cdot \right)=n_{\ell} \left(\cdot \right)_{,\ell} \, .
	\label{eq:12_L0}
\end{equation}
The solution of eq.(\ref{eq:equi_equa_SG_three_L0}) is
\begin{equation}
	\begin{split}
		v(x_{2}) = &
		\frac{\mu \, L_c^2 }{\mu_{\mbox{\tiny macro}}} \left(c_1 \, e^{-\frac{x_{2} \sqrt{\frac{2 \mu_{\mbox{\tiny macro}}}{\mu }}}{L_c}} + c_2 \, e^{\frac{x_{2} \sqrt{\frac{2 \mu_{\mbox{\tiny macro}}}{\mu }}}{L_c}}\right)
		+ c_4 \, x_{2}
		+ c_3 \, .
	\end{split}
	\label{eq:sol_v_SG_three_L0}
\end{equation}
After applying the boundary conditions eq.(\ref{eq:BC_SG_mono_L0}) it is possible to evaluate the displacement field solution of the cylindrical bending problem.
The classical bending moment, the higher-order bending moment, and energy (per unit area d$x_1$d$x_3$) definitions are reported in the following eq.(\ref{eq:sigm_ene_dimensionless_SG_three_L0}) (where $\boldsymbol{n} = \boldsymbol{e}_{1}$ in this case):

\begin{align}
	M_{\mbox{c}} (\boldsymbol{\kappa})
	=&
	\displaystyle\int\limits_{-h/2}^{h/2}
	\Big[\left(\widetilde{\sigma}_{jk} \, n_{j} - n_{i} \, n_{j} \, D\left(\mathfrak{m}_{ijk}\right)-2n_{j} \, D_{i}\left(\mathfrak{m}_{ijk}\right)\right) \,   n_{k}\Big] \, x_2 \, 
	\mbox{d}x_{2}
	=
	\frac{h^3}{12}
	\,
	2 \mu_{\mbox{\tiny macro}}
	\,
	\boldsymbol{\kappa}
	\, ,
	\label{eq:sigm_ene_dimensionless_SG_three_L0}
	\\*
	M_{\mbox{m}}(\boldsymbol{\kappa})
	=&
	\displaystyle\int\limits_{-h/2}^{h/2}
	\langle \Big(\underbrace{\boldsymbol{\mathfrak{m}} \cdot \boldsymbol{e}_1}_{\mathbb{R}^{3\times3}}\Big) \boldsymbol{e}_1 , \boldsymbol{e}_2 \rangle \,
	- \langle \Big(\underbrace{\boldsymbol{\mathfrak{m}} \cdot \boldsymbol{e}_2}_{\mathbb{R}^{3\times3}}\Big) \boldsymbol{e}_1 , \boldsymbol{e}_1 \rangle \,
	- \langle \Big(\underbrace{\boldsymbol{\mathfrak{m}} \cdot \boldsymbol{e}_1}_{\mathbb{R}^{3\times3}}\Big) \boldsymbol{e}_2 , \boldsymbol{e}_1 \rangle \,
	\mbox{d}x_{2}
	=
	\, \frac{h^3}{12}
	\left[
	36 \, \mu \, \left(\frac{L_c}{h}\right)^2
	\right]
	\boldsymbol{\kappa}
	\, ,
	\notag
	\\*
	W_{\mbox{tot}} (\boldsymbol{\kappa})
	=&
	\displaystyle\int\limits_{-h/2}^{+h/2} W \left( \boldsymbol{\mbox{D}u},\boldsymbol{\mbox{D}^2 u} \right) \, \mbox{d}x_{2}
	=
	\frac{1}{2}
	\frac{h^3}{12}
	\left[
	2 \mu_{\mbox{\tiny macro}}
	+ 36 \, \mu \, \left(\frac{L_c}{h}\right)^2
	\right]
	\boldsymbol{\kappa}^2
	\, .
	\notag
\end{align}
The plot of the bending moments and the strain energy divided by $\frac{h^3}{12}\boldsymbol{\kappa}$ and $\frac{1}{2}\frac{h^3}{12}\boldsymbol{\kappa}^2$, respectively, while changing $L_c$ is shown in Fig.~\ref{fig:all_plot_SG_three_L0}.
\subsection{One curvature parameter and arbitrary Poisson's ratio}
Substituting the ansatz eq.(\ref{eq:ansatz_SG_three}) in eq.(\ref{eq:equiMic_SG_three}) while choosing $a_1=a_2=1, a_3=\frac{3}{2}$, the equilibrium equation results to be:
\begin{equation}
- \mu \, L_c^2 \, v^{(4)}(x_{2}) + \lambda _{\tiny \mbox{macro}} \, \left(v''(x_{2}) - \boldsymbol{\kappa} \right) + 2 \mu _{\tiny \mbox{macro}} \, v''(x_{2}) = 0 \, .
\label{eq:equi_equa_SG_mono}
\end{equation}

\begin{figure}[H]
\centering
\includegraphics[width=0.5\linewidth]{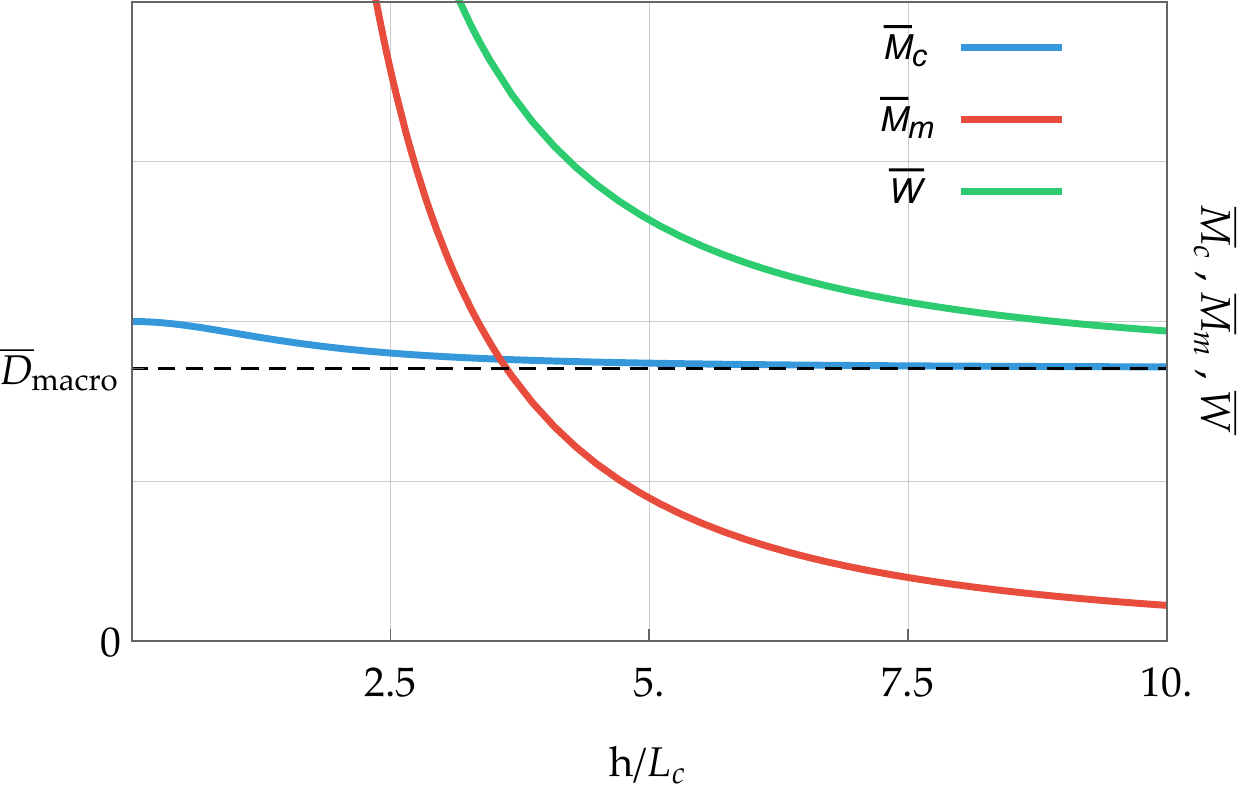}
\caption{(\textbf{Second gradient model}, one curvature parameter) Bending moments and energy while varying $L_c$. Observe that the bending stiffness is unbounded as $L_c \to \infty$ ($h\to 0$). This is a major difference with respect to the relaxed micromorphic model. The values of the parameters used are: $\mu _{\tiny \mbox{macro}} = 1$, $\lambda _{\tiny \mbox{micro}} = 1$, $\mu = 1$.}
\label{fig:all_plot_SG_mono}
\end{figure}

The solution of eq.(\ref{eq:equi_equa_SG_mono}) is
\begin{align}
v(x_{2}) = & \, \,
\frac{2 \, \mu \, L_c^2 }{\lambda_{\mbox{\tiny macro}} + 2 \mu_{\mbox{\tiny macro}}} \left(c_1 \, e^{-\frac{x_{2} \sqrt{\frac{\lambda_{\mbox{\tiny macro}}+2 \mu_{\mbox{\tiny macro}}}{\mu }}}{L_c}} + c_2 \, e^{\frac{x_{2} \sqrt{\frac{\lambda_{\mbox{\tiny macro}}+2 \mu_{\mbox{\tiny macro}}}{\mu }}}{L_c}}\right)
\label{eq:sol_v_SG_mono}
\\
&+ \frac{\lambda_{\mbox{\tiny macro}}}{\lambda_{\mbox{\tiny macro}}+2 \mu_{\mbox{\tiny macro}}} \, \frac{x_{2}^2}{2} \, \boldsymbol{\kappa}
+ c_4 \, x_{2}
+ c_3 \, .
\notag
\end{align}
After applying the boundary conditions (are reported here just the non zero terms), in the classical Mindlin formulation \cite{mindlin1964micro}, at the upper and lower surface (free surface)
\begin{align}
\widetilde{t}_{k} (x_2 = \pm \, h/2) &= 
\pm \, \left(\widetilde{\sigma}_{jk} \, n_{j} - n_{i} \, n_{j} \, D\left(\mathfrak{m}_{ijk}\right)-2n_{j} \, D_{i}\left(\mathfrak{m}_{ijk}\right)\right) = 
0 \, ,
\label{eq:BC_SG_mono}
\\
\eta_{k}(x_2 = \pm \, h/2) &= 
\pm \, \left(n_{i} \, n_{j} \, \mathfrak{m}_{ijk}\right) =
0 \, ,
\notag
\end{align}
where $\widetilde{\boldsymbol{\sigma}} = 2 \mu_{\mbox{\tiny macro}} \, \mbox{sym}\,\boldsymbol{\mbox{D}u}$, the third-order moment stress tensor $\boldsymbol{\mathfrak{m}}=\mu \, L_c \, \boldsymbol{\mbox{D}^2 u}$, $\boldsymbol{n} = \boldsymbol{e}_2$ is the normal to the upper or lower surface and
\begin{equation}
D_{j}\left(\cdot \right) =\left(\delta_{jl}-n_{j}n_{\ell}\right)\left(\cdot \right)_{,\ell}
\, ,
\qquad
D \left(\cdot \right)=n_{\ell} \left(\cdot \right)_{,\ell} \, ,
\label{eq:12_mono}
\end{equation}
it is possible to evaluate the displacement field solution of the cylindrical bending problem.
The classical bending moment, the higher-order bending moment, and energy (per unit area d$x_1$d$x_3$) definitions are reported in the following eq.(\ref{eq:sigm_ene_dimensionless_SG_three}) (where $\boldsymbol{n} = \boldsymbol{e}_{1}$ in this case)

\begin{align}
M_{\mbox{c}} (\boldsymbol{\kappa})
=&
\displaystyle\int\limits_{-h/2}^{h/2}
\Big[\left(\widetilde{\sigma}_{jk} \, n_{j} - n_{i} \, n_{j} \, D\left(\mathfrak{m}_{ijk}\right)-2n_{j} \, D_{i}\left(\mathfrak{m}_{ijk}\right)\right) \,   n_{k}\Big] \, x_2 \, 
\mbox{d}x_{2}
\notag
\\*
=& \, 
\frac{h^3}{12}
\left[
\frac{4 \mu_{\mbox{\tiny macro}} \left(\lambda_{\mbox{\tiny macro}}+\mu_{\mbox{\tiny macro}}\right)}{\lambda_{\mbox{\tiny macro}}+2 \mu_{\mbox{\tiny macro}}}
+ \frac{12 \, \mu \, \lambda_{\mbox{\tiny macro}}^2}{\left(\lambda_{\mbox{\tiny macro}}+2 \mu _{\mbox{\tiny macro}}\right){}^2} \left(\frac{L_c}{h}\right)^2
\right.
\notag
\\*
&
\hspace{23pt}
\left.
- \frac{24 \, \mu \, \lambda_{\mbox{\tiny macro}}^2}{\left(\lambda_{\mbox{\tiny macro}}+2 \mu _{\mbox{\tiny macro}}\right){}^2}
\left(\frac{L_c}{h}\right)^3
\sqrt{\frac{\mu }{\lambda_{\mbox{\tiny macro}}+2 \mu_{\mbox{\tiny macro}}}}
\tanh \left(\frac{h}{2 L_c} \sqrt{\frac{\lambda_{\mbox{\tiny macro}}+2 \mu_{\mbox{\tiny macro}}}{\mu }}\right)
\right]
\boldsymbol{\kappa}
\, ,
\notag
\\*
M_{\mbox{m}}(\boldsymbol{\kappa})
=&
\displaystyle\int\limits_{-h/2}^{h/2}
\langle \Big(\underbrace{\boldsymbol{\mathfrak{m}} \cdot \boldsymbol{e}_1}_{\mathbb{R}^{3\times3}}\Big) \boldsymbol{e}_1 , \boldsymbol{e}_2 \rangle \,
- \langle \Big(\underbrace{\boldsymbol{\mathfrak{m}} \cdot \boldsymbol{e}_2}_{\mathbb{R}^{3\times3}}\Big) \boldsymbol{e}_1 , \boldsymbol{e}_1 \rangle \,
- \langle \Big(\underbrace{\boldsymbol{\mathfrak{m}} \cdot \boldsymbol{e}_1}_{\mathbb{R}^{3\times3}}\Big) \boldsymbol{e}_2 , \boldsymbol{e}_1 \rangle \,
\mbox{d}x_{2}
=
\, \frac{h^3}{12}
\left[
36 \, \mu \, \left(\frac{L_c}{h}\right)^2
\right]
\boldsymbol{\kappa}
\, ,
\notag
\\*
W_{\mbox{tot}} (\boldsymbol{\kappa})
=&
\displaystyle\int\limits_{-h/2}^{+h/2} W \left( \boldsymbol{\mbox{D}u},\boldsymbol{\mbox{D}^2 u} \right) \, \mbox{d}x_{2}
\label{eq:sigm_ene_dimensionless_SG_mono}
\\*
=&
\frac{1}{2}
\frac{h^3}{12}
\left[
\frac{4 \mu_{\mbox{\tiny macro}} \left(\lambda_{\mbox{\tiny macro}}+\mu_{\mbox{\tiny macro}}\right)}{\lambda_{\mbox{\tiny macro}}+2 \mu_{\mbox{\tiny macro}}}
+ 12 \mu \, \frac{3 \left(\lambda_{\mbox{\tiny macro}}+2 \mu_{\mbox{\tiny macro}}\right){}^2+\lambda_{\mbox{\tiny macro}}^2}{\left(\lambda_{\mbox{\tiny macro}}+2 \mu_{\mbox{\tiny macro}}\right){}^2} \left(\frac{L_c}{h}\right)^2
\right.
\notag
\\*
&
\hspace{28pt}
\left.
- \frac{24 \, \mu \, \lambda_{\mbox{\tiny macro}}^2}{\left(\lambda_{\mbox{\tiny macro}}+2 \mu _{\mbox{\tiny macro}}\right){}^2}
\left(\frac{L_c}{h}\right)^3
\sqrt{\frac{\mu }{\lambda_{\mbox{\tiny macro}}+2 \mu_{\mbox{\tiny macro}}}}
\tanh \left(\frac{h}{2 L_c} \sqrt{\frac{\lambda_{\mbox{\tiny macro}}+2 \mu_{\mbox{\tiny macro}}}{\mu }}\right)
\right]
\boldsymbol{\kappa}^2
\, .
\notag
\end{align}
As always,
$
\frac{\mbox{d}}{\mbox{d}\boldsymbol{\kappa}}W_{\mbox{tot}}(\boldsymbol{\kappa}) = M_{\mbox{c}} (\boldsymbol{\kappa}) + M_{\mbox{m}} (\boldsymbol{\kappa})
= M_{\mbox{c}} (\boldsymbol{\kappa}) \, .
$
The plot of the bending moments and the strain energy divided by $\frac{h^3}{12}\boldsymbol{\kappa}$ and $\frac{1}{2}\frac{h^3}{12}\boldsymbol{\kappa}^2$, respectively, while changing $L_c$ is shown in Fig.~\ref{fig:all_plot_SG_mono}.
In Fig.~\ref{fig:P11_SG} we show the plot of $\left(\boldsymbol{\mbox{D}u}\right)_{11}$ across the thickness:
\begin{figure}[H]
	\centering
	\includegraphics[width=0.5\linewidth]{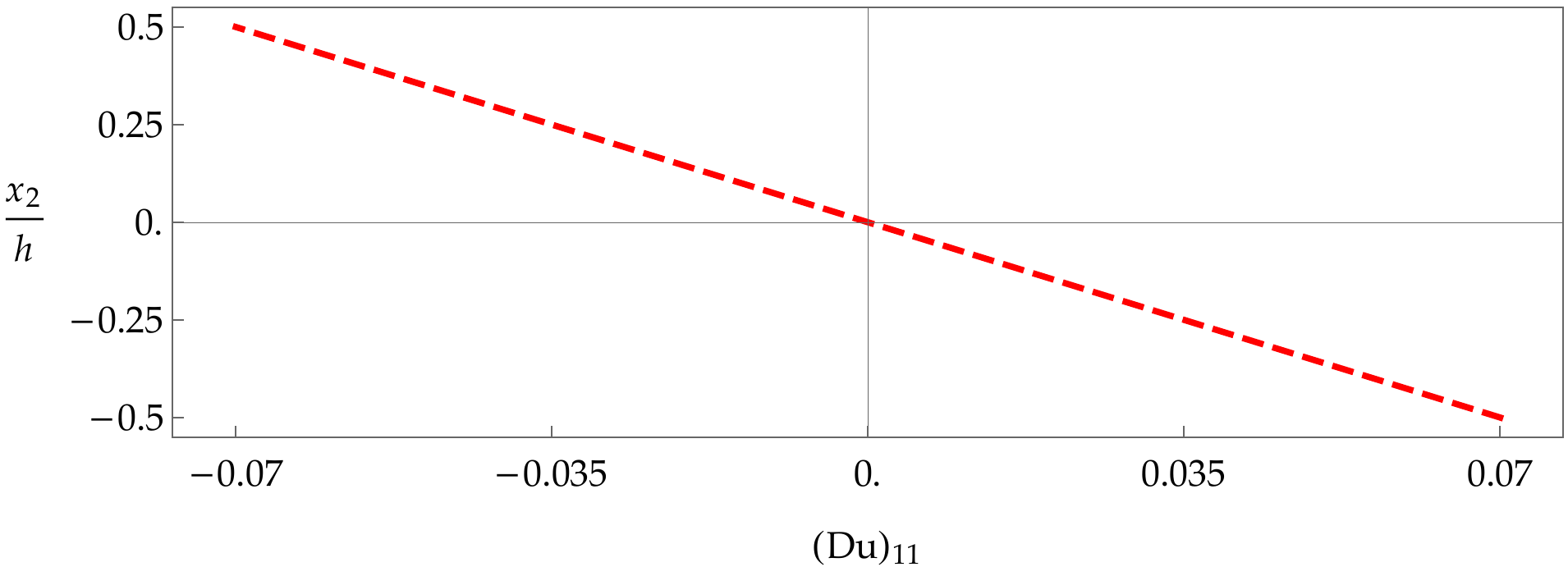}
	\caption{(\textbf{Second gradient model}, one curvature parameter) Plot of $(\boldsymbol{\mbox{D}u})_{11}$ across the thickness.
	}
	\label{fig:P11_SG}
\end{figure}
\subsection{Full isotropic curvature and arbitrary Poisson's ratios $\nu_{\tiny \mbox{macro}}$}
\begin{figure}[H]
	\centering
	\includegraphics[width=0.5\linewidth]{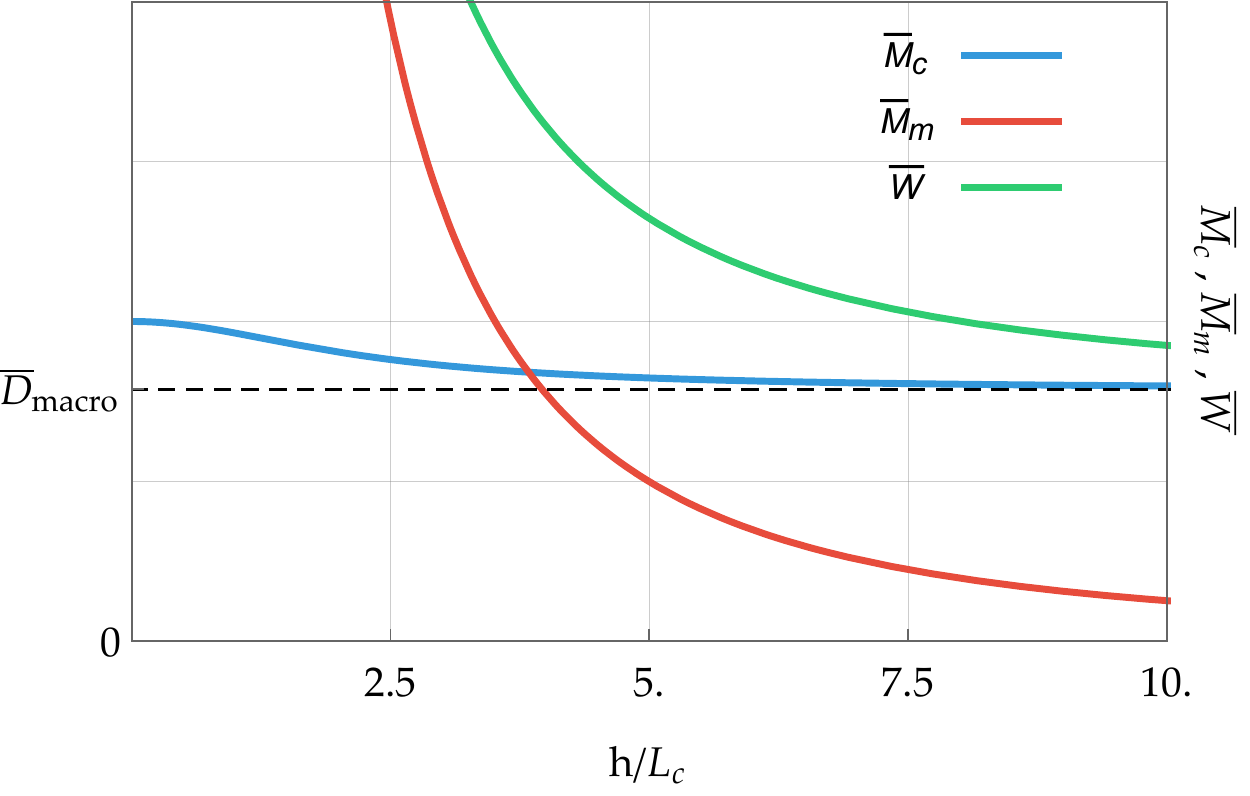}
	\caption{(\textbf{Second gradient model}, general case) Bending moments and energy while varying $L_c$. Observe that the bending stiffness is unbounded as $L_c \to \infty$ ($h\to 0$). This is a major difference with respect the relaxed micromorphic model. The values of the parameters used are: $\mu _{\tiny \mbox{macro}} = 1$, $\lambda _{\tiny \mbox{macro}} = 1$, $\mu = 1$, $a _1 = 2$, $a _2 = 1$, $a _3 = 1/2$.}
	\label{fig:all_plot_SG_three}
\end{figure}

The boundary conditions (see the Appendix~\ref{app:second_gradient_eq_eqa}) at the upper and lower surface (free surface) are (for the complete formulation in the classical notation see \cite{mindlin1964micro})
\begin{align}
	\widetilde{\boldsymbol{t}} (x_2 = \pm \, h/2) =& 
	\pm \, \widetilde{\boldsymbol{\sigma}} \cdot \boldsymbol{e}_2
	- \mu \, L_c^2 \,
	\displaystyle\sum_{i=1}^{3}
	\partial_{x_i}
	\boldsymbol{m}_i
	\cdot
	\boldsymbol{e}_2
	-
	\mu \, L_c^2 \,
	\displaystyle\sum_{i=1}^{3}
	\Big(
	\boldsymbol{m}_i \, (\boldsymbol{e}_{2})_{i}
	\Big)
	\cdot \boldsymbol{D^{\tau}}
	= \boldsymbol{0}
	\, ,
	\label{eq:BC_SG_three}
	\\
	\boldsymbol{\eta}(x_2 = \pm \, h/2) =& 
	\pm \, 
	\displaystyle\sum_{i=1}^{3}
	\left(
	\boldsymbol{m}_{i} \, (\boldsymbol{e}_{2})_{i}
	\right)
	\boldsymbol{e}_2
	=
	\boldsymbol{0}
	\, ,
	\notag
\end{align}
where $\widetilde{\boldsymbol{\sigma}} = 2 \mu_{\mbox{\tiny macro}} \, \mbox{sym}\,\boldsymbol{\mbox{D}u}
$ $+$ $ \lambda_{\tiny \mbox{macro}} \mbox{tr} \left(\boldsymbol{\mbox{D}u}\right) \boldsymbol{\mathbbm{1}}$,
$
\boldsymbol{m}_i=
\mu \, L_c^2 \,
\displaystyle\sum_{i=1}^{3}
\Big(
a_1 \mbox{dev} \, \mbox{sym} \left( \partial_{x_i}\boldsymbol{\mbox{D}u} \right)
$ $+$ $  a_2 \, \mbox{skew} \left( \partial_{x_i}\boldsymbol{\mbox{D}u} \right)
$ $+$ $ \frac{2}{9} \, a_3 \, \mbox{tr} \left( \partial_{x_i}\boldsymbol{\mbox{D}u} \right) \boldsymbol{\mathbbm{1}}
\Big)$, $i=1,2,3$, is the higher-order stress tensor,
$(\boldsymbol{e}_{2})_{i}$ is the scalar \textit{i}th component of the unit vector $\boldsymbol{e}_{2}$,
and
\begin{equation}
	\boldsymbol{D^{\tau}} = \left( \boldsymbol{\mathbbm{1}} - \boldsymbol{n} \otimes \boldsymbol{n}\right) \cdot (\partial_{x_1},\partial_{x_2},\partial_{x_3})
	\, .
	\label{eq:12}
\end{equation}
Substituting the ansatz eq.(\ref{eq:ansatz_SG_three}) in eq.(\ref{eq:equiMic_SG_three}) while choosing $a_1=a_2=1, a_3=\frac{3}{2}$, the equilibrium equation is
\begin{equation}
	-\frac{2}{9} \, \mu \, L_c^2 \, \left( 3 a_{1} + a_{3} \right) \, v^{(4)}(x_2) + \lambda _{\mbox{\tiny macro}} \, \left(v''(x_2) - \boldsymbol{\kappa} \right) + 2 \mu _{\mbox{\tiny macro}} \, v''(x_2) = 0 \, ,
	\label{eq:equi_equa_SG_three}
\end{equation}
and the solution of eq.(\ref{eq:equi_equa_SG_three}) is
\begin{align}
	v(x_{2}) =& \, \,
	\frac{1}{9} \, \frac{L_c^2}{f_1^2} \,
	\left(
	c_1 \, e^{-\frac{3 x_{2}}{L_c} f_1 }
	+
	c_2 \, e^{ \frac{3 x_{2}}{L_c} f_1 }
	\right)
	+ \frac{\lambda _{\mbox{\tiny macro}}}{2 \left(\lambda _{\mbox{\tiny macro}} + 2 \mu _{\mbox{\tiny macro}}\right)} \, \boldsymbol{\kappa}  \, x_{2}^2 
	+ c_4 \, x_2
	+ c_3 \, ,
	\label{eq:sol_v_SG_three}
	\\
	f_1 :=& \sqrt{\frac{\lambda _{\mbox{\tiny macro}}+2 \mu _{\mbox{\tiny macro}}}{2 \mu  (3 a_{1} + a_{3})}} \, .
	\notag
\end{align}
After applying the boundary conditions at the upper and lower surface (free surface)
it is possible to evaluate the displacement field solution of the cylindrical bending problem and then proceed to calculate the classical bending moment, the higher-order bending moment, and energy (per unit area d$x_1$d$x_3$) definitions which are reported in the following eq.(\ref{eq:sigm_ene_dimensionless_SG_three}):

\begin{align}
	M_{\mbox{c}} (\boldsymbol{\kappa})
	=&
	\displaystyle\int\limits_{-h/2}^{h/2}
	\left[
	\left(
	\widetilde{\boldsymbol{\sigma}} \cdot \boldsymbol{e}_1
	- \mu \, L_c^2 \,
	\displaystyle\sum_{i=1}^{3}
	\partial_{x_i}
	\boldsymbol{m}_i
	\,
	\boldsymbol{e}_1
	-
	\mu \, L_c^2 \,
	\displaystyle\sum_{i=1}^{3}
	\Big(
	\boldsymbol{m}_i \, (\boldsymbol{e}_{1})_{i}
	\Big)
	\cdot D^{\tau}
	\right) \,   \boldsymbol{e}_1
	\right]
	\, x_2 \, 
	\mbox{d}x_{2}
	\notag
	\\*
	=& \, 
	\frac{h^3}{12}
	\left[
	\frac{4 \mu_{\mbox{\tiny macro}} \left(\lambda_{\mbox{\tiny macro}}+\mu_{\mbox{\tiny macro}}\right)}{\lambda_{\mbox{\tiny macro}}+2 \mu_{\mbox{\tiny macro}}}
	- \frac{2}{3} \lambda _{\mbox{\tiny macro}} \left( \frac{4 \mu _{\mbox{\tiny macro}}}{\lambda _{\mbox{\tiny macro}}+2 \mu _{\mbox{\tiny macro}}}-\frac{9 a_{1}}{3 a_{1}+a_{3}} \right) \frac{1}{f_{1}^2} \left(\frac{L_c}{h}\right)^2
	\right.
	\notag
	\\*
	&
	\hspace{0.75cm}
	\left.
	- 
	\frac{4 \lambda _{\mbox{\tiny macro}}}{\lambda _{\mbox{\tiny macro}}+2 \mu _{\mbox{\tiny macro}}}
	\frac{1}{f_{1}^3}
	\left(\frac{a_{1} \lambda _{\mbox{\tiny macro}}}{3 a_{1}+a_{3}}+\frac{2 (3 a_{1}-2 a_{3}) \mu _{\mbox{\tiny macro}}}{9 (3 a_{1}+a_{3})}\right)
	\left(\frac{L_c}{h}\right)^3
	\right]
	\boldsymbol{\kappa}
	\, ,
	\notag
	\\
	M_{\mbox{m}}(\boldsymbol{\kappa})
	=&
	\displaystyle\int\limits_{-h/2}^{h/2}
	\langle
	\Big(\underbrace{\displaystyle\sum_{i=1}^{3}
		\left( \boldsymbol{m}_{i} \, (\boldsymbol{e}_{1})_{i} \right)}_{\mathbb{R}^{3\times3}}\Big) \boldsymbol{e}_1 , \boldsymbol{e}_2
	\rangle
	-
	\langle
	\Big(\underbrace{\displaystyle\sum_{i=1}^{3}
		\left( \boldsymbol{m}_{i} \, (\boldsymbol{e}_{2})_{i} \right)}_{\mathbb{R}^{3\times3}}\Big) \boldsymbol{e}_1 , \boldsymbol{e}_1
	\rangle
	-
	\langle
	\Big(\underbrace{\displaystyle\sum_{i=1}^{3}
		\left( \boldsymbol{m}_{i} \, (\boldsymbol{e}_{1})_{i} \right)}_{\mathbb{R}^{3\times3}}\Big) \boldsymbol{e}_2 , \boldsymbol{e}_1
	\rangle
	\,
	\mbox{d}x_{2}
	\notag
	\\*
	=&
	\frac{h^3}{12}
	\left[
	\frac{2}{3} \lambda _{\mbox{\tiny macro}}
	\left(
	\frac{4 \mu _{\mbox{\tiny macro}}}{\lambda _{\mbox{\tiny macro}}}
	+\frac{36 a_{2} \mu _{\mbox{\tiny macro}}}{(3 a_{1}+a_{3}) \lambda _{\mbox{\tiny macro}}}
	+\frac{9 (a_{1}+2 a_{2})}{3 a_{1}+a_{3}}
	\right)
	\frac{1}{f_{1}^2} \left(\frac{L_c}{h}\right)^2
	\right.
	\label{eq:sigm_ene_dimensionless_SG_three}
	\\*
	&
	\hspace{0.75cm}
	- 
	\frac{4 \lambda _{\mbox{\tiny macro}}}{\lambda _{\mbox{\tiny macro}}+2 \mu _{\mbox{\tiny macro}}}
	\frac{1}{f_{1}^3}
	\left(
	\frac{a_{1} (3 a_{1}-2 a_{3}) \left(\lambda _{\mbox{\tiny macro}}+2 \mu _{\mbox{\tiny macro}}\right)}{2 (3 a_{1}+a_{3})^2}
	\right.
	\notag
	\\*
	&
	\hspace{0.75cm}
	\left.\left.
	+\frac{(3 a_{1}-2 a_{3})^2 \mu _{\mbox{\tiny macro}} \left(\lambda _{\mbox{\tiny macro}}+2 \mu _{\mbox{\tiny macro}}\right)}{9 (3 a_{1}+a_{3})^2 \lambda _{\mbox{\tiny macro}}}
	\right)
	\left(\frac{L_c}{h}\right)^3
	\right]
	\boldsymbol{\kappa}
	\, ,
	\notag
	\\
	W_{\mbox{tot}} (\boldsymbol{\kappa})
	=&
	\displaystyle\int\limits_{-h/2}^{+h/2} W \left( \boldsymbol{\mbox{D}u},\boldsymbol{\mbox{D}^2 u} \right) \, \mbox{d}x_{2}
	\notag
	\\*
	=&
	\frac{1}{2}
	\frac{h^3}{12}
	\left[
	\frac{4 \mu_{\mbox{\tiny macro}} \left(\lambda_{\mbox{\tiny macro}}+\mu_{\mbox{\tiny macro}}\right)}{\lambda_{\mbox{\tiny macro}}+2 \mu_{\mbox{\tiny macro}}}
	+ \frac{2}{3} \lambda _{\mbox{\tiny macro}}
	\left(
	\frac{8 \mu_{\mbox{\tiny micro}}^2}{2 \lambda_{\mbox{\tiny micro}} \mu_{\mbox{\tiny micro}}+\lambda_{\mbox{\tiny micro}}^2}
	\right.
	\right.
	\notag
	\\*
	&
	\hspace{0.75cm}
	\left.
	\left.
	+\frac{36 a_{2} \mu _{\mbox{\tiny macro}}}{(3 a_{1}+a_{3}) \lambda _{\mbox{\tiny macro}}}
	+\frac{18 (a_{1}+a_{2})}{3 a_{1}+a_{3}}
	\right)
	\frac{1}{f_{1}^2} \left(\frac{L_c}{h}\right)^2
	\right.
	\notag
	\\*
	&
	\hspace{0.75cm}
	\left.
	- \frac{4 \lambda _{\mbox{\tiny macro}}}{\lambda _{\mbox{\tiny macro}}+2 \mu _{\mbox{\tiny macro}}}
	\frac{1}{f_{1}^3}
	\left(
	\frac{\left(2 (3 a_{1}-2 a_{3}) \mu _{\mbox{\tiny macro}}+9 a_{1} \lambda _{\mbox{\tiny macro}}\right){}^2}{18 (3 a_{1}+a_{3})^2 \lambda _{\mbox{\tiny macro}}}
	\right)
	\left(\frac{L_c}{h}\right)^3
	\right]
	\boldsymbol{\kappa}^2
	\, .
	\notag
\end{align}
As before
$
\frac{\mbox{d}}{\mbox{d}\boldsymbol{\kappa}}W_{\mbox{tot}}(\boldsymbol{\kappa}) = M_{\mbox{c}} (\boldsymbol{\kappa}) + M_{\mbox{m}} (\boldsymbol{\kappa})
= M_{\mbox{c}} (\boldsymbol{\kappa}) \, .
$
The plot of the bending moments and the strain energy divided by $\frac{h^3}{12}\boldsymbol{\kappa}$ and $\frac{1}{2}\frac{h^3}{12}\boldsymbol{\kappa}^2$, respectively, while changing $L_c$ is shown in Fig.~\ref{fig:all_plot_SG_three}.
\section{Summary and conclusions}
The present contribution presents the ansatz for solving the problem of pure cylindrical bending of elastic micromorphic continua. This ansatz is used to derive the solutions of different subclasses of micromorphic continua like the full micromorphic theory, microstrain theory and relaxed micromorphic theory with different approaches for the micro-curvature terms. The limiting case of very thin specimens (compared to the intrinsic length) is investigated and it is pointed out which theories yield bounded values of the flexural stiffness, thus providing hints on the choice of the type of theory and the respective parameters. Furthermore, the provided analytical solutions show the sensitivity of the flexural stiffness with respect to the constitutive parameters and thus offer a puzzle stone to identify these parameters. Finally, the analytical solutions are valuable benchmarks for numerical solution methods like FEM.

{\scriptsize
\paragraph{{\scriptsize Acknowledgements.}}
Angela Madeo and Gianluca Rizzi acknowledge funding from the French Research Agency ANR, “METASMART” (ANR-17CE08-0006).
Angela Madeo and Gianluca Rizzi acknowledge support from IDEXLYON in the framework of the “Programme Investissement d’Avenir” ANR-16-IDEX-0005.
Patrizio Neff acknowledges support in the framework of the DFG-Priority Programme 2256 ``Variational Methods for Predicting Complex Phenomena in Engineering Structures and Materials", Neff 902/10-1, Project-No. 440935806.
}

\let\oldbibliography\thebibliography
\renewcommand{\thebibliography}[1]{%
\oldbibliography{#1}%
\setlength{\itemsep}{2.5pt}%
}

\begin{scriptsize}
\bibliographystyle{plain}
\bibliography{Biblio_Bend_ArXiv}

\begin{thebibliography}{10}

\bibitem{altenbach2009linear}
H.~Altenbach and V.A. Eremeyev.
\newblock On the linear theory of micropolar plates.
\newblock {\em Zeitschrift für angewandte Mathematik und Mechanik},
  89(4):242--256, 2009.

\bibitem{altenbach2010generalized}
J.~Altenbach, H.~Altenbach, and V.A. Eremeyev.
\newblock On generalized {C}osserat-type theories of plates and shells: a short
  review and bibliography.
\newblock {\em Archive of Applied Mechanics}, 80(1):73--92, 2010.

\bibitem{arroyo2005continuum}
M.~Arroyo and T.~Belytschko.
\newblock Continuum mechanics modeling and simulation of carbon nanotubes.
\newblock {\em Meccanica}, 40(4-6):455--469, 2005.

\bibitem{barbagallo2016transparent}
G.~Barbagallo, A.~Madeo, M.V. {d}'Agostino, R.~Abreu, I.D. Ghiba, and P.~Neff.
\newblock {Transparent anisotropy for the relaxed micromorphic model:
  macroscopic consistency conditions and long wave length asymptotics}.
\newblock {\em International Journal of Solids and Structures}, 120:7--30,
  2017.

\bibitem{brcic2013estimation}
M.~Brcic, M.~Canadija, and J.~Brnic.
\newblock Estimation of material properties of nanocomposite structures.
\newblock {\em Meccanica}, 48(9):2209--2220, 2013.

\bibitem{corigliano2005chip}
A.~Corigliano, F.~Cacchione, B.~De~Masi, and C.~Riva.
\newblock On-chip electrostatically actuated bending tests for the mechanical
  characterization of polysilicon at the micro scale.
\newblock {\em Meccanica}, 40(4-6):485--503, 2005.

\bibitem{Cowin1983}
S.C. Cowin and J.W. Nunziato.
\newblock Linear elastic materials with voids.
\newblock {\em J. Elasticity.}, 13(2):125--147, 1983.

\bibitem{d2019effective}
M.V. d’Agostino, G.~Barbagallo, I.D. Ghiba, B.~Eidel, P.~Neff, and A.~Madeo.
\newblock Effective description of anisotropic wave dispersion in mechanical
  band-gap metamaterials via the relaxed micromorphic model.
\newblock {\em Journal of Elasticity}, 39:299--329, 2020.

\bibitem{DeCicco1997}
S.~De~Cicco and L.~Nappa.
\newblock Torsion and flexure of microstretch elastic circular cylinders.
\newblock {\em International Journal of Engineering Science}, 35(6):573--583,
  1997.

\bibitem{dell2009generalized}
F.~Dell'Isola, G.~Sciarra, and S.~Vidoli.
\newblock Generalized {H}ooke's law for isotropic second gradient materials.
\newblock {\em Proceedings of the Royal Society A: Mathematical, Physical and
  Engineering Sciences}, 465(2107):2177--2196, 2009.

\bibitem{Forest2018}
S.~Forest.
\newblock Micromorphic approach to materials with internal length.
\newblock In {\em Encyclopedia of Continuum Mechanics}, pages 1--11. Springer,
  Berlin, Heidelberg, 2018.

\bibitem{Forest2019}
S.~Forest.
\newblock Micromorphic approach to gradient plasticity and damage.
\newblock In {\em Handbook of Nonlocal Continuum Mechanics for Materials and
  Structures}, pages 499--546. Springer International Publishing, 2019.

\bibitem{forest2006nonlinear}
S.~Forest and R.~Sievert.
\newblock Nonlinear microstrain theories.
\newblock {\em International Journal of Solids and Structures},
  43(24):7224--7245, 2006.

\bibitem{gauthier1975quest}
R.D. Gauthier and W.E. Jahsman.
\newblock A quest for micropolar elastic constants.
\newblock {\em Journal of Applied Mechanics}, 42(2):369--374, 1975.

\bibitem{ghiba2017variant}
I.D. Ghiba, P.~Neff, A.~Madeo, and I.~M{\"u}nch.
\newblock A variant of the linear isotropic indeterminate couple-stress model
  with symmetric local force-stress, symmetric nonlocal force-stress, symmetric
  couple-stresses and orthogonal boundary conditions.
\newblock {\em Mathematics and Mechanics of Solids}, 22(6):1221--1266, 2017.

\bibitem{hadjesfandiari2011couple}
A.R. Hadjesfandiari and G.F. Dargush.
\newblock Couple stress theory for solids.
\newblock {\em International Journal of Solids and Structures},
  48(18):2496--2510, 2011.

\bibitem{hadjesfandiari2016pure}
A.R. Hadjesfandiari, A.~Hajesfandiari, and G.~F. Dargush.
\newblock Pure plate bending in couple stress theories.
\newblock {\em arXiv preprint arXiv:1606.02954}, 2016.

\bibitem{hutter2016application}
G.~H{\"u}tter.
\newblock Application of a microstrain continuum to size effects in bending and
  torsion of foams.
\newblock {\em International Journal of Engineering Science}, 101:81--91, 2016.

\bibitem{hutter2015micromorphic}
G.~H{\"u}tter, U.~M{\"u}hlich, and M.~Kuna.
\newblock Micromorphic homogenization of a porous medium: elastic behavior and
  quasi-brittle damage.
\newblock {\em Continuum Mechanics and Thermodynamics}, 27(6):1059--1072, 2015.

\bibitem{Iesan1971}
D.~Ieşan.
\newblock Torsion of micropolar elastic beams.
\newblock {\em International Journal of Engineering Science}, 9(11):1047--1060,
  1971.

\bibitem{Iesan1994}
D.~Ieşan and L.~Nappa.
\newblock {Saint-Venant's} problem for microstretch elastic solids.
\newblock {\em International Journal of Engineering Science}, 32(2):229--236,
  1994.

\bibitem{lakes1998elastic}
R.~Lakes.
\newblock Elastic freedom in cellular solids and composite materials.
\newblock In {\em Mathematics of Multiscale Materials}, pages 129--153.
  Springer, 1998.

\bibitem{Lakes1983}
R.S. Lakes.
\newblock Size effects and micromechanics of a porous solid.
\newblock {\em Journal of Materials Science}, 18(9):2572--2580, 1983.

\bibitem{lakes1995experimental}
R.S. Lakes.
\newblock Experimental methods for study of {C}osserat elastic solids and other
  generalized elastic continua.
\newblock {\em Continuum Models for Materials with Microstructure}, 70:1--25,
  1995.

\bibitem{Lakes2015}
R.S. Lakes and W.J. Drugan.
\newblock Bending of a {C}osserat elastic bar of square cross section: Theory
  and experiment.
\newblock {\em Journal of Applied Mechanics}, 82(9):091002, 2015.

\bibitem{lewintan2020korn}
P.~Lewintan, S.~M{\"u}ller, and P.~Neff.
\newblock Korn inequalities for incompatible tensor fields in three space
  dimensions with conformally invariant dislocation energy.
\newblock {\em arXiv preprint arXiv:2011.10573}, 2020.

\bibitem{lurie2018bending}
S.~Lurie, Y.~Solyaev, A.~Volkov, and D.~Volkov-Bogorodskiy.
\newblock Bending problems in the theory of elastic materials with voids and
  surface effects.
\newblock {\em Mathematics and Mechanics of Solids}, 23(5):787--804, 2018.

\bibitem{madeo2016new}
A.~Madeo, I.D. Ghiba, P.~Neff, and I.~M{\"{u}}nch.
\newblock A new view on boundary conditions in the
  {G}rioli--{K}oiter--{M}indlin--{T}oupin indeterminate couple stress model.
\newblock {\em European Journal of Mechanics-A/Solids}, 59:294--322, 2016.

\bibitem{mindlin1964micro}
R.D. Mindlin.
\newblock Micro-structure in linear elasticity.
\newblock {\em Archive for Rational Mechanics and Analysis}, 16(1):51--78,
  1964.

\bibitem{munch2018rotational}
I.~M{\"u}nch and P.~Neff.
\newblock Rotational invariance conditions in elasticity, gradient elasticity
  and its connection to isotropy.
\newblock {\em Mathematics and Mechanics of Solids}, 23(1):3--42, 2018.

\bibitem{munch2017modified}
I.~M{\"u}nch, P.~Neff, A.~Madeo, and I.D. Ghiba.
\newblock The modified indeterminate couple stress model: Why {Y}ang et al.'s
  arguments motivating a symmetric couple stress tensor contain a gap and why
  the couple stress tensor may be chosen symmetric nevertheless.
\newblock {\em Zeitschrift f{\"u}r Angewandte Mathematik und Mechanik},
  97(12):1524--1554, 2017.

\bibitem{neff2004material}
P.~Neff.
\newblock {On material constants for micromorphic continua}.
\newblock In {\em Trends in Applications of Mathematics to Mechanics, STAMM
  Proceedings, Seeheim}, pages 337--348. Shaker--Verlag, 2004.

\bibitem{neff2019identification}
P.~Neff, B.~Eidel, M.V. d’Agostino, and A.~Madeo.
\newblock Identification of scale-independent material parameters in the
  relaxed micromorphic model through model-adapted first order homogenization.
\newblock {\em Journal of Elasticity}, 139:269–298, 2020.

\bibitem{neff2015correct}
P.~Neff, I.D. Ghiba, A.~Madeo, and I.~M{\"u}nch.
\newblock Correct traction boundary conditions in the indeterminate couple
  stress model.
\newblock {\em arXiv preprint arXiv:1504.00448}, 2015.

\bibitem{neff2014unifying}
P.~Neff, I.D. Ghiba, A.~Madeo, L.~Placidi, and G.~Rosi.
\newblock A unifying perspective: the relaxed linear micromorphic continuum.
\newblock {\em Continuum Mechanics and Thermodynamics}, 26(5):639--681, 2014.

\bibitem{neff2009new}
P.~Neff and J.~Jeong.
\newblock {A new paradigm: the linear isotropic {C}osserat model with
  conformally invariant curvature energy}.
\newblock {\em Zeitschrift f{\"{u}}r Angewandte Mathematik und Mechanik},
  89(2):107--122, 2009.

\bibitem{neff2010stable}
P.~Neff, J.~Jeong, and A.~Fischle.
\newblock Stable identification of linear isotropic {C}osserat parameters:
  bounded stiffness in bending and torsion implies conformal invariance of
  curvature.
\newblock {\em Acta Mechanica}, 211(3-4):237--249, 2010.

\bibitem{Park1987}
H.C. Park and R.S. Lakes.
\newblock Torsion of a micropolar elastic prism of square cross-section.
\newblock {\em International Journal of Solids and Structures}, 23(4):485--503,
  1987.

\bibitem{renda2020geometric}
F.~Renda, C.~Armanini, V.~Lebastard, F.~Candelier, and F.~Boyer.
\newblock A geometric variable-strain approach for static modeling of soft
  manipulators with tendon and fluidic actuation.
\newblock {\em IEEE Robotics and Automation Letters}, 5(3):4006--4013, 2020.

\bibitem{rizzi2020shear}
G.~Rizzi, G.~H{\"u}tter, A.~Madeo, and P.~Neff.
\newblock Analytical solutions of the simple shear problem for micromorphic
  models and other generalized continua.
\newblock {\em to appear in Archive of Applied Mechanics}, 2021.

\bibitem{rueger2019cosserat}
Z.~Rueger, C.S. Ha, and R.S. Lakes.
\newblock Cosserat elastic lattices.
\newblock {\em Meccanica}, 54(13):1983--1999, 2019.

\bibitem{shaat2018reduced}
M.~Shaat.
\newblock A reduced micromorphic model for multiscale materials and its
  applications in wave propagation.
\newblock {\em Composite Structures}, 201:446--454, 2018.

\bibitem{Taliercio2010}
A.~Taliercio.
\newblock Torsion of micropolar hollow circular cylinders.
\newblock {\em Mechanics Research Communications}, 37(4):406--411, 2010.

\bibitem{tekouglu2008size}
C.~Teko{\u{g}}lu and P.R. Onck.
\newblock Size effects in two-dimensional {V}oronoi foams: a comparison between
  generalized continua and discrete models.
\newblock {\em Journal of the Mechanics and Physics of Solids},
  56(12):3541--3564, 2008.

\bibitem{Waseem2013}
A.~Waseem, A.J. Beveridge, M.A. Wheel, and D.H. Nash.
\newblock The influence of void size on the micropolar constitutive properties
  of model heterogeneous materials.
\newblock {\em European Journal of Mechanics-A/Solids}, 40:148--157, 2013.

\bibitem{Yang1982}
J.F.C. Yang and R.S. Lakes.
\newblock Experimental study of micropolar and couple stress elasticity in
  compact bone in bending.
\newblock {\em Journal of Biomechanics}, 15(2):91--98, 1982.

\bibitem{zhang2017application}
L.~Zhang, Binbin L., S.~Zhou, B.~Wang, and Y.~Xue.
\newblock An application of a size-dependent model on microplate with elastic
  medium based on strain gradient elasticity theory.
\newblock {\em Meccanica}, 52(1-2):251--262, 2017.

\end{thebibliography}
\end{scriptsize}


\begin{footnotesize}
\appendix
\section{Appendix}
With the ansatz eq.(\ref{eq:ansatz_RM}) the following relations hold:
\begin{align}
\mbox{dev} \, \mbox{Curl} \boldsymbol{P} &= 
\mbox{Curl} \boldsymbol{P} = 
\left(
\begin{array}{ccc}
0 & 0 & -\boldsymbol{\kappa} -P_{11}'(x_{2}) \\
0 & 0 & 0 \\
P_{33}'(x_{2}) & 0 & 0 \\
\end{array}
\right) \, ,
\qquad
\mbox{tr}\left(\mbox{Curl} \boldsymbol{P}\right) =
\boldsymbol{0} \, ,
\notag
\\ 
\mbox{skew} \, \mbox{Curl} \boldsymbol{P} &= 
\frac{1}{2}
\left(
\begin{array}{ccc}
0 & 0 & -\boldsymbol{\kappa} -P_{11}'(x_{2})-P_{33}'(x_{2}) \\
0 & 0 & 0 \\
\boldsymbol{\kappa} +P_{11}'(x_{2})+P_{33}'(x_{2}) & 0 & 0 \\
\end{array}
\right) \, ,
\label{eq:curl_part}
\\ 
\mbox{dev} \, \mbox{sym} \, \mbox{Curl} \boldsymbol{P} &= 
\mbox{sym} \, \mbox{Curl} \boldsymbol{P} = 
\frac{1}{2}
\left(
\begin{array}{ccc}
0 & 0 & P_{33}'(x_{2}) - \boldsymbol{\kappa} -P_{11}'(x_{2}) \\
0 & 0 & 0 \\
P_{33}'(x_{2}) -\boldsymbol{\kappa} -P_{11}'(x_{2}) & 0 & 0 \\
\end{array}
\right) \, .
\notag
\end{align}
Given eq.(\ref{eq:curl_part}) we observe
\begin{equation}
\begin{array}{l}
\left \lVert \mbox{Curl}\boldsymbol{P} \right \rVert^2 =
\left \lVert \mbox{dev} \, \mbox{Curl}\boldsymbol{P} \right \rVert^2
\, ,
\quad
\left \lVert \mbox{sym} \, \mbox{Curl}\boldsymbol{P} \right \rVert^2 =
\left \lVert \mbox{dev} \, \mbox{sym} \, \mbox{Curl}\boldsymbol{P} \right \rVert^2
\end{array}
\end{equation}
and also
\begin{align}
&
\pushleft{\quad a_1\left \lVert \mbox{dev} \, \mbox{sym} \, \mbox{Curl}\boldsymbol{P} \right \rVert^2 + 
a_2\left \lVert \mbox{skew} \, \mbox{Curl}\boldsymbol{P} \right \rVert^2 + 
\frac{a_3}{3} \, \mbox{tr}^2 \left( \mbox{Curl}\boldsymbol{P} \right)}
\notag
\\*
&
\hspace{2.5cm}
=\frac{1}{2} \left(a_{1} \left(\boldsymbol{\kappa} +P_{11}'(x_{2})-P_{33}'(x_{2})\right)^2+a_{2} \left(\boldsymbol{\kappa} +P_{11}'(x_{2})+P_{33}'(x_{2})\right)^2\right)
\\*
&
\hspace{2.5cm}
=\frac{1}{2} \left((a_{1}+a_{2}) \left(\left(P_{11}'(x_2) + \boldsymbol{\kappa}\right)^2+P_{33}'(x_2)^2\right)-2 (a_{1}-a_{2}) P_{33}'(x_2) \left(\boldsymbol{\kappa} +P_{11}'(x_2)\right)\right) \, .
\notag
\end{align}
\section{Generalized Neumann boundary conditions for the relaxed micromorphic model and the Cosserat model}
\label{app:neumann_coss}
Partial integration for the matrix-Curl operator can be written as
\begin{align}
\displaystyle\int\limits_{\partial \Omega} \left\langle \boldsymbol{P} \times \boldsymbol{\nu} , \boldsymbol{Q} \right\rangle \, \mbox{dS} &=
\displaystyle\int\limits_{\Omega} \left\langle \mbox{Curl} \, \boldsymbol{P} , \boldsymbol{Q} \right\rangle
- \left\langle \boldsymbol{P} , \mbox{Curl} \, \boldsymbol{Q} \right\rangle \mbox{d}x \, ,
\\
\displaystyle\int\limits_{\partial \Omega} \left\langle \boldsymbol{P} \times \boldsymbol{\nu} ,\mbox{Curl} \, \boldsymbol{A} \right\rangle  \, \mbox{dS} &=
\displaystyle\int\limits_{\Omega} \left\langle \mbox{Curl} \, \boldsymbol{P} ,\mbox{Curl} \, \boldsymbol{A} \right\rangle
- \left\langle \boldsymbol{P} , \mbox{Curl} \, \mbox{Curl} \, \boldsymbol{A} \right\rangle \mbox{d}x \, ,
\notag
\end{align}
where $\boldsymbol{P,Q} \in \mbox{C}^{1} \left(\overline{\Omega},\mathbb{R}^{3\times3}\right)$ are sufficiently smooth square $3\times3$ matrix fields and $\nu$ is the outward unit normal vector to $\partial\Omega$.
Inserting $\boldsymbol{P}= \delta \boldsymbol{A}$ and $\boldsymbol{Q}= \mbox{Curl} \boldsymbol{A}$ for $\boldsymbol{A} \in \mbox{C}^{1} \left(\overline{\Omega},\mathfrak{so}(3)\right)$, for a test field $\delta \boldsymbol{A}$ and argument $\boldsymbol{A} \in \mathfrak{so}(3)$, we obtain
\begin{equation}
\displaystyle\int\limits_{\partial \Omega} \left\langle \delta \boldsymbol{A} \times \boldsymbol{\nu} ,\mbox{Curl} \, \boldsymbol{A} \right\rangle \, \mbox{dS} =
\displaystyle\int\limits_{\Omega} \left\langle \mbox{Curl} \, \delta \boldsymbol{A} ,\mbox{Curl} \, \boldsymbol{A} \right\rangle 
- \left\langle \delta \boldsymbol{A} , \mbox{Curl} \, \mbox{Curl} \, \boldsymbol{A} \right\rangle \mbox{d}x \, .
\end{equation}
The scalar-product on the left-hand side is interpreted row-wise.
Making use of the permutation properties of the scalar product, namely 
\begin{equation}
\langle \boldsymbol{a} \times \boldsymbol{\nu}, \boldsymbol{b}\rangle_{\mathbb{R}^3} = 
\mbox{det}\left[\boldsymbol{b},\boldsymbol{a},\boldsymbol{\nu}\right] = 
-\mbox{det}\left[\boldsymbol{a},\boldsymbol{b},\boldsymbol{\nu}\right] = 
-\langle \boldsymbol{a} , \boldsymbol{b} \times \boldsymbol{\nu} \rangle
\end{equation}
we arrive at
\begin{equation}
\displaystyle\int\limits_{\partial \Omega} \left\langle \delta \boldsymbol{A} ,\mbox{Curl} \, \boldsymbol{A} \times \boldsymbol{\nu} \right\rangle \, \mbox{dS} =
\displaystyle\int\limits_{\Omega} \left\langle \delta \boldsymbol{A} , \mbox{Curl} \, \mbox{Curl} \, \boldsymbol{A} \right\rangle
- \left\langle \mbox{Curl} \, \delta \boldsymbol{A} ,\mbox{Curl} \, \boldsymbol{A} \right\rangle \mbox{d}x \, .
\end{equation}
Since $\delta \boldsymbol{A} \in \mbox{C}^{1} \left(\overline{\Omega},\mathfrak{so}(3)\right)$ this gives equivalently 
\begin{equation}
\displaystyle\int\limits_{\partial \Omega} \left\langle \delta \boldsymbol{A} ,\mbox{skew} \left(\mbox{Curl} \, \boldsymbol{A} \times \boldsymbol{\nu}\right) \right\rangle \, \mbox{dS} =
\displaystyle\int\limits_{\Omega} \left\langle \delta \boldsymbol{A} , \mbox{skew} \left( \mbox{Curl} \, \mbox{Curl} \, \boldsymbol{A} \right) \right\rangle 
- \left\langle \mbox{Curl} \, \delta \boldsymbol{A} , \mbox{Curl} \, \boldsymbol{A} \right\rangle \mbox{d}x \, .
\end{equation}
Replacing $\boldsymbol{Q} = \boldsymbol{m} = \mbox{Curl} \, \boldsymbol{A}$ yields
\begin{equation}
\displaystyle\int\limits_{\partial \Omega} \left\langle \delta \boldsymbol{A} ,\mbox{skew} \left(\boldsymbol{m} \times \boldsymbol{\nu}\right) \right\rangle \, \mbox{dS} =
\displaystyle\int\limits_{\Omega} \left\langle \delta \boldsymbol{A} , \mbox{skew} \left( \mbox{Curl} \, \boldsymbol{m} \right) \right\rangle 
-  \left\langle \mbox{Curl} \, \delta \boldsymbol{A} , \boldsymbol{m} \right\rangle \mbox{d}x \, ,
\end{equation}
where the appropriate localization shows $\mbox{skew} \left(\boldsymbol{m} \times \boldsymbol{\nu}\right)\Bigr|_{\Gamma}=0$ since $\delta \boldsymbol{A}$ is arbitrary on $\Gamma$.
\section{The Lie-algebra $\mathfrak{so}(3)$, the 3D-Curl on $\mathfrak{so}(3)$ and Nye's relation}
\label{Sec:appendix_3D-Curl}

Given $\boldsymbol{A} \in \mathfrak{so}(3)$ 

\begin{equation}
\boldsymbol{A} = 
\left(
\begin{array}{ccc}
0  & -a_3 &  a_2 \\
a_3 &   0  & -a_1 \\
-a_2 &  a_1 &   0 
\end{array}
\right)
\label{eq:appC_A}
\end{equation}
the operator axl: $\mathfrak{so}(3)\in \mathbb{R}^3$ is introduced
\begin{equation}
\mbox{axl}
\left(
\begin{array}{ccc}
0  & -a_3 &  a_2 \\
a_3 &   0  & -a_1 \\
-a_2 &  a_1 &   0 
\end{array}
\right)
:=
\left(
\begin{array}{ccc}
a_1 \\
a_2 \\
a_3 
\end{array}
\right),
\qquad
\boldsymbol{A} \cdot \boldsymbol{v} = \left(\mbox{axl} \boldsymbol{A}\right) \times \boldsymbol{v},
\quad
\forall \boldsymbol{v} \in \mathbb{R}^3 \, .
\label{eq:appC_axl_A}
\end{equation}
Given the definition eqs.~(\ref{eq:appC_A})-(\ref{eq:appC_axl_A}), the following identities hold (\textbf{Nye's relation})
\begin{equation}
- \mbox{Curl} \boldsymbol{A} = 
\left( \boldsymbol{\mbox{D}} \mbox{axl} \boldsymbol{A} \right)^T 
- \mbox{tr}\left[\left( \boldsymbol{\mbox{D}} \mbox{axl} \boldsymbol{A} \right)^T\right] \cdot \boldsymbol{\mathbbm{1}},
\quad
\boldsymbol{\mbox{D}} \mbox{axl} \boldsymbol{A}   = 
- \left(\mbox{Curl} \boldsymbol{A}\right)^T
+ \frac{1}{2} \mbox{tr}\left[\left( \mbox{Curl} \boldsymbol{A} \right)^T\right] \cdot \boldsymbol{\mathbbm{1}} \, .
\label{eq:appC_ident}
\end{equation}

If we now have $\boldsymbol{A} = \mbox{skew} \, \boldsymbol{\mbox{D}u}$ it is possible to show that 
\begin{equation}
- \left( \mbox{Curl} \, \mbox{skew} \, \boldsymbol{\mbox{D}u} \right)^T= 
\boldsymbol{\mbox{D}} \, \mbox{axl}  \left(\mbox{skew} \, \boldsymbol{\mbox{D}u}\right)  ,
\qquad\qquad
\frac{1}{2} \, \mbox{curl} \boldsymbol{u} = 
\mbox{axl}\left( \mbox{skew} \, \boldsymbol{\mbox{D}u}\right),
\label{eq:nye_2}
\end{equation}
which leads to the following identity for the full Curl,
$\left\lVert \mbox{Curl} \, \mbox{skew} \, \boldsymbol{\mbox{D}u} \right\rVert^{2}
=
\frac{1}{4} \left\lVert \boldsymbol{\mbox{D}} \mbox{curl} \, \boldsymbol{u} \right\rVert^{2}.
$

It is also highlighted here that, thanks to eq.(\ref{eq:nye_2})$_1$ this relation holds
$
\mbox{Curl} \, \mbox{skew} \, \boldsymbol{\mbox{D}u}  = 
-  \mbox{Curl} \, \mbox{sym} \, \boldsymbol{\mbox{D}u}
$
(see \cite{ghiba2017variant}),
since 
$
0 = \mbox{Curl} \boldsymbol{\mbox{D}u}=\mbox{Curl} \, \mbox{skew} \, \boldsymbol{\mbox{D}u}  +  \mbox{Curl} \, \mbox{sym} \, \boldsymbol{\mbox{D}u}
$
, which implies that choosing the symmetric or the skew-symmetric part of the gradient of the displacement does not make any difference besides a sign.
Moreover $\mbox{tr}(\mbox{Curl} \boldsymbol{S})=0$ for any symmetric matrix $\boldsymbol{S} \in \mbox{Sym}(3)$
\section{Cylindrical bending for the isotropic Cosserat continuum with classical notation}
\label{sec:Cos_Class}
In \cite{neff2009new} (eq.(2.2)) there is the correspondence between the isotropic Cosserat model with rotation vector and the Curl representation in dislocation format.
Both have three curvature parameters and the identification is given by
\begin{equation}
\alpha = \frac{1}{3}\left(4a_3-a_1\right),
\qquad\qquad \beta = \frac{a_1-a_2}{2},
\qquad\qquad \gamma = \frac{a_1+a_2}{2} \, .
\label{eq:coeff_coss_class}
\end{equation}
Setting $\vartheta := \mbox{axl} (\boldsymbol{A})$ and taking into account eqs.~(\ref{eq:appC_axl_A})--(\ref{eq:coeff_coss_class}), the expression of the strain energy for the isotropic Cosserat continuum can be equivalently written as:
\begin{align}
W \left(\boldsymbol{\mbox{D}u}, \boldsymbol{A},\mbox{Curl}\,\boldsymbol{A}\right) = &
\, \mu_{\mbox{\tiny macro}} \left\lVert \mbox{sym} \, \boldsymbol{\mbox{D}u} \right\rVert^{2}
+ \frac{\lambda_{\mbox{\tiny macro}}}{2} \mbox{tr}^2 \left(\boldsymbol{\mbox{D}u} \right) 
+ \mu_{c} \left\lVert \mbox{skew} \, \boldsymbol{\mbox{D}u} - \boldsymbol{A} \right\rVert^{2}
\notag
\\
&+ \frac{\mu \, L_c^2}{2}
\underbrace{
	\left(
	a_1 \, \left \lVert \mbox{dev} \, \mbox{sym} \, \mbox{Curl} \, \boldsymbol{A}\right \rVert^2 \, 
	+ a_2 \, \left \lVert \mbox{skew} \, \mbox{Curl} \, \boldsymbol{A}\right \rVert^2 \, 
	+ \frac{a_3}{3} \, \mbox{tr}^2 \left(\mbox{Curl} \, \boldsymbol{A} \right)
	\right)
}_{\mbox{dislocation tensor format}}
\label{eq:energy_Cos_classic}
\\
= W \left(\boldsymbol{\mbox{D}u}, \boldsymbol{\vartheta},\boldsymbol{\mbox{D} \vartheta}\right) = &
\, \mu_{\mbox{\tiny macro}} \left\lVert \mbox{sym} \, \boldsymbol{\mbox{D}u} \right\rVert^{2}
+ \frac{\lambda_{\mbox{\tiny macro}}}{2} \mbox{tr}^2 \left(\boldsymbol{\mbox{D}u} \right) 
+ \frac{\mu_{c}}{2} \left\lVert \mbox{curl} \boldsymbol{u} - 2\boldsymbol{\vartheta} \right\rVert^{2}
\notag
\\
&+ \frac{\mu \, L_c^2}{2}
\underbrace{
	\left(
	\alpha \, \mbox{tr}^2 \left(\boldsymbol{\mbox{D} \vartheta} \right)
	+ \beta \, \langle \boldsymbol{\mbox{D} \vartheta}^{T},\boldsymbol{\mbox{D} \vartheta} \rangle 
	+ \gamma \, \left \lVert \boldsymbol{\mbox{D} \vartheta} \right \rVert^2
	\right)
}_{\mbox{micro-rotation vector format}}  \, ,
\notag
\end{align}
since 
\begin{equation}
\begin{split}
\left\lVert \mbox{skew} \, \boldsymbol{\mbox{D}u} - \boldsymbol{A} \right\rVert^{2} =
2\left\lVert \mbox{axl}(\mbox{skew} \, \boldsymbol{\mbox{D}u} - \mbox{Anti}(\boldsymbol{\vartheta})) \right\rVert^{2} =
2\lVert \frac{1}{2}\mbox{curl} \boldsymbol{u} - \boldsymbol{\vartheta} \rVert^{2} =
\frac{1}{2}\left\lVert \mbox{curl} \boldsymbol{u} - 2\boldsymbol{\vartheta} \right\rVert^{2} \, .
\end{split}
\end{equation}
The equilibrium equations without body forces in the classical notation are now the following
\begin{align}
\mbox{Div}\left[
2\mu_{\mbox{\tiny macro}}\,\mbox{sym} \, \boldsymbol{\mbox{D}u} + \lambda_{\mbox{\tiny macro}} \mbox{tr} \left(\boldsymbol{\mbox{D}u} \right) \boldsymbol{\mathbbm{1}}
\right]
-\mu_{c} \, \mbox{curl} \left[\mbox{curl} \, \boldsymbol{u} - 2 \boldsymbol{\vartheta}\right]  \,  
&= \boldsymbol{0} \, ,
\label{eq:equi_Cos_classic}
\\*
\mu \, L_c^2 \, \mbox{Div}\left[
\alpha \, \mbox{tr} \left(\boldsymbol{\mbox{D} \vartheta} \right) \, \boldsymbol{\mathbbm{1}}
+ \beta \, \left(\boldsymbol{\mbox{D} \vartheta}\right)^{T} 
+ \gamma \, \boldsymbol{\mbox{D} \vartheta}  \, 
\right] 
+2\mu_c \, \left(\mbox{curl} \, \boldsymbol{u} - 2 \boldsymbol{\vartheta}\right)
&= \boldsymbol{0} \, .
\notag
\end{align}
The boundary conditions at the upper and lower surface (free surface) are 
\begin{align}
\boldsymbol{\widetilde{t}}(x_2 = \pm \, h/2) = 
\pm \, \boldsymbol{\widetilde{\sigma}}(x_2) \cdot \boldsymbol{e}_2 = 
\boldsymbol{0} \, ,
\qquad\qquad
\overline{\boldsymbol{\eta}}(x_2 = \pm \, h/2) = 
\pm \, \overline{\boldsymbol{m}} (x_2) \cdot \boldsymbol{e}_2 = 
\boldsymbol{0} \, ,
\label{eq:BC_Cos_classic_gen}
\end{align}
where $\boldsymbol{\widetilde{\sigma}} = 2\mu_{\mbox{\tiny macro}}\,\mbox{sym} \, \boldsymbol{\mbox{D}u} + \lambda_{\mbox{\tiny macro}} \mbox{tr} \left(\boldsymbol{\mbox{D}u} \right) \boldsymbol{\mathbbm{1}} + 2\mu_c \, \left(\mbox{skew} \, \boldsymbol{\mbox{D}u} - \mbox{Anti}(\boldsymbol{\vartheta})\right)$, $\boldsymbol{e}_2$ is the unit vector aligned to the $x_2$-direction, and the second-order moment stress tensor $\overline{\boldsymbol{m}} := \mu \, L_c^2 \, \left(\alpha \, \mbox{tr} \left(\boldsymbol{\mbox{D} \vartheta} \right) \, \boldsymbol{\mathbbm{1}}+ \beta \, \left(\boldsymbol{\mbox{D} \vartheta}\right)^{T} + \gamma \, \boldsymbol{\mbox{D} \vartheta}\right)$.
The relation to the higher-order stress tensor reported in Sect.~\ref{sec:Cos} is the following:
\begin{equation}
\mbox{dev} \left(\boldsymbol{m}^{T}\right) = -\mbox{dev} \left(\overline{\boldsymbol{m}}\right) \, ,
\qquad
\mbox{tr} \left(\boldsymbol{m}\right) = \frac{1}{2}\mbox{tr} \left(\overline{\boldsymbol{m}}\right) \, .
\label{eq:m_class_vs_m_curl}
\end{equation}

According to the reference system shown in Fig.~\ref{fig:intro}, the ansatz for the displacement field and the micro-rotation vector is
\begin{equation}
\boldsymbol{u}(x_1,x_2)=
\left(
\begin{array}{c}
-\boldsymbol{\kappa}_1 \, x_1 x_2 \\
v(x_2)+\frac{\boldsymbol{\kappa}_1  x_1^2}{2} \\
0 \\
\end{array}
\right) \, ,
\qquad\qquad
\boldsymbol{\vartheta}(x_1,x_2) =
\left(
\begin{array}{ccc}
0  \\
0 \\
- \boldsymbol{\kappa}_2 \, x_1 \\
\end{array}
\right) \, .
\label{eq:ansatz_Cos_classic}
\end{equation}
Substituting the ansatz eq.(\ref{eq:ansatz_Cos_classic}) in eq.(\ref{eq:equi_Cos_classic}) the equilibrium equations result in
\begin{equation}
\begin{split}
2 \mu _c (\boldsymbol{\kappa}_1-\boldsymbol{\kappa}_2 )-\boldsymbol{\kappa}_1 \lambda _e+\left(\lambda _e+2 \mu _e\right) v''(x_{2}) = 0 \, ,
\quad \, \,
2 x_1 \mu _c (\boldsymbol{\kappa}_2 -\boldsymbol{\kappa}_1) = 0 \, ,
\quad \, \,
2 x_1 \mu _c  (\boldsymbol{\kappa}_1 -\boldsymbol{\kappa}_2) = 0 \, ,
\end{split}
\label{eq:equi_equa_Cos_classic}
\end{equation}
which are exactly the same as the eq.(\ref{eq:equi_equa_Cos}) in Section~\ref{sec:Cos}. Since also the boundary conditions eq.(\ref{eq:BC_Cos_classic_gen}) are equivalent to the boundary condition eq.(\ref{eq:BC_Cos_gen}) in Section~\ref{sec:Cos}, further calculations are avoided.
It is nevertheless interesting to show the definition of the higher-order bending moment:
\begin{equation}
M_{\mbox{m}}(\boldsymbol{\kappa}) =
\displaystyle\int\limits_{-h/2}^{h/2}
\langle \overline{\boldsymbol{m}} \, \boldsymbol{e}_1 , \boldsymbol{e}_3 \rangle \,
\mbox{d}x_{2}
=
\displaystyle\int\limits_{-h/2}^{h/2}
\langle \left(\boldsymbol{m} \times \boldsymbol{e}_1 \right) \boldsymbol{e}_2 , \boldsymbol{e}_1 \rangle \,
\mbox{d}x_{2}
= 
h \,  \mu \, L_c^2 \, \gamma \, \boldsymbol{\kappa} \, ,
\label{eq:sigm_ene_dimensionless_Cos_classic}
\end{equation}
given that 
\begin{equation}
\langle \left(\boldsymbol{m} \times \boldsymbol{e}_1 \right) \boldsymbol{e}_2 , \boldsymbol{e}_1 \rangle = 
\langle \left(\boldsymbol{m}, \boldsymbol{e}_2 \times \boldsymbol{e}_1 \right), \boldsymbol{e}_1 \rangle = 
-\langle \boldsymbol{m} \, \boldsymbol{e}_3 , \boldsymbol{e}_1 \rangle = 
-\langle \boldsymbol{m}^T \, \boldsymbol{e}_1, \boldsymbol{e}_3 \rangle = 
-\langle \overline{\boldsymbol{m}} \, \boldsymbol{e}_1, \boldsymbol{e}_3 \rangle \, .
\end{equation}
\subsection{Lakes formula}
\label{app:lakes}
In order to connect ourselves to the existing literature, we provide the reader with an excerpt taken from Lakes \cite{lakes1998elastic} to which we compare our results:
``Cosserat solids may be characterized via size effects in rigidity. Exact analytical solutions
for size effects form the basis of a variety of experiments for the characterization of Cosserat solids.
For example, Gauthier and Jahsman \cite{gauthier1975quest} give $\Omega$, \textbf{the ratio of rigidity to its classical value}, for cylindrical bending of a plate
\begin{equation}
\Omega := \left(1 + 24 \frac{\ell_b^2 \, \left(1-\nu \right)}{h^2} \right),
\qquad\qquad
\ell_b = \sqrt{\frac{\gamma}{4\mu}}
\label{eq:Omega_lakes}
\end{equation}
with $h$ the plate thickness''.

In plate bending, the anticlastic curvature due to the Poisson's effect is constrained, in contrast to (classical) beam bending.
\begin{align}
M^{\mbox{\small Lakes}}_{\mbox{tot}} =& \frac{\widehat{E} \, J_{x_3}}{R} \left(1 + 24 \left(\frac{l_b}{h}\right)^2 \left(1-\nu \right)\right)
=
\frac{E}{1-\nu^2}\frac{h^3}{12}\frac{1}{R} \left(1 + 24 \frac{\gamma}{4\mu}\frac{1}{h^2} \left(1-\nu \right)\right)
\notag
\\*
=&
\frac{E}{1-\nu^2}\frac{h^3}{12}\frac{1}{R} \left(1 + 24 \gamma \, \frac{(1+\nu)}{2 \, E}\frac{1}{h^2} \left(1-\nu \right)\right)
=
\frac{E}{1-\nu^2}\frac{h^3}{12}\frac{1}{R} + \frac{\gamma \, h}{R} \, ,
\notag
\\
M_{\mbox{tot}} =&
\frac{4\mu _e \left(\lambda _e+\mu _e\right)}{\lambda _e+2 \mu _e} \, 
\frac{h^3 }{12} \, \boldsymbol{\kappa} 
+
h \, \widetilde{\gamma} \, \boldsymbol{\kappa} 
=
\frac{E}{1-\nu^2} \, 
\frac{h^3 }{12} \, \frac{1}{R}
+
\frac{\widetilde{\gamma} \, h}{R} \, ,
\\*
\boldsymbol{\kappa}=&\frac{1}{R} \, ,
\quad \,\,\,\,\,
\widetilde{\gamma} := \mu \, L_c^2 \, \gamma \, ,
\quad \,\,\,\,\,
\widehat{E} = \frac{E}{1-\nu^2} \, ,
\quad \,\,\,\,\,
J_{x_3} = \frac{h^3}{12} \, ,
\quad \,\,\,\,\,
h=2a \, ,
\quad \,\,\,\,\,
\mu = \frac{E}{2(1+\nu)} \, ,
\notag
\end{align}
where $\gamma = \frac{a_1 + a_2}{2}$.
It is possible to define a dimensionless bending moment by dividing by the classical Cauchy bending moment eq.(\ref{eq:sigm_ene_dimensionless_Cau})$_1$ obtaining \cite{gauthier1975quest,lakes1995experimental}
\begin{equation}
\Omega := \frac{M_{\mbox{tot}}}{M_{\mbox{c}}} 
= 1 + 24 \frac{\gamma}{4\mu}\frac{1}{h^2} \left(1-\nu \right)
= 1 + 24 \left(\frac{\ell_b}{h}\right)^2 \, \left(1-\nu \right)
\, ,
\label{eq:Omega_dimensionless}
\end{equation}
which coincide with eq.(\ref{eq:Omega_lakes}).
\section{Equilibrium equation and boundary conditions for the full micromorphic and micro-strain model}
\label{app:full_micro_and_strain_eq_eqa}
The only critical part in this calculus is connected to the used isotropic curvature expression
\begin{equation}
\begin{split}
W_{\tiny \mbox{curv}} = 
\frac{\mu \, L_c^2}{2} \,
\displaystyle\int_{\Omega}
\displaystyle\sum_{i=1}^{3}
\left(
a_1 \, \left\lVert \mbox{dev} \, \mbox{sym} \Big( \partial_{x_i} \boldsymbol{P} \Big) \right\rVert^2_{\mathbb{R}^{3\times 3}}
+ a_2 \, \left\lVert \mbox{skew} \Big( \partial_{x_i} \boldsymbol{P} \Big) \right\rVert^2_{\mathbb{R}^{3\times 3}}
+ \frac{2}{9} \, a_3 \, \mbox{tr}^2 \Big( \partial_{x_i} \boldsymbol{P} \Big)
\right) \, \mbox{d}x
\, .
\end{split}
\end{equation}
The first variation of $W_{\tiny \mbox{curv}}$ with respect to $\boldsymbol{P}$ is
\begin{align}
\delta W_{\tiny \mbox{curv}} 
= &
\,
\mu \, L_c^2 \,
\displaystyle\int_{\Omega}
\displaystyle\sum_{i=1}^{3}
\bigg(
a_1 \, \langle \mbox{dev} \, \mbox{sym} \left( \partial_{x_i} \boldsymbol{P} \right),\mbox{dev} \, \mbox{sym} \left( \partial_{x_i} \delta \boldsymbol{P} \right) \rangle
+ a_2 \, \langle \mbox{skew} \left( \partial_{x_i} \boldsymbol{P} \right),\mbox{skew} \left( \partial_{x_i} \delta \boldsymbol{P} \right) \rangle
\notag
\\*
&
\pushright{\left.
+ \frac{2}{9} \, a_3 \, \mbox{tr} \left( \partial_{x_i} \boldsymbol{P} \right) \, \mbox{tr} \left( \partial_{x_i} \delta \boldsymbol{P} \right)
\right) \, \mbox{d}x}
\notag
\\*
= &
\,
\mu \, L_c^2 \,
\displaystyle\int_{\Omega}
\displaystyle\sum_{i=1}^{3}
\bigg(
a_1 \, \langle \mbox{dev} \, \mbox{sym} \left( \partial_{x_i} \boldsymbol{P} \right) , \partial_{x_i} \delta \boldsymbol{P}\rangle
+ a_2 \, \langle \mbox{skew} \left( \partial_{x_i} \boldsymbol{P} \right) , \partial_{x_i} \delta \boldsymbol{P} \rangle
\label{eq:equiequa_micro_strain_varia_1}
\\*
&
\pushright{\left.
+ \frac{2}{9} \, a_3 \, \langle \mbox{tr} \left( \partial_{x_i} \boldsymbol{P} \right) \boldsymbol{\mathbbm{1}}
,
\partial_{x_i} \delta \boldsymbol{P} \rangle
\right) \, \mbox{d}x \, .}
\notag
\end{align}
The product rule implies that
\begin{align}
\displaystyle\sum_{i=1}^{3}
\partial_{x_i} \, \bigg(
a_1 \, \langle \mbox{dev} \, &\mbox{sym} \left( \partial_{x_i} \boldsymbol{P} \right)
,
\delta \boldsymbol{P} \rangle_{\mathbb{R}^{3\times 3}}
+ a_2 \, \langle \mbox{skew} \left( \partial_{x_i} \boldsymbol{P} \right)
,
\delta \boldsymbol{P} \rangle_{\mathbb{R}^{3\times 3}}
\notag
\\
&
\pushright{\left.
+ \frac{2}{9} \, a_3 \, \langle \mbox{tr} \left( \partial_{x_i} \boldsymbol{P} \right) \boldsymbol{\mathbbm{1}}
,
\delta \boldsymbol{P} \rangle_{\mathbb{R}^{3\times 3}}
\right) \, \mbox{d}x}
\notag
\\
=
&
\, \displaystyle\sum_{i=1}^{3}
\bigg(
a_1 \, \langle \partial_{x_i} \, \mbox{dev} \, \mbox{sym} \left( \partial_{x_i} \boldsymbol{P} \right)
,
\delta \boldsymbol{P} \rangle_{\mathbb{R}^{3\times 3}}
+ a_2 \, \langle \partial_{x_i} \, \mbox{skew} \left( \partial_{x_i} \boldsymbol{P} \right)
,
\delta \boldsymbol{P} \rangle_{\mathbb{R}^{3\times 3}}
\label{eq:equiequa_micro_strain_varia_1b}
\\
&
\pushright{\left.
+ \frac{2}{9} \, a_3 \, \langle \partial_{x_i} \, \mbox{tr} \left( \partial_{x_i} \boldsymbol{P} \right) \boldsymbol{\mathbbm{1}}
,
\delta \boldsymbol{P} \rangle_{\mathbb{R}^{3\times 3}}
\right) \, \mbox{d}x}
\notag
\\
&
+
\displaystyle\sum_{i=1}^{3}
\bigg(
a_1 \, \langle \mbox{dev} \, \mbox{sym} \left( \partial_{x_i} \boldsymbol{P} \right) , \partial_{x_i} \delta \boldsymbol{P}\rangle
+ a_2 \, \langle \mbox{skew} \left( \partial_{x_i} \boldsymbol{P} \right) , \partial_{x_i} \delta \boldsymbol{P} \rangle
\notag
\\
&
\pushright{
\left.
+ \frac{2}{9} \, a_3 \, \langle \mbox{tr} \left( \partial_{x_i} \boldsymbol{P} \right) \boldsymbol{\mathbbm{1}}
,
\partial_{x_i} \delta \boldsymbol{P} \rangle
\right) \, \mbox{d}x \, ,
}
\notag
\end{align}
thus $\delta W_{\tiny \mbox{curv}} $ can be written as
\begin{align}
\delta W_{\tiny \mbox{curv}} 
= &
\,
\mu \, L_c^2 \,
\displaystyle\int_{\Omega}
\displaystyle\sum_{i=1}^{3}
\partial_{x_i} \, \bigg(
a_1 \, \langle \mbox{dev} \, \mbox{sym} \left( \partial_{x_i} \boldsymbol{P} \right)
,
\delta \boldsymbol{P} \rangle_{\mathbb{R}^{3\times 3}}
+ a_2 \, \langle \mbox{skew} \left( \partial_{x_i} \boldsymbol{P} \right)
,
\delta \boldsymbol{P} \rangle_{\mathbb{R}^{3\times 3}}
\notag
\\*
&
\pushright{\left.
+ \frac{2}{9} \, a_3 \, \langle \mbox{tr} \left( \partial_{x_i} \boldsymbol{P} \right) \boldsymbol{\mathbbm{1}}
,
\delta \boldsymbol{P} \rangle_{\mathbb{R}^{3\times 3}}
\right) \, \mbox{d}x}
\notag
\\*
&
-\mu \, L_c^2 \,
\displaystyle\int_{\Omega}
\displaystyle\sum_{i=1}^{3}
\bigg(
a_1 \, \langle \partial_{x_i} \, \mbox{dev} \, \mbox{sym} \left( \partial_{x_i} \boldsymbol{P} \right)
,
\delta \boldsymbol{P} \rangle_{\mathbb{R}^{3\times 3}}
+ a_2 \, \langle \partial_{x_i} \, \mbox{skew} \left( \partial_{x_i} \boldsymbol{P} \right)
,
\delta \boldsymbol{P} \rangle_{\mathbb{R}^{3\times 3}}
\label{eq:equiequa_micro_strain_varia_2}
\\*
&
\pushright{\left.
+ \frac{2}{9} \, a_3 \, \langle \partial_{x_i} \, \mbox{tr} \left( \partial_{x_i} \boldsymbol{P} \right) \boldsymbol{\mathbbm{1}}
,
\delta \boldsymbol{P} \rangle_{\mathbb{R}^{3\times 3}}
\right) \, \mbox{d}x \, ,}
\notag
\end{align}
or, defining $\boldsymbol{m}_{i} = \mu \, L_c \, \Big(a_1 \, \mbox{dev} \, \mbox{sym} \left( \partial_{x_i} \boldsymbol{P} \right)$ $+$ $a_2 \, \mbox{skew} \left( \partial_{x_i} \boldsymbol{P} \right)$ $+$ $\frac{2}{9} \, a_3 \, \mbox{tr} \left( \partial_{x_i} \boldsymbol{P} \right) \boldsymbol{\mathbbm{1}}\Big)$, $i=1,2,3$, we can equivalently write
\begin{align}
\delta W_{\tiny \mbox{curv}} 
= &
\,
\displaystyle\int_{\Omega}
\displaystyle\sum_{i=1}^{3}
\partial_{x_i} \, 
\langle
\boldsymbol{m}_{i}
,
\delta \boldsymbol{P}
\rangle_{\mathbb{R}^{3\times 3}} \, \mbox{d}x
-
\displaystyle\int_{\Omega}
\displaystyle\sum_{i=1}^{3}
\langle
\partial_{x_i} \, \boldsymbol{m}_{i}
,
\delta \boldsymbol{P}
\rangle_{\mathbb{R}^{3\times 3}} \, \mbox{d}x
\label{eq:equiequa_micro_strain_varia_2b}
\\
=
&
\displaystyle\int_{\Omega}
\mbox{Div}
\left(
\begin{array}{c}
\langle \boldsymbol{m}_{1} , \delta \boldsymbol{P} \rangle_{\mathbb{R}^{3\times 3}}
\\
\langle \boldsymbol{m}_{2} , \delta \boldsymbol{P} \rangle_{\mathbb{R}^{3\times 3}}
\\
\langle \boldsymbol{m}_{3} , \delta \boldsymbol{P} \rangle_{\mathbb{R}^{3\times 3}}
\end{array}
\right)
\, \mbox{d}x
-
\displaystyle\int_{\Omega}
\displaystyle\sum_{i=1}^{3}
\langle
\partial_{x_i} \, \boldsymbol{m}_{i}
,
\delta \boldsymbol{P}
\rangle_{\mathbb{R}^{3\times 3}} \, \mbox{d}x \, .
\notag
\end{align}
After applying the divergence theorem to the first integral of eq.(\ref{eq:equiequa_micro_strain_varia_2}) and eq.(\ref{eq:equiequa_micro_strain_varia_2b}) respectively, we obtain

\begin{align}
\delta W_{\tiny \mbox{curv}} 
= &
\,
\mu \, L_c^2 \,
\displaystyle\int_{\partial \Omega}
\displaystyle\sum_{i=1}^{3}
\langle 
\left(
a_1 \, \mbox{dev} \, \mbox{sym} \left( \partial_{x_i} \boldsymbol{P} \right)
+ a_2 \, \mbox{skew} \left( \partial_{x_i} \boldsymbol{P} \right)
+ \frac{2}{9} \, a_3 \, \mbox{tr} \left( \partial_{x_i} \boldsymbol{P} \right) \boldsymbol{\mathbbm{1}}
\right)
\, n_{i}
,
\delta \boldsymbol{P}
\rangle_{\mathbb{R}^{3\times 3}} \, \mbox{d}S
\notag
\\
&
-\mu \, L_c^2 \,
\displaystyle\int_{\Omega}
\displaystyle\sum_{i=1}^{3}
\langle
a_1 \, \mbox{dev} \, \mbox{sym} \left( \partial_{x_i}^2 \boldsymbol{P} \right)
+ a_2 \, \mbox{skew} \left( \partial_{x_i}^2 \boldsymbol{P} \right)
+ \frac{2}{9} \, a_3 \, \mbox{tr} \left( \partial_{x_i}^2 \boldsymbol{P} \right) \boldsymbol{\mathbbm{1}}
,
\delta \boldsymbol{P} \rangle_{\mathbb{R}^{3\times 3}}
\, \mbox{d}x
\notag
\\
= &
\displaystyle\int_{\partial \Omega}
\langle \displaystyle\sum_{i=1}^{3} \boldsymbol{m}_{i} \, n_{i} , \delta \boldsymbol{P} \, \rangle_{\mathbb{R}^{3\times 3}} \mbox{d}S
-\mu \, L_c^2
\displaystyle\int_{\Omega}
\langle
a_1 \, \mbox{dev} \, \mbox{sym} \, \boldsymbol{\Delta P}
+ a_2 \, \mbox{skew} \, \boldsymbol{\Delta P}
\label{eq:equiequa_micro_strain_varia_3}
\\
&
\pushright{+ \frac{2}{9} \, a_3 \, \mbox{tr} \left( \boldsymbol{\Delta P} \right) \boldsymbol{\mathbbm{1}}
\, ,
\delta \boldsymbol{P}
\rangle_{\mathbb{R}^{3\times 3}} \, \mbox{d}x
\, ,}
\notag
\end{align}
or equivalently in the form of eq.(\ref{eq:equiequa_micro_strain_varia_2b})
\begin{align}
\delta W_{\tiny \mbox{curv}} 
=
&
\displaystyle\int_{\partial \Omega}
\langle
\left(
\begin{array}{c}
\langle \boldsymbol{m}_{1} , \delta \boldsymbol{P} \rangle_{\mathbb{R}^{3\times 3}}
\\
\langle \boldsymbol{m}_{2} , \delta \boldsymbol{P} \rangle_{\mathbb{R}^{3\times 3}}
\\
\langle \boldsymbol{m}_{3} , \delta \boldsymbol{P} \rangle_{\mathbb{R}^{3\times 3}}
\end{array}
\right)
,
\boldsymbol{n}
\rangle_{\mathbb{R}^{3}}
\, \mbox{d}S
-\mu \, L_c^2 \,
\displaystyle\int_{\Omega}
\displaystyle\sum_{i=1}^{3}
\langle
a_1 \, \mbox{dev} \, \mbox{sym} \left( \boldsymbol{\Delta} \boldsymbol{P} \right)
+ a_2 \, \mbox{skew} \left( \boldsymbol{\Delta} \boldsymbol{P} \right)
\hspace{1.1cm}
\label{eq:equiequa_micro_strain_varia_2c}
\\
&
\pushright{
+ \frac{2}{9} \, a_3 \, \mbox{tr} \left( \boldsymbol{\Delta} \boldsymbol{P} \right) \boldsymbol{\mathbbm{1}}
,
\delta \boldsymbol{P} \rangle_{\mathbb{R}^{3\times 3}}
\, \mbox{d}x
\, .
}
\notag
\end{align}
Since $\delta \boldsymbol{P}$ is arbitrary in $\Omega$, and on $\partial \Omega$, and since $\delta W_{\tiny \mbox{curv}} = 0$ for a stationary point, we are now in a position to write the curvature terms of the equilibrium equation and the boundary condition without body forces and external load for the micromorphic model
\begin{empheq}[right=\empheqrbrace]{align}	
&\mbox{in} \, \Omega&
-\mu \, L_c^2
\left[
a_1 \, \mbox{dev} \, \mbox{sym} \, \underset{\in \, \mathbb{R}^{3\times 3}}{\boldsymbol{\Delta P}}
+ a_2 \, \mbox{skew} \, \underset{\in \, \mathbb{R}^{3\times 3}}{\boldsymbol{\Delta P}}
+ \frac{2}{9} \, a_3 \, \mbox{tr} \left( \underset{\in \, \mathbb{R}^{3\times 3}}{\boldsymbol{\Delta P}} \right) \boldsymbol{\mathbbm{1}}
\right] &= \boldsymbol{0}_{\mathbb{R}^{3\times3}} \, ,
\label{eq:equiequa_full_micro_varia}
\\
&\mbox{on} \, \partial \Omega&
\displaystyle\sum_{i=1}^{3} \boldsymbol{m}_{i} \, n_{i} &= \boldsymbol{0}_{\mathbb{R}^{3\times3}}
 \, ,
\notag
\end{empheq}
If $\boldsymbol{P}=\boldsymbol{S}$ is symmetric, then $\delta \boldsymbol{S} \in \mbox{Sym}(3)$ as well and, like in the micro-strain model, eq.(\ref{eq:equiequa_full_micro_varia}) turns into
\begin{empheq}[right=\empheqrbrace]{align}	
& \mbox{on} \, \Omega&
-\mu \, L_c^2
\left[
a_1 \, \mbox{dev} \, \mbox{sym} \, \underset{\in \, \mathbb{R}^{3\times 3}}{\boldsymbol{\Delta S}}
+ \frac{2}{9} \, a_3 \, \mbox{tr} \left( \underset{\in \, \mathbb{R}^{3\times 3}}{\boldsymbol{\Delta S}} \right) \boldsymbol{\mathbbm{1}}
\right] &= \boldsymbol{0}_{\mathbb{R}^{3\times3}} \, ,
\label{eq:equiequa_micro_strain_varia}
\\
& \mbox{on} \, \partial \Omega&
\mbox{sym}\left(\displaystyle\sum_{i=1}^{3} \boldsymbol{m}_{i} \, n_{i}\right) &= \boldsymbol{0}_{\mathbb{R}^{3\times3}}
 \, .
\notag
\end{empheq}
\section{Equilibrium equation for the second gradient elastic model}
\label{app:second_gradient_eq_eqa}
\begin{align}
W =
&
\displaystyle\int_{\Omega}
\mu_{\mbox{\tiny macro}} \left\lVert \mbox{sym} \, \boldsymbol{\mbox{D}u} \right\rVert^{2}_{\mathbb{R}^{3\times 3}}
+
\frac{\lambda_{\mbox{\tiny macro}}}{2} \mbox{tr}^2 \left(\boldsymbol{\mbox{D}u} \right)
\\
&
+
\frac{\mu \, L_c^2}{2} \,
\displaystyle\sum_{i=1}^{3}
\left(
a_1 \, \left\lVert \mbox{dev} \, \mbox{sym} \Big( \partial_{x_i}\boldsymbol{\mbox{D}u} \Big) \right\rVert^2_{\mathbb{R}^{3\times3}}
+ a_2 \, \left\lVert \mbox{skew} \Big( \partial_{x_i}\boldsymbol{\mbox{D}u} \Big) \right\rVert^2_{\mathbb{R}^{3\times3}}
+ \frac{2}{9} \, a_3 \, \mbox{tr}^2 \Big( \partial_{x_i}\boldsymbol{\mbox{D}u} \Big)
\right) \mbox{d}x
\, .
\notag
\end{align}
It is clear that the proposed curvature energy is isotropic \cite{munch2018rotational}.
However, the expression is not the most general isotropic curvature term, which would have 5 independent constants (see also \cite{mindlin1964micro,dell2009generalized} and Appendix~\ref{app:coeff}).
The first variation with respect to $\boldsymbol{u}$ is
\begin{align}
\delta W
= &
\displaystyle\int_{\Omega}
2\mu_{\mbox{\tiny macro}} \, \langle \mbox{sym} \, \boldsymbol{\mbox{D}u} ,  \mbox{sym} \, \boldsymbol{\mbox{D}} \delta \boldsymbol{u} \rangle_{\mathbb{R}^{3\times 3}}
+
\lambda_{\mbox{\tiny macro}} \, \mbox{tr} \left(\boldsymbol{\mbox{D}u} \right) \mbox{tr} \left(\boldsymbol{\mbox{D}} \delta \boldsymbol{u} \right)  \, \mbox{d}x
\notag
\\
&
+
\displaystyle\int_{\Omega}
\mu \, L_c^2 \,
\displaystyle\sum_{i=1}^{3}
\bigg(
a_1 \, \langle \mbox{dev} \, \mbox{sym} \left( \partial_{x_i}\boldsymbol{\mbox{D}u} \right),\mbox{dev} \, \mbox{sym} \left( \partial_{x_i} \boldsymbol{\mbox{D}} \delta \boldsymbol{u} \right) \rangle_{\mathbb{R}^{3\times 3}}
\notag
\\
&
\pushright{
+ a_2 \, \langle \mbox{skew} \left( \partial_{x_i}\boldsymbol{\mbox{D}u} \right),\mbox{skew} \left( \partial_{x_i} \boldsymbol{\mbox{D}} \delta \boldsymbol{u} \right) \rangle_{\mathbb{R}^{3\times 3}}
+ \frac{2}{9} \, a_3 \, \mbox{tr} \left( \partial_{x_i}\boldsymbol{\mbox{D}u} \right) \, \mbox{tr} \left( \partial_{x_i} \boldsymbol{\mbox{D}} \delta \boldsymbol{u} \right)
\bigg) \, \mbox{d}x
}
\notag
\\
=
&
\displaystyle\int_{\Omega}
2\mu_{\mbox{\tiny macro}} \, \langle \mbox{sym} \, \boldsymbol{\mbox{D}u} , \boldsymbol{\mbox{D}} \delta \boldsymbol{u} \rangle
+
\lambda_{\mbox{\tiny macro}} \, \langle \mbox{tr} \left(\boldsymbol{\mbox{D}u} \right) \boldsymbol{\mathbbm{1}} , \boldsymbol{\mbox{D}} \delta \boldsymbol{u} \rangle  \, \mbox{d}x
\label{eq:first_var_SG_1}
\\
&
+
\displaystyle\int_{\Omega}
\mu \, L_c^2 \,
\displaystyle\sum_{i=1}^{3}
\bigg(
a_1 \, \langle \mbox{dev} \, \mbox{sym} \left( \partial_{x_i}\boldsymbol{\mbox{D}u} \right), \partial_{x_i} \boldsymbol{\mbox{D}} \delta \boldsymbol{u} \rangle_{\mathbb{R}^{3\times 3}}
+ a_2 \, \langle \mbox{skew} \left( \partial_{x_i}\boldsymbol{\mbox{D}u} \right),\partial_{x_i} \boldsymbol{\mbox{D}} \delta \boldsymbol{u}\rangle_{\mathbb{R}^{3\times 3}}
\notag
\\
&
\hspace{8cm}
+ \frac{2}{9} \, a_3 \, \langle \mbox{tr} \left( \partial_{x_i} \boldsymbol{\mbox{D}u} \right) \boldsymbol{\mathbbm{1}}
,
\partial_{x_i} \, \boldsymbol{\mbox{D}} \delta \boldsymbol{u} \rangle_{\mathbb{R}^{3\times 3}}
\bigg) \, \mbox{d}x
\, .
\notag
\end{align}
The product rule implies that (here reported just for the curvature part of the energy as an example)
\begin{align}
\displaystyle\sum_{i=1}^{3}
&\partial_{x_i} \, \bigg(
a_1 \, \langle \mbox{dev} \, \mbox{sym} \left( \partial_{x_i} \boldsymbol{\mbox{D}u} \right)
,
\boldsymbol{\mbox{D}} \delta \boldsymbol{u} \rangle_{\mathbb{R}^{3\times 3}}
+ a_2 \, \langle \mbox{skew} \left( \partial_{x_i} \boldsymbol{\mbox{D}u} \right)
,
\boldsymbol{\mbox{D}} \delta \boldsymbol{u} \rangle_{\mathbb{R}^{3\times 3}}
\notag
\\*
&
\pushright{\left.
+ \frac{2}{9} \, a_3 \, \langle \mbox{tr} \left( \partial_{x_i} \boldsymbol{\mbox{D}u} \right) \boldsymbol{\mathbbm{1}}
,
\boldsymbol{\mbox{D}} \delta \boldsymbol{u} \rangle_{\mathbb{R}^{3\times 3}}
\right) \, \mbox{d}x}
\notag
\\*
=
&
\, \displaystyle\sum_{i=1}^{3}
\bigg(
a_1 \, \langle \partial_{x_i} \, \mbox{dev} \, \mbox{sym} \left( \partial_{x_i} \boldsymbol{\mbox{D}u} \right)
,
\boldsymbol{\mbox{D}} \delta \boldsymbol{u} \rangle_{\mathbb{R}^{3\times 3}}
+ a_2 \, \langle \partial_{x_i} \, \mbox{skew} \left( \partial_{x_i} \boldsymbol{\mbox{D}u} \right)
,
\boldsymbol{\mbox{D}} \delta \boldsymbol{u} \rangle_{\mathbb{R}^{3\times 3}}
\label{eq:first_var_SG_1b}
\\*
&
\pushright{\left.
+ \frac{2}{9} \, a_3 \, \langle \partial_{x_i} \, \mbox{tr} \left( \partial_{x_i} \boldsymbol{\mbox{D}u} \right) \boldsymbol{\mathbbm{1}}
,
\boldsymbol{\mbox{D}} \delta \boldsymbol{u} \rangle_{\mathbb{R}^{3\times 3}}
\right) \, \mbox{d}x}
\notag
\\*
&
+ \displaystyle\sum_{i=1}^{3}
\bigg(
a_1 \, \langle \mbox{dev} \, \mbox{sym} \left( \partial_{x_i} \boldsymbol{\mbox{D}u} \right) , \partial_{x_i} \boldsymbol{\mbox{D}} \delta \boldsymbol{u}\rangle_{\mathbb{R}^{3\times 3}}
+ a_2 \, \langle \mbox{skew} \left( \partial_{x_i} \boldsymbol{\mbox{D}u} \right) , \partial_{x_i} \boldsymbol{\mbox{D}} \delta \boldsymbol{u} \rangle_{\mathbb{R}^{3\times 3}}
\notag
\\*
&
\pushright{
\left.
+ \frac{2}{9} \, a_3 \, \langle \mbox{tr} \left( \partial_{x_i} \boldsymbol{\mbox{D}u} \right) \boldsymbol{\mathbbm{1}}
,
\partial_{x_i} \boldsymbol{\mbox{D}} \delta \boldsymbol{u} \rangle_{\mathbb{R}^{3\times 3}}
\right) \, \mbox{d}x \, .
}
\notag
\end{align}
Since
$
\mbox{div} \left(\boldsymbol{S}^T_1 \, \delta \boldsymbol{u}\right) =
\langle
\mbox{Div} \, \boldsymbol{S}_1
,
\delta \boldsymbol{u}
\rangle_{\mathbb{R}^3}
+
\langle
\boldsymbol{S}_1
,
\nabla \delta \boldsymbol{u}
\rangle_{\mathbb{R}^{3\times3}}
$,
we can express $\delta W$ as
\begin{align}
\delta W
= &
\displaystyle\int_{\Omega}
\mbox{div}\bigg[
\Big(
2\mu_{\mbox{\tiny macro}} \, \mbox{sym} \, \boldsymbol{\mbox{D}u} +
\lambda_{\mbox{\tiny macro}} \, \mbox{tr} \left(\boldsymbol{\mbox{D}u}\right) \boldsymbol{\mathbbm{1}}
\Big)^T
\,
\delta \boldsymbol{u}
\bigg]
\, \mbox{d}x
- \displaystyle\int_{\Omega}
\langle \mbox{Div}\left[2\mu_{\mbox{\tiny macro}} \, \mbox{sym} \, \boldsymbol{\mbox{D}u} +
\lambda_{\mbox{\tiny macro}} \, \mbox{tr} \left(\boldsymbol{\mbox{D}u}\right) \boldsymbol{\mathbbm{1}} \right] , \delta \boldsymbol{u} \rangle_{\mathbb{R}^{3}}
\, \mbox{d}x
\notag
\\*
&
+
\displaystyle\int_{\Omega}
\mu \, L_c^2 \,
\displaystyle\sum_{i=1}^{3}
\partial_{x_i} \bigg(
\langle a_1 \, \mbox{dev} \, \mbox{sym} \left( \partial_{x_i}\boldsymbol{\mbox{D}u} \right)
+ a_2 \, \mbox{skew} \left( \partial_{x_i}\boldsymbol{\mbox{D}u} \right)
+ \frac{2}{9} \, a_3 \, \mbox{tr} \left( \partial_{x_i} \boldsymbol{\mbox{D}u} \right) \boldsymbol{\mathbbm{1}}
,
\boldsymbol{\mbox{D}} \delta \boldsymbol{u} \rangle_{\mathbb{R}^{3\times 3}}
\bigg) \, \mbox{d}x
\label{eq:first_var_SG_2}
\\*
&
-
\displaystyle\int_{\Omega}
\mu \, L_c^2 \,
\displaystyle\sum_{i=1}^{3}
\bigg(
\langle a_1 \, \mbox{dev} \, \mbox{sym} \left( \partial_{x_i}^2\boldsymbol{\mbox{D}u} \right)
+ a_2 \, \mbox{skew} \left( \partial_{x_i}^2\boldsymbol{\mbox{D}u} \right)
+ \frac{2}{9} \, a_3 \, \mbox{tr} \left( \partial_{x_i}^2 \boldsymbol{\mbox{D}u} \right) \boldsymbol{\mathbbm{1}}
,
\boldsymbol{\mbox{D}} \delta \boldsymbol{u} \rangle_{\mathbb{R}^{3\times 3}}
\bigg) \, \mbox{d}x
\, ,
\notag
\end{align}
or equivalently
\begin{align}
\delta W
= &
\displaystyle\int_{\Omega}
\mbox{div}\bigg[
\Big(
2\mu_{\mbox{\tiny macro}} \, \mbox{sym} \, \boldsymbol{\mbox{D}u} +
\lambda_{\mbox{\tiny macro}} \, \mbox{tr} \left(\boldsymbol{\mbox{D}u}\right) \boldsymbol{\mathbbm{1}}
\Big)
\,
\delta \boldsymbol{u}
\bigg]
\, \mbox{d}x
- \displaystyle\int_{\Omega}
\langle \mbox{Div}\left[2\mu_{\mbox{\tiny macro}} \, \mbox{sym} \, \boldsymbol{\mbox{D}u} +
\lambda_{\mbox{\tiny macro}} \, \mbox{tr} \left(\boldsymbol{\mbox{D}u}\right) \boldsymbol{\mathbbm{1}} \right] , \delta \boldsymbol{u} \rangle_{\mathbb{R}^{3}}
\, \mbox{d}x
\notag
\\
&
+
\displaystyle\int_{\Omega}
\mbox{div}
\left(
\begin{array}{c}
\langle \widetilde{\boldsymbol{m}}_1, \boldsymbol{\mbox{D}} \delta \boldsymbol{u} \rangle_{\mathbb{R}^{3\times 3}} \\
\langle \widetilde{\boldsymbol{m}}_2, \boldsymbol{\mbox{D}} \delta \boldsymbol{u} \rangle_{\mathbb{R}^{3\times 3}} \\
\langle \widetilde{\boldsymbol{m}}_3, \boldsymbol{\mbox{D}} \delta \boldsymbol{u} \rangle_{\mathbb{R}^{3\times 3}}
\end{array}
\right)
\, \mbox{d}x
\label{eq:first_var_SG_2b}
\\
&
-
\displaystyle\int_{\Omega}
\mu \, L_c^2 \,
\displaystyle\sum_{i=1}^{3}
\bigg(
\langle a_1 \, \mbox{dev} \, \mbox{sym} \left( \partial_{x_i}^2\boldsymbol{\mbox{D}u} \right)
+ a_2 \, \mbox{skew} \left( \partial_{x_i}^2\boldsymbol{\mbox{D}u} \right)
+ \frac{2}{9} \, a_3 \, \mbox{tr} \left( \partial_{x_i}^2 \boldsymbol{\mbox{D}u} \right) \boldsymbol{\mathbbm{1}}
,
\boldsymbol{\mbox{D}} \delta \boldsymbol{u} \rangle_{\mathbb{R}^{3\times 3}}
\bigg) \, \mbox{d}x
\, ,
\notag
\end{align}
where $m_{i} = \mu \, L_c^2 \,
\sum_{i=1}^{3}
\bigg(
\langle a_1 \, \mbox{dev} \, \mbox{sym} \left( \partial_{x_i}\boldsymbol{\mbox{D}u} \right)
+ a_2 \, \mbox{skew} \left( \partial_{x_i}\boldsymbol{\mbox{D}u} \right)
+ \frac{2}{9} \, a_3 \, \mbox{tr} \left( \partial_{x_i} \boldsymbol{\mbox{D}u} \right) \boldsymbol{\mathbbm{1}}
,
\boldsymbol{\mbox{D}} \delta \boldsymbol{u} \rangle_{\mathbb{R}^{3\times 3}}
\bigg)$.
After applying the divergence theorem to the first and the third term of eq.(\ref{eq:first_var_SG_2}) we obtain
\begin{align}
\delta W
=
&
\displaystyle\int_{\partial \Omega}
\langle
\Big(
2\mu_{\mbox{\tiny macro}} \, \mbox{sym} \, \boldsymbol{\mbox{D}u}
+ \lambda_{\mbox{\tiny macro}} \, \mbox{tr} \left(\boldsymbol{\mbox{D}u}\right) \boldsymbol{\mathbbm{1}}
\Big)^T
\delta \boldsymbol{u}
,
\boldsymbol{n}
\rangle _{\mathbb{R}^{3}}
\, \mbox{d}S
- \displaystyle\int_{\Omega}
\langle \mbox{Div}\left[2\mu_{\mbox{\tiny macro}} \, \mbox{sym} \, \boldsymbol{\mbox{D}u}
+ \lambda_{\mbox{\tiny macro}} \, \mbox{tr} \left(\boldsymbol{\mbox{D}u}\right) \boldsymbol{\mathbbm{1}}\right] , \delta \boldsymbol{u} \rangle_{\mathbb{R}^{3}}
\, \mbox{d}x
\notag
\\
&
+ \displaystyle\int_{\partial \Omega}
\mu \, L_c^2 \,
\displaystyle\sum_{i=1}^{3}
\langle 
\Big(
a_1 \, \mbox{dev} \, \mbox{sym} \left( \partial_{x_i}\boldsymbol{\mbox{D}u} \right)
+ a_2 \, \mbox{skew} \left( \partial_{x_i}\boldsymbol{\mbox{D}u} \right)
+ \frac{2}{9} \, a_3 \, \mbox{tr} \left( \partial_{x_i}\boldsymbol{\mbox{D}u} \right) \boldsymbol{\mathbbm{1}}
\Big)
,
\boldsymbol{\mbox{D}} \, \delta \boldsymbol{u}
\rangle_{\mathbb{R}^{3\times 3}}
\, n_{i}
\, \mbox{d}S
\label{eq:first_var_SG_3}
\\
&
- \displaystyle\int_{\Omega}
\mu \, L_c^2 \,
\displaystyle\sum_{i=1}^{3}
\langle
a_1 \, \mbox{dev} \, \mbox{sym} \left( \partial_{x_i}^2\boldsymbol{\mbox{D}u} \right)
+ a_2 \, \mbox{skew} \left( \partial_{x_i}^2\boldsymbol{\mbox{D}u} \right)
+ \frac{2}{9} \, a_3 \, \mbox{tr} \left( \partial_{x_i}^2\boldsymbol{\mbox{D}u} \right) \boldsymbol{\mathbbm{1}}, \boldsymbol{\mbox{D}} \delta \boldsymbol{u}
\rangle_{\mathbb{R}^{3\times 3}}
\, \mbox{d}x
\, .
\notag
\end{align}
or equivalently
\begin{align}
\delta W
=
&
\displaystyle\int_{\partial \Omega}
\langle
\Big(
2\mu_{\mbox{\tiny macro}} \, \mbox{sym} \, \boldsymbol{\mbox{D}u}
+ \lambda_{\mbox{\tiny macro}} \, \mbox{tr} \left(\boldsymbol{\mbox{D}u}\right) \boldsymbol{\mathbbm{1}}
\Big)
\delta \boldsymbol{u}
,
\boldsymbol{n}
\rangle _{\mathbb{R}^{3}}
\, \mbox{d}S
- \displaystyle\int_{\Omega}
\langle \mbox{Div}\left[2\mu_{\mbox{\tiny macro}} \, \mbox{sym} \, \boldsymbol{\mbox{D}u}
+ \lambda_{\mbox{\tiny macro}} \, \mbox{tr} \left(\boldsymbol{\mbox{D}u}\right) \boldsymbol{\mathbbm{1}}\right] , \delta \boldsymbol{u} \rangle_{\mathbb{R}^{3}}
\, \mbox{d}x
\notag
\\
&
+ \displaystyle\int_{\partial \Omega}
\left(
\langle
\begin{array}{c}
\langle \widetilde{\boldsymbol{m}}_1, \boldsymbol{\mbox{D}} \delta \boldsymbol{u} \rangle_{\mathbb{R}^{3\times 3}} \\
\langle \widetilde{\boldsymbol{m}}_1, \boldsymbol{\mbox{D}} \delta \boldsymbol{u} \rangle_{\mathbb{R}^{3\times 3}} \\
\langle \widetilde{\boldsymbol{m}}_1, \boldsymbol{\mbox{D}} \delta \boldsymbol{u} \rangle_{\mathbb{R}^{3\times 3}}
\end{array}
\right)
,
\boldsymbol{n}
\rangle_{\mathbb{R}^{3}}
\, \mbox{d}S
\label{eq:first_var_SG_3b}
\\
&
- \displaystyle\int_{\Omega}
\mu \, L_c^2 \,
\displaystyle\sum_{i=1}^{3}
\langle
a_1 \, \mbox{dev} \, \mbox{sym} \left( \partial_{x_i}^2\boldsymbol{\mbox{D}u} \right)
+ a_2 \, \mbox{skew} \left( \partial_{x_i}^2\boldsymbol{\mbox{D}u} \right)
+ \frac{2}{9} \, a_3 \, \mbox{tr} \left( \partial_{x_i}^2\boldsymbol{\mbox{D}u} \right) \boldsymbol{\mathbbm{1}}, \boldsymbol{\mbox{D}} \delta \boldsymbol{u}
\rangle_{\mathbb{R}^{3\times 3}}
\, \mbox{d}x
\, .
\notag
\end{align}
Integrating by part the fourth term of eq.(\ref{eq:first_var_SG_3}) it is possible to write
\begin{align}
&
- \displaystyle\int_{\Omega}
\mu \, L_c^2 \,
\displaystyle\sum_{i=1}^{3}
\langle
a_1 \, \mbox{dev} \, \mbox{sym} \left( \partial_{x_i}^2\boldsymbol{\mbox{D}u} \right)
+ a_2 \, \mbox{skew} \left( \partial_{x_i}^2\boldsymbol{\mbox{D}u} \right)
+ \frac{2}{9} \, a_3 \, \mbox{tr} \left( \partial_{x_i}^2\boldsymbol{\mbox{D}u} \right) \boldsymbol{\mathbbm{1}}, \boldsymbol{\mbox{D}} \delta \boldsymbol{u}
\rangle_{\mathbb{R}^{3\times 3}}
\, \mbox{d}x
\notag
\\*
=
&
- \displaystyle\int_{\Omega}
\mbox{div}
\left[
\Big(
\mu \, L_c^2 \,
\displaystyle\sum_{i=1}^{3}
\langle
a_1 \, \mbox{dev} \, \mbox{sym} \left( \partial_{x_i}^2\boldsymbol{\mbox{D}u} \right)
+ a_2 \, \mbox{skew} \left( \partial_{x_i}^2\boldsymbol{\mbox{D}u} \right)
+ \frac{2}{9} \, a_3 \, \mbox{tr} \left( \partial_{x_i}^2\boldsymbol{\mbox{D}u} \right) \boldsymbol{\mathbbm{1}}
\Big)^T
\delta \boldsymbol{u}
\right]
\, \mbox{d}x
\label{eq:first_var_SG_4}
\\*
&
+ \displaystyle\int_{\Omega}
\langle
\mbox{Div}
\left[
\mu \, L_c^2 \,
\displaystyle\sum_{i=1}^{3}
\left(
a_1 \, \mbox{dev} \, \mbox{sym} \left( \partial_{x_i}^2\boldsymbol{\mbox{D}u} \right)
+ a_2 \, \mbox{skew} \left( \partial_{x_i}^2\boldsymbol{\mbox{D}u} \right)
+ \frac{2}{9} \, a_3 \, \mbox{tr} \left( \partial_{x_i}^2\boldsymbol{\mbox{D}u} \right) \boldsymbol{\mathbbm{1}}
\right)
\right]
,
\delta \boldsymbol{u}
\rangle_{\mathbb{R}^{3}}
\, \mbox{d}x
\notag
\end{align}
and using the divergence theorem on the first term of eq.(\ref{eq:first_var_SG_4}) we have
\begin{align}
&
- \displaystyle\int_{\Omega}
\mbox{div}
\left[
\Big(
\mu \, L_c^2 \,
\displaystyle\sum_{i=1}^{3}
a_1 \, \mbox{dev} \, \mbox{sym} \left( \partial_{x_i}^2\boldsymbol{\mbox{D}u} \right)
+ a_2 \, \mbox{skew} \left( \partial_{x_i}^2\boldsymbol{\mbox{D}u} \right)
+ \frac{2}{9} \, a_3 \, \mbox{tr} \left( \partial_{x_i}^2\boldsymbol{\mbox{D}u} \right) \boldsymbol{\mathbbm{1}}
\Big)^T
\delta \boldsymbol{u}
\right]
\, \mbox{d}x
\label{eq:first_var_SG_5}
\\
=
&
- \displaystyle\int_{\partial \Omega}
\mu \, L_c^2 \,
\displaystyle\sum_{i=1}^{3}
\langle
\left(
a_1 \, \mbox{dev} \, \mbox{sym} \left( \partial_{x_i}^2\boldsymbol{\mbox{D}u} \right)
+ a_2 \, \mbox{skew} \left( \partial_{x_i}^2\boldsymbol{\mbox{D}u} \right)
+ \frac{2}{9} \, a_3 \, \mbox{tr} \left( \partial_{x_i}^2\boldsymbol{\mbox{D}u} \right) \boldsymbol{\mathbbm{1}}
\right)
\,
\boldsymbol{n}
,
\delta \boldsymbol{u}
\rangle_{\mathbb{R}^{3}}
\, \mbox{d}S
\notag
\end{align}
Now, the only term that needs further attention is the third one in eq.(\ref{eq:first_var_SG_3}), since only the normal component to the surface of $\boldsymbol{\mbox{D}} \delta \boldsymbol{u}$ is independent with respect to $\delta \boldsymbol{u}$.
We introduce these three operators 
\begin{equation}
\boldsymbol{D^{\tau}} = \left( \boldsymbol{\mathbbm{1}} - \boldsymbol{n} \otimes \boldsymbol{n}\right) \cdot (\partial_{x_1},\partial_{x_2},\partial_{x_3}) \, ,
\qquad
\boldsymbol{D^n_a} = \left(\boldsymbol{n} \otimes \boldsymbol{n}\right) \cdot (\partial_{x_1},\partial_{x_2},\partial_{x_3}) \, ,
\qquad
D^n_b = (\partial_{x_1},\partial_{x_2},\partial_{x_3}) \cdot \boldsymbol{n} \, ,
\label{eq:deco_grad_SG}
\end{equation}
which, using the fact that $\left(\boldsymbol{\mbox{D} \delta u}\right)^T = (\partial_{x_1},\partial_{x_2},\partial_{x_3}) \otimes \boldsymbol{\delta u}$ and that $\boldsymbol{D^{\tau}} + \boldsymbol{D^n_a} = (\partial_{x_1},\partial_{x_2},\partial_{x_3})$, allows us to write the third term in eq.(\ref{eq:first_var_SG_3}) as
\begin{align}
&
\displaystyle\int_{\partial \Omega}
\mu \, L_c^2 \,
\displaystyle\sum_{i=1}^{3}
\langle 
\Big(a_1 \, \mbox{dev} \, \mbox{sym} \left( \partial_{x_i}\boldsymbol{\mbox{D}u} \right)
+ a_2 \, \mbox{skew} \left( \partial_{x_i}\boldsymbol{\mbox{D}u} \right)
+ \frac{2}{9} \, a_3 \, \mbox{tr} \left( \partial_{x_i}\boldsymbol{\mbox{D}u} \right) \boldsymbol{\mathbbm{1}}\Big)  n_{i}, \boldsymbol{\mbox{D}} \, \delta \boldsymbol{u} \rangle_{\mathbb{R}^{3\times 3}}
\, \mbox{d}S
\notag
\\
\hspace{-2.5cm}
=
&
\displaystyle\int_{\partial \Omega}
\mu \, L_c^2 \,
\displaystyle\sum_{i=1}^{3}
\langle 
\Big(a_1 \, \mbox{dev} \, \mbox{sym} \left( \partial_{x_i}\boldsymbol{\mbox{D}u} \right)
+ a_2 \, \mbox{skew} \left( \partial_{x_i}\boldsymbol{\mbox{D}u} \right)
+ \frac{2}{9} \, a_3 \, \mbox{tr} \left( \partial_{x_i}\boldsymbol{\mbox{D}u} \right) \boldsymbol{\mathbbm{1}}\Big)  n_{i}
,
\delta \boldsymbol{u} \otimes \boldsymbol{D^{\tau}} \rangle_{\mathbb{R}^{3\times 3}}
\, \mbox{d}S
\notag
\\
&
+
\displaystyle\int_{\partial \Omega}
\mu \, L_c^2 \,
\displaystyle\sum_{i=1}^{3}
\langle 
\Big(a_1 \, \mbox{dev} \, \mbox{sym} \left( \partial_{x_i}\boldsymbol{\mbox{D}u} \right)
+ a_2 \, \mbox{skew} \left( \partial_{x_i}\boldsymbol{\mbox{D}u} \right)
+ \frac{2}{9} \, a_3 \, \mbox{tr} \left( \partial_{x_i}\boldsymbol{\mbox{D}u} \right) \boldsymbol{\mathbbm{1}}\Big)  n_{i}
,
\delta \boldsymbol{u} \otimes \boldsymbol{D^n_a} \rangle_{\mathbb{R}^{3}}
\, \mbox{d}S
\notag
\\
\hspace{-2.5cm}
=
&
\displaystyle\int_{\partial \Omega}
\overbrace{\mu \, L_c^2 \,
\displaystyle\sum_{i=1}^{3}
\langle 
\Big(a_1 \, \mbox{dev} \, \mbox{sym} \left( \partial_{x_i}\boldsymbol{\mbox{D}u} \right)
+ a_2 \, \mbox{skew} \left( \partial_{x_i}\boldsymbol{\mbox{D}u} \right)
+ \frac{2}{9} \, a_3 \, \mbox{tr} \left( \partial_{x_i}\boldsymbol{\mbox{D}u} \right) \boldsymbol{\mathbbm{1}}\Big)  n_{i}}^{\textstyle B_{ij} \mathstrut}
,
\delta \boldsymbol{u} \otimes \boldsymbol{D^{\tau}} \rangle_{\mathbb{R}^{3\times 3}}
\, \mbox{d}S
\label{eq:first_var_SG_6}
\\
&
+
\displaystyle\int_{\partial \Omega}
\mu \, L_c^2 \,
\displaystyle\sum_{i=1}^{3}
\langle 
\left(
\Big(a_1 \, \mbox{dev} \, \mbox{sym} \left( \partial_{x_i}\boldsymbol{\mbox{D}u} \right)
+ a_2 \, \mbox{skew} \left( \partial_{x_i}\boldsymbol{\mbox{D}u} \right)
+ \frac{2}{9} \, a_3 \, \mbox{tr} \left( \partial_{x_i}\boldsymbol{\mbox{D}u} \right) \boldsymbol{\mathbbm{1}}\Big)  n_{i}
\right)
\boldsymbol{n}
,
D^n_b \delta \boldsymbol{u} \rangle_{\mathbb{R}^{3}}
\, \mbox{d}S
\, .
\notag
\end{align}
While the term involving the normal derivative of the virtual displacement ($D^n_b \boldsymbol{\partial u} = \boldsymbol{\mbox{D}u} \cdot \boldsymbol{n}$) is independent of $\boldsymbol{\delta u}$, the term involving the tangential projection of $\boldsymbol{\mbox{D} \delta u}$ ($\boldsymbol{\delta u} \otimes \boldsymbol{D^{\tau}} = \boldsymbol{\mbox{D} \delta u} \cdot \left( \boldsymbol{\mathbbm{1}} - \boldsymbol{n} \otimes \boldsymbol{n}\right) $) is not, and must be further manipulated.
It is known (\cite{madeo2016new,mindlin1964micro}) that this tangential term can be still manipulated integrating by parts and using the surface divergence theorem so implying (see eq.(3.5) in \cite{madeo2016new})
\begin{equation}
	\int\limits_{\partial \Omega} \langle \boldsymbol{B}, \boldsymbol{\delta u} \otimes \boldsymbol{D^{\tau}} \rangle =
	- \int\limits_{\partial \Omega} \langle \boldsymbol{\mbox{D}} \left[\boldsymbol{B} \cdot \left( \boldsymbol{\mathbbm{1}} - \boldsymbol{n} \otimes \boldsymbol{n}\right)\right] \cdot \left( \boldsymbol{\mathbbm{1}} - \boldsymbol{n} \otimes \boldsymbol{n}\right)
	, \boldsymbol{\delta u}
	 \rangle
	+
	\int\limits_{\partial\partial \Omega}
	\llbracket
	\langle
	\boldsymbol{B} \nu,\boldsymbol{\delta u}
	\rangle
	\rrbracket
	\, ,
\end{equation}
where $\nu$ is the normal to $\partial\partial \Omega$.
In our bending problem there is no $\partial\partial \Omega$ and the normal is constant, so that the preceding equation reduce to

\begin{align}
	\int\limits_{\partial \Omega} \langle \boldsymbol{B}, \boldsymbol{\delta u} \otimes \boldsymbol{D^{\tau}} \rangle &=
	- \int\limits_{\partial \Omega}
	\langle
	\boldsymbol{\mbox{D}} \boldsymbol{B} \cdot \left( \boldsymbol{\mathbbm{1}} - \boldsymbol{n} \otimes \boldsymbol{n}\right) \cdot \left( \boldsymbol{\mathbbm{1}} - \boldsymbol{n} \otimes \boldsymbol{n}\right)
	,
	\boldsymbol{\delta u}
	\rangle
	\\
	&
	=
	- \int\limits_{\partial \Omega}
	\langle
	\boldsymbol{\mbox{D}} \boldsymbol{B} \cdot \left( \boldsymbol{\mathbbm{1}} - \boldsymbol{n} \otimes \boldsymbol{n}\right)
	,
	\boldsymbol{\delta u}
	\rangle
	=
	- \int\limits_{\partial \Omega}
	\langle
	\boldsymbol{B} \cdot \boldsymbol{D^{\tau}}
	,
	\boldsymbol{\delta u}
	\rangle
	\, .
	\notag
\end{align}
Hence, in our particular case, the first term in eq.(\ref{eq:first_var_SG_6}) can be written as:
\begin{align}
&\displaystyle\int_{\partial \Omega}
\mu \, L_c^2 \,
\displaystyle\sum_{i=1}^{3}
\langle 
\Big(a_1 \, \mbox{dev} \, \mbox{sym} \left( \partial_{x_i}\boldsymbol{\mbox{D}u} \right)
+ a_2 \, \mbox{skew} \left( \partial_{x_i}\boldsymbol{\mbox{D}u} \right)
+ \frac{2}{9} \, a_3 \, \mbox{tr} \left( \partial_{x_i}\boldsymbol{\mbox{D}u} \right) \boldsymbol{\mathbbm{1}}\Big)  n_{i}
,
\delta \boldsymbol{u} \otimes \boldsymbol{D^{\tau}} \rangle_{\mathbb{R}^{3\times 3}}
\, \mbox{d}S
\label{eq:first_var_SG_7}
\\
=
&
- \displaystyle\int_{\partial \Omega}
\mu \, L_c^2 \,
\displaystyle\sum_{i=1}^{3}
\langle 
\Big[
\Big(
a_1 \, \mbox{dev} \, \mbox{sym} \left( \partial_{x_i}\boldsymbol{\mbox{D}u} \right)
+ a_2 \, \mbox{skew} \left( \partial_{x_i}\boldsymbol{\mbox{D}u} \right)
+ \frac{2}{9} \, a_3 \, \mbox{tr} \left( \partial_{x_i}\boldsymbol{\mbox{D}u} \right) \boldsymbol{\mathbbm{1}}
\Big)  n_{i}
\Big]
\cdot \boldsymbol{D^{\tau}}
,
\delta \boldsymbol{u} \rangle_{\mathbb{R}^{3}}
\, \mbox{d}S \, .
\notag
\end{align}
It is now possible to write all together
\begin{align}
\delta W
=
&
- \displaystyle\int_{\Omega}
\langle
\mbox{Div}
\left[
2\mu_{\mbox{\tiny macro}} \, \mbox{sym} \, \boldsymbol{\mbox{D}u}
+ \lambda_{\mbox{\tiny macro}} \, \mbox{tr} \left(\boldsymbol{\mbox{D}u}\right) \boldsymbol{\mathbbm{1}}
\right]
, \delta \boldsymbol{u}
\rangle_{\mathbb{R}^{3}}
\, \mbox{d}x
\notag
\\*
&
+ \displaystyle\int_{\Omega}
\langle
\mbox{Div}
\left[
\mu \, L_c^2 \,
\displaystyle\sum_{i=1}^{3}
\left(
a_1 \, \mbox{dev} \, \mbox{sym} \left( \partial_{x_i}^2\boldsymbol{\mbox{D}u} \right)
+ a_2 \, \mbox{skew} \left( \partial_{x_i}^2\boldsymbol{\mbox{D}u} \right)
+ \frac{2}{9} \, a_3 \, \mbox{tr} \left( \partial_{x_i}^2\boldsymbol{\mbox{D}u} \right) \boldsymbol{\mathbbm{1}}
\right)
\right]
,
\delta \boldsymbol{u}
\rangle_{\mathbb{R}^{3}}
\mbox{d}x
\notag
\\*
&
+\displaystyle\int_{\partial \Omega}
\langle
\Big(2\mu_{\mbox{\tiny macro}} \, \mbox{sym} \, \boldsymbol{\mbox{D}u}
+ \lambda_{\mbox{\tiny macro}} \, \mbox{tr} \left(\boldsymbol{\mbox{D}u}\right) \boldsymbol{\mathbbm{1}}\Big) \, \boldsymbol{n}, \delta \boldsymbol{u}
\rangle_{\mathbb{R}^{3}}
\, \mbox{d}S
\label{eq:first_var_SG_8}
\\*
&
- \displaystyle\int_{\partial \Omega}
\mu \, L_c^2 \,
\displaystyle\sum_{i=1}^{3}
\langle
\left(
a_1 \, \mbox{dev} \, \mbox{sym} \left( \partial_{x_i}^2\boldsymbol{\mbox{D}u} \right)
+ a_2 \, \mbox{skew} \left( \partial_{x_i}^2\boldsymbol{\mbox{D}u} \right)
+ \frac{2}{9} \, a_3 \, \mbox{tr} \left( \partial_{x_i}^2\boldsymbol{\mbox{D}u} \right) \boldsymbol{\mathbbm{1}}
\right)
\,
\boldsymbol{n}
,
\delta \boldsymbol{u}
\rangle_{\mathbb{R}^{3}}
\, \mbox{d}S
\notag
\\*
&
- \displaystyle\int_{\partial \Omega}
\mu \, L_c^2 \,
\displaystyle\sum_{i=1}^{3}
\langle 
\Big[
\Big(
a_1 \, \mbox{dev} \, \mbox{sym} \left( \partial_{x_i}\boldsymbol{\mbox{D}u} \right)
+ a_2 \, \mbox{skew} \left( \partial_{x_i}\boldsymbol{\mbox{D}u} \right)
+ \frac{2}{9} \, a_3 \, \mbox{tr} \left( \partial_{x_i}\boldsymbol{\mbox{D}u} \right) \boldsymbol{\mathbbm{1}}
\Big)  n_{i}
\Big]
\cdot \boldsymbol{D^{\tau}}
,
\delta \boldsymbol{u} \rangle_{\mathbb{R}^{3}}
\, \mbox{d}S
\notag
\\*
&
+
\displaystyle\int_{\partial \Omega}
\mu \, L_c^2 \,
\displaystyle\sum_{i=1}^{3}
\langle 
\left(
\Big(a_1 \, \mbox{dev} \, \mbox{sym} \left( \partial_{x_i}\boldsymbol{\mbox{D}u} \right)
+ a_2 \, \mbox{skew} \left( \partial_{x_i}\boldsymbol{\mbox{D}u} \right)
+ \frac{2}{9} \, a_3 \, \mbox{tr} \left( \partial_{x_i}\boldsymbol{\mbox{D}u} \right) \boldsymbol{\mathbbm{1}}\Big)  n_{i}
\right)
\boldsymbol{n}
,
D^n \delta \boldsymbol{u} \rangle_{\mathbb{R}^{3}}
\, \mbox{d}S
\, ,
\notag
\end{align}
or in an equivalent way
\begin{align}
\delta W_{\tiny \mbox{curv}}
=
&
- \displaystyle\int_{\Omega}
\langle
\mbox{Div}
\left(
\widetilde{\boldsymbol{\sigma}}
\right)
,
\delta \boldsymbol{u}
\rangle_{\mathbb{R}^{3}}
\, \mbox{d}x
+ \displaystyle\int_{\Omega}
\langle
\mbox{Div}
\left[
\mbox{Div}
\left(
\begin{array}{c}
\boldsymbol{m}_1
\\
\boldsymbol{m}_2
\\
\boldsymbol{m}_3
\end{array}
\right)
\right]
,
\delta \boldsymbol{u}
\rangle_{\mathbb{R}^{3}}
\label{eq:first_var_SG_9}
\\*
&
+\displaystyle\int_{\partial \Omega}
\langle \widetilde{\boldsymbol{\sigma}} \, \boldsymbol{n}, \delta \boldsymbol{u} \rangle_{\mathbb{R}^{3}} \, \mbox{d}x
- \displaystyle\int_{\partial \Omega}
\langle 
\mbox{Div}
\left[
\left(
\begin{array}{c}
\boldsymbol{m}_1
\\
\boldsymbol{m}_2
\\
\boldsymbol{m}_3
\end{array}
\right)
\right]
\,
\boldsymbol{n}
,
\delta \boldsymbol{u}
\rangle_{\mathbb{R}^{3}}
\, \mbox{d}x
- \displaystyle\int_{\partial \Omega}
\langle 
\mbox{Div}
\left[
\left(
\begin{array}{c}
\boldsymbol{m}_1 \, n_{1}
\\
\boldsymbol{m}_2 \, n_{2}
\\
\boldsymbol{m}_3 \, n_{3}
\end{array}
\right)
\right]
\cdot \boldsymbol{D^{\tau}}
,
\delta \boldsymbol{u} \rangle_{\mathbb{R}^{3}} \, \mbox{d}x
\notag
\\*
&
+
\displaystyle\int_{\partial \Omega}
\langle 
\left(
\begin{array}{c}
\boldsymbol{m}_1
\\
\boldsymbol{m}_2
\\
\boldsymbol{m}_3
\end{array}
\right) \cdot
\boldsymbol{n} \cdot 
\boldsymbol{n}
,
D^n \delta \boldsymbol{u} \rangle_{\mathbb{R}^{3}} \, \mbox{d}x
\, .
\end{align}
where 
$ \boldsymbol{m}_{i} = 
\mu \, L_c^2 \, \sum_{i=1}^{3} 
\Big(
a_1 \, \mbox{dev} \, \mbox{sym} \left( \partial_{x_i}\boldsymbol{\mbox{D}u} \right)
$ $+$ $ a_2 \, \mbox{skew} \left( \partial_{x_i}\boldsymbol{\mbox{D}u} \right)
$ $+$ $ \frac{2}{9} \, a_3 \, \mbox{tr} \left( \partial_{x_i}\boldsymbol{\mbox{D}u} \right) \boldsymbol{\mathbbm{1}}
\Big)$
and $\widetilde{\boldsymbol{\sigma}} =
2\mu_{\mbox{\tiny macro}} \, \mbox{sym} \, \boldsymbol{\mbox{D}u}
$ $+$ $ \lambda_{\mbox{\tiny macro}} \, \mbox{tr} \left(\boldsymbol{\mbox{D}u}\right) \boldsymbol{\mathbbm{1}}$.

We are now in a position to write the equilibrium equation and the boundary conditions without body forces and external load for the strain gradient model for the cylindrical plate bending since the only variations that remain are $\delta \boldsymbol{u}$ and $D^n_b \delta \boldsymbol{u}$ which are independent with respect each other:
\begin{align}
\mbox{Div} \left[\widetilde{\sigma} - \mu \, L_c^2 \,
\displaystyle\sum_{i=1}^{3}
\left(
a_1 \, \mbox{dev} \, \mbox{sym} \left( \partial_{x_i}^2\boldsymbol{\mbox{D}u} \right)
+ a_2 \, \mbox{skew} \left( \partial_{x_i}^2\boldsymbol{\mbox{D}u} \right)
+ \frac{2}{9} \, a_3 \, \mbox{tr} \left( \partial_{x_i}^2\boldsymbol{\mbox{D}u} \right) \boldsymbol{\mathbbm{1}}
\right)  \right] = \boldsymbol{0}_{\mathbb{R}^{3}}
& \qquad \mbox{in} \, \Omega
\notag
\\
\Big(2\mu_{\mbox{\tiny macro}} \, \mbox{sym} \, \boldsymbol{\mbox{D}u}
+ \lambda_{\mbox{\tiny macro}} \, \mbox{tr} \left(\boldsymbol{\mbox{D}u}\right) \boldsymbol{\mathbbm{1}}\Big) \, \boldsymbol{n} 
\hspace{1.1cm}
\notag
\\
-
\mu \, L_c^2 \,
\displaystyle\sum_{i=1}^{3}
\left(
a_1 \, \mbox{dev} \, \mbox{sym} \left( \partial_{x_i}^2\boldsymbol{\mbox{D}u} \right)
+ a_2 \, \mbox{skew} \left( \partial_{x_i}^2\boldsymbol{\mbox{D}u} \right)
+ \frac{2}{9} \, a_3 \, \mbox{tr} \left( \partial_{x_i}^2\boldsymbol{\mbox{D}u} \right) \boldsymbol{\mathbbm{1}}
\right)
\,
\boldsymbol{n} 
\hspace{1.1cm}
\label{eq:first_var_SG_Equi_equa}
\\
-
\mu \, L_c^2 \,
\displaystyle\sum_{i=1}^{3}
\Big[
\Big(
a_1 \, \mbox{dev} \, \mbox{sym} \left( \partial_{x_i}\boldsymbol{\mbox{D}u} \right)
+ a_2 \, \mbox{skew} \left( \partial_{x_i}\boldsymbol{\mbox{D}u} \right)
+ \frac{2}{9} \, a_3 \, \mbox{tr} \left( \partial_{x_i}\boldsymbol{\mbox{D}u} \right) \boldsymbol{\mathbbm{1}}
\Big)  n_{i}
\Big]
\cdot \boldsymbol{D^{\tau}}
= \boldsymbol{0}_{\mathbb{R}^{3}}
& \qquad \mbox{on} \, \partial \Omega
\notag
\\
\left(
\mu \, L_c^2 \,
\displaystyle\sum_{i=1}^{3}
\Big(a_1 \, \mbox{dev} \, \mbox{sym} \left( \partial_{x_i}\boldsymbol{\mbox{D}u} \right)
+ a_2 \, \mbox{skew} \left( \partial_{x_i}\boldsymbol{\mbox{D}u} \right)
+ \frac{2}{9} \, a_3 \, \mbox{tr} \left( \partial_{x_i}\boldsymbol{\mbox{D}u} \right) \boldsymbol{\mathbbm{1}}\Big)  n_{i}
\right)
\boldsymbol{n}
=
\boldsymbol{0}_{\mathbb{R}^{3}}
& \qquad \mbox{on} \, \partial \Omega
\notag
\end{align}
\section{Relations of parameters in terms of the classical Mindlin-Eringen formulation}
\label{app:coeff}
\subsection{Reduced Mindlin-Eringen formulation in terms of the classical Mindlin-Eringen formulation}
The curvature energy in eq.(\ref{eq:energy_MM_Mind_neff}) can be represented in index form as 
\begin{align}
W_{\tiny \mbox{curv}}
&=
\frac{\mu \, L_c^2}{2}
\Bigg[
a_1
\Bigg(
\frac{1}{2} \left( P_{jk} + P_{kj} \right)-\frac{1}{3}\delta_{jk}P_{mm}
\Bigg)_{,i}
\Bigg(
\frac{1}{2} \left( P_{jk} + P_{kj} \right)-\frac{1}{3}\delta_{jk}P_{mm}
\Bigg)_{,i}
\notag
\\
&
\qquad\qquad
+ a_2
\Bigg(
\frac{1}{2} \left( P_{jk} - P_{kj} \right)
\Bigg)_{,i}
\Bigg(
\frac{1}{2} \left( P_{jk} - P_{kj} \right)
\Bigg)_{,i}
+ \frac{2}{9}a_3
\delta_{jk}P_{mm,i}\delta_{jk}P_{nn,i}
\Bigg]
\notag
\\
&=
\frac{\mu \, L_c^2}{2}
\Bigg[
a_1
\Bigg(
\frac{1}{2} \left( P_{jk,i} + P_{kj,i} \right)-\frac{1}{3}\delta_{jk}P_{mm,i}
\Bigg)
\Bigg(
\frac{1}{2} \left( P_{jk,i} + P_{kj,i} \right)-\frac{1}{3}\delta_{jk}P_{mm,i}
\Bigg)
\notag
\\
&
\qquad\qquad
+ a_2
\Bigg(
\frac{1}{2} \left( P_{jk,i} - P_{kj,i} \right)
\Bigg)
\Bigg(
\frac{1}{2} \left( P_{jk,i} - P_{kj,i} \right)
\Bigg)
+ \frac{2}{3}a_3
P_{mm,i}P_{nn,i}
\Bigg]
\notag
\\
&=
\frac{\mu \, L_c^2}{2}
\Bigg[
a_1
\Bigg(
\frac{1}{4} \left( P_{jk,i} + P_{kj,i} \right)\left( P_{jk,i} + P_{kj,i} \right)
-\frac{1}{3} \left( P_{jk,i} + P_{kj,i} \right) \delta_{jk}P_{mm,i}
+\frac{1}{9}\delta_{jk}P_{mm,i}\delta_{jk}P_{nn,i}
\Bigg)
\notag
\\
&
\qquad\qquad
+ a_2
\frac{1}{4} 
\Bigg(
P_{jk,i}P_{jk,i} -2 P_{jk,i} P_{kj,i} + P_{kj,i}P_{kj,i}
\Bigg)
+ \frac{2}{3}a_3
P_{mm,i}P_{nn,i}
\Bigg]
\notag
\\
&=
\frac{\mu \, L_c^2}{2}
\Bigg[
a_1
\Bigg(
\frac{1}{4} \left( P_{jk,i} + P_{kj,i} \right)\left( P_{jk,i} + P_{kj,i} \right)
-\frac{2}{3} P_{nn,i}P_{mm,i}
+\frac{1}{3}P_{mm,i}P_{nn,i}
\Bigg)
\\
&
\qquad\qquad
+ a_2
\frac{1}{4} 
\Bigg(
P_{jk,i}P_{jk,i} -2 P_{jk,i} P_{kj,i} + P_{kj,i}P_{kj,i}
\Bigg)
+ \frac{2}{3}a_3
P_{mm,i}P_{nn,i}
\Bigg]
\notag
\\
&=
\frac{\mu \, L_c^2}{2}
\Bigg[
\frac{a_1}{4}
\left( P_{jk,i}P_{jk,i} + 2 P_{jk,i}P_{kj,i} + P_{kj,i}P_{kj,i} \right)
-
\frac{a_1}{3}
P_{nn,i}P_{mm,i}
\notag
\\
&
\qquad\qquad
+ \frac{a_2}{2} 
\Bigg(
P_{jk,i}P_{jk,i} - P_{jk,i} P_{kj,i}
\Bigg)
+ \frac{2}{3}a_3
P_{mm,i}P_{nn,i}
\Bigg]
\notag
\\
&=
\frac{\mu \, L_c^2}{2}
\Bigg[
\frac{a_1}{2}
\left( P_{jk,i}P_{jk,i} + P_{jk,i}P_{kj,i} \right)
+ \frac{a_2}{2} 
\Bigg(
P_{jk,i}P_{jk,i} - P_{jk,i} P_{kj,i}
\Bigg)
+ \frac{2a_3-a_1}{3}
P_{mm,i}P_{nn,i}
\Bigg]
\notag
\\
&=
\frac{\mu \, L_c^2}{2}
\Bigg[
\frac{a_1 + a_2}{2}
P_{jk,i}P_{jk,i}
+ \frac{a_1 - a_2}{2} 
P_{jk,i} P_{kj,i}
+ \frac{2a_3-a_1}{3}
P_{mm,i}P_{nn,i}
\Bigg]
\notag
\\
&=
\mu \, L_c^2 \,
\frac{a_1 + a_2}{4}
\chi_{ijk} \, \chi_{ijk}
+ 
\mu \, L_c^2 \,
\frac{a_1 - a_2}{4} 
\chi_{ijk} \, \chi_{ikj}
+ 
\mu \, L_c^2 \,
\frac{2a_3-a_1}{6}
\chi_{imm}\chi_{inn}
\, ,
\notag
\end{align}
which, compared to eq.(\ref{eq:energy_MM_Mind}), gives the following relations between the curvature parameters
\begin{align}
\widehat{a}_{1,2,3,5,8,11,14,15} &= 0 \, ,
\qquad \widehat{a}_{4} = \mu \, L_c^2 \frac{2a_3-a_1}{3} \, ,
\qquad \widehat{a}_{10} = \mu \, L_c^2 \frac{a_1 + a_2}{2} \, ,
\qquad \widehat{a}_{13} = \mu \, L_c^2 \frac{a_1 - a_2}{2}
\, .
\end{align}
\subsection{Relaxed micromorphic formulation in terms of the classical Mindlin-Eringen formulation}
Given the following expressions
\begin{align}
\mbox{sym Curl}\boldsymbol{P} &= \frac{1}{2}
\left(
\epsilon_{jnp}P_{in,p}
+ \epsilon_{inp}P_{jn,p}
\right)
e_{i} \otimes e_{j} \, ,
\notag
\\
\mbox{tr Curl}\boldsymbol{P} &= \epsilon_{knp}P_{kn,p} \, ,
\\
\mbox{dev sym Curl}\boldsymbol{P} &= \frac{1}{2}
\left(
\epsilon_{jnp}P_{in,p}
+ \epsilon_{inp}P_{jn,p}
- \frac{2}{3} \epsilon_{knp}P_{kn,p}\delta_{ij}
\right)
e_{i} \otimes e_{j} \, ,
\notag
\\
\mbox{skew Curl}\boldsymbol{P} &= \frac{1}{2}
\left(
\epsilon_{jnp}P_{in,p}
- \epsilon_{inp}P_{jn,p}
\right)
e_{i} \otimes e_{j} \, ,
\notag
\end{align}
it is possible to represent each single quadratic curvature term in eq.(\ref{eq:energy_RM}) in index form as 
\begin{align}
\rVert &\mbox{dev sym Curl}\boldsymbol{P} \lVert^2 = 
\notag
\\
&
\frac{1}{4}
\left(
\epsilon_{jnp} P_{in,p}
+ \epsilon_{inp} P_{jn,p}
- \frac{2}{3} \epsilon_{knp} P_{kn,p}\delta_{ij}
\right)
\left(
\epsilon_{jmq} P_{im,q}
+ \epsilon_{imq} P_{jm,q}
- \frac{2}{3} \epsilon_{qrs} P_{qr,s}\delta_{ij}
\right)
\notag
\\
=&
\frac{1}{2}
\left(
\epsilon_{jnp} P_{in,p}
+ \epsilon_{inp} P_{jn,p}
- \frac{2}{3} \epsilon_{knp} P_{kn,p}\delta_{ij}
\right)
\epsilon_{jmq} P_{im,q}
\notag
\\
=&
\frac{1}{2}
\left(
\epsilon_{jnp} P_{in,p}
+ \epsilon_{inp} P_{jn,p}\epsilon_{jmq} P_{im,q}
- \frac{2}{3} \epsilon_{knp} P_{kn,p}\delta_{ij}\epsilon_{jmq} P_{im,q}
\right)
\notag
\\
=&
\frac{1}{2}
\left(
\delta_{mn}\delta_{pq}-\delta_{nq}\delta_{mp}
\right)
P_{in,p}P_{im,q}
+\frac{1}{2}
\begin{vmatrix}
	\delta_{ij} & \delta_{im} & \delta_{iq} \\ 
	\delta_{nj} & \delta_{nm} & \delta_{nq} \\ 
	\delta_{pj} & \delta_{pm} & \delta_{pq} 
\end{vmatrix}
P_{jn,p}P_{im,q}
-\frac{1}{3}
\begin{vmatrix}
	\delta_{ki} & \delta_{km} & \delta_{kq} \\ 
	\delta_{ni} & \delta_{nm} & \delta_{nq} \\ 
	\delta_{pi} & \delta_{pm} & \delta_{pq} 
\end{vmatrix}
P_{kn,p}P_{im,q}
\notag
\\
=& \,
\frac{1}{2}
\left(
P_{im,p}P_{im,p} - P_{in,m}P_{im,n}
\right)
\\
&
+\frac{1}{2}
\left(
 \delta_{ij}\delta_{nm}\delta_{pq}
+\delta_{im}\delta_{nq}\delta_{pj}
+\delta_{iq}\delta_{nj}\delta_{pm}
-\delta_{pj}\delta_{nm}\delta_{iq}
-\delta_{pm}\delta_{nq}\delta_{ij}
-\delta_{pq}\delta_{nj}\delta_{im}
\right)P_{jn,p}P_{im,q}
\notag
\\
&
-\frac{1}{3}
\left(
 \delta_{ki}\delta_{nm}\delta_{pq}
+\delta_{km}\delta_{nq}\delta_{pi}
+\delta_{kq}\delta_{ni}\delta_{pm}
-\delta_{pi}\delta_{nm}\delta_{kq}
-\delta_{pm}\delta_{nq}\delta_{ki}
-\delta_{pq}\delta_{ni}\delta_{km}
\right)P_{kn,p}P_{im,q}
\notag
\\
=& \,
\frac{1}{2}
\left(
P_{im,p}P_{im,p} - P_{in,m}P_{im,n}
\right)
\notag
\\
&
+\frac{1}{2}
\left(
 P_{im,p}P_{im,p}
+P_{jn,j}P_{ii,n}
+P_{jj,m}P_{im,i}
-P_{jm,j}P_{im,i}
-P_{in,m}P_{im,n}
-P_{jj,p}P_{ii,p}
\right)
\notag
\\
&
-\frac{1}{3}
\left(
P_{im,p}P_{im,p}
+P_{kn,i}P_{ik,n}
+P_{ki,m}P_{im,k}
-P_{km,i}P_{im,k}
-P_{in,m}P_{im,n}
-P_{ki,p}P_{ik,p}
\right)
\notag
\\
=& \,
  \frac{2}{3} P_{im,p}P_{im,p}
- \frac{2}{3} P_{in,m}P_{im,n}
+ \frac{1}{3} P_{km,i}P_{im,k}
+ \frac{1}{3} P_{ki,p}P_{ik,p}
- \frac{1}{3} P_{kn,i}P_{ik,n}
- \frac{1}{3} P_{ki,m}P_{im,k}
\notag
\\
&
+ \frac{1}{2} P_{jn,j}P_{ii,n}
+ \frac{1}{2} P_{jj,m}P_{im,i}
- \frac{1}{2} P_{jm,j}P_{im,i}
- \frac{1}{2} P_{jj,p}P_{ii,p}
\notag
\\
=& \,
  \frac{2}{3} \, \chi_{ijk}\chi_{ijk}
- \frac{2}{3} \, \chi_{ijk}\chi_{kji}
+ \frac{1}{3} \, \chi_{ijk}\chi_{jik}
+ \frac{1}{3} \, \chi_{ijk}\chi_{ikj}
- \frac{1}{3} \, \chi_{ijk}\chi_{kij}
- \frac{1}{3} \, \chi_{ijk}\chi_{jki}
+              \chi_{iik}\chi_{kjj}
\notag
\\
&
- \frac{1}{2} \, \chi_{jjk}\chi_{iik}
- \frac{1}{2} \, \chi_{ijj}\chi_{ikk}
 \, ,
\notag
\end{align}
\begin{align}
\rVert &\mbox{skew Curl}\boldsymbol{P} \lVert^2 = 
\frac{1}{4}
\left(
\epsilon_{jnp} P_{in,p}
- \epsilon_{inp} P_{jn,p}
\right)
\left(
\epsilon_{jmq} P_{im,q}
- \epsilon_{imq} P_{jm,q}
\right)
\notag
\\
=&
\frac{1}{2}
\left(
  \epsilon_{jnp}P_{in,p}
- \epsilon_{inp}P_{jn,p}
\right)
\epsilon_{jmq}P_{im,q}
=
\frac{1}{2}
\left(
\epsilon_{jnp} \epsilon_{jmq} P_{in,p} P_{im,q}
- \epsilon_{inp} \epsilon_{jmq} P_{jn,p} P_{im,q}
\right)
\notag
\\
=&
\frac{1}{2}
\left(
\delta_{mn}\delta_{pq} - \delta_{nq}\delta_{mp}
\right)
P_{in,p} P_{im,q}
-\frac{1}{2}
\begin{vmatrix}
	\delta_{ij} & \delta_{im} & \delta_{iq} \\ 
	\delta_{nj} & \delta_{nm} & \delta_{nq} \\ 
	\delta_{pj} & \delta_{pm} & \delta_{pq} 
\end{vmatrix}
P_{jn,p} P_{im,q}
\notag
\\
=&
\frac{1}{2}
\left(
\delta_{mn}\delta_{pq} - \delta_{nq}\delta_{mp}
\right)
P_{in,p} P_{im,q}
\\
&
-
\frac{1}{2}
\left(
  \delta_{ij}\delta_{nm}\delta_{pq}
+ \delta_{im}\delta_{nq}\delta_{pj}
+ \delta_{iq}\delta_{nj}\delta_{pm}
- \delta_{pj}\delta_{nm}\delta_{iq}
- \delta_{pm}\delta_{nq}\delta_{ij}
- \delta_{pq}\delta_{nj}\delta_{im}
\right)
P_{jn,p} P_{im,q}
\notag
\\
=&
\frac{1}{2}
\left(
  P_{im,p} P_{im,p}
- P_{in,m} P_{im,n}
- P_{im,p} P_{im,p}
- P_{jn,j} P_{ii,n}
- P_{jj,m} P_{im,i}
+ P_{jm,j} P_{im,i}
\right.
\notag
\\*
&
\quad
\left.
+ P_{in,m} P_{im,n}
+ P_{jj,p} P_{ii,p}
\right)
\notag
\\
=&
\frac{1}{2}
\left(
  \chi_{ijk} \, \chi_{ijk}
- \chi_{ijk} \, \chi_{kji}
- \chi_{ijk} \, \chi_{ijk}
+ \chi_{ijk} \, \chi_{kji}
- \chi_{iik} \, \chi_{kjj}
- \chi_{iik} \, \chi_{kjj}
+ \chi_{iik} \, \chi_{jjk}
+ \chi_{ijj} \, \chi_{ikk}
\right)
\notag
\\
=&
\frac{1}{2}
\left(
- 2\chi_{iik} \, \chi_{kjj}
+  \chi_{iik} \, \chi_{jjk}
+  \chi_{ijj} \, \chi_{ikk}
\right)
 \, ,
\notag
\end{align}
\begin{align}
\mbox{tr}^2 \mbox{ Curl}\boldsymbol{P} =& \, 
\epsilon_{inp} P_{in,p} \epsilon_{jmq} P_{jm,q}
\notag
\\*
=&
\left(
  \delta_{ij}\delta_{nm}\delta_{pq}
+ \delta_{im}\delta_{nq}\delta_{pj}
+ \delta_{iq}\delta_{nj}\delta_{pm}
- \delta_{pj}\delta_{nm}\delta_{iq}
- \delta_{pm}\delta_{nq}\delta_{ij}
- \delta_{pq}\delta_{nj}\delta_{im}
\right)
P_{in,p} P_{jm,q}
\notag
\\*
=&
  P_{im,p} P_{im,p}
+ P_{in,j} P_{ji,n}
+ P_{ij,m} P_{jm,i}
- P_{im,j} P_{jm,i}
- P_{in,m} P_{im,n}
- P_{ij,p} P_{ji,p}
\\*
=&
  \chi_{ijk} \, \chi_{ijk}
+ \chi_{ijk} \, \chi_{jki}
+ \chi_{ijk} \, \chi_{kij}
- \chi_{ijk} \, \chi_{jik}
- \chi_{ijk} \, \chi_{kji}
- \chi_{ijk} \, \chi_{ikj}
 \, .
\notag
\end{align}
Inserting all the expressions written before in eq.(\ref{eq:energy_RM}) we have
\begin{align}
W_{\tiny \mbox{curv}}=&
\,
\frac{\mu \, L_c^2}{2}a_1
\Bigg(
  \frac{2}{3}\chi_{ijk}\chi_{ijk}
- \frac{2}{3}\chi_{ijk}\chi_{kji}
+ \frac{1}{3}\chi_{ijk}\chi_{jik}
+ \frac{1}{3}\chi_{ijk}\chi_{ikj}
\notag
\\
&
\qquad\qquad
- \frac{1}{3}\chi_{ijk}\chi_{kij}
- \frac{1}{3}\chi_{ijk}\chi_{jki}
+             \chi_{iik}\chi_{kjj}
- \frac{1}{2}\chi_{jjk}\chi_{iik}
- \frac{1}{2}\chi_{ijj}\chi_{ikk}
\Bigg)
\notag
\\
&
+
\frac{\mu \, L_c^2}{2}a_2
\Bigg(
\frac{1}{2}
\left(
- 2\chi_{iik}\chi_{kjj}
+ \chi_{iik}\chi_{jjk}
+ \chi_{ijj}\chi_{ikk}
\right)
\Bigg)
\notag
\\
&
+
\frac{\mu \, L_c^2}{6}a_3
\Bigg(
\chi_{ijk}\chi_{ijk}
+ \chi_{ijk}\chi_{jki}
+ \chi_{ijk}\chi_{kij}
- \chi_{ijk}\chi_{jik}
- \chi_{ijk}\chi_{kji}
- \chi_{ijk}\chi_{ikj}
\Bigg)
\notag
\\
=& \, 
  \frac{\mu \, L_c^2 \, a_1}{3}\chi_{ijk}\chi_{ijk}
- \frac{\mu \, L_c^2 \, a_1}{3}\chi_{ijk}\chi_{kji}
+ \frac{\mu \, L_c^2 \, a_1}{6}\chi_{ijk}\chi_{jik}
+ \frac{\mu \, L_c^2 \, a_1}{6}\chi_{ijk}\chi_{ikj}
\\
&
- \frac{\mu \, L_c^2 \, a_1}{6}\chi_{ijk}\chi_{kij}
- \frac{\mu \, L_c^2 \, a_1}{6}\chi_{ijk}\chi_{jki}
+ \frac{\mu \, L_c^2 \, a_1}{2}\chi_{iik}\chi_{kjj}
- \frac{\mu \, L_c^2 \, a_1}{4}\chi_{jjk}\chi_{iik}
\notag
\\
&
- \frac{\mu \, L_c^2 \, a_1}{4}\chi_{ijj}\chi_{ikk}
- \frac{\mu \, L_c^2 \, a_2}{2}\chi_{iik}\chi_{kjj}
+ \frac{\mu \, L_c^2 \, a_2}{4}\chi_{iik}\chi_{jjk}
+ \frac{\mu \, L_c^2 \, a_2}{4}\chi_{ijj}\chi_{ikk}
\notag
\\
&
+ \frac{\mu \, L_c^2 \, a_3}{6}\chi_{ijk}\chi_{ijk}
+ \frac{\mu \, L_c^2 \, a_3}{6}\chi_{ijk}\chi_{jki}
+ \frac{\mu \, L_c^2 \, a_3}{6}\chi_{ijk}\chi_{kij}
- \frac{\mu \, L_c^2 \, a_3}{6}\chi_{ijk}\chi_{jik}
\notag
\\
&
- \frac{\mu \, L_c^2 \, a_3}{6}\chi_{ijk}\chi_{kji}
- \frac{\mu \, L_c^2 \, a_3}{6}\chi_{ijk}\chi_{ikj}
\notag
\\
=& \, 
  \frac{\mu \, L_c^2 \, \left(2a_1 + a_3\right)}{6}\chi_{ijk}\chi_{ijk}
- \frac{\mu \, L_c^2 \, \left(2a_1 + a_3\right)}{6}\chi_{ijk}\chi_{kji}
+ \frac{\mu \, L_c^2 \, \left( a_1 - a_3\right)}{6}\chi_{ijk}\chi_{jik}
\notag
\\
&
+ \frac{\mu \, L_c^2 \, \left( a_1 - a_3\right)}{6}\chi_{ijk}\chi_{ikj}
+ \frac{\mu \, L_c^2 \, \left( a_3 - a_1\right)}{6}\chi_{ijk}\chi_{kij}
+ \frac{\mu \, L_c^2 \, \left( a_3 - a_1\right)}{6}\chi_{ijk}\chi_{jki}
\notag
\\
&
+ \frac{\mu \, L_c^2 \, \left(2a_1 - a_2\right)}{4}\chi_{iik}\chi_{kjj}
+ \frac{\mu \, L_c^2 \, \left( a_2 - a_1\right)}{4}\chi_{iik}\chi_{jjk}
+ \frac{\mu \, L_c^2 \, \left( a_2 - a_1\right)}{4}\chi_{ijj}\chi_{ikk}
\notag
\end{align}
which compared to eq.(\ref{eq:energy_MM_Mind}) gives the following relations between the curvature parameters
\begin{align}
\quad
\widehat{a}_{1} &= \mu \, L_c^2 \, \frac{2a_1 - a_2}{4} \, ,
\quad
\widehat{a}_{2} = \widehat{a}_{5} = \widehat{a}_{8} = 0 \, ,
\quad
\widehat{a}_{3} = \widehat{a}_{4} = \mu \, L_c^2 \, \frac{a_2 - a_1}{2} \, ,
\quad
\widehat{a}_{10} = \mu \, L_c^2 \, \frac{2a_1 + a_3}{3} \, ,
\\
\widehat{a}_{11} &= \mu \, L_c^2 \, \frac{a_3 - a_1}{6} \, ,
\quad
\widehat{a}_{13} =\widehat{a}_{14} = \mu \, L_c^2 \, \frac{a_1 - a_3}{3} \, ,
\quad
\widehat{a}_{15} = - \mu \, L_c^2 \, \frac{2a_1 + a_3}{3}
\, .
\notag
\end{align}

The most simple curvature $\frac{\mu \, L_c^2}{2} \rVert \mbox{Curl} \boldsymbol{P} \lVert^2$ corresponds to $a_1 = a_2 = a_3 = 1$ and is therefore given by $\widehat{a}_{2,3,4,5,8,11,13,14}=0$, $\widehat{a}_{1}=\frac{\mu \, L_c}{4}$, $\widehat{a}_{10} = -\widehat{a}_{15} = \mu \, L_c$.
\subsection{Second Gradient formulation in terms of the classical Mindlin-Eringen formulation}
Given that $\chi_{ijk} := u_{k,ij} = \chi_{jik}$, the curvature energy in eq.(\ref{eq:energy_SG}) can be represented in index form as 
\begin{align}
	&W_{\tiny \mbox{curv}}
	\notag
	\\
	&=
	\frac{\mu \, L_c^2}{2}
	\Bigg[
	a_1
	\Bigg(
	\frac{1}{2} \left( u_{i,jk} + u_{j,ik} \right)-\frac{1}{3}\delta_{ij}u_{m,mk}
	\Bigg)
	u_{i,jk}
	+
	\frac{1}{2} a_2
	\left( u_{i,jk} - u_{j,ik} \right)
	u_{i,jk}
	+ \frac{2}{9}a_3 \, 
	u_{m,mk} \, u_{i,ik}
	\Bigg]
	\notag
	\\
	&=
	\frac{\mu \, L_c^2}{2}
	\Bigg[
	a_1
	\Bigg(
	\frac{1}{2} \left( \chi_{jki} + \chi_{ikj} \right)-\frac{1}{3}\delta_{ij}\chi_{kmm}
	\Bigg)
	\chi_{jki}
	+
	\frac{1}{2} a_2
	\left( \chi_{jki} - \chi_{ikj} \right)
	\chi_{jki}
	+ \frac{2}{9}a_3 \, 
	\chi_{kmm} \, \chi_{kii}
	\Bigg]
	\\
	&=
	\frac{\mu \, L_c^2}{2}
	\Bigg[
	  \frac{1}{2} a_1 \, \chi_{jki} \, \chi_{jki}
	+ \frac{1}{2} a_1 \, \chi_{ikj} \, \chi_{jki}
	- \frac{1}{3} a_1 \, \chi_{kmm} \, \chi_{jki}
	+ \frac{1}{3} a_2 \, \chi_{jki} \, \chi_{jki}
	+ \frac{1}{3} a_2 \, \chi_{ikj} \, \chi_{jki}
	+ \frac{2}{9} a_3 \, \chi_{kmm} \, \chi_{kii}
	\Bigg]
	\notag
	\\
	&=
	\frac{\mu \, L_c^2}{2}
	\Bigg[
	  \frac{a_1 + a_2}{2} \, \chi_{jki} \, \chi_{jki}
	+ \frac{a_1 - a_2}{2} \, \chi_{ikj} \, \chi_{jki}
	- \frac{1}{3} a_1 \, \chi_{kmm} \, \chi_{jki}
	+ \frac{2}{9} a_3 \, \chi_{kmm} \, \chi_{kii}
	\Bigg]
	\notag
	\\
	&=
	\frac{\mu \, L_c^2}{2}
	\Bigg[
	  \frac{a_1 + a_2}{2} \, \chi_{ijk} \, \chi_{ijk}
	+ \frac{a_1 - a_2}{2} \, \chi_{ijk} \, \chi_{kji}
	+ \frac{2a_3 - 3a_1}{9} a_3 \, \chi_{ijj} \, \chi_{ikk}
	\Bigg]
	\, ,
	\notag
\end{align}
which compared to eq.(\ref{eq:energy_SG_gen}) gives the following relations between the curvature parameters
\begin{align}
	\widehat{a}_{1,3} &= 0 \, ,
	\qquad\qquad \widehat{a}_{2} = \mu \, L_c^2 \frac{2a_3-3a_1}{18} \, ,
	\qquad\qquad \widehat{a}_{4} = \mu \, L_c^2 \frac{a_1 + a_2}{4} \, ,
	\qquad\qquad \widehat{a}_{5} = \mu \, L_c^2 \frac{a_1 - a_2}{4}
	\, .
\end{align}
\end{footnotesize}
%
%
%
\end{document}